\documentclass[11pt, oneside, letterpaper]{book}

\usepackage{amsthm}
\usepackage{amssymb, mathrsfs, amsthm, amsmath, graphicx}
\usepackage{stmaryrd}
\usepackage[all]{xy}
\usepackage{enumerate}
\usepackage{array}
\usepackage{changepage}

\usepackage[usenames,dvipsnames]{xcolor}
\usepackage[pagebackref]{hyperref}
\hypersetup{
	linktocpage,
    colorlinks,
    citecolor=black,
    filecolor=black,
    linkcolor=black,
    urlcolor=black,
}

\changepage{2cm}{2cm}{-1cm}{-1cm}{}{-1cm}{}{}{}

\usepackage[small, md, sc]{titlesec}
\usepackage{lettrine} 
\usepackage{fancyhdr} 
\usepackage[palatino]{quotchap}
\definecolor{chaptergrey}{rgb}{0,0,0}
\usepackage{titling}
\usepackage{setspace} 
\usepackage{booktabs} 
\usepackage[labelfont={sf,bf,small},textfont={sf,small},justification=RaggedRight,margin=0pt, figurewithin=section, tablewithin=section]{caption}
\renewcommand{\AA}{\mathbb{A}}
\newcommand{\A}{\mathcal{A}}
\newcommand{\C}{\mathcal{C}}

\newcommand{\E}{\mathcal{E}}
\newcommand{\uE}{\underline{\E}}

\newcommand{\F}{\mathcal{F}}

\newcommand{\G}{\mathcal{G}}
\newcommand{\GG}{\mathbb{G}}
\renewcommand{\H}{\mathsf{H}}
\newcommand{\HH}{\mathbb{H}}

\newcommand{\m}{{m}} 
\newcommand{\M}{\mathscr{M}}

\newcommand{\OO}{\mathcal{O}}

\renewcommand{\P}{\mathcal{P}}
\newcommand{\PP}{\mathbb{P}}

\newcommand{\QQ}{\mathbb{Q}}

\renewcommand{\S}{\mathcal{S}}
\newcommand{\SH}{\mathsf{SH}}

\newcommand{\T}{\mathscr{T}}

\newcommand{\W}{\mathcal{W}}
\newcommand{\X}{\mathcal{X}}
\newcommand{\Y}{\mathcal{Y}}
\newcommand{\Z}{\mathcal{Z}}
\newcommand{\ZZ}{\mathbb{Z}}
\newcommand{\zll}{{\mathbb{Z}_{(\ell)}}}
\newcommand{\zpi}{\mathbb{Z}[\tfrac{1}{p}]}
\newcommand{\ozpi}{[\tfrac{1}{p}]}
\newcommand{\HZ}{\mathsf{H}\mathbb{Z}}
\newcommand{\HZl}{{\mathsf{H}\mathbb{Z}_{(\ell)}}}
\newcommand{\HZpi}{{\mathsf{H}\mathbb{Z}[\tfrac{1}{p}]}}

\newcommand{\KH}{\mathsf{KH}}

\renewcommand{\mp}[1]{M(#1)[\tfrac{1}{p}]}
\newcommand{\ml}[1]{M(#1)_{(\ell)}}
\newcommand{\mcp}[1]{M^c(#1)[\tfrac{1}{p}]}

\newcommand{\dua}{\mathcal{D}}
\newcommand{\dis}{cdd}
\newcommand{\id}{\mathrm{id}}
\newcommand{\im}{\mathrm{im}}
\newcommand{\ldh}{\ensuremath{\ell dh}}
\newcommand{\fpsl}{\ensuremath{\mathrm{fps}\ell'}}
\newcommand{\fps}{\ensuremath{\mathrm{fps}}}

\newcommand{\un}{\ensuremath{1\!\!1}}

\newcommand{\Tr}{\mathrm{Tr}}
\newcommand{\Ext}{\mathrm{Ext}}
\newcommand{\sing}{\mathrm{sing}}

\newcommand{\corr}{\ensuremath{\bullet\!\!\to}}
\renewcommand{\mod}{\ensuremath{\textrm{-mod}}}
\newcommand{\eff}{\ensuremath{^{{eff}}}}
\newcommand{\gm}{\ensuremath{_{{gm}}}}

\newtheorem{theo}{Theorem}[section]
\newtheorem*{theoNo}{Theorem}

\newtheorem{coro}[theo]{Corollary}
\newtheorem{lemm}[theo]{Lemma}

\newtheorem{prop}[theo]{Proposition}

\theoremstyle{definition}
\newtheorem{defi}[theo]{Definition}
\newtheorem*{defiNo}{Definition}
\newtheorem{rema}[theo]{Remark}
\newtheorem{exam}[theo]{Example}

\DeclareMathOperator{\hocolim}{hocolim}

\DeclareMathOperator{\ihom}{\underline{hom}}
\DeclareMathOperator{\length}{length}
\DeclareMathOperator{\rs}{rs}

\renewcommand{\maketitle}{ 
	\begin{onehalfspacing}
	\thispagestyle{empty}
	\begin{center}
	\normalsize \textbf{UNIVERSIT{\'E} PARIS 13 - Institut Galil{\'e}e \\
	Laboratoire Analyse, G{\'e}om{\'e}trie et Applications, UMR 7539} \\
	\vspace{30pt}
	\large \textbf{TH{\`E}SE} \\
	\normalsize pour obtenir le grade de \\
	\large \textbf{DOCTEUR DE L'UNIVERSIT{\'E} PARIS 13} \\
	\normalsize Discipline : \textbf{Math{\'e}matiques} \\
	\vspace{30pt}
	\normalsize pr{\'e}sent{\'e}e et soutenue publiquement par \\
	\large \textbf{Shane KELLY} \\
	\normalsize le 19 octobre 2012 \\
	\vspace{30pt}
	\LARGE \textbf{Triangulated categories of motives in positive characteristic} \\
	\vspace{30pt}
	\normalsize \textbf{Directeurs de th{\`e}se} \\
\begin{table}[h]
\begin{tabular}{  >{\raggedright}p{7cm}  >{\raggedleft}p{7cm}  }
M. Denis-Charles CISINSKI & Universit{\'e} Paul Sabatier \tabularnewline
M. Amnon NEEMAN & Australian National University
\end{tabular}
\end{table}
	\normalsize \textbf{Rapporteurs} \\
\begin{table}[h]
\begin{tabular}{  >{\raggedright}p{7cm}  >{\raggedleft}p{7cm}  }
M. Paul-Arne {\O}STV{\AE}R  & Universitetet i Oslo \tabularnewline
M. Charles WEIBEL & Rutgers University
\end{tabular}
\end{table}
	\normalsize \textbf{Jury} \\
\begin{table}[h]
\begin{tabular}{  >{\raggedright}p{7cm}  >{\raggedleft}p{7cm}  }
M. Yves ANDR{\'E}  &  Universit{\'e} de Paris 6, CNRS \tabularnewline
M. Pascal BOYER &  Universit{\'e} Paris 13 \tabularnewline
M. Denis-Charles CISINSKI &  Universit{\'e} Paul Sabatier \tabularnewline
M. Fr{\'e}d{\'e}ric D{\'E}GLISE & ENS Lyon, CNRS \tabularnewline
M. Ofer GABBER & IHES, CNRS \tabularnewline
M. Bruno KAHN & Universit{\'e} Paris 7, CNRS \tabularnewline
M. J{\"o}rg WILDESHAUS & Universit{\'e} Paris 13 \tabularnewline
\end{tabular}
\end{table}	
	\end{center}
\end{onehalfspacing}
\singlespacing
\clearpage\mbox{}\thispagestyle{empty}\clearpage
}

\newenvironment{acknowledgements}
{\clearpage\thispagestyle{empty}\null\vfill\begin{center}\bfseries Acknowledgements\end{center}\small}
{\normalsize\vfill\null\clearpage\mbox{}\thispagestyle{empty}\clearpage }

\title{Triangulated categories of motives in positive characteristic}
\author{Shane Kelly}

\begin{document}

\maketitle

\begin{acknowledgements}

I would like to express my gratitude to Denis-Charles Cisinski for being the ideal supervisor. Incredibly generous with his time, ideas, encouragement, and indefatigable optimism, I am very fortunate to have had the opportunity to work with him.

I would like to thank Amnon Neeman for cosupervising my thesis, for his confidence and support.


I would like to thank the incredible generosity of the ANU Vice-Chancellor's Scholarship for providing me with enough funding to live in the mathematical hub that is Paris, attend numerous conferences, and the independance I was given to use the travel money as a I saw fit.

Chuck Weibel and Paul-Arne {\O}˜stv{\ae}r accepted the ungratifying task of \emph{rapporteur}. I am very appreciative for their close reading and feedback.

I am honored to have Yves Andr{\'e}, Pascal Boyer, Ofer Gabber, Bruno Kahn, and J{\"o}rg Wildeshaus in my jury. I thank all of them.

I thank Fr{\'e}d{\'e}ric D{\'e}glise for his course ``Complexes Motiviques'', for the groupe de travail ``L'hypoth{\`e}se de Riemann d'apr{\`e}s Deligne'', but mostly for the many useful and interesting discussions that he afforded me on many topics, especially relative cycles and trace morphisms for $K$-theory.

Over the course of preparing this work, I profited from useful discussions with Marc Levine, Thomas Geisser, Chuck Weibel, Jo{\"e}l Riou, Mark Walker, Aravind Asok,  Pablo Pelaez, Simon Pepin Lehalleur, Marc Hoyois, James Wallbridge, Valentina Sala and many others who I thank for their time and encouragement.

For comments on earlier versions of this work and for presentation suggestions I am grateful to Chuck Weibel, Thomas Geisser, Ariyan Javanpeykar, Giuseppe Ancona, Nicola Mazzari, Brad Drew, and Javier Fres{\`a}n.

I thank Shuji Saito for his hospitality at Tokyo Institute of Technology where Chapter~\ref{chap:applications} was written. I am also grateful for the opportunity I was given to speak in the Number Theory Seminar at The University of Tokyo. During that stay I also had many useful discussions with Thomas Geisser and I thank him for the invitation to speak at The University of Nagoya.

I am grateful to Vicky Hoskins for her hospitality at the University of Oxford where Section~\ref{sec:tracesToTransfers} was completed.

I am obliged for the opportunity to speak at the University of Zurich in the Oberseminar in Algebraic Geometry which resulted in my discovery of a much more streamlined proof of the $K$-theory vanishing.

I thank Yolande Jimenez and Isabelle Barbotin at Paris XIII, and Kelly Wicks and Alison Anderson at ANU for their incredible administrative efficiency.

I also thank Tristan Buckmaster, Lashi Bandara, and Hemanth Saratchandran for their encouragement and support.


\end{acknowledgements}

\tableofcontents

\chapter{Introduction}


\section{Main results}

\lettrine[lines=2, findent=3pt, nindent=0pt]{S}{ome} might consider the following theorems as the principal applications of the material in this thesis.

\begin{theoNo}[{Theorem~\ref{theo:suslin}, cf. \cite[Theorem 4.2, Corollary 4.3]{Sus00}}]
Let $X$ be an equidimensional quasi-projective scheme over an algebraically closed field $k$. Let $i \geq d = \dim X$ and suppose that $m$ is coprime to the characteristic of $k$. Then
\[ CH^i(X, n; \ZZ/m) \cong H_c^{2(d - i) + n}(X, \ZZ/m(d - i))^\# \]
where $H_c$ is the {\'e}tale cohomology with compact supports. If the scheme $X$ is smooth then this formula simplifies to $CH^i(X, n; \ZZ/m) \cong H_{\acute{e}t}^{2i - n}(X, \ZZ/m(i))$.
\end{theoNo}

\begin{theoNo}[{Theorem~\ref{theo:FV82}, Theorem~\ref{theo:FV83}, cf. \cite[Theorem 8.2, Theorem 8.3]{FV}}]
Let $k$ be a perfect field of exponential characteristic $p$, let $U$ be a smooth scheme of pure dimension $n$ over $k$, and let $X, Y$ be separated schemes of finite type over $k$. We denote by $A_{r, i}(Y, X)$ the bivariant cycle cohomology of \cite[Definition 4.3]{FV}. There are canonical isomorphisms
\[ A_{r,i}(Y \times U, X)\ozpi \cong A_{r + n, i}(Y, X \times U)\ozpi. \]
We also have the following properties:
\begin{enumerate}
 \item (Homotopy invariance) The pull-back homomorphism $z_{equi}(X, r) \to z_{equi}(X \times \AA^1, r + 1)$ induces for any $i \in \ZZ$ an isomorphism
\[ A_{r,i}(Y, X)\ozpi \stackrel{\sim}{\to} A_{r + 1, i}(Y, X \times \AA^1)\ozpi. \]
 \item (Suspension) Let
\[ p: X \times \PP^1 \to X \]
\[ i: X \to X \times \PP^1 \]
be the natural projection and closed embedding. Then the morphism
\[ i_* \oplus p^*: z_{equi}(X, r + 1) \oplus z_{equi}(X, r) \to z_{equi}(X \times \PP^1, r + 1) \]
induces an isomorphism
\[ A_{r + 1, i}(Y, X)\ozpi \oplus A_{r, i}(Y, X)\ozpi \stackrel{\sim}{\to} A_{r + 1, i}(Y, X \times \PP^1)\ozpi. \]
 \item (Cosuspension) There are canonical isomorphisms:
\[ A_{r,i}(Y \times \PP^1, X)\ozpi \stackrel{\sim}{\to} A_{r + 1, i}(Y, X)\ozpi \oplus A_{r,i}(Y, X)\ozpi. \]
 \item (Gysin) Let $Z \subset U$ be a closed immersion of smooth schemes everywhere of codimension $c$ in $U$. Then there is a canonical long exact sequence of abelian groups of the form
\[ \dots A_{r + c, i}(Z, X)\ozpi \to A_{r,i}(U, X)\ozpi \to A_{r, i}(U - Z, X)\ozpi \]
\[ \to A_{r + c, i - 1}(Z, X)\ozpi \to \dots \]
\end{enumerate}
\end{theoNo}

\begin{theoNo}[{cf. \cite[Corollary 3.5.5, 4.1.4, 4.1.6, Theorem 4.3.7]{Voev00}}]
Let $k$ be a perfect field of exponential characteristic $p$.
\begin{enumerate}
 \item (Lemma~\ref{lemm:geometricMotive}, Lemma~\ref{lemm:geometricCompactMotive}) The subcategory $DM_{gm}(k, \zpi) \subset DM(k, \zpi)$ contains the objects $\mp{X} = \underline{C}_*(c_{equi}(X/k, 0))$ and $\mcp{X} = \underline{C}_*(z_{equi}(X/k, 0))$ and for any separated scheme $X$ of finite type over $k$.
 \item (Proposition~\ref{prop:generatedBysmproj}) $DM_{gm}(k, \zpi)$ is generated by $\mp{X}$ for $X$ smooth and projective over $k$.
 \item $DM\gm(k, \zpi)$ has an internal hom.
 \item (Theorem~\ref{theo:dualityDMp}) Denoting $A^* = \ihom_{DM\gm(k, \zpi)}(A, \zpi)$ one has:
\begin{enumerate}
 \item For any object $A$ in $DM\gm(k, \zpi)$ the canonical morphism $A \to (A^*)^*$ is an isomorphism.
 \item For any pair of objects $A, B$ of $DM\gm(k, \zpi)$ there are canonical morphisms
\[ (A \otimes B)^* \cong A^* \otimes B^* \]
\[ \ihom_{DM\gm(k, \zpi)}(A, B) \cong A^* \otimes B. \]
 \item For a smooth scheme $X$ of pure dimension $n$ over $k$ one has canonical isomorphisms
\[ \mp{X}^* \cong \mcp{X}(-n)[-2n] \]
\[ \mcp{X}^* \cong \mp{X}(-n)[-2n]. \]
\end{enumerate}
\end{enumerate}
\end{theoNo}

In \cite[2.9]{Wei80} Weibel asks if $K_n(X) = 0$ for $n < - \dim X$ for every noetherian scheme $X$ where $K_n$ is the $K$-theory of Bass-Thomason-Trobaugh. This question was answered in the affirmative in \cite{CHSW} for schemes essentially of finite type over a field of characteristic zero. Assuming strong resolution of singularities, it is also answered in the affirmative in \cite{GH10} for schemes essentially of finite type over a field of positive characteristic. Both of these proofs compare $K$-theory with cyclic homology, and then use a cdh descent argument. The main new ingredients in the following theorem are the representability of homotopy invariant algebraic $K$-theory (\cite{Wei89}) in the Morel-Voevodsky stable homotopy category (\cite{Cis13}) and a theorem of Gabber (a weak version is mentioned further down in this introduction).

\begin{theoNo}[{Theorem~\ref{theo:KtheoryVanishing}}]
Let $X$ be a quasi-excellent noetherian scheme and $p$ a prime that is nilpotent on $X$. Then $K_n(X) \otimes \zpi = 0$ for $n < - \dim X$.
\end{theoNo}

The first three theorems above appear in \cite{FV}, \cite{Voev00}, and \cite{Sus00} under the assumption of resolution of singularities. The resolution of singularities assumption is applied through the theorem \cite[Theorem 4.1.2]{Voev00} of Voevodsky. Our main technical result is the following $\zpi$-linear version of this theorem.

\begin{theoNo}[{Theorem~\ref{theo:MainResultRos}, cf. \cite[Theorem 4.1.2]{Voev00}}]
Let $k$ be a perfect field of exponential characteristic $p$. Suppose that $F$ is a presheaf with transfers on the category $Sch(k)$ of separated $k$-schemes of finite type such that $F_{cdh} \otimes \zpi = 0$. Then $(\underline{C}_*(F)_{Nis} \otimes \zpi)|_{Sm(k)}$ is quasi-isomorphic to zero as a complex of Nisnevich sheaves on $Sm(k)$.
\end{theoNo}

With this $\zpi$-linear version version of \cite[Theorem 4.1.2]{Voev00}, Suslin's proof of Theorem~\ref{theo:suslin} goes through unchanged. Similarly, $\zpi$ versions of the material in \cite{FV} and \cite{Voev00} that previously assumed resolution of singularities are now all valid with minor changes to some arguments.

To apply resolution of singularities in his work on algebraic cycle cohomology theories, Voevodsky introduced the cdh topology. This is an enlargement of the Nisnevich topology so that the proper birational morphisms coming from resolution of singularities may be used as covers. We will use the following theorem of Gabber as a replacement for resolution of singularities.

\begin{theoNo}[Gabber {\cite[1.3]{Ill09} or \cite[Theorem 3, Theorem 3.2.1]{ILO}}]
Let $X$ be a separated scheme of finite type over a perfect field $k$ and $\ell$ a prime distinct from the characteristic of $k$. There exists a smooth quasi-projective $k$ scheme $Y$, and a $k$-morphism $f: Y \to X$ such that
\begin{enumerate}
 \item $f$ is proper, surjective, and sends every generic point to a generic point, and
 \item for each generic point $\xi$ of $X$ there is a unique point $\eta$ of $Y$ over it, and $[k(\eta):k(\xi)]$ is finite of degree prime to $\ell$.
\end{enumerate}
\end{theoNo}

To apply this theorem of Gabber, we need to enlarge the Nisnevich topology further. Hence, we search a topology which contains the Nisnevich topology and for which we can use Gabber's theorem to show that every scheme of finite type admits a smooth quasi-projective covering. There are many possible choices, some making the proofs easier than others. The following definition could be considered as the first new contribution of this thesis.

\begin{defiNo}[{Definition~\ref{defi:ltopology}}]
Let $\ell$ be a prime and $S$ a scheme. The \emph{$\ldh$ topology} is the Grothendieck topology on the category $Sch(S)$ of schemes of finite type over $S$ generated by the cdh topology, and the pretopology for which the covers are singletons $\{ Y \stackrel{f}{\to} X \}$ containing a finite flat surjective morphism of constant degree prime to $\ell$, such that $f_*\OO_Y$ is a globally free $\OO_X$-module.
\end{defiNo}

We will call the latter pretopology the \emph{$\fpsl$ pretopology}.

We have defined the $\ldh$ topology in such a way that proofs may be reduced to a cdh part, and an $\fpsl$ part. To structure the proofs dealing with $\fpsl$ part, we formalise a notion of presheaf with traces. In our work the class $\P$ in the following definition will always consist of finite flat surjective morphisms.

\begin{defiNo}[{Definition~\ref{defi:traces}}]
A presheaf with traces $(F, \S, \A, \Tr, \P)$ is an additive functor $F: \S^{op} \to \A$ from a category of schemes $\S$ to an additive category $\A$, together with a class $\P$ of morphisms of $\S$, and a morphism $\Tr_f: F(Y) \to F(X)$ for every morphism $f \in \P$. The morphisms $\Tr$ are required to satisfy the following axioms.
\begin{enumerate}
 \item[(Add)] For morphisms $f_1: Y_1 \to X_1$ and $f_2: Y_2 \to X_2$ in $\P$ we have
\[ \Tr_{f_1 \amalg f_2} = \Tr_{f_1} \oplus \Tr_{f_2}. \]
 \item[(Fon)] For morphisms $W \stackrel{g}{\to} Y \stackrel{f}{\to} X$ in $\P$ we have 
\[ \Tr_f \Tr_g = \Tr_{fg} \textrm{ and } \Tr_{\id_X} = \id_{F(X)}. \]
 \item[(CdB)] For every cartesian square in $\S$
\[\xymatrix{
Y \times_X W \ar[r]^-g \ar[d]_q & W \ar[d]^p \\
Y  \ar[r]_f & X
} \]
such that $f, g \in \P$ we have
\[ F(p)\Tr_f = \Tr_gF(q). \]
 \item[(Deg)] For every finite flat surjective morphism $f: Y \to X$ in $\P$ such that $f_*\OO_Y$ is a globally free $\OO_X$ module we have 
\[ \Tr_f F(f) = \deg f \cdot \id_{F(X)}. \]
\end{enumerate}
\end{defiNo}

It falls straight out of the definition that if $(F, Sch(S), \zll\mod, \Tr, \P)$ is a presheaf with traces such and $\P$ contains each finite flat surjective morphism of degree prime to $\ell$, then $F$ is an $\fpsl$ sheaf (Lemma~\ref{lemm:fpslAcyclic}). Hence, a cdh sheaf with traces of that form is an $\ldh$ sheaf.

\section{Outline}

We outline now the proof of Theorem~\ref{theo:MainResultRos}. It suffices to show $\underline{C}_*(F)_{Nis} \otimes \zll$ quasi-isomorphic to zero for each $\ell \neq p$. As the $\ldh$ topology is finer than the cdh topology, clearly we can assume that $F_{\ldh} \otimes \zll = 0$. As $F_{\ldh} \otimes \zll$ is the image of $F$ in the derived category of $\ldh$ sheaves of $\zll$-modules, and $\underline{C}_*(F)_{Nis} \otimes \zll$ is the image of $F$ in $DM\eff(k, \zll)$, our result will follow if we can find a factorisation 
\[ D(Shv_{Nis}(SmCor(k), \zll)) \to D(Shv_{\ldh}(SmCor(k), \zll)) \to DM\eff(k, \zll). \]
We can describe the derived category of Nisnevich (resp. $\ldh$) sheaves as the derived category of presheaves with Nisnevich (resp. $\ldh$) hypercovers inverted. To obtain such a factorisation then, it suffices to show that every $\ldh$ hypercover in $DM(k)$ is isomorphic to the scheme it covers. The functor $Sm \to DM(k, \zll)$ factors through the category of modules over the ring spectrum $\HZl$ in $\SH(k)$ where $\HZl$ is the object of the Morel-Voevodsky stable homotopy category $\SH(k)$ that represents motivic cohomology with $\zll$ coefficients. So we have converted the problem into showing that each $\ldh$ hypercover in this category of modules is isomorphic to the scheme it covers, or equivalently, that every such object of $\SH(k)$ satisfies $\ldh$-descent.

This descent problem is broken up into a cdh part, and a $\fpsl$ part. The cdh part is taken care of by a theorem of Cisinski (\cite[3.7]{Cis13}) which applies work of Ayoub (\cite[Corollary 1.7.18]{Ayo07}) and Voevodsky (\cite{Voe10a}, \cite{Voe10b}) to show that every object of $\SH(k)$ satisfies cdh descent. For the $\fpsl$ part we define a suitable notion of what it means for an object of $\SH(k)$ to have a structure of traces (Definition~\ref{defi:sectionWithTraces}), and show that certain objects in $\SH(k)$ have such a structure of traces (Corollary~\ref{coro:tracesOnHZmod}). To show that cdh descent and a structure of traces implies $\ldh$ descent, we compare the cdh and $\ldh$ descent spectral sequences and use the following theorem, which is a weak summary of the results of Chapter~\ref{chap:comparison}.

\begin{theoNo}[{Theorem~\ref{theo:cdhlprimeComparison}}]
Let $k$ be a perfect field and $\ell$ a prime that is invertible in $k$. Let $F$ be a presheaf of $\zll$-modules with traces on $Sch(k)$, such that
\begin{enumerate}
 \item $F(X) \to F(X_{red})$ is an isomorphism for every $X \in Sch(k)$,
 \item $F(X) \to F(\AA^1_X)$ is an isomorphism for every $X \in Sm(k)$, and
 \item $F|_{Sm(k)}$ has a structure of transfers,
\end{enumerate}
then for every $n \in \ZZ_{\geq 0}$ and every $X \in Sch(S)$, the canonical morphism
\[ H^n_{cdh}(X, F_{cdh}) \to H^n_{\ldh}(X, F_{\ldh}) \]
is an isomorphism.
\end{theoNo}

We actually prove Theorem~\ref{theo:cdhlprimeComparison} with much weaker hypotheses (see Theorem~\ref{theo:bigComparisonTheorem}). The statement above is designed to be applied to the homotopy presheaves of an oriented $\zll$ local object of $\SH(k)$.

To obtain a structure of traces on $\HZl$, which is the hard part, we proceed as follows. We can define a structure of traces on $\KH$, the object representing homotopy invariant algebraic $K$-theory, fairly naturally (see \cite{Wei89} for homotopy invariant algebraic $K$-theory and \cite{Cis13} for its representability in the Morel-Voevodsky stable homotopy category). We then notice that $\HZ$ is the zero slice of $\KH$ due to work of Levine (\cite[Theorems 6.4.2 and 9.0.3]{Lev08}). Traces on $\HZ$ now follows from the following theorem, which is the main goal of Chapter~\ref{chap:slices}.

\begin{theoNo}[{Proposition~\ref{prop:tracesOnSlices}}]
Suppose that $k$ is a perfect field of exponential characteristic $p$. If $\E \in \SH(k)$ is a $\zpi$-local object with a structure of traces, then the slices $s_q\E$ have a canonical structure of traces.
\end{theoNo}

We will now give an outline of this thesis.

In Chapter~\ref{chap:cycles} we present a part of the Suslin-Voevodsky theory of relative cycles \cite{SV}. Instead of defining the presheaves of relative cycles $c_{equi}(X / S, 0)$ gradually via the presheaves $Cycl(X/S, 0)$ as is the usual treatment, we present a definition of them via a universal property. We then show that they exist using a reworking of the usual construction that hopefully is more accessible for a novice to the theory.

One thing worth mentioning is that there is a small error in \cite{SV} that we correct. In particular, \cite[Corollary 3.2.4]{SV} is not true if $S$ is not reduced at its generic points. As a consequence, we lose the claim made in \cite[Corollary 3.3.11]{SV} that $cycl: \ZZ PropHilb(X/S, 0) \to c_{equi}(X/S, 0)$ is a natural transformation. To see that the naturality breaks, it suffices to consider the morphism $S_{red} \to S$ for any $S$ which is not reduced at its generic points, and $Z = X = S$. Note that if we restrict to the category of reduced noetherian schemes there is no problem, and so this does not affect correspondences between normal schemes at all.

We end this chapter with a similar exposition of the category $Cor(S)$ of correspondences, but cite the original article for the hard technical work.

In Chapter~\ref{chap:comparison} we present our definition of the $\ldh$-topology and compare it to the cdh topology. The main technical result of this chapter is Theorem~\ref{theo:cdhlprimeComparison} which has already been mentioned. We outline briefly the steps involved in getting there. The comparison of cohomologies is fairly straight-forward if we are working with presheaves with transfers. We show that the cohomologies of cdh/$\ldh$ sheaves with transfers can be calculated as Ext's in the category of presheaves with transfers, and then the result follows immediately from the easy facts that every presheaf with transfers is a presheaf with traces (Lemma~\ref{lemm:transfersImpliesTraces}), and every presheaf of $\zll$ modules with traces is a sheaf for the topology generated by finite flat surjective morphims of degree prime to $\ell$ (Lemma~\ref{lemm:fpslAcyclic}). Hence, the categories of cdh sheaves of $\zll$ modules with transfers and $\ldh$ sheaves of $\zll$ modules with transfers are equivalent (in fact, they are equal).

To get to transfers we show that under certain conditions a structure of traces on a presheaf $F$ induces a structure of traces on the cdh associated sheaf (Proposition~\ref{prop:tracesForcdh}), and a nice enough cdh sheaf with traces has a canonical structure of transfers (Theorem~\ref{theo:cdhTriImpliesTransfers}). The latter is straight-forward using Raynaud-Gruson's platification theorem (Theorem~\ref{theo:platification}) to convert every correspondence into a sum of compositions of morphisms of schemes  and ``transposes'' of finite flat morphisms.

Pushing the structure of traces through the cdh sheafification is harder. For this we introduce the notion of a Gersten presheaf (Definition~\ref{defi:gersten}) and a topology that we call the completely decomposed discrete topology or cdd topology (Definition~\ref{defi:cd}). A Gersten presheaf is a presheaf which satisfies some analogue of Gertsen's sequence for algebraic $K$-theory. The most important property of the cdd topology is that the cdd associated sheaf $F_{cdd}$ of a presheaf $F$ satisfies $F_{cdd}(X) = \prod_{x \in X} F(x)$ where the product is over the points of $X$ of every codimension. The Gersten exact sequence implies that if $F$ is a presheaf of $\zll$ modules then we have a sequence of monomorphims $F \to F_{cdh} \to F_{\ldh} \to F_{cdd}$. We show that a structure of traces on $F$ passes to a structure of traces on $F_{cdd}$ (Theorem~\ref{theo:tracesForcd}), give a criterion for a section to be in the image of $F_{cdh} \to F_{cdd}$, and show that the trace morphisms of $F_{cdd}$ preserve this criterion. Hence, the structure of traces on $F_{cdd}$ induces a structure of traces on $F_{cdh}$ (Theorem~\ref{theo:tracesForcd}). For an explanation of why the cdd topology arises quite naturally for us see Remark~\ref{rema:whyCdd}.

In Chapter~\ref{chap:slices} we shift focus to the Morel-Voevodsky stable homotopy category. The idea is that we can define trace morphisms in algebraic $K$-theory quite easily, and motivic cohomology is a graded piece of algebraic $K$-theory, so we might be able to descend the algebraic $K$-theory trace morphisms to motivic cohomology. In the context of $\SH$, this involves a study of the slice filtration. We begin the chapter by translating some work of Pelaez on the functoriality of the slice filtration into Ayoub's language of a stable homotopy 2-functors (cf. Theorem~\ref{theo:ayo} and Remark~\ref{rema:stableHomotopy2functor}), which makes it easier to study the functoriality of the slice filtration. The main theorem of Pelaez that we use is Theorem~\ref{theo:criteria} which gives criteria for a triangulated functor to preserve the slices of an object. We show that the functors we are interested in satisfy his criteria (Theorem~\ref{theo:pullbackAll}, Proposition~\ref{prop:ffslices}).

We then define what it means for an object $\E \in \SH(S)$ to have a structure of traces (Definition~\ref{defi:sectionWithTraces}), and use the material we have developed to show that a structure of traces on an object induces a canonical structure of traces on its slices. This is Proposition~\ref{prop:tracesOnSlices} stated above. We also mention that a structure of traces on an object induces a structure of traces on its homotopy presheaves (Lemma~\ref{lemm:tracesImpliesTraces}), that a structure of traces on a ring spectrum induces a structure of traces on each of its free modules (Proposition~\ref{prop:tracesOnProducts}), and that structures of traces are preserved morphisms of 2-functors which commute with the right adjoints (Lemma~\ref{lemm:localisationOfTraces}).

In Chapter~\ref{chap:applications} we apply all the previous material. We begin by showing that the object $\KH$ representing algebraic $K$-theory in $\SH(k)$ has a structure of traces (Proposition~\ref{prop:ktheoryTraces}), and that the object $\HZ$ representing motivic cohomology has what we have called a weak structure of smooth traces (Definition~\ref{defi:wsst}, Proposition~\ref{prop:weakTracesOnHZ}). We show that cdh descent plus a structure of traces implies $\ldh$ descent (Theorem~\ref{theo:ldescentInSH}). We have already mentioned that every object in $\SH$ satisfies cdh descent (\cite[3.7]{Cis13}), and so end up with the result that objects of the form $\HZl \wedge M$ satisfy $\ldh$ descent. It follows that every $\HZl$ module satisfies $\ldh$ descent. In particular, every smooth $\ldh$ hypercover in the category of $\HZl$ modules is isomorphic to the scheme that it covers. We apply this in the way outlined above to obtain Theorem~\ref{theo:MainResultRos}. We recall some parts of \cite{Sus00} and show how our Theorem~\ref{theo:MainResultRos} implies Theorem~\ref{theo:suslin}.

Finally, we discuss the conjecture of Weibel mentioned above about vanishing of algebraic $K$-theory in sufficiently low degrees. 

\section{Notation and conventions}

All schemes will be separated unless otherwise stated. Associated to a scheme $S$ we consider the following categories.
\begin{enumerate}
 \item[] $Sch(S)$ the category of schemes of finite type over $S$.
 \item[] $Sm(S)$ the full subcategory of $Sch(S)$ whose objects are smooth $S$-schemes.
 \item[] $Reg(S)$ the full subcategory of $Sch(S)$ whose objects are regular $S$-schemes.
 \item[] $QProj(S)$ the full subcategory of $Sch(S)$ which are quasi-projective.
 \item[] $EssSch(S)$ the category of $S$-schemes that are inverse limits of left filtering systems in $Sch(S)$ in which the transition morphisms are all affine open immersions.
 \item[] $EssQProj(S)$ the category of $S$-schemes that are inverse limits of left filtering systems in $QProj(S)$ in which the transition morphisms are all affine open immersions.
\end{enumerate}

\chapter{Relative cycles} \label{chap:cycles}


\section{Introduction}

\lettrine[lines=2, findent=3pt, nindent=0pt]{T}{he} goal of this chapter is to give a construction of the presheaves of relative cycles $c_{equi}(X/S, 0)$ of Suslin-Voevodsky \cite{SV}. The construction of $c_{equi}(X/S, r)$, $c(X/S, r)$, $z_{equi}(X/S, r)$, and $z(X/S, r)$, is analogous. The culmination of the first four sections is Theorem~\ref{theo:defcequi} which suggests a definition of the presheaf $c_{equi}(X/S, 0)$ as the unique presheaf $F$ satisfying:
\begin{enumerate}
 \item[(Gen)] $F(T)$ is a subgroup of the free abelian group generated by the points $z$ of $X \times_S T$ such that $\overline{\{z\}} \to T$ is finite and dominates an irreducible component of $T$, where $\overline{\{z\}}$ is the closure of $z$ in $X \times_S T$.
 \item[(Red)] If $i: T_{red} \to T$ is the canonical inclusion, then $F(i)$ is the morphism induced by the canonical identification of the points of $X \times_S T$ with the points of $X \times_S T_{red}$.
 \item[(Pla)] If $\sum n_iz_i \in F(T)$, $k$ is a field, and $\iota: Spec(k) \to T$ is a $k$-point of $T$ such that the image of $i$ is in the flat locus of $\amalg \overline{\{z_i\}} \to T$, then 
\[ F(\iota)\Z = \sum n_i m_{ij} w_{ij} \]
where the $w_{ij}$ are the (generic) points of $k \times_T \overline{\{z_i\}}$ and $m_{ij} = \length \OO_{k \times_T \overline{\{z_i\}}, w_{ij}}$.
 \item[(Uni)] Any other presheaf possessing the above three properties is a subpresheaf of $F$.
\end{enumerate}

A definition of this form was clearly known to Suslin and Voevodsky (see for example the beginning of Section 2 of \cite{FV}) and the reader familiar with their theory will fail to be surprised by it.

Such a definition has the advantage that it takes less than half a page to write down and the pullbacks for any morphism $f: T \to S$ for the presheaves $c_{equi}(X/S, 0)$ can be calculated using these axioms and the platification theorem (reproduced as Theorem~\ref{theo:platification}). In the original article \cite{SV} the definition of $c_{equi}(X/S, 0)$ appears on page 36 (actually the 27th page of the article) and everything preceding it is more or less necessary to arrive at that definition. There is also a criterion for a formal sum to belong to the subgroup $c_{equi}(X/S, 0)$ which can be stated using morphisms calculated from the axioms above. Thus, if desiring to do so, a reader could potentially develop a working knowledge of these presheaves without having to wade through the construction that proves they exist. 

The idea behind the $c_{equi}(X/S, 0)$ is that these relative cycles should be finite sums $\Z = \sum n_i z_i$ of points $z$ of $X$ that lie over generic points of $S$, and such that $\overline{\{z\}} \to S$ is a finite morphism. The free abelian group generated by such points will play an important r{\^o}le and we denote it $c_{equi}^{nai}(X/S, 0)$. This is an adaptation of the notation $c_{equi}(X/S, 0)$ where ``equi'' refers to the requirement that the morphisms $\overline{\{z\}} \to S$ are equidimensional, and the $0$ to the fact that they are of relative dimension zero. We have added ``nai'' to indicate that these free abelian groups are what one might na{\"i}vely expect to be the groups of relative cycles. The problem is that with the definitions of pullbacks $c_{equi}^{nai}(X/S, 0) \to c_{equi}^{nai}(S' \times_S X/S', 0)$ that we want (associated to a morphism $f: S' \to S$) these groups don't form presheaves (see Example~\ref{exam:fonFailure}). The solution is to keep the pullbacks that we like for certain kinds of $\Z$ and $f$, and then jettison any cycles that don't respect the induced functoriality.

Our particular choices of pullbacks that appear in the above axioms determine all of the other pullbacks uniquely (this is the content of Proposition~\ref{prop:sameMorphisms}), and so the groups $c_{equi}(X/S, 0)$ are then defined as the largest collection of subgroups of the free abelian groups $c_{equi}^{nai}(X/S, 0)$ that forms a presheaf with these chosen pullbacks. For a precise definition, see Definition~\ref{defi:cequi}.

To prove that presheaves $c_{equi}(X/S, 0)$ satisfying the above axioms exist we construct them. Our exposition of this construction is very strongly influenced by the original article \cite{SV}, however we deviate mildly from the way in which Suslin-Voevodsky present the material. Our goals were expositional: to introduce as little notation as possible, and to try and avoid any definition whose motivation wasn't immediately obvious on a first reading. We also decided to avoid the use of fat points (\cite[Definition 3.1.1]{SV}) to see if this could be done, but the concept we replace them with -- good factorisations (Definition~\ref{defi:goodFactorisation}) -- is more or less equivalent (cf. \cite[Proposition 3.1.5]{SV}) and morally our proofs are the same as theirs.

We remark that there is a small fixable problem in \cite{SV} due to nilpotents. In particular, \cite[Corollary 3.2.4]{SV} is not true if $S$ is not reduced at its generic points. As a consequence, we lose the claim made in \cite[Corollary 3.3.11]{SV} that $cycl: \ZZ PropHilb(X/S, 0) \to c_{equi}(X/S, 0)$ is a natural transformation. To see that the naturality breaks, it suffices to consider the morphism $S_{red} \to S$ for any $S$ which is not reduced at its generic points, and $Z = X = S$. Note that if we restrict to the category of reduced noetherian schemes there is no problem, and so this does not affect correspondences between normal schemes at all.

We propose a way of fixing this by using a slightly different version of their $cycl$. If $f: X \to S$ is a morphism of finite type and $Z \to X$ a closed subscheme that is flat over $S$, they define the cycle associated to $Z$ as $\sum n_i z_i$ with $z_i$ the generic points of $Z$ and $n_i = \length \OO_{Z, z_i}$. We propose, however to take $n_i = \length \OO_{f(z_i) \times_S Z, z_i}$. This altered definition does not affect the presheaves $c_{equi}(X/S, 0)$ at all. The reader can check in \cite[Proposition 3.1.5]{SV} and \cite[Theorem 3.3.1]{SV} that $c_{equi}(X/S, 0) \to c_{equi}(S_{red} \times_S X/S_{red}, 0)$ is implicitly forced to be the isomorphism induced by the canonical identification of the points of $X$ with the points of $S_{red} \times_S X$. Hence the $c_{equi}(X/S, 0)$ are completely determined by their values on reduced schemes. It is easily checked that with out new choice of $cycl$, the morphisms $cycl: \ZZ PropHilb(X/S, 0) \to c_{equi}(X/S, 0)$ are natural transformations of presheaves (see Proposition~\ref{prop:cyclesOfFlatSubschemes}).

Lastly, we mention that Ivorra \cite{Ivo05} (published as \cite{Ivo11}) has a produced an extremely readable version of Suslin-Voevodsky's \cite{SV} from which we learn't a lot. It is unclear how he treats the problem of nilpotents we mention above as his version $[-]$ of Suslin-Voevodsky's $cycl$ mentioned above is not defined. His application is to regular schemes and so this poses no serious problem to him. There is also an extension of the theory in development by Cisinski-D{\'e}glise. A preliminary version appears in \cite{CD}. The idea is that cycles (i.e., a scheme equipped with a formal sum of its points) should be the objects of a category in their own right. This category is equipped with a relative product. The Suslin-Voevodsky pullback, as well as the Suslin-Voevodsky product of relative cycles, are recovered special cases of this relative product.

\emph{Index.} As a guide to the reader for what is to come, and as a reference, we collect here the notation we introduce. As we already mentioned we tried to keep this as minimal as possible, and wherever we could to use notation that already existed in the literature.

\begin{enumerate}
 \item[] $c_{equi}^{nai}(X/S, 0)$. Definition~\ref{defi:C}. This we introduce as we find it clearer than the $C(X/S, 0)$ in \cite{Ivo05}. Suslin-Voevodsky don't have a notation for these groups.
 \item[] $cycl_{X/S}$. Definition~\ref{defi:cycl}. This is a version of the $cycl_X$ from \cite{SV}, but our version is modified to adjust for nilpotents.
 \item[] $f^*_{nai}$. Definition~\ref{defi:naivePullback}. This is the na{\"i}ve pullback that we might expect. In many cases, it is indeed the correct pullback (i.e., it agrees with $f^\circledast$).
 \item[] $(\iota, p)$ Definition~\ref{defi:goodFactorisation}. This is our analogue of the Suslin-Voevodsky fat points.
 \item[] $(\iota, p)^*$. Definition~\ref{defi:goodFactorisation}. This is our analogue of the pullback along a fat point of Suslin-Voevodsky.
 \item[] $c_{equi}(X/S, 0)$. Definition~\ref{defi:cequi}. These are the subgroups of $c_{equi}^{nai}(X/S, 0)$ that behave well with respect to the pullback.
 \item[] $f^\circledast$. Definition~\ref{defi:cyclePullbacks}. The pullbacks  of the presheaves $c_{equi}(X/S, 0)$. This notation is from \cite{Ivo05} and replaces the clunky $cycl(f)$ of \cite{SV}, which incidentally is in conflict with their $cycl$ which is mentioned above.
\end{enumerate}

\section{First definitions}

In this section we define the free abelian groups $c_{equi}^{nai}(X/S, 0)$ which contain the groups $c_{equi}(X/S, 0)$. We define a na{\"i}ve pullback $f^*_{nai}$ for the groups $c_{equi}^{nai}(X/S, 0)$, give an example of why these pullbacks don't equip these groups with the structure of a presheaf, and prove some properties about them that we will need. We then give our version of the Suslin-Voevodsky fat points, which we call good factorisations. We define the pullback $(\iota, p)^*$ along a good factorisation and show that good factorisations always exist (up to field extension). In certain cases the pullbacks $f^*_{nai}$ and $(\iota, p)^*$ agree with the pullbacks $f^\circledast$ of $c_{equi}(X/S, 0)$ (see Lemma~\ref{lemm:agreesWithOldPullbacks} for a precise statement) and so these definitions can also be regarded as calculations.

We begin with the free abelian groups that will contain our relative cycle groups.

\begin{defi} \label{defi:C}
Suppose that $f: X \to S$ is a scheme of finite type over a noetherian base scheme $S$. We define $c_{equi}^{nai}(X/S, 0)$ to be the free abelian group generated by the points $z \in X$ such that $\overline{\{z\}} \to S$ is finite, and dominates an irreducible component of $S$. That is, $z$ is in one of the generic fibres of $X \to S$.
\end{defi}

We will most often come across elements of $c_{equi}^{nai}(X/S, 0)$ using the following definition. The notation $(-)^{(0)}$ indicates points of codimension zero as usual.

\begin{defi} \label{defi:cycl}
With the notation as in Definition~\ref{defi:C} suppose that $Z$ is a closed subscheme of $X$ which is finite over $S$. We define
\[ cycl_{X/S}(Z) = \sum_{\substack{z_i \in Z^{(0)} \textrm{ s.t.} \\ f(z_i) \in S^{(0)}}} n_i z_i \] 
where $n_i$ is the length of the local ring of the point $z_i$ in its fibre. That is, $n_i = \length \OO_{f(z_i) \times_S Z, z_i}$. We will sometimes omit the subscript and just write $cycl$ if the morphism $X \to S$ is clear from the context.
\end{defi}

\begin{rema}
This differs from the $cycl_X(Z)$ defined in \cite{SV} as their coefficients are $n_i = \length \OO_{Z, z_i}$. Our choice of definition for $cycl_{X/S}$ is a proposed fix for the problem mentioned in the introduction that they don't actually get a morphism of presheaves $\ZZ PropHilb(X/S, 0) \to c_{equi}(X/S, 0)$ over all non-reduced schemes.

We also added the hypothesis that the sum only counts those points that lie over generic points to assure that our cycle is in $c_{equi}^{nai}(X/S, 0)$, but this is just to avoid introducing another notation for the free abelian group generated by all the points of $X$.
\end{rema}

\begin{rema} \label{rema:cyclNoNilpotence}
From our definition it follows that, in the notation of the definition, we have $cycl_{X/S}(Z) = cycl_{S_{red} \times_S X/S_{red}}(S_{red} \times_S Z)$ via the canonical identification $c_{equi}^{nai}(X/S, 0) \cong c_{equi}^{nai}(S_{red} \times_S X/S_{red}, 0)$. This is not true of the Suslin-Voevodsky $cycl_X$.
\end{rema}

We define now the obvious pull-back morphism for the $c_{equi}^{nai}(X/S, 0)$. When we restrict to relative cycles, these will end up being the actual pullbacks $f^\circledast$ in certain settings so this definition can also be seen as an explicit calculation of certain examples of $f^\circledast\Z$. For a precise description of some cases when $f^\circledast\Z = f^*_{nai}\Z$ see Lemma~\ref{lemm:agreesWithOldPullbacks}.

\begin{defi} \label{defi:naivePullback}
Suppose that $f: T \to S$ is a morphism between noetherian schemes and $X \to S$ a morphism of finite type. We define a morphism
\[ f^*_{nai}: c_{equi}^{nai}(X/S, 0) \to c_{equi}^{nai}(X \times_S T/T, 0) \]
by
\[ f^*_{nai}(\sum n_iz_i) = \sum n_i cycl_{X \times_S T/T}(T \times_S \overline{\{z_i\}}). \]
More explicitly, we have $f^*_{nai}\Z = \sum n_i m_{ij} w_{ij}$ where $w_{ij}$ are the generic points of $\overline{\{z_i\}} \times_S T$ that lie over generic points of $T$, and $m_{ij}$ are the lengths of their local rings $\length \OO_{t_{ij} \times_S \overline{\{z_i\}}, w_{ij}}$ in their fibres over $T$ (the point $t_{ij}$ is the image of $w_{ij}$ in $T$).
\end{defi}

\begin{exam} \label{exam:fonFailure}
The pullback defined above does not equip the groups $c_{equi}^{nai}(X/S, 0)$ with a structure of presheaf. Consider $S = S_1 \cup S_2$ to be the union of two affine lines $S_1 \cong \AA^1, S_2 \cong \AA^1$ joined at a closed point $s = S_1 \cap S_2$. Let $\eta_1$ be the generic point of $S_1$ so we get an element $\eta \in c_{equi}^{nai}(S/S, 0)$. Consider the inclusion $\iota_2: s \to S_1 \amalg S_2$ of the point $s$ into $S_2$, and the canonical morphism $p: S_1 \amalg S_2 \to S$. Now we have $(p i_2)^*_{nai}\eta = s \in c_{equi}^{nai}(s/s, 0)$ but $i^*_{2, nai}p^*_{nai}\eta = 0 \in c_{equi}^{nai}(s/s, 0)$.
\end{exam}

\begin{rema}
The example above suggests that the problem occurs when we have multiple choices of branches, and this is indeed the case. Notably, if $S$ is regular, then we have equality $c_{equi}(X/S, 0) = c_{equi}^{nai}(X/S, 0)$ (see \cite[Corollary 3.4.6]{SV}). Our definition of a good factorisation (and the Suslin-Voevodsky idea of a fat point) can be thought of as a choice of branch.
\end{rema}

\begin{exam}
If $f: S' \to S$ is a birational morphism\footnote{We recall that a \emph{birational} morphism $f: T \to S$ is a morphism which sends every generic point of $T$ to a unique generic point of $S$, every generic point of $S$ is in the image, and the field extensions induced on generic points are all trivial. In particular $S_{red} \to S$ is birational.} then for some $\Z = \sum n_i z_i \in c_{equi}^{nai}(X/S, 0)$, the na{\"i}ve pull-back is just $f^*_{nai}\Z = \sum n_i z_i'$ where $z_i'$ is $z_i$ seen as a point of $S' \times_S X$ via the canonical identification of the generic fibres of $X \to S$ and $S' \times_S X \to S'$.
\end{exam}

Before the following easy lemma, we recall the following definition from Cisinski-D{\'e}glise \cite{CD}. An earlier version of \cite{ILO} calls these morphisms ``horizontal'', and the current version calls them ``maximally dominant''.

\begin{defi} \label{defi:psuedodominant}
A morphism of schemes $f: Y \to X$ is said to be \emph{pseudo-dominant} if every generic point of $Y$ is sent to a generic point of $X$.
\end{defi}


\begin{lemm} \label{lemm:fonBir}
Suppose that $S_3 \stackrel{g}{\to} S_2 \stackrel{f}{\to} S_1$ are morphisms between noetherian schemes, suppose $g$ is pseudo-dominant, and let $X \to S$ be a morphism of finite type. Then
\[ g^*_{nai}f^*_{nai} = (fg)^*_{nai}. \]
\end{lemm}

\begin{rema}
Example~\ref{exam:fonFailure} shows that this is not true if we remove the hypothesis that $g$ is pseudo-dominant, even if we add the hypothesis that $f$ is birational.
\end{rema}

\begin{proof}
First note that if $\iota$ is the inclusion of the generic points of a scheme, then $\iota^*_{nai}$ is injective. Now due to the commutative square
\[ \xymatrix{
\amalg_{s' \in S_{3}^{(0)}} s' \ar[r] \ar[d] & S_3 \ar[d] \\
\amalg_{s \in S_{2}^{(0)}} s \ar[r] & S_2
} \]
and this injectivity, it suffices to consider the two cases (i) when $S_3$ and $S_2$ are reduced of dimension zero, and (ii) when $g$ is the inclusion of a subset of the generic points of $S_2$. 

Consider the case (i). We can assume that $S_3$ and $S_2$ are actually integral of dimension zero. Let $z \in X$ be a point over a generic point of $S_1$ such that $\overline{\{z\}} \to S_1$ is finite. Suppose that $w_i$ are the generic points of $S_2 \times_{S_1} \overline{\{z\}}$ and $v_{ij}$ are the generic points of $S_3 \times_{S_2} \overline{\{w_i\}}$ and set
\[ \begin{split}
 m_i &= \length \OO_{S_2 \times_{S_1} \overline{\{z\}}, w_i} \\
 n_{ij} &= \length \OO_{S_3 \times_{S_2} \overline{\{w_i\}}, v_{ij}} \\
 \ell_{ij} &= \length \OO_{S_3 \times_{S_1} \overline{\{z\}}, v_{ij}}
\end{split} \]
so that we have
\[ \begin{split}
 f^*_{nai}z &= \sum m_i w_i \\
 g^*_{nai}f^*_{nai}z &= \sum m_i n_{ij} v_{ij} \\
 (fg)^*_{nai}z &= \sum \ell_{ij}v_{ij}
\end{split} \]
Hence, it suffices to show that we have $m_i n_{ij} = \ell_{ij}$. This is precisely what Lemma~\ref{lemm:appendix3} says.

The case (ii) follows straight from the definition of $(-)^*_{nai}$.
\end{proof}

\begin{lemm} \label{lemm:naiInjective}
Suppose that $f: T \to S$ is a morphism between noetherian schemes and $X \to S$ a morphism of finite type. If $f$ is dominant, then $f^*_{nai}$ is injective.
\end{lemm}

\begin{proof}
Suppose that $s_i$ are the generic points of $S$ and for each $i$ choose a generic point $t_i$ of $T$ which is over it. We find the commutative square
\[ \xymatrix{
\amalg t_i \ar[r] \ar[d]_q & T \ar[d] \\
\amalg s_i \ar[r]_p & S
} \]
and by the functoriality given in Lemma~\ref{lemm:fonBir} it suffices to show that $p^*_{nai}$ and $q^*_{nai}$ are injective. In both these cases, the injectivity is clear from the definitions.
\end{proof}

The following theorem is a cut down version of \cite[3.2.2]{SV} with more or less the same proof. There is a small mistake in the proof of \cite[3.2.2]{SV}. Using their notation, in their final case they claim $\eta'$ is the only point over $\tau'$ which is not always true -- consider the case when $\tau'$ and $\eta$ are the same non-trivial finite separable field extension of $\tau$. We don't reproduce their error.

\begin{theo} \label{theo:naiFon}
Suppose that $T \to S$ is a morphism of noetherian schemes and $X \to S$ is a morphism of finite type. Let $\sum n_i Z_i$ be a finite sum of closed subschemes of $X$ that are finite and flat over $S$. Then for $\sum n_i cycl_{T \times_S X / T}(T \times_S Z_i)$ to be zero in $c_{equi}^{nai}(T \times_S X / T, 0)$ it is sufficient that $\sum n_i cycl_{X/S}(Z_i)$ is zero in $c_{equi}^{nai}(X/S, 0)$.

\end{theo}

\begin{proof}
 \emph{Reduction to $T$ integral dimension zero, and $S$ local reduced.} When $T \to S$ is birational, the generic fibres of $X \to S$ and $T \times_S X \to T$ are canonical isomorphic. Via this identification, we have the equality $\sum n_i cycl_{T \times_S X / T}(T \times_S Z_i) = \sum n_i cycl_{X/S}(Z_i)$. Therefore we have the stronger statement that $\sum n_i cycl_{X/S}(Z_i)$ is zero if and only if $\sum n_i cycl_{T \times_S X / T}(T \times_S Z_i)$ is zero. Hence, we can replace $S$ by $S_{red}$, and we can replace $T$ by the disjoint union of its generic points. To show $\sum n_i cycl_{T \times_S X / T}(T \times_S Z_i)$ is zero it is enough to consider each generic point of $T$ separately. So we assume that $T$ is an integral scheme of dimension zero. Without affecting any of the multiplicities, we can replace $S$ with any subscheme that contains the generic points of $S$, and the image $s$ of $T$. For example, the disjoint union of the localisation of $S$ at $s$, and any generic points not contained in this localisation. The generic points not involved in the localisation of $S$ at $s$ do not affect $\sum n_i cycl_{T \times_S X / T}(T \times_S Z_i)$ in any way, and so we can forget them. That is, we assume $S$ is a reduced local scheme and the image of $T$ its closed point.

 \emph{The case where $S$ and $T$ are both integral dimension zero.} Without losing any information we can assume that $X = \cup Z_i$. Moreover, since it suffices to consider each connected component of $X$ one at a time, we can assume $X$ has a unique point $x$. If $y_j$ are the points of $T \times_S X$ then Lemma~\ref{lemm:appendix3} says that 
\begin{equation} \label{equa:hilbertPullback1}
\length \OO_{Z_i, x} \length \OO_{T \times_S x, y_j} = \length \OO_{T \times_S Z_i, y_j}
\end{equation}
for each $i, j$. By definition we have
\begin{equation} \label{equa:hilbertPullback2}
\sum_i n_i cycl_{X/S}(Z_i) = \left ( \sum_i n_i \length \OO_{Z_i, x} \right ) x
\end{equation}
and
\begin{equation} \label{equa:hilbertPullback3}
\sum_{i} n_i cycl_{T \times_S X / T}(T \times_S Z_i) = \sum_j \left ( \sum_i n_i \length \OO_{T \times_S Z_i, y_j} \right ) y_j.
\end{equation}
Multiplying Equation~(\ref{equa:hilbertPullback2}) by $\length \OO_{T \times_S x, y_j}$, using the substitution given by Equation~(\ref{equa:hilbertPullback1}), and comparing it with Equation~(\ref{equa:hilbertPullback3}), we see that in this case we actually have the stronger $\sum_i n_i cycl_{X/S}(Z_i)$ is zero if and only if the sum \mbox{$\sum_{i} n_i cycl_{T \times_S X / T}(T \times_S Z_i)$} is zero.

Notice that the reduction above, together with the dimension zero case, answers our question when $T \to S$ is pseudo-dominant. In particular, when $T \to S$ is flat.

\emph{The case when $S$ is local henselian and $T$ is the closed point of $S$.} Now we return to the case when $S$ was a reduced local ring, and suppose that $T$ is the inclusion of the closed point $s$ of $S$. Since we know the theorem is true for flat morphisms, we can replace $S$ by its henselisation. In this case since $X = \cup Z_i$ is finite over $S$, the scheme $X$ is a disjoint union of local schemes. It suffices to consider each connected component by itself, and so we can assume that $X$ is finite and local over $S$. This means that there is a unique $x \in X$ over the closed point $s \in S$, and that the $Z_i \to S$ are of constant degree $d_i$. In this case, we must show that $\sum n_i \length \OO_{s \times_S Z_i, x} = 0$. Since $\length \OO_{s \times_S Z_i, x} = d_i \cdot [k(x):k(s)]$ (Lemma~\ref{lemm:appendix1}), it is enough to show that $\sum n_i d_i = 0$. 

Consider a generic point $\eta \in S$ and the generic points of $X$ that lie over it. By Lemma~\ref{lemm:appendix1} and the fact that $S$ is reduced we know that
\[ d_i = \sum_ {\xi \in Z_i^{(0)}} [k(\xi):k(\eta)] \length \OO_{Z_i,\xi} \]
and so to show $\sum n_i d_i = 0$ it is enough to show that 
\[ \sum_i n_i \sum_ {\xi \in Z_i^{(0)}} [k(\xi):k(\eta)] \length \OO_{Z_i,\xi} = 0. \]
Interchanging the summands, we rewrite this sum as
\[ \begin{split}
\sum_i n_i \sum_ {\xi \in Z_i^{(0)}} [k(\xi):k(\eta)] \length \OO_{Z_i,\xi} &= \sum_{\xi \in X^{(0)}} \left ( \sum_{\substack{Z_i \textrm{ s.t. } \\ \xi \in Z_i}} n_i [k(\xi):k(\eta)] \length \OO_{Z_i,\xi} \right ) \\
&= \sum_{\xi \in X^{(0)}} [k(\xi):k(\eta)]  \sum_{\substack{Z_i \textrm{ s.t. } \\ \xi \in Z_i}} n_i \length \OO_{Z_i,\xi} 
\end{split} \]
and we see that it is enough to show that for each $\xi$ we have $\sum_{\substack{Z_i \textrm{ s.t. } \\ \xi \in Z_i}} \length \OO_{Z_i, \xi} = 0$. But since $S$ is reduced, this is equivalent to the statement $\sum n_i cycl_{X/S}(Z_i) = 0$. Hence, the result is true in this case.

\emph{The case $T$ integral dimension zero, and $S$ local reduced.} We have seen that the theorem holds when $T \to S$ is a flat morphism so we can replace $S$ by its henselisation at the closed point. Now we factor the morphism as $T \to s \to S$ where $s$ is the closed point of $S$ and we have already considered these two cases.
\end{proof}

\begin{coro} \label{coro:naiFon}
Suppose that $S_3 \stackrel{g}{\to} S_2 \stackrel{f}{\to} S_1$ is a pair of composable morphisms of noetherian schemes, $X \to S_1$ is a morphism of finite type, and $\Z = \sum n_i z_i \in c_{equi}^{nai}(X/S_1, 0)$. Let $\W = \sum m_j w_j = f^*_{nai}\Z$. We suppose that the image of the generic points of $S_3$ (resp. $S_2$) is in the flat locus of $\amalg \overline{\{w_j\}} \to S_2$ (resp. $\amalg \overline{\{z_i\}} \to S_1$). Then
\[ g^*_{nai}f^*_{nai}\Z = (fg)^*_{nai}\Z. \]
\end{coro}

\begin{proof}
We can assume that $\Z = z$ consists of a single point with coefficient one. Since we are concerned only with phenomena that occur over generic points, we can replace $S_2$ (resp. $S_3$) with any open subset that contains the image of $g$ (resp. $f$). Hence, we can assume that $\amalg \overline{\{w_j\}} \to S_2$ (resp. $\overline{\{z\}} \to S_1$) is flat. We must show that $cycl(S_3 \times_{S_1} \overline{\{z\}}) = \sum m_j cycl(S_3 \times_{S_2} \overline{\{w_j\}})$. By Theorem~\ref{theo:naiFon} this will follow if $cycl(S_2 \times_{S_1} \overline{\{z\}}) = \sum m_j cycl(\overline{\{w_j\}})$. But this was the definition of the $m_j, w_j$.
\end{proof}

Finally we introduce a pullback that is closely related to the pullback along a fat point discussed in \cite{SV}. We will see later on that $(\iota, p)^* = (p \iota)^\circledast$ (see Lemma~\ref{lemm:agreesWithOldPullbacks}) so again, this definition can be considered as a calculation.

\begin{defi} \label{defi:goodFactorisation}
Suppose that $S$ is a noetherian scheme and $X \to S$ is a morphism if finite type. Suppose that $\Z = \sum n_i z_i \in c_{equi}^{nai}(X /S, 0)$. Let $Spec(k) \to S$ be a $k$ point of $S$ with $k$ a field. A \emph{good factorisation} of $Spec(k) \to S$ with respect to $\Z$ is a factorisation of the form 
\[ Spec(k) \stackrel{\iota}{\to} S' \stackrel{p}{\to} S \]
such that
\begin{enumerate}
 \item $p$ is proper and birational, and
 \item considering the $z_i$ as points of $S' \times_S X$ via the canonical identification of the generic fibres of $S' \times_S X \to S'$ and $X \to S$, the morphisms $\overline{\{z_i\}} \to S'$ are flat. 
\end{enumerate}
We define the pullback of $\Z$ along such a good factorisation as
\[ (\iota, p)^*\Z = \iota^*_{nai}p^*_{nai}\Z. \]
\end{defi}

We will construct good factorisations using the following theorem. We use the statement from \cite[Theorem 2.2.2]{SV}.

\begin{theo}[{\cite{RG71}}] \label{theo:platification}
Let $p: X \to S$ be a morphism of noetherian schemes and $U$ an open subscheme in $S$ such that $p$ is flat over $U$. Then there exists a closed subscheme $Z$ in $S$ such that $U \cap Z = \varnothing$ and the proper transform of $X$ with respect to the blow-up $Bl_Z S \to S$ with centre in $Z$ is flat over $S$.
\end{theo}

\begin{lemm} \label{lemm:existanceGoodFactor}
Suppose that $S$ is a noetherian scheme and $X \to S$ is a morphism if finite type. Suppose that $\Z = \sum n_i z_i \in c_{equi}^{nai}(X /S, 0)$. Let $Spec(k) \to S$ be a $k$ point of $S$. 

Then there exists a finite extension $L$ of $k$ such that the induced $L$ point $Spec(L) \to S$ has a good factorisation with respect to $\Z$.
\end{lemm}

\begin{proof}
The platification theorem (Theorem~\ref{theo:platification}) gives the existence of a blow-up $S' \to S_{red}$ of $S_{red}$ such that the strict transform of the morphism $\amalg \overline{\{z_i\}} \to S_{red}$ is flat. The composition $S' \to S$ is proper and birational and satisfies the necessary flatness condition for the $z_i$. Since $S' \to S$ is a morphism of finite type, for every point $s \in S$ there exists a point $s' \to S'$ such that $[k(s'):k(s)]$ is finite. Hence, there exists a finite extension $L$ of $k$ such that the induced $L$ point $Spec(L) \to S$ factors through $S'$, i.e., we have found a good factorisation.
\end{proof}

\section{Presheaves of relative cycles}

We now come to our precise description of the properties we wish our presheaves $c_{equi}(X/S, 0)$ to have. There are various other choices that give the same presheaves but we have chosen these.

\begin{defi} \label{defi:presheafOfRelativeCycles}
Suppose that $S$ is a noetherian scheme, $X \to S$ a morphism of finite type and $F$ a presheaf on the category of noetherian schemes over $S$. We will say that $F$ is a \emph{presheaf of relative cycles} if the following conditions are satisfied:
\begin{enumerate}
 \item[(Gen)] $F(T)$ is a subgroup of $c_{equi}^{nai}(X \times_S T / T, 0)$.
 \item[(Red)] If $\Z \in F(T)$ and if $i: T_{red} \to T$ is the canonical inclusion then
\[ F(i)\Z = i^*_{nai}\Z. \]
 \item[(Pla)] If $\sum n_iz_i \in F(T)$ and $\iota: Spec(k) \to T$ is a $k$-point of $T$ (with $k$ a field) such that the image of $i$ is in the flat locus of $\amalg \overline{\{z_i\}} \to T$, then 
\[ F(\iota)\Z = i^*_{nai}\Z. \]
\end{enumerate}
\end{defi}

The following lemma contains properties that we will use shortly.

\begin{lemm} \label{lemm:domBlow}
Suppose that $F$ is a presheaf of relative cycles and $f: T_2 \to T_1$ a morphism of noetherian $S$ schemes.
\begin{enumerate}
 \item If $f$ is dominant then $F(f)$ is injective, and
 \item if $f$ is birational then $F(f) =f^*_{nai}$.
\end{enumerate}
\end{lemm}

\begin{proof}
For the first statement, it suffices to consider the cases (i) when $f: T_2 \to T_1$ is the inclusion of the generic points and (ii) when $f: Spec(L) \to Spec(k)$ is a field extension. In the first case, $f$ factors through $(T_1)_{red}$, and so the result follows from (Red) and (Pla). The second follows from (Pla) and Lemma~\ref{lemm:naiInjective}.

Now suppose that $f$ is birational. We have a commutative square
\[ \xymatrix{
\amalg \tau_i \ar[r] \ar[d] & (T_1)_{red} \ar[d] \\
T_2 \ar[r] & T_1
} \]
and so the result follows from the case when $f$ is dominant, (Pla), and (Red).
\end{proof}

The following proposition shows that our axioms completely determine the pullback morphisms. It follows that the class of presheaves of relative cycles (associated to the same $X / S$) is partially ordered by inclusion. In particular, it makes sense to speak of a potential maximal element.

\begin{prop} \label{prop:sameMorphisms}
Suppose that $F_1$ and $F_2$ are two presheaves of relative cycles (associated to the same $X / S$), suppose that $f: T_2 \to T_1$ is a morphism between noetherian $S$ schemes. Then for any formal sum $\Z \in F_1(T_1) \cap F_2(T_1)$ that is in both presheaves, we have $F_1(f)\Z = F_2(f)\Z$.
\end{prop}

\begin{proof}
The morphism $f$ induces a morphism $f_{red}: (T_2)_{red} \to (T_1)_{red}$ and so due to the axiom (Red) it suffices to consider the case when $T_1$ and $T_2$ are reduced. Let $\Z = \sum n_i z_i$. Since $T_1$ is reduced, by the platification theorem (Theorem~\ref{theo:platification}) there exists a blow-up of $T_1$ with nowhere dense centre such that the proper transform of $\amalg \overline{\{z_i\}} \to T_1$ is flat. We construct the following commutative diagram
\[ \xymatrix{
\amalg Spec(k_i) \ar[rr] \ar[d] && \tilde{T}_1 \ar[d] \\
\amalg \tau_i \ar[r] & T_2 \ar[r] & T_1
} \]
where the $\tau_i$ are the generic points of $T_2$ and $k_i / k(\tau_i)$ is field extension such that $Spec(k_i) \to T_1$ lifts through the blow-up $\tilde{T}_1 \to T_1$. Since $F_j(T_2) \to F_j(\amalg \tau_i) \to F_j(\amalg Spec(k_i))$ is injective for $j = 1, 2$ it suffices to show that $F_1$ and $F_2$ agree $\amalg Spec(k_i) \to \tilde{T}_1$ and $\tilde{T}_1 \to T_1$. The latter is given to us by Lemma~\ref{lemm:domBlow} and the former is (Pla) since the $\overline{\{z_i\}} \to \tilde{T}_1$ are flat.
\end{proof}

Lastly, we show that any presheaf of relative cycles satisfies two important properties that we will use to define the pullbacks $f^\circledast$.

\begin{prop} \label{prop:cIsprc}
Suppose that $F$ is a presheaf of relative cycles, $T$ is a noetherian $S$ scheme, $\Z = \sum n_i z_i \in F(T)$ is a section. Then we have the following properties.
\begin{enumerate}
 \item For any field $k$, any $k$ point $Spec(k) \to T$, and any pair of good factorisations $(\iota_1, p_1)$, $(\iota_2, p_2)$ we have
\[ (\iota_1, p_1)^*\Z = (\iota_2, p_2)^*\Z \]
 \item For any field $k$, any $k$-point $Spec(k) \to S$ with image $s \in S$ and induced morphism $q: Spec(k) \to s$, and any good factorisation $(\iota, p)$ with respect to $\Z$,
\[ \xymatrix{
Spec(k) \ar[d]_q \ar[r]^\iota & S' \ar[d]^p \\
s \ar[r] & S
} \]
there exists a unique $\Z' \in F(s)$ such that
\[ q^*_{nai}\Z' = (\iota, p)^*\Z. \]
\end{enumerate}
\end{prop}

\begin{proof}
This is a direct consequence of functoriality, the axioms (Red), (Pla), (Gen), and Lemma~\ref{lemm:domBlow}.
\end{proof}

\section{The groups $c_{equi}(X/S, 0)$ and the pull-backs $f^\circledast$}

We make the following definition with two motivations. The first is Proposition~\ref{prop:cIsprc} : if we wish the axioms to hold, then clearly we need these properties. The second is our choice of definition of the pullbacks $f^\circledast$ (see Definition~\ref{defi:cyclePullbacks} and Theorem~\ref{theo:pullback}). These two properties are what we will use to define the pullbacks.

\begin{defi} \label{defi:cequi}
Suppose that $S$ is a noetherian scheme and $X \to S$ a morphism of finite type. We define $c_{equi}(X/S, 0)$ to be the subgroup of $c_{equi}^{nai}(X/S, 0)$ of formal sums $\Z = \sum n_iz_i$ which have the properties of Proposition~\ref{prop:cIsprc}. That is:
\begin{enumerate}
 \item For every field $k$, every $k$-point $Spec(k) \to S$ of $S$, and every pair of good factorisations $(\iota_1, p_1), (\iota_2, p_2)$ with respect to $\Z$ we have
\[ (\iota_1, p_1)^*\Z = (\iota_2, p_2)^*\Z. \]
 \item For any field $k$, any $k$-point $Spec(k) \to S$ with image $s \in S$ and induced morphism $q: Spec(k) \to s$, and any good factorisation $(\iota, p)$ with respect to $\Z$,
\[ \xymatrix{
Spec(k) \ar[d]_q \ar[r]^\iota & S' \ar[d]^p \\
s \ar[r] & S
} \]
 there exists a unique $\W \in c_{equi}^{nai}(s \times_S X / s, 0)$ such that
\[ q^*_{nai}\W = (\iota, p)^*\Z. \]
\end{enumerate}
\end{defi}

The following proposition shows that the second condition can actually be made much weaker.

\begin{prop} \label{prop:alternativeConditionsForcequi}
Suppose that $S$ is a noetherian scheme and $X \to S$ a morphism of finite type, and $\Z = \sum n_iz_i \in c_{equi}^{nai}(X/S, 0)$. The following conditions are equivalent.
\begin{enumerate}
 \item Condition (1) from Definition~\ref{defi:cequi}.
 \item The same condition, except that for each $s$ we only need to consider an algebraic closure of $k(s)$, and we only need to find one $p: S' \to S$.
\end{enumerate}
More explicitly:
\begin{enumerate}
 \item For every field $k$, every $k$-point $Spec(k) \to S$ of $S$, and every pair of good factorisations $(\iota_1, p_1), (\iota_2, p_2)$ with respect to $\Z$ we have
\[ (\iota_1, p_1)^*\Z = (\iota_2, p_2)^*\Z. \]
 \item For every point $s \in S$, and algebraic closure $\Omega$ of $k(s)$ with induced $\Omega$ point $Spec(\Omega) \to S$, there exists a good factorisation $(Spec(\Omega) \stackrel{\iota}{\to} S', S' \stackrel{p}{\to} S)$ of $Spec(\Omega) \to S$ with respect to $\Z$ such for any other factorisation $\iota': Spec(\Omega) \to S'$ we have
\[ (\iota, p)^*\Z = (\iota', p)^*\Z. \]
\end{enumerate}
The commutative diagram for the second conditions is:
\[ \xymatrix{
Spec(\Omega) \ar@<0.5ex>[r]^\iota \ar@<-0.5ex>[r]_{\iota'} \ar[d] & S' \ar[d]^p \\
s \ar[r] & S
} \]
\end{prop}

\begin{proof}
Clearly the first condition implies the second (c.f Lemma~\ref{lemm:existanceGoodFactor}). So suppose that the second condition is satisfied.

We wish to show that the first condition is true. Suppose that $k'$ is a field, $Spec(k') \to S$ is a $k'$ point with target $s$, and $(\phi_1, S'_1 \stackrel{q_1}{\to} S), (\phi_2, S_2' \stackrel{q_2}{\to} S)$ is a pair of good factorisations of $Spec(k') \to S$ with respect to $\Z$. By the definition of pullback with respect to a good factorisation, if $\psi: Spec(L) \to Spec(k')$ is any field extension, then for $j = 1, 2$ we have
\[ \psi^*_{nai}(\phi_j, q_j)^*\Z = \psi^*_{nai}\phi^*_{j, nai} q^*_{nai} \Z \stackrel{\ref{coro:naiFon}}{=} (\phi_j\psi)^*_{nai} q^*_{nai} \Z = (\phi_j\psi, q)^*\Z \]
so since $\psi^*_{nai}$ is injective (Lemma~\ref{lemm:naiInjective}) we see that $(\phi_1\psi, p)^*\Z = (\phi_2\psi, p)^*\Z$ if and only if $(\phi_1, p)^*\Z = (\phi_2, p)^*\Z$. So we can assume that $k' = \Omega$ is an algebraic closure of $k(s)$.

Since the morphisms $q_1: S'_1 \to S$ and $q_2: S'_2 \to S$ are proper and birational, there exists a proper birational morphism $q_3: S'_3 \to S$ which factors through both $q_1$ and $q_2$ and such that any factorisation $Spec(\Omega) \to S_3' \to S$ is a good factorisation.\footnote{Let $U \subset S$ be a dense open subset over which $(q_1)_{red}$ and $(q_2)_{red}$ are both isomorphisms, and let $S'_3$ be the closure of the pre-image of $U$ in $S'_1 \times_S S'_2$. If the morphisms $\overline{\{z_i'\}} \to S'_3$ are not flat (where $\Z = \sum n_i z_i$ and $z_i'$ is $z_i$ seen as a point of $S'_3 \times_S X$) then the platification theorem (Theorem~\ref{theo:platification}) gives a blow-up of $S'_3$ with nowhere dense centre such that the proper transforms $\overline{\{z_i'\}}^\sim \to {S'_3}^\sim$ are flat. We then replace $S'_3$ with ${S'_3}^\sim$.} The morphisms $Spec(\Omega) \to S'_j$ factor through $S'_3$ for $j = 1, 2$. Let $\phi_j': Spec(\Omega) \to S'_3$ be the resulting morphisms.
\[ \xymatrix{
Spec(\Omega) \ar@<0.5ex>[dr]^{\phi'_1} \ar@<-0.5ex>[dr]_{\phi'_2} \ar@/^12pt/[drr]^{\phi_1} \ar@/_12pt/[ddr]_{\phi_2} \\
& S'_3 \ar[r]^{r_1} \ar[d]_{r_2} \ar[dr]^{q_3} & S'_1 \ar[d]^{q_1} \\
& S'_2 \ar[r]_{q_2} & S
} \]
Now for $j = 1, 2$ we have
\begin{equation} \label{equa:twoFatPoints}
\begin{split}
(q_j, \phi_j)^*\Z &= \phi^*_{j,nai} q^*_{j, nai}\Z \stackrel{\ref{coro:naiFon}}{=} \phi'^*_{j,nai} r^*_{j, nai}q^*_{j, nai}\Z \\
&\stackrel{\ref{lemm:fonBir}}{=} \phi'^*_{j,nai} (q_3)^*_{nai}\Z \\
&= (q_3, \phi'_j)^*\Z
\end{split}
\end{equation}
So we have reduced to showing that $(\phi'_1, q_3)^*\Z = (\phi'_2, q_3)^*\Z$.

Now use the same argument to build the following diagram
\[ \xymatrix{
Spec(\Omega) \ar@<0.5ex>[dr]^{\phi''_j} \ar@<-0.5ex>[dr]_{\iota'} \ar@/^12pt/[drr]^{\phi'_j} \ar@/_12pt/[ddr]_{\iota} \\
& S'_4 \ar[r]^{r_3} \ar[d]_{r_4} \ar[dr]^h & S'_3 \ar[d]^{q_3} \\
& S' \ar[r]_{p} & S
} \]
where $j = 1$ or $2$ and $r_3, r_4$ are birational and proper. The same argument as in Equation~\ref{equa:twoFatPoints} shows that $(\iota, p)^*\Z = (\iota', h)^*\Z$ and 
$(\phi'_j, q_3)^*\Z = (\phi''_j, h)^*\Z$. The second condition now says that $(\phi''_j, h)^*\Z = (\iota', h)^*\Z$ for $i = 1, 2$ and hence, $(\iota, p)^*\Z = (\phi'_j, q_3)^*\Z$ for $j = 1, 2$ and so $(\phi'_1, q_3)^*\Z = (\phi'_2, q_3)^*\Z$.
\end{proof}

The following theorem gives our definition of the $f^\circledast$ (see Definition~\ref{defi:cyclePullbacks}). It is morally equivalent to the definition given by Suslin-Voevodsky which is described before \cite[Lemma 3.3.9]{SV}.

\begin{theo}[{cf. \cite[Theorem 3.3.1]{SV}}] \label{theo:pullback}
Suppose that $S$ is a noetherian scheme, $X \to S$ a morphism of finite type and $\Z  \in c_{equi}(X / S, 0)$. Let $f: T \to S$ be a morphism of noetherian schemes. There exists a commutative diagram
\[ \xymatrix{
\amalg Spec(\Omega_j) \ar[rr]^{\iota = \sum \iota_j} \ar[d]_{q = \amalg q_i} && S' \ar[d]^p \\
\amalg \tau_j \ar[r]_t & T \ar[r]_f & S
} \]
where
\begin{enumerate}
 \item $t: \amalg \tau_i \to T$ is the inclusion of the generic points of $T$,
 \item $\Omega_i$ are algebraic closures of the $k(\tau_i)$,
 \item $(\iota_j, p)$ is a good factorisation of $Spec(\Omega_j) \to S$ with respect to $\Z$ \end{enumerate}
More importantly, there also exists a unique cycle $\W \in c_{equi}^{nai}(T \times_S X / T, 0)$ such that
\[ (t q)^*_{nai}\W = (\iota, p)^*\Z \]
and this $\W$ is uniquely determined by $f$ and $\Z$.
\end{theo}

\begin{proof}
\emph{Existence of the diagram.} The diagram exists by Lemma~\ref{lemm:existanceGoodFactor}.

\emph{Uniqueness.} This is clear since $(t q)^*_{nai}$ is injective (Lemma~\ref{lemm:naiInjective}).

\emph{Construction of $\W$ and membership in $c_{equi}^{nai}(T \times_S X / T, 0)$.} By the second axiom in Definition~\ref{defi:cequi} there exists a unique cycle $\W' = \sum m_k w_k' \in c_{equi}^{nai}((\amalg \tau) \times_S X / (\amalg \tau_i), 0)$ such that $q^*_{nai}\W' = (\iota, p)^*\Z$. Since $\amalg \tau_i \to T$ is birational, we can consider the $w_k'$ as points $w_k$ in $T \times_S X$ that lie over generic points of $T$. If the morphisms $\overline{\{w_k\}} \to T$ are finite, then sum $\W = n_i w_i$ belongs to $c_{equi}^{nai}(T \times_S X / T, 0)$ and satisfies $(t q)^*_{nai}\W = (\iota, p)^*\Z$.

It is enough to show that $w_k \in \cup \overline{\{z_i\}}$ since then $\overline{\{w_k\}}$ is a closed subscheme of $T \times_S (\cup \overline{\{z_i\}})$ which is finite over $T$ because each $\overline{\{z_i\}}$ is finite over $S$. 

By the definition of $q^*_{nai}$ and $t^*_{nai}$ there is a point $w'_k \in (\amalg Spec(\Omega_j)) \times_S X$ that maps to $w_k$. By the definition of $\iota^*_{nai}$ this $w_k'$ is a generic point of some $(\amalg Spec(\Omega_j)) \times_{S'} \overline{\{z_i'\}}$ where $z_i'$ is the point $z_i$ thought of as a point of $S' \times_S X$. This means that $w'_k$ is mapped inside one of the $\overline{\{z_i'\}}$. Clearly, $z_i'$ is mapped to $z_i$ by $p$ and so $\overline{\{z_i'\}}$ is contained in $S' \times_S \overline{\{z_i\}}$ and therefore the image of $\iota(w_k') \in S' \times_S \overline{\{z_i\}}$. This means that $p\iota(w_k') \in \overline{\{z_i\}}$ and so since $f(w_k) = p\iota(w_k')$ we are done.

\emph{Independence of choices}. A second choice of $q, \iota, p$ gives a commutative diagram
\[ \xymatrix{
\amalg Spec(\Omega_j) \ar[rr]^{\iota'} \ar@/_24pt/[dd]_{q'} \ar[d]^{\alpha} \ar[drr]^{\iota''} && S'' \ar@/^12pt/[dd]^{p'} \\
\amalg Spec(\Omega_j) \ar[rr]^{\iota} \ar[d]_{q} && S' \ar[d]^p \\
\amalg \tau_j \ar[r]_t & T \ar[r]_f & S
} \]
where $\alpha$ is an isomorphism. The cycle $\W$ coming from the choice $q, \iota, p$ satisfies the criterion for the choice $q', \iota'', p$ since
\[ (t q')^*_{nai}\W \stackrel{\ref{coro:naiFon}}{=} \alpha^*_{nai}(t q)^*_{nai}\W =  \alpha^*_{nai}(\iota, p)^*\Z = \alpha^*_{nai}\iota^*_{nai}p^*_{nai} \stackrel{\ref{coro:naiFon}}{=} \iota''^*_{nai}p^*_{nai} = (\iota'', p)^*\Z. \]
Moreover, by our assumption that $\Z \in c_{equi}(X/S, 0)$ we have $(\iota', p')\Z = (\iota'', p)\Z$. So $\W$ satisfies the criterion for the choice $q', \iota', p'$. So independence of the choices follows from uniqueness.
\end{proof}

\begin{prop} \label{prop:pullbackIsACorr}
The cycle $\W \in c_{equi}^{nai}(T \times_S X / T, 0)$ described in Theorem~\ref{theo:pullback} is in fact in $c_{equi}(T \times_S X / T, 0)$.
\end{prop}

\begin{proof}
Continuing with the notations from Theorem~\ref{theo:pullback}, suppose that $\Omega'$ is an algebraically closed field, $Spec(\Omega') \to T$ is an $\Omega'$ point of $T$ and $(\iota', T' \stackrel{p'}{\to} T)$ is a good factorisation of this $\Omega'$ point. Since $\Omega'$ is algebraically closed, the composition $f p' \iota'$ admits a lifting $\phi: Spec(\Omega') \to S'$. The diagram is the following.
\[ \xymatrix{
Spec(\Omega') \ar[r]_{\iota'} \ar@{-->}@/^12pt/[rr]^{\phi} & T' \ar[d]_{p'} & S' \ar[d]^p \\
& T \ar[r]_f & S
} \]
 We claim that
\begin{equation} \label{equa:pullbackIsACorr}
(\iota', p')^* \W = \phi^*_{nai}p^*_{nai}\Z.
\end{equation}
Since $\Z \in c_{equi}(X/S, 0)$ if this equality holds then it follows from the definition of $c_{equi}(X/S, 0)$ that $\W \in c_{equi}(T \times_S X / T, 0)$.

We claim that there exists a commutative diagram
\[ \xymatrix{
& W' \ar[r]^c & W \ar[d]^b \ar[dr]^a \\
Spec(\Omega') \ar[rr]_{\iota'} \ar[ur]^d && T' \ar[d]_{p'} & S' \ar[d]^p \\
&& T \ar[r]_f & S
} \]
such that
\begin{enumerate}
 \item $W$ is integral and the generic point of $W$ hits a generic point of $T'$, and the induced field extension is finite, and
 \item $c$ induces an isomorphism over a dense open subscheme of $W$, and 
 \item if we write $c^*_{bir}b^*_{nai}p'^*_{bir}(\W) = \sum \ell_i x_i$ then the $\overline{\{x_i\}}$ are flat over $W'$. 
\end{enumerate}
Notice that with these hypotheses, Lemma~\ref{lemm:fonBir} and Corollary~\ref{coro:naiFon} imply that
\begin{equation} \label{equa:bcd}
(bcd)^*_{nai} = b^*_{nai}c^*_{nai}d^*_{nai} \qquad \textrm{ and } \qquad (acd)^*_{nai} = a^*_{nai}c^*_{nai}d^*_{nai}
\end{equation}

To find such a diagram, consider the composition $T' \times_S S' \to T' \to T$. Choose a generalisation $\tau$ of $\iota'(Spec(\Omega')) \in T'$. Since this composition $T' \times_S S' \to T' \to T$ is finite type and surjective, there is a point $\tau' \in T' \times_S S'$ in the pre-image of $\tau$ such that the induced field extension is finite. We set $W = \overline{\{\tau'\}}$. This gives us $a, b$ and a factorisation of $\iota'$ through $W$. Now we use the platification theorem (Theorem~\ref{theo:platification}) to find a blow-up $c: W' \to W$ of $W$ such that the strict transform of $\amalg \overline{\{x_i\}} \to W$ is flat. Since $W' \to W$ is surjective of finite type, every point of $W$ has a point over it such that the induced field extension is finite, hence the morphism $d$. 

To prove the equality (\ref{equa:pullbackIsACorr}) we will show
\begin{equation} \label{equa:pullbackIsACorr2}
b^*_{nai}p'^*_{nai} \W = a^*_{nai}p^*_{nai}\Z.
\end{equation}
The equality (\ref{equa:pullbackIsACorr}) will then follow from (\ref{equa:bcd}).

Let $w$ be the generic point of $W$ and $\Sigma$ be an algebraic closure of $k(w)$ with induced morphism $\theta: Spec(\Sigma) \to W$. Since $\Sigma$ is also an algebraic closure of the field of functions of a generic point of $T$ the definition of $\W$ says that we have $(b\theta)^*_{nai}p'^*_{nai} \W = (a\theta)^*_{nai}p^*_{nai}\Z$. It follows now from Corollary~\ref{coro:naiFon} that we have the equality~\ref{equa:pullbackIsACorr2}.
\end{proof}

\begin{defi} \label{defi:cyclePullbacks}
Suppose that $S$ is a noetherian scheme, $X \to S$ a morphism of finite type and $\Z \in c_{equi}(X / S, 0)$. Let $f: T \to S$ be a morphism of noetherian schemes. We define
\[ f^\circledast\Z = \W \in c_{equi}(T \times_S X / T, 0) \]
where $\W$ is the cycle given by Theorem~\ref{theo:pullback} (cf. Proposition~\ref{prop:pullbackIsACorr} as well). By the uniqueness of $\W$, there is an induced homomorphism of abelian groups
\[ f^\circledast: c_{equi}(X / S, 0) \to c_{equi}(T \times_S X / T, 0). \]
\end{defi}

\begin{lemm} \label{lemm:agreesWithOldPullbacks}
Suppose that $S$ is a noetherian scheme, $X \to S$ a morphism of finite type and $\Z = \sum n_i z_i \in c_{equi}(X / S, 0)$. Let $f: T \to S$ be a morphism of noetherian schemes. 
\begin{enumerate}
 \item If the images of the generic points of $T$ are in the flat locus of each $\overline{\{z_i\}} \to S$ then
\[ f^\circledast\Z = f^*_{nai}\Z. \]
 \item If $f$ is pseudo-dominant, then
\[ f^\circledast\Z = f^*_{nai}\Z. \]
 \item If $k$ is a field, $f: Spec(k) \to S$ is a $k$ point and $(\iota', p)$ is a good factorisation with respect to $\Z$ then
\[ f^\circledast\Z = (\iota', p)^*\Z. \]
\end{enumerate}
\end{lemm}

\begin{proof}
We use the notation of Theorem~\ref{theo:pullback}.
\begin{enumerate}
 \item By the the platification theorem (Theorem~\ref{theo:platification}) we can find a proper birational morphism $S' \to S$ that is an isomorphism over the flat locus of $\overline{\{z_i\}} \to S$. Consequently, we have the following commutative diagram
\[ \xymatrix{
\amalg Spec(\Omega_j) \ar[rr]^{\iota} \ar[d]_{q} && S' \ar[d]^p \\
\amalg \tau_j \ar[r]_t \ar[urr]^{t'} & T \ar[r]_f & S
} \]
We then have
\[ \iota^*_{nai}p^*_{nai}\Z \stackrel{\ref{lemm:fonBir}}{=} q^*_{nai}t'^*_{nai}p^*_{nai}\Z \stackrel{\ref{coro:naiFon}}{=} q^*_{nai}(pt')^*_{nai}\Z \stackrel{\ref{lemm:fonBir}}{=} q^*_{nai}t^*_{nai}f^*_{nai}\Z \]
and so $f^*_{nai}\Z$ satisfies the criterion defining $f^\circledast\Z$.

 \item In this case $\iota$ (and of course $q$ and $t$) are pseudo-dominant as well, and so it follows from Lemma~\ref{lemm:fonBir}.

 \item Our diagram is
\[ \xymatrix{
Spec(\Omega) \ar[r]^\iota \ar[d]_q & S' \ar[d]^p \\
Spec(k) \ar[r]_f \ar[ur]^{\iota'} & S 
} \]
and we have
\[ q^*_{nai}f^\circledast\Z \stackrel{def}{=} \iota^*_{nai}p^*_{nai}\Z \stackrel{\ref{coro:naiFon}}{=} q^*_{nai}\iota'^*_{nai}p^*_{nai}\Z \]
so the claim follows from the fact that $q^*_{nai}$ is injective (Lemma~\ref{lemm:naiInjective}).
\end{enumerate}
\end{proof}

\begin{lemm} \label{lemm:compositionPullback}
Suppose that $S$ is a noetherian scheme, $X \to S$ a morphism of finite type and $\Z = \sum n_i z_i \in c_{equi}(X / S, 0)$. Let $U \stackrel{g}{\to} T \stackrel{f}{\to} S$ be a pair of composable morphisms of noetherian schemes. Then
\[ g^\circledast f^\circledast\Z = (fg)^\circledast\Z. \]
\end{lemm}

\begin{proof}
We use Lemma~\ref{lemm:agreesWithOldPullbacks}. It follows directly from the definition that the result is true if $g$ is of the form $g: \amalg Spec(\Omega_j) \to T$ where $\Omega_j$ are algebraic closures of the function fields $k(\tau_j)$ at the generic points $\tau_j$ of $T$. So it suffices now to consider the case where $U$ is of the form $Spec(\Omega)$ (but not necessarily hitting a generic point of $T$). In this situation however, the result follows immediately from the claim (\ref{equa:pullbackIsACorr}) in the proof of Proposition~\ref{prop:pullbackIsACorr}.
\end{proof}

\begin{theo} \label{theo:defcequi}
Suppose that $S$ is a noetherian scheme and $X \to S$ a morphism of finite type. Then the groups $c_{equi}(- \times_S X / -, 0)$ form a presheaf of relative cycles. Moreover, every other presheaf of relative cycles associated to $X/S$ is a subpresheaf of this presheaf.
\end{theo}

\begin{proof}
We have $c_{equi}(T \times_S X / T, 0) \subset c_{equi}^{nai}(T \times_S X / T, 0)$ by definition. We have seen that the $c_{equi}(- \times_S X / -, 0)$ with the morphisms $(-)^\circledast$ are a presheaf (Lemma~\ref{lemm:compositionPullback}) and that they satisfy the two properties asked of a presheaf of relative cycles (Lemma~\ref{lemm:agreesWithOldPullbacks}). Moreover, if $F$ is a presheaf of relative cycles (associated to a morphism of finite type $X \to S$) then we have also seen that the elements $\Z \in F(T)$ satisfy the properties asked of an element of $c_{equi}(T \times_S X / T, 0)$ (Proposition~\ref{prop:cIsprc}). Proposition~\ref{prop:sameMorphisms} tells us that $F$ is a subpresheaf of $c_{equi}(- \times_S X / -, 0)$.
\end{proof}

\begin{defi} \label{defi:presheaf}
Suppose $X \to S$ is a morphism of finite type with $S$ a noetherian scheme. We abusively use $c_{equi}(X /S, 0)$ to also denote the presheaf of relative cycles $c_{equi}(- \times_S X / -, 0)$.
\end{defi}




\section{The category of correspondences}

We discuss now the category $Cor(S)$ of correspondences (cf. \cite{SV}, \cite{CD}, \cite{Ivo05}, \cite[p.141]{FV}). As for relative presheaves, we define $Cor(S)$ by means of a universal property. We give a short proof of its existence but for the hardest part -- the construction of the correspondence homomorphisms of \cite[Section 3.7]{SV}, and the fact that the induced composition in $Cor(S)$ is associative -- we cite the literature. We give an explicit expression for various compositions of correspondences, and also show that $Cor(S)$ satisfies analogues of the axioms for a presheaf with traces that we will define later.

We begin with an easy corollary of Theorem~\ref{theo:naiFon}.

\begin{prop} \label{prop:cyclesOfFlatSubschemes}
Let $S$ be a noetherian scheme, $X \to S$ a morphism of finite type, and $Z \subset X$ a closed subscheme such that $Z \to S$ is flat and finite. Then $cycl_{X/S}(Z) \in c_{equi}(X / S, 0)$ and if $f: T \to S$ is an morphism of noetherian schemes then $f^\circledast cycl_{X/S}(Z) = cycl_{T \times_S X / T}(T \times_S Z)$.
\end{prop}

\begin{proof}
Firstly, notice that in the case when $f$ is birational, we have $f^*_{nai} cycl_{X/S}(Z) = cycl_{T \times_S X / T}(T \times_S Z)$.

Secondly, suppose that $S$ is reduced, and that $T$ is the spectrum of a field whose image in $S$ is in the flat locus of $\amalg \overline{\{z_i\}} \to S$ where the $z_i$ are the generic points of $Z$. Suppose further that $Z_i \to S$ is flat where the $Z_i$ are the irreducible components of $Z$ with their reduced structure. We claim that $f^*_{nai}cycl(Z) = cycl(T \times_S Z)$ in this case. Let $\Z = Z - \sum n_i Z_i$ where $n_i = \length \OO_{Z, z_i}$. We have $cycl(\Z) = 0$ and so $cycl(k \times_S \Z) = 0$ by Theorem~\ref{theo:naiFon}. Consequently, by linearity, it suffices to consider the case when $Z$ is integral. But this case follows immediately from the definition of $f_{nai}^*$.

Now that we have these two facts, the statement that $cycl_{X/S}(Z) \in c_{equi}(X / S, 0)$ is a direct consequence of Definition~\ref{defi:cequi}(1) and Proposition~\ref{prop:alternativeConditionsForcequi}(2). The statement $f^\circledast cycl_{X/S}(Z) = cycl_{T \times_S X / T}(T \times_S Z)$ follows for the same reasons from the definition of $f^\circledast$ (see Theorem~\ref{theo:pullback}).
\end{proof}

\begin{defi} \label{defi:graphAndTranspose}
Let $f: Y \to X$ be a morphism in $Sch(S)$ and $\Gamma_f \subset Y \times_S X$ the closed subscheme that is its graph. We define
\[ [f] = cycl_{Y \times_S X/Y}({\Gamma_f}_{red}) \in c_{equi}^{nai}(Y \times_S X/ Y, 0). \]
If $f$ is finite flat then we define
\[ {[^tf]} = cycl_{X \times_S X/X}({^t\Gamma_f}) \in c_{equi}^{nai}(X \times_S Y / X, 0) \]
where ${^t\Gamma_f} \subset X \times_S Y$ is the closed subscheme corresponding to $\Gamma_f \subset Y \times_S X$. Explicitly, we have $[f] = \sum z_i$ and ${[^tf]} = \sum n_i z_i$ where the $z_i$ are the generic points of $Y$ (seen as points of $X \times_S Y$ or $Y \times_S X$) and $n_i = \length \OO_{f(z_i) \times_X Y, z_i}$.
\end{defi}

\begin{rema}
Notice that $[f]$ has no coefficients, even when $Y$ and $X$ are non-reduced, whereas in ${[^tf]}$ we have taken care to include the multiplicities of the generic points of $Y$ in their fibres. We insist that this is necessary to make the theory work.
\end{rema}

\begin{lemm}
For every morphism $f: Y \to X$ in $Sch(S)$ the formal sum $[f]$ lies in $c_{equi}(Y \times_S X / Y, 0)$. If $f$ is finite flat, then $[^tf]$ lies in $c_{equi}(X \times_S Y / X, 0)$.
\end{lemm}

\begin{proof}
For ${[^tf]}$, we have already proven in Proposition~\ref{prop:cyclesOfFlatSubschemes} that formal sums of the form $cycl_{X \times_S Y / X}(Z)$ lie in $c_{equi}(X \times_S Y / X, 0)$ for closed subschemes $Z \subset X \times_S Y$ that are flat and finite over $X$.

For $[f]$, it is clear from the definition that $c_{equi}(Y \times_S X / Y, 0) \to c_{equi}(Y_{red} \times_S X / Y_{red}, 0)$ is an isomorphism and so it suffices to consider the case when $Y$ is reduced. But then $(\Gamma_f)_{red} \to Y$ is flat and finite (it is an isomorphism) and so $[f] \in c_{equi}(Y \times_S X / Y, 0)$ for the same reasons as ${[^tf]}$.
\end{proof}

\begin{defi} \label{defi:compositionCorrespondances}
Suppose that $S$ is a noetherian scheme and $W, Y, X$ are three $S$-schemes of finite type. We define a bilinear morphism
\[ - \circ - : c_{equi}(Y \times_S W / Y, 0) \otimes c_{equi}^{nai}(X \times_S Y / X, 0) \to c_{equi}^{nai}(X \times_S W / X, 0) \]
as follows. Let $\beta \in c_{equi}(Y \times_S W / Y, 0)$ and $\alpha = n_i z_i \in c_{equi}^{nai}(X \times_S Y / X, 0)$. Let $Z_i = \overline{\{z_i\}}$ and let $\iota_i:  Z_i \to Y$ be the canonical morphisms. Then we define
\begin{equation} \label{equa:cycleComposition}
\beta \circ \alpha = \sum n_i m_{ij} d_{ij} w_{ij}'
\end{equation}
where $\iota_i^\circledast \beta = \sum m_{ij} w_{ij} \in c_{equi}(Z_i \times_S W / Z_i, 0)$, the $w_{ij}'$ are the images of the $w_{ij}$ in $X \times_S W$ under the canonical (finite) morphism $Z_i \times_S W \to X \times_S W$ and $d_{ij} = [k(w_{ij}):k(w_{ij}')]$. 
\end{defi}

The following theorem we cite from the literature.

\begin{theo}[{\cite[Theorem 3.7.3]{SV}, \cite[Section 2.1.1]{Ivo05}}] \label{theo:cycleComposition}
The morphism $- \circ -$ of Definition~\ref{defi:compositionCorrespondances} satisfies the following properties.
\begin{enumerate}
 \item If $\alpha \in c_{equi}(X \times_S Y / X, 0), \beta \in c_{equi}(Y \times_S W / Y, 0)$ then $\beta \circ \alpha \in c_{equi}(X \times_S W / X, 0)$.
 \item Suppose $V, W, X, Y$ are four $S$-schemes of finite type and $\alpha \in c_{equi}(V \times_S W / V, 0), \beta \in c_{equi}(W \times_S X / W, 0), \gamma \in c_{equi}(X \times_S Y / X, 0)$. Then $(\gamma \circ \beta) \circ \alpha = \gamma \circ (\beta \circ \alpha)$.
\end{enumerate}
\end{theo}

\begin{prop} \label{prop:variousCompositions}
The morphism $- \circ -$ of Definition~\ref{defi:compositionCorrespondances} satisfies the following properties.
\begin{enumerate}
 \item If $f: X \to Y$ is a morphism in $Sch(S)$ and $\beta \in c_{equi}(Y \times_S W / Y, 0)$ for some $W \in Sch(S)$ then 
\[ \beta \circ [f] = f^\circledast\beta. \]

 \item Suppose $f: V \to X$ is a finite flat surjective morphism in $Sch(S)$ and $\beta \in c_{equi}(V \times_S W / V, 0)$. Then
\begin{equation} \label{equa:compTrans}
\beta \circ [^tf] = \sum n_j m_j d_j w_j'
\end{equation}
where $\beta = \sum m_j w_j$, the points $x_j, v_j, w_j'$ are the respective images of $w_j$ in $X$, $V$, and $X \times_S W$ via the obvious morphisms, $d_j = [k(w_j):k(w_j')]$, and $n_j = \length \OO_{x_j \times_X V, v_j}$.

 \item Suppose $\alpha = n_iz_i \in c_{equi}(X \times_S Y / X, 0)$ and $g: Y \to W$ is a morphism in $Sch(S)$. Then
\[ [g] \circ \alpha = \sum n_i d_i w_i \]
where $w_i = (X \times_S g)(z_i)$ and $d_i = [k(z_i): k(w_i)]$.
\end{enumerate}
\end{prop}

\begin{proof}
\begin{enumerate}
 \item This is straight-forward from our explicit description. Notice, that in our case, in the definition of the composition the closed integral subschemes $Z_i$ are canonically isomorphic to the irreducible components $X_i$ of $X$ and the morphisms $Z_i \to X$ are the compositions $X_i \to X \to Y$. Consequently, the morphisms $Z_i \times_S W \to X \times_S W$ are closed immersions and so the $d_{ij}$ are all 1. The result follows from the fact that $c_{equi}(X \times_S Y / X, 0) \to \oplus c_{equi}(X_i \times_S Y / X_i, 0)$ is the obvious morphism.

 \item As everything happens generically, we can replace $X$ by any dense open subscheme without changing the result. Hence, shrinking $X$ and using additivity, we are permitted to assume that $V$ and $X$ each have a unique irreducible component. In this case the $\alpha$ of Definition~\ref{defi:compositionCorrespondances} is $[^tf] = n v$ where $n = \length \OO_{x \times_X V, v}$ and $x, v$ are the generic points of $X$ and $V$ respectively. In the notation of the definition there is a unique $Z_i$ and the morphism $Z_i \to Y$ is isomorphic to the canonical morphism $V_{red} \to V$. Since $c_{equi}(V \times_S W / V, 0) \to \oplus c_{equi}(V_{red} \times_S W / V_{red}, 0)$ is the obvious morphism, the $\iota_i^\circledast \beta$ of the definition is $\sum m_j w_j$ now considered as an element of $c_{equi}(V_{red} \times_S W / V_{red}, 0)$. Finally, the $d_{ij}$ of the definition, of which we have only one, is $d = [k(v):k(x)]$. So the $\beta \circ \alpha = \sum n_i m_{ij} d_{ij} w_{ij}'$ of Equation~\ref{equa:cycleComposition} is, in our case, $\sum n m_j d w_j$ (note that we have no need of indices on $n$ and $d$ because we have assumed $V$ and $X$ are irreducible).

 \item First consider the case that $Y$ is reduced. By definition $[g]$ is $cycl_{Y \times_S W / Y}(\Gamma_{g_{red}})$ but since $Y$ is reduced this is just $cycl_{Y \times_S W / Y}(\Gamma_g)$. Moreover, $\Gamma_g$ is canonically isomorphic to $Y$. Let $\iota_i: Z_i \to Y$ and $\beta$ be as in the definition of $- \circ -$ (so $\beta = cycl_{Y \times_S W / Y}(\Gamma_g)$). By Proposition~\ref{prop:cyclesOfFlatSubschemes} we see that $\iota_i^\circledast \beta = cycl_{Z_i \times_S W / Z_i}(\Gamma_{(Z_i \to W)})$. The result is now clear from the explicit formula in the definition of $- \circ -$.

Now we remove the assumption that $Y$ is reduced. Notice that as $\alpha$ doesn't actually depend on the ambient scheme $X \times_S Y$, it also defines an element $\alpha'$ of ${c_{equi}(X \times_S Y_{red} / X, 0)}$. Moreover, using the case when $Y$ was reduced we can write $\alpha = [i] \circ \alpha'$ where $i: Y_{red} \to Y$ is the canonical morphism. So now using the fact that $- \circ -$ is associative, it suffices to prove that $[g] \circ [i] = [gi]$. This follows from the first part as $[gi] = i^\circledast[g] = i^*_{nai}[g]$ since $i$ is birational.
\end{enumerate}
\end{proof}

\begin{prop} \label{prop:corHasTraces}
The morphism $- \circ -$ of Definition~\ref{defi:compositionCorrespondances} satisfies the following properties.
\begin{enumerate}
 \item \emph{Functoriality.} For finite flat surjective morphisms $W \stackrel{g}{\to} Y \stackrel{f}{\to} X$ we have 
\[ {[^tf]} \circ {[^tg]} = {[^tfg]} . \]

 \item \emph{Base-change.} For every cartesian square
\begin{equation} \label{equa:cartSquare}
 \xymatrix{
Y \times_X W \ar[r]^-g \ar[d]_q & W \ar[d]^p \\
Y  \ar[r]_f & X
} \end{equation}
such that $f$ is finite flat surjective we have
\[ {[^tf]} \circ [p] = [q] \circ {[^tg]}. \]

 \item \emph{Degree.} For every finite flat surjective morphism $f: Y \to X$ of constant degree $d$ we have 
\[ [f] \circ {[^tf]} = d \cdot [id_X]. \]

 \item \emph{Triangles.} Consider a commutative triangle of schemes with $f, g$ finite flat surjective and $X$ integral.
\[ \xymatrix{ Y' \ar[rr]^h \ar[dr]_g && Y \ar[dl]^f \\
& X
} \]
\begin{enumerate}
 \item Suppose that the scheme $Y'$ is the disjoint union of the integral components $Y_i'$ of $Y$, and $h$ is the canonical morphism \mbox{$Y' = \amalg Y_i' \to Y$}. Then 
\[ [^tf] = \sum m_i [h_i] \circ [^tg_{i}] \]
where $h_i, g_i$ are the restrictions to $Y_i'$ and $m_i = \length \OO_{Y, \eta_i}$ with $\eta_i$ the generic point of $Y_i'$.
 \item Forgetting the hypotheses of (a), suppose now that all schemes are integral. Then 
\[ \tfrac{\deg g}{\deg f} [^tf] = [h] \circ [^tg]. \]
\end{enumerate}
\end{enumerate}
\end{prop}

\begin{proof}
\emph{Functoriality}. Let $\xi_i$ be the generic points of $X$, let $\eta_{ij}$ be the generic points of $Y$ (over $\xi_i$) and let $\omega_{ijk}$ be the generic points of $W$ (over $\eta_{ij}$). By definition ${[^tfg]}$ is $\sum_{ijk} \length \OO_{\xi \times_X W, \omega_{ijk}} \xi_i$. Using (\ref{equa:compTrans}) we calculate ${[^tg]}{[^tf]}$ as
\[ \sum_{ijk} \length \OO_{\eta_{ij} \times_Y W, \omega_{ijk}} \length \OO_{\xi_{i} \times_X Y, \omega_{ij}} [k(\omega_{ijk}) : k(\eta_{ij})] \xi_i. \]
So we must show that for each $i$ we have
\[ \sum_{jk} \length \OO_{\eta_{ij} \times_Y W, \omega_{ijk}} \length \OO_{\xi_{i} \times_X Y, \omega_{ij}} [k(\omega_{ijk}) : k(\eta_{ij})] = \sum_{jk} \length \OO_{\xi \times_X W, \omega_{ijk}}. \]
This is done in Lemma~\ref{lemm:appendix3}.

\emph{Base change formula}. This follows directly from Proposition~\ref{prop:cyclesOfFlatSubschemes}. Let $V = Y \times_X W$. We have
\[ {[^tf]}[p] = p^\circledast({[^tf]}) =  p^\circledast(cycl_{X \times_S Y / X}(Y)) = cycl_{W \times_S Y / W}(V) = \sum n_iv_i \in c_{equi}(W \times_S Y / W, 0) \]
where $v_i$ are the generic points of $V$ and $n_i = \length \OO_{g(v_i) \times_W V, v_i}$. We also have $[q] = \sum v_i \in c_{equi}(V \times_S Y / V, 0)$ and using the formula (\ref{equa:compTrans}) we find that $[q]{[^tg]} = \sum n_i v_i \in c_{equi}(W \times_S Y / W, 0)$ as $V = W \times_X Y \to W \times_S Y$ is a closed immersion.

\emph{Degree formula}. Suppose $\xi_i$ are the generic points of $X$ and $\eta_{ij}$ the generic points of $Y$ with $\eta_{ij}$ over $\xi_i$. Still using (\ref{equa:compTrans}) we calculate $[f] {[^tf]}$ as $\sum_{ij} \length \OO_{\xi_i \times_X Y, \eta_{ij}} [k(\eta_{ij}) : k(\xi_i)] \xi_i$. For each $\xi_i$, the degree of $\xi_i \times_X Y \to \xi_i$ is $d$, and this is equal to $\sum_j \length \OO_{\xi_i \times_X Y, \eta_{ij}} [k(\eta_{ij}) : k(\xi_i)]$, hence the degree formula.

\emph{Triangles}. Both equalities follow directly from Proposition~\ref{prop:variousCompositions}.
\end{proof}

\begin{theo} \label{theo:categoryOfCorrespondances}
For each noetherian separated scheme $S$ there exists a unique category $Cor(S)$ with the following properties.
\begin{enumerate}
 \item[(Ob)] There is a one-to-one correspondence between the objects of $Cor(S)$ and the objects of $Sch(S)$. If $X \in Sch(S)$ we denote the corresponding object in $Cor(S)$ by $[X]$.
 \item[(Mor)] Consider $X, Y \in Sch(S)$. Then $\hom_{Cor(S)}([X], [Y])$ is a subgroup of the free abelian group generated by the points $z$ of $X \times_S Y$ such that the canonical morphism $\overline{\{z\}} \to X$ is finite and dominates an irreducible component of $X$.
 \item[(Gra)] Let $f: Y \to X$ be a morphism in $Sch(S)$. Then $[f] \in \hom_{Cor(S)}([X], [Y])$. Furthermore, if $f$ is finite and flat then $[^tf] \in \hom_{Cor(S)}([Y], [X])$.
 \item[(Com1a)] If $f: X_{red} \to X$ is the canonical inclusion then the morphism $\hom_{Cor(S)}([X], [Y]) \to \hom_{Cor(S)}([X_{red}], [Y])$ induced by composition with $[f]$ is the obvious one coming from the canonical identification of the points of $X \times_S Y$ and the points of $X_{red} \times_S Y$.
 \item[(Com1b)] If $\alpha = \sum n_i z_i \in \hom_{Cor(S)}([X], [Y])$, $k$ is a field, and $\iota: k \to X$ is a $k$-point of $X$ such that the image is in the flat locus of $\amalg \overline{\{z_i\}} \to X$ then $\alpha \circ [\iota] = \sum n_i m_{ij} w_{ij}$ where the $w_{ij}$ are the generic points of $k \times_Y \overline{\{z_i\}}$ and $m_{ij} = \length \OO_{k \times_Y \overline{\{z_i\}}, w_{ij}}$.
 \item[(Com2)] Suppose that we have $f: V \to X$ a finite flat surjective morphism in $Sch(S)$ and that $\alpha \in \hom_{Cor(S)}([V], [Y])$. Then
\[ \alpha \circ [^tf] = \sum n_i m_i d_i y_i' \]
where $\alpha = \sum n_i y_i$, the points $x_i, v_i, y_i'$ are the respective images of $y_i$ in $X$, $V$, and $X \times Y$ via the obvious morphisms, $d_i = [k(y_i):k(y_i')]$, and $m_i = \length \OO_{x_i \times_X V, v_i}$.
 \item[(Uni)] Any other category satisfying the above axioms is a subcategory of $Cor(S)$.
\end{enumerate}
Moreover, the composition in this unique category $Cor(S)$ is the one given in Definition~\ref{defi:compositionCorrespondances}.
\end{theo}


\begin{rema}
In light of Definition~\ref{defi:presheafOfRelativeCycles} and Theorem~\ref{theo:defcequi} we could have replaced the axioms (Mor), (Gra), (Com1a), (Com1b), (Uni) with the two axioms (Mor$'$) and (Com1$'$). This would have given the following list.
\begin{enumerate}
 \item[(Ob)] There is a one-to-one correspondence between the objects of $Cor(S)$ and the objects of $Sch(S)$. If $X \in Sch(S)$ we denote the corresponding object in $Cor(S)$ by $[X]$.
 \item[(Mor$'$)] Consider $X, Y \in Sch(S)$. Then $\hom_{Cor(S)}([X], [Y]) = c_{equi}(X \times_S Y / X, 0)$.
 \item[(Com1$'$)] If $f: X \to Y$ is a morphism in $Sch(S)$ and $\alpha \in \hom_{Cor(S)}([Y], [W])$ for some $W \in Sch(S)$ then $\alpha \circ [f] = f^\circledast\alpha$.
 \item[(Com2)] Suppose that we have $f: V \to X$ a finite flat surjective morphism in $Sch(S)$ and that $\alpha \in \hom_{Cor(S)}([V], [Y])$. Then
\[ \alpha \circ [^tf] = \sum n_i m_i d_i y_i' \]
where $\alpha = \sum n_i y_i$, the points $x_i, v_i, y_i'$ are the respective images of $y_i$ in $X$, $V$, and $X \times Y$ via the obvious morphisms, $d_i = [k(y_i):k(y_i')]$, and $m_i = \length \OO_{x_i \times_X V, v_i}$.
\end{enumerate}
We chose the statement in the theorem because there is no explicit reference to presheaves of relative cycles.
\end{rema}

\begin{proof}
We begin with uniqueness. Since the objects and the morphisms are completely described, it suffices to show that the composition is determined (Com1$'$) and (Com2). Let $\alpha \in \hom_{Cor(S)}([X], [Y])$ and $\beta \in \hom_{Cor(S)}([Y], [W])$ and suppose that $\circ$ and $\circ'$ are two different compositions. Since pullback along a birational morphism is injective (Lemma~\ref{lemm:domBlow} for example), to show that $\beta \circ \alpha = \beta \circ' \alpha$ it suffices to show that $f^\circledast (\beta \circ \alpha) = f^\circledast(\beta \circ' \alpha)$ for some birational $f: X' \to X$. Let $\alpha = n_i z_i$. The platification theorem (Theorem~\ref{theo:platification}) provides a birational morphism $f: X' \to X$ such that the proper transforms of the $\overline{\{z_i\}} \to X$ are flat over $X'$. Let $f^\circledast\alpha = \sum n_i z'_i$ and let $g_i: \overline{\{z_i'\}} \to Y$ and $h_i: \overline{\{z_i'\}} \to X'$ be the canonical morphisms. Then we have
\[ f^\circledast (\beta \circ \alpha) = \beta \circ \alpha \circ f = \beta \circ [g_i] \circ [^th_i] \]
and similarly for $\circ'$. Due to (Com1$'$) and (Com2) the cycles $\beta \circ [g_i] \circ [^th_i]$ and $\beta \circ' [g_i] \circ' [^th_i]$ are equal. Therefore $\beta \circ \alpha = \beta \circ' \alpha$.

Now for existence. The majority of the difficulty of the proof of existence is contained in Theorem~\ref{theo:cycleComposition}. Since we are admitting this, it remains to show that the composition has identities, and satisfies (Com1$'$) and (Com2). These all follow from Proposition~\ref{prop:variousCompositions}.
\end{proof}

\begin{defi} \label{defi:categoryOfCorrespondances}
The category $Cor(S)$ of Theorem~\ref{theo:categoryOfCorrespondances} is call the \emph{category of correspondences}. The \emph{category of smooth correspondences} $SmCor(S)$ is the full subcategory of $Cor(S)$ whose objects are smooth schemes over $S$.
\end{defi}


\chapter{Comparison of cdh and $\ldh$ sheafification and cohomology} \label{chap:comparison}
\section{Introduction}

\lettrine[lines=2, findent=3pt, nindent=0pt]{I}{n} this chapter we introduce the $\ldh$-topology (Definition~\ref{defi:ltopology}) where $\ell$ denotes a prime. We compare the cdh and $\ldh$ sheafifications and cohomologies. The idea of the cdh topology is to enlarge the Nisnevich topology enough so that the morphisms coming from resolution of singularities may be used as covers. Similarly, the idea of the $\ldh$ topology is that it should be an enlargement of the Nisnevich topology so that morphisms given by a theorem of Gabber on alterations (Theorem~\ref{theo:gabberGlobal} or Theorem~\ref{theo:gabberLocal}) may be used as covers.

In Section~\ref{sec:ltopology} we begin the chapter by introducing our definition of the \emph{$\ldh$ topology}. Our definition (Definition~\ref{defi:ltopology}) -- equivalent to many others\footnote{While our definition is equivalent to many other possible definitions, it is different from the topology of $\ell'$-alterations appearing in \cite{ILO} and \cite{Ill09}. This is because they work with a category of reduced finitely horizontal schemes (i.e., every generic point is sent to a generic point of the base, and the induced field extension is finite) while we work with a more general category. Ours is a ``global'' version of theirs which is ``local'' where ``local'' and ``global'' are in the sense of resolution of singularities.} -- is inspired directly by the techniques that we will use to study it. Namely, it is generated in some sense by the cdh topology, and a topology we refer to as the {\fpsl topology} (fini-plat-surjectif-premier-{\`a} -$\ell$). In shorthand we could write ``cdh $+ $\fpsl$ = \ldh$''. In this section we also convert Gabber's Theorem into the form that we will apply it in: every nice scheme admits an $\ldh$ cover with regular source (Corollary~\ref{coro:regularlCover}).

The literature abounds with techniques to work with the cdh topology, and so we are left with the study of the $\fpsl$ topology. Our main tool here is the concept of a \emph{presheaf with traces} which we introduce in Definition~\ref{defi:traces}. A presheaf with traces is a presheaf which in addition to being a contravariant functor, also has a covariant functoriality for morphisms that are finite flat and surjective, and furthermore satisfies a change-of-base and degree formula. It falls straight out of our definition that every presheaf of $\zll$ modules with traces is an acyclic $\fpsl$ sheaf (Lemma~\ref{lemm:fpslAcyclic}).

In Section~\ref{sec:comparisonCohomology} we show that if we have a cdh sheaf of $\zll$-modules with transfers, then the cdh and $\ldh$ cohomologies agree. We do this using the technique of \cite[3.1.8]{Voev00}, that is, we claim that these cohomologies can be calculated using Ext's in the categories of sheaves with transfers (Proposition~\ref{prop:318}, Proposition~\ref{prop:318l}). This comes down to proving an acyclicity result, which we do in a more general context (Proposition~\ref{prop:sigmaTauAcyclic}). Accepting that we can use Ext's to calculate the cohomologies, since every presheaf with transfers is a presheaf with traces (Lemma~\ref{lemm:transfersImpliesTraces}), and hence an $\fpsl$ sheaf (Lemma~\ref{lemm:fpslAcyclic}), the categories of cdh and $\ldh$ sheaves with transfers are equivalent (Corollary~\ref{coro:equivalencecdhltransfers}) and so we deduce that the cohomologies agree.

In Section~\ref{sec:cd} we introduce a topology which will help us study the cdh and $\ldh$ associated sheaves of a presheaf with traces. The idea is to embed $F_{cdh}$ and $F_{\ldh}$ into a larger presheaf, and then descend properties of this larger presheaf to $F_{cdh}$ and $F_{\ldh}$. The larger presheaf that we use is the sheafification $F_{\dis}$ for a Grothendieck topology that we call the \emph{completely decomposed discrete topology} or cdd topology (Definition~\ref{defi:cd}). The most important property of the cdd topology is that the cdd associated sheaf $F_{cdd}$ of a presheaf $F$ satisfies $F_{cdd}(X) = \prod_{x \in X} F(X)$ where the product is over the points of $X$ of every codimension. For an explanation of why the cdd topology arises quite naturally for us see Remark~\ref{rema:whyCdd}. In this section we prove that if $F$ has a structure of traces, then there is a canonical induced structure of traces on $F_{\dis}$ (Theorem~\ref{theo:tracesForcd}) and moreover, this structure satisfies some particularly important properties (\ref{prop:sepImpliesTrin}).

In Section~\ref{sec:gersten} we introduce the concept of a \emph{Gersten presheaf} (Definition~\ref{defi:gersten}) which is a presheaf satisfying an analogue of the first part of the Gersten sequence in $K$-theory. We prove that if $F$ is a presheaf of $\zll$-modules with traces that satisfies the very first part of the Gersten sequence, then the cdh and $\ldh$ associated sheaves are isomorphic (Corollary~\ref{coro:cdhliso}). We deduce this in a convoluted way (see the diagram in the proof) from $F_{\ldh}$ being a subsheaf of $F_{\dis}$. We also use the Gersten property to find a criteria for a section of $F_{\dis}$ to belong to the image of $F_{\ldh}$ (and hence the image of $F_{cdh}$ since $F_{cdh} \cong F_{\ldh}$) and show that the trace morphisms of $F_{\dis}$ preserve this property. This implies that $F_{cdh}$ has a structure of traces (Proposition~\ref{prop:tracesForcdh}).

In Section~\ref{sec:tracesToTransfers} we prove Theorem~\ref{theo:cdhTriImpliesTransfers} which says that every cdh sheaf with traces that satisfies two additional properties has a canonical structure of transfers. We use the principle that every correspondence can be decomposed (locally for the cdh topology) into a formal sum of compositions of ``traces'', and morphisms of schemes (Lemma~\ref{lemm:platificationCorr}).

For the convenience of the reader let us make a small index here.

\begin{enumerate}
 \item[] (Definition~\ref{defi:topologies}) The $\fpsl$ and $\ldh$ topologies.
 \item[] (Definition~\ref{defi:traces}) Presheaf with traces, properties (Fon), (CdB), and (Deg).
 \item[] (Definition~\ref{defi:tri}) Properties (Tri1), (Tri2), (Tri1)$_{\leq n}$, (Tri2)$_{\leq n}$.
 \item[] (Definition~\ref{defi:refinable}) A refinable topology.
 \item[] (Definition~\ref{defi:cd}) The discrete topology.
 \item[] (Definition~\ref{defi:gersten}) Gersten presheaf.
 \item[] (Definition~\ref{defi:fn}) Correspondences of the form (FN).
\end{enumerate}

\begin{enumerate}
 \item[] \emph{Throughout this chapter we will state at the beginning of each section what class of schemes the results of that section hold for. In general, everything is true for the category of separated schemes essentially of finite type over a base scheme $S$ which is a quasi-excellent separated noetherian scheme. By essentially of finite type, we mean an inverse limit of a left filtering system of schemes of finite type, for which each of the transition morphisms is an affine open immersion.}
\end{enumerate}


\section{The \texorpdfstring{$\ldh$}{ldh} topology} \label{sec:ltopology}

In this section we present the definition of the $\ldh$ topology that we will use (Definition~\ref{defi:ltopology}). We state Gabber's theorem in some original versions (Theorem~\ref{theo:gabberLocal}, Theorem~\ref{theo:gabberGlobal}) and the version that we will use (Corollary~\ref{coro:regularlCover}).

Recall that if $\{ U_i \to X \}_{i \in I}$ is a finite family of morphisms, a refinement is a finite family of morphisms $\{ V_j \to X \}_{j \in J}$ such that for each $j \in J$ there is an $i_j \in I$ and a factorisation $V_j \to U_{i_j} \to X$. The reader not familiar with the cdh topology can find it in \cite{SV00}.

\begin{defi} \label{defi:topologies}
Let $\ell \in \ZZ$ be a prime.
\begin{enumerate}
 \item We will call an \emph{\fpsl cover} (fini-plat-surjectif-premier-{\`a}-$\ell$) a singleton set $\{f: U \to X\}$ containing a morphism $f$ that is finite flat surjective and globally free of degree prime to $\ell$. That is, $f_*\OO_U$ is a free $\OO_X$-module of rank prime to $\ell$.
 \item An \emph{$\ldh$ cover} is a finite family of morphisms of finite type $\{U_i \to X\}$ such that there exists a refinement of the form $\{ V_j' \to V_j \to X \}$ where $\{ V_j \to X \}$ is a cdh cover and $\{ V_j' \to V_j \}$ are $\fpsl$ covers.
\end{enumerate}
\end{defi}

\begin{rema}
Note that we can assume the $V_j, V_j'$ are affine as the Zariski topology is coarser than the cdh topology.
\end{rema}

For a pretopology $\tau$, we observe the usual abuse of terminology and refer to a morphism $Y \to X$ as a $\tau$ cover if $\{ Y \to X \}$ is a $\tau$ cover. The standard reference for the Nisnevich topology is \cite{Nis89} where it is referred to as the cd topology.

\begin{lemm} \label{lemm:NisFpslSwap}
Suppose that $Y' \to Y$ is a Nisnevich cover and $Y \to X$ a flat finite surjective morphism of constant degree (not necessarily globally free). Then there exists a Nisnevich cover $X' \to X$ such that $Y \times_X X' \to Y$ refines $Y' \to Y$, and $Y \times_X X' \to X'$ is an $\fpsl$ cover.
\end{lemm}

\begin{rema} \label{rema:triangleProblem}
This lemma is false if we replace the Nisnevich topology by the proper cdh topology. For example let $Y$ to be a non-normal curve, $Y \to X$ any flat finite morphism to a normal curve $X$, and $Y' \to Y$ the normalisation. Clearly $Y' \to Y$ doesn't split, but every proper cdh cover of $X$ is refinable by the trivial cover (this is true of any regular excellent scheme of dimension one). This failing is an obstacle to passing a structure of traces for a presheaf to its cdh sheafification (cf. Proposition~\ref{prop:tracesForNis}).
\end{rema}

\begin{proof}
If $X$ is henselian, then $Y$ is also henselian and $Y' \to Y$ splits. So we can take $X' = X$. If not then for every point $x \in X$ we consider the pullback along the henselisation $\ ^hx \to X$. The result now follows from the limit arguments in \cite[Section 8]{EGAIV3} and the description of the henselisation as a suitable limit of {\'e}tale neighbourhoods.
\end{proof}

The following proposition shows that the $\ldh$ covers as we have defined them form a pretopology in the sense of \cite[II.1.3]{SGA41}.

\begin{prop} \label{prop:cdhlswap}
Let $X$ be a noetherian scheme and suppose that $Y \to X$ is an $\fpsl$ morphism and $\{ U_i \to Y \}_{i \in I}$ is a cdh cover. Then there exists a cdh cover $\{ V_j \to X \}_{j \in J}$ and a set of $\fpsl$ morphisms $V'_j \to V_j$ such that $\{ V'_j \to X \}_{j \in J}$ refines $\{ U_i \to Y \to X \}_{i \in I}$.
\end{prop}

\begin{proof}
It suffices to consider the case when the cardinality of $I$ is one (replace $\{ U_i \to Y \}_{i \in I}$ by $\{ \amalg_{i \in I} U_i \to Y \}$). Recall that every cdh cover $U \to Y$ admits a refinement of the form $U'' \to U' \to Y$ where $U'' \to U'$ is a Nisnevich cover and $U' \to Y$ is a proper morphism which is a cdh cover (\cite[12.28]{MVW} or \cite[5.9]{SV00}). We have already treated the Nisnevich case in Lemma~\ref{lemm:NisFpslSwap} so it suffices to treat the proper cdh case.

We prove by noetherian induction that if we have $U \to Y \to X$ with $U \to Y$ proper cdh and $Y \to X$ $\fpsl$ then there exists $V' \to V \to X$ such that $V \to X$ is proper cdh, $V' \to V$ is $\fpsl$ and the composition $V' \to X$ factors through the composition $U \to Y$. Suppose that this statement is true for all proper closed subschemes of $X$. Indeed, by the inductive hypothesis, it suffices to prove that in the situation just mentioned we have the morphisms and properties just mentioned but with $V \to X$ proper and birational instead of proper cdh.

We can assume that $X$ is reduced, and even integral since the inclusion of the irreducible components is a proper birational morphism. Since $U \to Y$ is completely decomposed (i.e., the pullback along each $y \in Y$ admits a section), by replacing $U$ with an appropriate disjoint union of closed irreducible subschemes of $U$ we can assume that $U_{red} \to Y_{red}$ an isomorphism over a dense open subscheme of $Y$. If $\eta_i$ are the generic points of $Y$ and $m_i$ the lengths of their local rings, then the degree of $Y \to X$ is $\sum m_i [k(\eta_i): k(\xi)]$ where $\xi$ is the generic point of $X$. Since $\ell$ doesn't divide $\sum m_i [k(\eta_i): k(\xi)]$, there is some $i$ for which $\ell$ doesn't divide $[k(\eta_i): k(\xi)]$. By the platification theorem \cite{RG71} there exists a blowup of $X' \to X$ with nowhere dense centre such that the strict transform of $\overline{\{ \eta_i \}} \to X$ is flat, and hence finite flat surjective of degree prime to $\ell$. That is, we can assume $Y$ is reduced and even integral.

To recap, we have reduced to the case where $U$ is reduced, $U \to Y$ is an isomorphism over a dense open subscheme of $Y$, and $Y$ and $X$ are integral. In particular, the composition $U \to X$ is generically an $\fpsl$ cover. Using again the platification theorem, this time applied to the composition $U \to Y \to X$, we can find a blowup of $X$ with nowhere dense centre such that the strict transform of $U \to Y \to X$ is flat. Since it is generically an $\fpsl$ cover, flatness implies that it is actually an $\fpsl$ cover. So we are done.
\end{proof}

\begin{defi} \label{defi:ltopology}
The \emph{$\ldh$ pretopology} on a category of schemes is the pretopology for which the covers are $\ldh$ covers.
\end{defi}

\begin{rema}
Our choice of definition of an $\ldh$ topology is motivated by the following two ideas. Firstly, the theorem of Gabber (Theorem~\ref{theo:gabberLocal}) should provide the existence of regular $\ldh$ covers (or smooth depending on the context). Secondly, we want to make use of the vast literature available on the cdh topology. That is, we want to be able to reduce statements about the $\ldh$ topology, to statements about the cdh topology and statements about the $\fpsl$ topology. This way we only need to deal with the $\fpsl$ topology. This we usually do using a structure of traces -- cf. Lemma~\ref{lemm:fpslAcyclic}.

The name $\ldh$ is an acronym for $\ell$-decomposed $h$-topology (see \cite[Definition 3.1.2]{Voev96} fo the $h$-topology). For any set of primes $L$ one can define an $L$-decomposed morphism as a morphism $Y \to X$ such that for every point $x \in X$ there exists a point $y \in Y$ over $x$ such that no element of $L$ divides $[k(y):k(x)]$. We recover the notion of a completely decomposed morphism as a $\mathbb{P}$-decomposed morphism where $\mathbb{P}$ is the set of all primes. The cdh topology on a category of noetherian schemes is generated by the Nisnevich topology and covers which are proper and completely decomposed. Similarly, if we consider the pretopology generated by Nisnevich covers and proper $\{\ell\}$-decomposed morphisms, we obtain a pretopology which gives the same sheaves as our $\ldh$ topology. Hence, in some sense, the $\ldh$ topology is a legitimate generalisation of the cdh topology.

While we are discussing etymology, we mention the following counterexample. The na{\"i}ve reader may suspect that the cdh topology is equivalent to the topology obtained from the pretopology whose covers are $h$ covers that are completely decomposed. This is false. Let $k$ be a field, $Y = Spec(k[x, x^{-1}]) \amalg Spec(k[x^{\frac{1}{2}}])$ and $X = Spec(k[x])$ and let $Y \to X$ be the obvious morphism. This morphism is flat (and is therefore an $h$ cover), and is completely decomposed. However, the corresponding morphism of representable presheaves on $Sch(k)$ is not surjective for the cdh topology. If it were surjective, then there would exist a morphism $U \to Y$ in $Sch(k)$ such that the composition $U \to X$ is a cdh cover (the section $\id_X$ of $\hom(-, X)_{cdh}$ over $X$ would lift to $\hom(-, Y)_{cdh}$). Since $k[[x]]$ is a complete discrete valuation ring, every cdh cover of $Spec(k[[x]])$ admits a section. In particular, the canonical morphism $Spec(k[[x]]) \to X$ factors through $U \to X$, and pulling back $Y \to X$ along $ Spec(k[[x]]) \to X$, this implies that $Spec(k((x))) \amalg Spec(k[[x^{\frac{1}{2}}]]) \to Spec(k[[x]])$ admits a section, which is impossible.
\end{rema}

\begin{rema}
Our $\ldh$ pretopology differs from the topology of $\ell'$ alterations of \cite{ILO} (see the beginning of Section III.3) principally because the underlying categories are different -- they use the category denoted $alt/S$ (\cite[Definition 1.2.2]{ILO}). Their category $alt/S$ consists of reduced schemes $f: T \to S$ that are of finite type, surjective, and psuedo-dominant (Definition~\ref{defi:psuedodominant}) over $S$, and such that for every generic point $t$ of $T$ the extension $k(t)/k(f(t))$ is finite. Their topology of $\ell'$ alterations satisfies the following property (\cite[Theorem 3.2.1]{ILO}). If $X$ is irreducible and quasi-excellent, then every covering family for the topology of $\ell'$-alterations has a refinement of the form $\{ V_i \to Y \to X \}$ such that $Y$ is integral, $Y \to X$ is proper surjective of generic degree prime to $\ell$, and $\{V_i \to Y \}$ is a Nisnevich cover. Our pretopology is in some way a ``global'' version of their ``local'' pretopology where global and local are in the resolution of singularities sense.
\end{rema}

\begin{rema}
As with the cdh pretopology, we do not get an $\ldh$ pretopology on the category of smooth schemes $Sm(S)$ over some base $S$ as there are not enough fibre products. As with the cdh pretopology we do however get an induced topology.
\end{rema}

\begin{defi}
The $\ldh$ topology on $Sm(S)$ is the topology for which the covering sieves of a scheme $X$ are sieves $R \subseteq h_X$ that contain a sieve of the form $\im(h_{\amalg U_i} \to h_X)$ for some $\ldh$ cover $\{U_i \to X\}$.
\end{defi}

We now reproduce two versions of a theorem of Gabber. We follow them with a corollary which converts them into a form that we will use. For a statement and an outline of the proof of Gabber's Theorem of see \cite{Ill09}, or \cite{Gab05}. There is also a book in preparation \cite{ILO}.

\begin{theo}[{\cite[Theorem 2, Theorem 3.2.1]{ILO}}] \label{theo:gabberLocal}
Let $X$ be a noetherian quasi-excellent scheme, let $\ell$ be a prime number invertible on $X$. There exists a finite family of morphisms $\{U_i \to X \}_{i \in I}$ with each $U_i$ regular, and a refinement of the form $\{ V_j \to Y \to X \}_{j \in J}$ such that
\begin{enumerate}
 \item $\{V_j \to Y \}$ is a Nisnevich cover,
 \item $Y$ is locally integral,
 \item $Y \to X$ is proper and surjective, and
 \item for each generic point $\xi$ of $X$ there is a unique point $\eta$ of $Y$ over it, and $[k(\eta):k(\xi)]$ is finite of degree prime to $\ell$.
\end{enumerate}
\end{theo}

\begin{theo}[Gabber {\cite[1.3]{Ill09} or \cite[Theorem 3, Theorem 3.2.1]{ILO}}] \label{theo:gabberGlobal}
Let $X$ be a separated scheme of finite type over a perfect field $k$ and $\ell$ a prime distinct from the characteristic of $k$. There exists a smooth quasi-projective $k$ scheme $Y$, and a $k$-morphism $f: Y \to X$ such that
\begin{enumerate}
 \item $f$ is proper, surjective, pseudo-dominant (Definition~\ref{defi:psuedodominant}), and
 \item for each generic point $\xi$ of $X$ there is a unique point $\eta$ of $Y$ over it, and $[k(\eta):k(\xi)]$ is finite of degree prime to $\ell$.
\end{enumerate}
\end{theo}

\begin{coro} \label{coro:regularlCover}
Let $X$ be a scheme and $\ell$ a prime number invertible on $X$. If $X$ is noetherian and quasi-excellent then there exists an $\ldh$ cover $\{ U_i \to X \}$ of $X$ such that each $U_i$ is regular. If $X$ happens to be separated of finite type over a perfect field $k$, then there exists such a cover with each $U_i$ smooth and quasi-projective over $k$.
\end{coro}

\begin{proof}
The proof in both cases is the same so we give it only once. Let $X$ be noetherian and quasi-excellent. We proceed by noetherian induction. Suppose that the result is true for all proper closed subschemes of $X$. We can assume that $X$ is integral since the set of inclusions of irreducible components is a cdh cover. Let $\{U_i \to X \}_{i \in I}$ and $\{ V_j \to Y \to X \}_{j \in J}$ be as in the statement of Theorem~\ref{theo:gabberLocal} (or in the second case, just the $Y \to X$ from Theorem~\ref{theo:gabberGlobal}). We must show that the latter has a refinement which is a composition of $\fpsl$ and cdh covers.  By the platification theorem (Theorem~\ref{theo:platification}) there exists a blowup with nowhere dense centre $X' \to X$ such that the proper transform $Y' \to X'$ of $Y \to X$ is finite flat surjective morphism of constant degree (but not necessarily globally free). Let $\{ V'_i \to Y' \}$ be the pullback of the Nisnevich cover $\{ V_i \to Y \}$.
\[ \xymatrix{
X''' \ar[d]_{\textrm{\footnotesize \fpsl}} \ar[r] & \ar[r] V' \ar[rr] \ar[d]^{\textrm{\footnotesize Nis}} && V \ar[d]^{\textrm{\footnotesize Nis}} \\
X'' \ar[dr]_{\textrm{\footnotesize Nis}} & Y' \ar[rr] \ar[d]^{\textrm{\footnotesize fps }(\ell, \deg) = 1} && Y \ar[d]^{\textrm{\footnotesize prop. surj. gen. fin. } (\ell, \deg) = 1} \\
& X' \ar[rr]_{\textrm{\footnotesize blowup}} && X
} \]
By Lemma~\ref{lemm:NisFpslSwap} there exists a finite set of morphisms of the form $\{ X'''_j \to X''_j \to X' \}$ such that $\{X''_j \to X' \}$ is a Nisnevich cover and each $X_j''' \to X_j''$ is an $\fpsl$ cover, and furthermore, $\{ X'''_j \to X''_j \to X' \}$  is a refinement of $\{ V'_i \to Y' \to X' \}$. If $Z \in X$ is a closed subscheme such that $X' \to X$ is an isomorphism outside of $Z$, then $\{Z \to X, X''_j \to X' \to X \}$ is a cdh cover. By the inductive hypothesis, there exists an $\ldh$ cover $\{Z_k' \to Z\}_{k \in K}$ of $Z$ with each $Z_k'$ regular (or in the second case, quasi-projective and smooth over $k$). Hence, $\{Z_k' \to X\}_{k \in K} \cup \{ U_i \to X\}_{i \in I}$ is a finite family of morphisms with regular (resp. smooth quasi-projective) sources, such that there exists a refinement which is a composition of a cdh cover and $\fpsl$ covers as in the definition of an $\ldh$ cover.
\end{proof}


\section{Presheaves with traces} \label{sec:presheavesWithTraces}

In this section we present our definitions of a presheaf with traces (Definition~\ref{defi:traces}), and a presheaf with transfers (Definition~\ref{defi:presheafWithTransfers}). 
We mention that a presheaf with transfers is a presheaf with traces (Lemma~\ref{lemm:transfersImpliesTraces}).


\subsection{Presheaves with traces}

\begin{defi} \label{defi:traces}
A \emph{presheaf with traces} $(F, \S, \A, \Tr, \P)$ is an additive functor $F: \S^{op} \to \A$ from a category of schemes $\S$ to an additive category $\A$, together with a class $\P$ of morphisms of $\S$, and a morphism $\Tr_f: F(Y) \to F(X)$ for every morphism $f \in \P$. The morphisms $\Tr$ are required to satisfy the following axioms.
\begin{enumerate}
 \item[(Add)] For morphisms $f_1: Y_1 \to X_1$ and $f_2: Y_2 \to X_2$ in $\P$ we have
\[ \Tr_{f_1 \amalg f_2} = \Tr_{f_1} \oplus \Tr_{f_2}. \]
 \item[(Fon)] For morphisms $W \stackrel{g}{\to} Y \stackrel{f}{\to} X$ in $\P$ we have 
\[ \Tr_f \Tr_g = \Tr_{fg} \textrm{ and } \Tr_{\id_X} = \id_{F(X)}. \]
 \item[(CdB)] For every cartesian square in $\S$
\begin{equation} \label{equa:cartSquareDefiTraces}
\xymatrix{
Y \times_X W \ar[r]^-g \ar[d]_q & W \ar[d]^p \\
Y  \ar[r]_f & X
}
\end{equation}
such that $f, g \in \P$ we have
\[ F(p)\Tr_f = \Tr_gF(q). \]
 \item[(Deg)] For every finite flat surjective morphism $f: Y \to X$ in $\P$ such that $f_*\OO_Y$ is a globally free $\OO_X$ module we have 
\[ \Tr_f F(f) = \deg f \cdot \id_{F(X)}. \]
\end{enumerate}
Sometimes we will just denote a presheaf with traces by $F$ if the rest of the data is already established. In this thesis $\P$ will always be the class of finite flat surjective morphisms in $\S$. We will denote this class by $\S^\fps$.

A \emph{morphism of presheaves with traces} $(F, \S, \A, \Tr, \P) \to (G, \S, \A, \Tr, \P)$ is a morphism of the underlying presheaves $F \to G$ such that for every $f \in \P$ the square
\[ \xymatrix{
F(Y) \ar[r] \ar[d]_{\Tr_f} & G(Y) \ar[d]^{\Tr_f} \\
F(X) \ar[r] & G(X)
} \]
is commutative.
\end{defi}


\begin{exam} \label{exam:presheavesWithTraces}
Here are some examples of presheaves with traces.
\begin{enumerate}
 \item Suppose $\S$ is any category of schemes, $\A$ is any additive category and $F$ is any constant additive sheaf. Then since every finite flat surjective morphism Zariski locally satisfies the hypotheses of (Deg), there is a unique structure $\Tr$ such that $(F, \S, \A, \Tr, \S^\fps)$ is a presheaf with traces. This is because if $f: Y \to X$ is a morphism satisfying the hypotheses of (Deg) between connected schemes then $F(f)$ is an isomorphism, and so (Deg) requires that $\Tr_f = \deg f \cdot F(f)^{-1}$. The other axioms are straight-forward.

 \item The presheaves $\OO_X^*$ and $\OO_X$ (represented by the group schemes $\GG_m$ and $\AA^1$) have canonical structures of traces induced by the determinant and trace of matrices. More explicitly if $Spec\ B \to Spec\ A$ is a morphism of affine schemes and there is an isomorphism of $A$ algebras $B \cong \oplus_{i = 1}^d A$ then there is an induced morphism $B \to M_d(A)$ of $B$ into the ring of $d$ by $d$ matrices with coefficients in $A$ (induced by right or left multiplication of $B$ on itself). Then the determinant and trace define group homomorphisms $(B^*, *) \to (A^*, *)$ and $(B, +) \to (A, +)$. It can be checked that these morphisms are independent of the chosen isomorphism $B \cong \oplus_{i = 1}^d A$ and glue to give a structure of traces on non-affine schemes.

 \item The example described above is a special case of a more general phenomena. On the category of quasi-projective normal schemes, any presheaf represented by an algebraic group has transfers, and hence a structure of traces (any presheaf with transfers has a structure of traces; this is mentioned further down the list). 


 \item We might like to say that $\fpsl$ sheaves have traces using a similar tactic to \cite[Section 5]{SV96} to define traces using pseudo-Galois covers. However, when passing to a normal extension, we lose control of the degree and cannot ensure it stays prime to $\ell$. The converse is true: a presheaf of $\zll$ modules with traces is an $\fpsl$ sheaf (Lemma~\ref{lemm:fpslAcyclic}).

 \item Algebraic $K$-theory and homotopy invariant algebraic $K$-theory have structures of traces due to the constructions being functorial with respect to biWaldhausan categories (cf. the proof of Proposition~\ref{prop:ktheoryTraces}).

 \item We will see that any presheaf with transfers (Definition~\ref{defi:presheafWithTransfers}) has a canonical structure of presheaf with traces (Lemma~\ref{lemm:transfersImpliesTraces}).

 \item If $F$ is a presheaf with traces, its Nisnevich and {\'e}tale sheafifications have a canonical structure of traces (Proposition~\ref{prop:tracesForNis}, Lemma~\ref{lemm:NisFpslSwap}). We will also see conditions on $F$ under which the discrete sheafification (Definition~\ref{defi:cd}) has a canonical structure of traces, the associated cdh separated presheaf has a canonical structure of traces (Proposition~\ref{prop:tauSepHasTraces}), and the associated cdh and $\ldh$ sheaves have a structure of traces (Proposition~\ref{prop:tracesForcdh}, Corollary~\ref{coro:cdhliso}). The latter is one of the main results of this paper.

 \item Let $S$ be a noetherian scheme, $Sch(S)$ the category of $S$-schemes of finite type, $\A$ an additive category with small colimits, $EssSch(S)$ the category of schemes essentially of finite type. We remind the reader that when we say \emph{essentially}, we are talking about limits of left filtering systems in which the transition morphisms are affine open immersions. It is a standard application of the results in \cite[Section 8]{EGAIV3} that if $F: Sch(S) \to \A$ is a presheaf with traces then $F$ gives rise canonically to a presheaf with traces on $EssSch(S)$.

\end{enumerate}
\end{exam}

We eventually want to find a criteria for when a structure of traces on a presheaf induces a structure of traces on the cdh sheafification (this is achieved in Proposition~\ref{prop:tracesForcdh}). The following proposition, applicable in the case $\tau = $ Nisnevich, is a first step in this direction.

\begin{prop} \label{prop:tracesForNis}
Let $(F, \S, \A, \Tr, \P)$ be a presheaf with traces and suppose that $\A$ is an abelian category and the class $\P$ is closed under fibre products. That is, if $f: Y \to X \in \P$ then so is $W \times_X f$ for any $W \to X$ in $\S$. Now suppose that $\tau$ is a pretopology on $\S$ such that 
\begin{enumerate}
 \item[] for every morphism $f: Y \to X \in \P$ and every $\tau$ cover $V \to Y$ there exists a $\tau$ cover $U \to X$ such that $Y \times_X U \to Y$ is a refinement of $V \to Y$.
\end{enumerate}
Then there is a unique class of morphisms $\Tr^{\tau}$ such that $(F_\tau, \S, \A, \Tr^\tau, \P)$ is a presheaf with traces and such that the canonical morphism $F \to F_\tau$ is a morphism of presheaves with traces.
\end{prop}

\begin{proof}
Let $f: Y \to X$ be a morphism in $\P$, $U \to X$ a $\tau$ cover of $X$, and $s \in \ker(F(Y \times_X U) \to F((Y \times_X U) \times_Y (Y \times_X U))$. We claim that there is a unique element $t \in F_\tau(X)$ such that the restriction to $F_\tau(U)$ agrees with the image of $\Tr_{(f \times_X U)}s$. Indeed, if follows immediately from (CdB) and the isomorphism $(Y \times_X U) \times_Y (Y \times_X U) \cong Y \times_X (U \times_X U)$ that $\Tr_{(f \times_X U)}s$ is a cocycle and so it descends to a unique element of $F_\tau(X)$.

Now for every element $s$ of $F_\tau(Y)$ there exists a $\tau$ cover $V \to Y$ so that $s|_V$ is in the image of $F \to F_\tau$. By our hypothesis, we can assume $V$ is of the form $Y \times_X U \to Y$ for some $\tau$ cover $U \to X$. By what we have just shown, we have a corresponding element in $F_\tau(X)$ that is independent of the choice of $U$. Hence a morphism $\Tr^\tau_f: F_\tau(Y) \to F_\tau(X)$.

The axioms (Fon) and (Deg) follow immediately from the way we have defined the morphisms $\Tr^\tau$. It is also immediate from the definition that these are compatible with $F \to F_\tau$, and are the only possible such choice. For (CdB) it is enough to draw the appropriate cube and do the diagram chase.
\end{proof}

We have cause to discuss two further properties that might be satisfied by a presheaf with traces. In the case of a cdh sheaf, these two properties bridge the gap between a structure of traces and a structure of transfers (cf. Lemma~\ref{lemm:transfersImpliesTraces}, Theorem~\ref{theo:cdhTriImpliesTransfers}). They deal with commutative triangles:

\begin{equation} \label{equa:tri}
 \xymatrix{
Y' \ar[dr]_g \ar[rr]^h && Y \ar[dl]^f \\
& X
}
\end{equation}

\begin{defi} \label{defi:tri}
Suppose that we have a commutative triangle (\ref{equa:tri}) as above and $(F, \S, \A, \Tr, \P)$ a presheaf with traces. We define the following two properties.
\begin{enumerate}
 \item[(Tri1)$_{\leq d}$] Suppose that in the commutative triangle (\ref{equa:tri}) the scheme $X$ is integral of dimension $\leq d$, the scheme $Y'$ is the disjoint union of the integral components $Y_i'$ of $Y$, and $h$ is the canonical morphism \mbox{$Y' = \amalg Y_i' \to Y$}, and the morphisms $f, g_i$ are in $\P$ where $h_i, g_i$ are the restrictions of $h, g$ to $Y_i'$. Then 
\[ \Tr_f = \sum m_i \Tr_{g_{i}}F(h_i) \]
where $m_i = \length \OO_{Y, \eta_i}$ with $\eta_i$ the generic point of $Y_i'$.
 \item[(Tri2)$_{\leq d}$] Forgetting the hypotheses of (Tri1)$_{\leq d}$, suppose that in the commutative triangle (\ref{equa:tri}) the morphisms $f$ and $g$ are in $\P$, and all the schemes $X, Y, Y'$ are integral of dimension $\leq d$. Then 
\[ \tfrac{\deg g}{\deg f} \Tr_f = \Tr_gF(h). \]
 \item[] We will use just (Tri1) and (Tri2) if we require these axioms without restriction on the dimension.
\end{enumerate}
\end{defi}

\begin{rema} \label{rema:otherPresheavesWithTraces}
\begin{enumerate}
 \item We will almost always only ask for (Tri1)$_{\leq 0}$. This is because we will end up using presheaves that are separated for the cdd topology (Definition~\ref{defi:cd}) and for such presheaves (Tri2)$_{\leq n}$ is true for all $n$ and (Tri1)$_{\leq n}$ for all $n$ is implied by (Tri1)$_{\leq 0}$ (Proposition~\ref{prop:sepImpliesTrin}). We will soon give a criteria under which (Tri1)$_{\leq 0}$ is satisfied (Lemma~\ref{lemm:Tr2dimZero}). 

 \item We will see below that if $F$ is a presheaf with transfers then $F$ is a presheaf with traces that satisfies (Tri1) and (Tri2) (Lemma~\ref{lemm:transfersImpliesTraces}). We will prove that conversely if $F$ is a cdh sheaf with traces that satisfies (Tri1) and (Tri2) then $F$ has a canonical structure of presheaf with transfers (Theorem~\ref{theo:cdhTriImpliesTransfers}). 




 \item The morphisms on algebraic $K$-theory described in Example~\ref{exam:presheavesWithTraces} do not satisfy (Tri1) before we sheafify it. This is for a similar reason to the fact that algebraic $K$-theory does not have transfers \cite[Section 3.4]{Vctpt}. For example, let $X$ be a projective line and choose a closed point $x$. Let $Y'$ be two disjoint copies of $X$, and suppose that $Y$ has two irreducible components, each isomorphic to $X$, and the intersection of these two irreducible components is the chosen point $x$. The morphisms $f, g, h$ are the obvious ones. Both $f$ and $g$ are finite flat and surjective, and $h$ is the inclusion of the integral components. However the class of $f_*\OO_X$ in $K_0(X)$ is different from that of $g_*\OO_Y$. This is the only example we know of a presheaf with traces that doesn't satisfy (Tri1). If we were to require (Deg) to hold for all finite flat surjective morphisms of constant degree (and not just globally free ones) we would lose this counter-example.

\end{enumerate}
\end{rema}

\begin{lemm} \label{lemm:Tr2dimZero}
Suppose $(F, \S, \A, \Tr, \S^\fps)$ is a presheaf with traces.
\begin{enumerate}
 \item (Tri2)$_{\leq 0}$ is always satisfied.
 \item Suppose for every finite morphism of schemes of dimension zero $Y \to X$, if $X \in \S$ then $Y \in \S$. The axiom (Tri1)$_{\leq 0}$ is satisfied if $F(Y) \to F(Y_{red})$ is an isomorphism for every $Y$, and for every field $k$ the exponential characteristic of $k$ is invertible in $F(Spec(k))$ .
\end{enumerate}
\end{lemm}

\begin{proof}
\begin{enumerate}
 \item We have $\Tr_gF(h) \stackrel{(Fon)}{=} \Tr_f \Tr_hF(h) \stackrel{(Deg)}{=} \deg h \cdot \Tr_f = \tfrac{\deg g}{\deg f} \Tr_f$.

 \item Consider a triangle (\ref{equa:tri}) of schemes of dimension zero with the hypotheses of (Tri1)$_{\leq 0}$. We can and do assume that $Y$ is connected, and so our schemes are of the form $Y = Spec(A)$, $X = Spec(k)$, and $Y' = Spec(L)$ where $k$ is a field, $A$ is a finite local $k$-algebra, and $L$ is the residue field of $A$.

Suppose for the moment that $L / k$ is purely inseparable. Since $p = \mathrm{char}(k)$ is invertible in $F(k)$ and $[L : k]$ is a power of $p$, the axiom (Deg) implies that $F(k) \to F(L)$ is injective. Consider the pull-back along $Spec(L) \to Spec(k)$ and the resulting diagram
\[ \xymatrix{
Spec(L \otimes_k L) \ar[rr]^h \ar[dr]_g && Spec(A \otimes_k L) \ar[dl]^f \\
& Spec(L)
} \]
By Lemma~\ref{lemm:degLength} applied to the original triangle and (CdB) applied to the pullback squares along $Spec(L) \to Spec(k)$ it suffices to prove that $\Tr_f = \tfrac{\deg f}{\deg g} \Tr_g F(h)$ due to the injectivity of $F(k) \to F(L)$.

Since purely inseparable extensions are universal monomorphisms \cite[3.5]{EGAI} the schemes in this diagram all have a unique point. Moreover, $g$ now admits a section $s$ and since $L = (A \otimes_k L)_{red} = (L \otimes_k L)_{red}$ the morphism $F(f)$ (resp. $F(g)$) is an isomorphism (by hypothesis) with inverse $F(hs)$ (resp. $F(s)$). We have
\[ \Tr_f = \Tr_f F(f) F(hs) \stackrel{(Deg)}{=} \deg f \cdot F(h s) = \deg f \cdot F(s) F(h)  \]
\[ \stackrel{(Deg)}{=} \tfrac{\deg f}{\deg g} \Tr_g F(g) F(s) F(h) = \tfrac{\deg f}{\deg g} \Tr_g F(h). \]

Now we remove the assumption that $L/k$ is purely inseparable. Let $k \subset K \subset L$ be a maximal separable subextension so that $K / k$ is separable and $L / K$ is purely inseparable. Let $B$ be the preimage of $K$ under the canonical $A \to L$. By Cohen's Structure Theorem for complete local rings \cite[28.J]{Mat}, the morphism $B \to K$ admits a section which is a $k$-morphism. Consequently, we have a commutative diagram
\[ \xymatrix@!=6pt{
Spec(L) \ar[rrrr] \ar[dr] &&&& Spec(A) \ar[ddll] \ar[dlll]\\
& Spec(K) \ar[dr] &&&& \\
&& Spec(k) 
} \]
and the result follows from the purely inseparable case.
\end{enumerate}
\end{proof}




\subsection{Presheaves with transfers}

\begin{defi} \label{defi:presheafWithTransfers}
Suppose $S$ is a noetherian scheme. A \emph{presheaf with transfers} is an additive presheaf on $Cor(S)$ (Definition~\ref{defi:categoryOfCorrespondances}). The category of presheaves with transfers is denoted $PreShv(Cor(S))$, and if $\Lambda$ is a ring, then the category of presheaves of $\Lambda$-modules with transfers is denoted $PreShv(Cor(S), \Lambda)$.

If $\tau$ is a Grothendieck topology on $Sch(S)$ then a \emph{$\tau$ sheaf with transfers} is a presheaf with transfers whose restriction to $Sch(S)$ is a $\tau$ sheaf. We have corresponding categories $Shv_\tau(Cor(S))$ and $Shv_\tau(Cor(S), \Lambda)$.

If $\S$ is a class of schemes in $Sch(S)$ with corresponding full subcategory $\C$ in $Cor(S)$, we make the analogous definitions of a presheaf with transfers on $\S$, presheaf of $\Lambda$-modules on $\S$, $\tau$ sheaf with transfers on $\S$, $\tau$ sheaf of $\Lambda$-modules on $\S$, with corresponding categories $PreShv(\C)$, $PreShv(\C, \Lambda)$, $Shv_\tau(\C)$, and $Shv_\tau(\C, \Lambda)$.
\end{defi}

\begin{defi} \label{defi:representableSheafWithTransfers}
For any object $[X] \in Cor(S)$ we denote the corresponding \emph{representable presheaf} with transfers by $L(X)$. The cdh (resp. Nisnevich) sheafification of $L(X)$ has a canonical structures of transfers (Theorem~\ref{theo:CD}, resp. \cite[Lemma 3.1.6]{Voev00}) and we denote this sheafification by $L_{cdh}(X)$ (resp. $L_{Nis}(X)$).
\end{defi}

\begin{lemm} \label{lemm:transfersImpliesTraces}
Every presheaf with transfers is a presheaf with traces that satisfies (Tri1) and (Tri2).
\end{lemm}

\begin{proof}
A direct consequence of Proposition~\ref{prop:corHasTraces}.
\end{proof}

\begin{rema}
Notice that since the presheaf with traces induced by a presheaf with transfers necessarily satisfies (Tri1) and (Tri2), any presheaf with traces whose structure is extendible to a structure of transfers necessarily satisfies these two properties.
\end{rema}


\section{Comparison of cdh and \texorpdfstring{$\ldh$}{ldh} cohomology} \label{sec:comparisonCohomology}

The goal of this subsection is Theorem~\ref{theo:lcdhcohomologyAgree} which says that the $\ldh$ and cdh cohomology of a cdh sheaf with transfers agree (Theorem~\ref{theo:lcdhcohomologyAgree}). We pass by an equivalence of the categories of cdh sheaves of $\zll$ modules with transfers and $\ldh$ sheaves of $\zll$ modules with transfers (Corollary~\ref{coro:equivalencecdhltransfers}).

\subsection{\v{C}ech cohomology and refinable topologies} \label{sec:cechRefine}

Most of the subsection is devoted to building a proof for Proposition~\ref{prop:sigmaTauAcyclic} which will be used to compare the cdh and $\ldh$ cohomology via Proposition~\ref{prop:318l}. In the Proposition~\ref{prop:compatSheaf} (which is independent from the rest of the subsection) we note some easily proved facts that we will need later.

\begin{enumerate}
 \item[] \emph{In this subsection we work with an essentially small category $\C$ which we will assume to be equipped with fibre products.}
\end{enumerate}

Our interest in \v{C}ech cohomology stems from the following well know lemma.

We say that a presheaf $F$ is \emph{acyclic} for a topology $\tau$ if $F = F_\tau$ and $H_\tau^n(-, F) = 0$ for all $n > 0$.

\begin{lemm}[{\cite[V.4.3]{SGA42} or \cite[III.2.11]{Mil80}}] \label{lemm:acyclicIsCechAcyclic}
A presheaf $F$ on $\C$ is acyclic for a topology $\tau$ if and only if its \v{C}ech cohomology groups vanish.
\end{lemm}

The following is another well-known lemma.

\begin{lemm}[{\cite[Exp. V 2.3.5]{SGA42}, \cite[III.2.1]{Mil80}, or \cite[1.4.3]{Art62}}] \label{lemm:homotopyGeneral}
Let $V / X, U/X$ be two $X$-objects in $\C$. Suppose that $F$ is a presheaf on $\C$. Then any two $X$ morphisms $V \rightrightarrows U$ induce the same morphism $\check{H}^n(U / X, F) \to \check{H}^n(V / X, F)$. Consequently, any $X$-morphism $V \to U$ that admits an $X$-section $V \leftarrow U$ induces isomorphisms $\check{H}^n(U / X, F) \stackrel{\sim}{\to} \check{H}^n(V / X, F)$.
\end{lemm}

\begin{proof}
Suppose that $f_0, f_1: V \rightrightarrows U$ are the morphisms. We construct a simplicial homotopy\footnote{See \cite[Exp Vbis 3.0.2]{SGA42} for the definition of a simplicial homotopy that we use.} between $cosk_0(f_0)$ and $cosk_0(f_1)$ in $\C$. This induces a homotopy between the associated morphisms of chain complexes (see \cite[Exp Vbis 3.0.2.3]{SGA42}). For each $\phi: [n] \to [1]$ we define $V^{\times_X n} \to U^{\times_X n}$ to be the map whose $i$th component is $f_{\phi(i)}$. It is easily checked that this is a homotopy $\Delta[1] \times cosk_0(f_0) \to cosk_0(f_1)$ between simplicial objects. Hence the associated morphisms of chain complexes are homotopic, and therefore induce the same morphisms on cohomology.

For the second statement, let $f: V \to U$ and $s: U \to V$ be $X$-morphisms such that $fs = \id_U$. By functoriality, $\check{H}^n(f / X, F)$ is a left inverse to $\check{H}^n(s / X, F)$ for every $n \geq 0$, and so it suffices to show that $\check{H}^n(s / X, F)$ is a left inverse to $\check{H}^n(f / X, F)$ for every $n \geq 0$. That is, we wish to see that $\check{H}^n(sf / X, F)$ is the identity for every $n \geq 0$. Applying the previous result to the two morphisms $sf, \id_V: V \rightrightarrows V$ shows this.
\end{proof}

\begin{lemm} \label{lemm:SS}
Let $V \stackrel{g}{\to} U \stackrel{f}{\to} X$ be a pair of composable morphisms. Then there exists a bisimplicial object $W_{p, q}$ such that the $p$th column is $cosk_0(V^{\times_X (p + 1)} \to U^{\times_X (p + 1)})$ and the $q$th row is $cosk_0(V^{\times_U (q + 1)} \to X)$. Notably, for any presheaf $F$ we get a first quadrant spectral sequence
\[ E^{p, q}_1 = \check{H}^{q}(V^{\times_X (p + 1)} / U^{\times_X (p + 1)}, F) \implies \check{H}^{p + q}(V / X, F) \]
\end{lemm}

\begin{proof}
We define the following objects.
\[ W_{(p - 1),(q - 1)} \stackrel{def}{=} \underbrace{V^{\times_X p} \underset{(U^{\times_X p})}{\times} \cdots \underset{(U^{\times_X p})}{\times} V^{\times_X p}}_{q \textrm{ times }} = \underbrace{V^{\times_U q} \underset{X}{\times} \cdots \underset{X}{\times} V^{\times_U q}}_{p \textrm{ times }}. \]
The object $W_{(p - 1),(q - 1)}$ is also the limit of a diagram that has $p \times q$ copies of $V$ with $p$ copies of $U$, and an $X$, and the $(i, j)$th $V$ has a morphism towards the $i$th $U$ and every $U$ and $V$ has a morphism towards $X$ (the morphisms being either $g, f$ or $fg$). Presented this way, there are obvious face $W_{p, q} \to W_{p + 1, q}, W_{p, q} \to W_{p, q + 1}$ and degeneracy $W_{p, q} \to W_{p - 1, q}, W_{p, q} \to W_{p, q - 1}$ morphisms coming from the projections and diagonals (in the $p$th and $q$th directions) and these are compatible in the sense that we get a bisimplicial object.


We consider the double complex associated to the bicosimplicial abelian group obtained by applying $F$ to the $W_{p, q}$. There are two associated spectral sequences, one from the filtration of the total complex by rows and one from the filtration of the total complex by columns. We start with the filtration for which the $E_0$ differentials are in the $p$ direction. The $E_1$ terms are $\check{H}^{p}(V^{\times_U {(q + 1)}} / X, F)$ and the $E_1$ differentials are induced by the differentials in the $q$ direction. Now every face morphism $\partial_i$ of $cosk_0(V/U)$ has a section, namely the degeneracy $\sigma_i$. So by Lemma~\ref{lemm:homotopyGeneral} the morphisms induced by these face morphisms $H^p(\partial_i): H^p(cosk_0(V^{\times_U {q - 1}} \to X)) \to H^p(cosk_0(V^{\times_U q} \to X))$ are all the same isomorphism. Hence, the differentials $d = \sum_{i = 0}^q (-1)^i H^p(\partial_i)$ are zero if $q$ is odd and an isomorphism if $q$ is even. Consequently, on the $E_2$ sheet everything is zero except the bottom row $E_2^{0, q} \cong H^p(cosk_0(V^{\times_U 1} \to X))$ and so we see that the cohomology of the total complex is $H^p(cosk_0(V^{\times_U {1}} \to X)) = \check{H}^{p}(V / X, F)$. Considering the other filtration gives the $E_1$ terms in the statement of the result.
\end{proof}

\begin{defi} \label{defi:refinable}
Let $\tau$ be a Grothendieck pretopology and $\rho, \sigma$ two classes of $\tau$ covers. We say that $\tau$ is \emph{$\stackrel{\rho}{\to}\stackrel{\sigma}{\to}$ refinable} if every $\tau$ cover admits a refinement of the form $\{ V_{ij} \stackrel{r_{ij}}{\to} U_i \stackrel{s_i}{\to} X \}$ such that $\{V_{ij} \to U_i\} \in \rho$ and $\{U_i \to X\} \in \sigma$.
\end{defi}

\begin{rema}
Suppose that $\tau$ is a Grothendieck pretopology such that for every cover $\{W_i \to X\}_{i \in I}$ the morphism  $\amalg_{i \in I} W_i \to X$ exists in the category that we are working with. Then to show that $\tau$ is $\stackrel{\rho}{\to}\stackrel{\sigma}{\to}$ refinable it is enough to consider $\tau$ covers that contain a single morphism, since $\{W_i \to X\}_{i \in I}$ is refinable if and only if $\{ \amalg_{i \in I} W_i \to X \}$ is.
\end{rema}

\begin{exam} \label{exam:refinable}
\begin{enumerate}
 \item In Lemma~\ref{lemm:NisFpslSwap} we have seen that the pretopology generated by Nis and $\fpsl$ is $\stackrel{\fpsl}{\to}\stackrel{Nis}{\to}$ refinable.
 \item By definition the $\ldh$ pretopology is $\stackrel{\fpsl}{\to}\stackrel{cdh}{\to}$ refinable. Since Zariski covers are cdh covers, we can even restrict to the class of $\fpsl$ covers of affine schemes.
 \item The cdh pretopology is $\stackrel{Nis}{\to}\stackrel{cdp}{\to}$ refinable where cdp is the class of cdh covers $\{U_i \to X \}$ such that each morphism $U_i \to X$ is proper \cite[12.28]{MVW} or \cite[5.9]{SV00} (the proof they give works over any noetherian base scheme).
 \item In ``Homology of schemes I'' \cite{Voev96} Voevodsky shows that over noetherian excellent schemes the h pretopology is $\stackrel{Zar}{\to}\stackrel{ps}{\to}$ refinable (and therefore $\stackrel{et}{\to}\stackrel{ps}{\to}$ refinable) where ps is the class of proper, surjective morphisms. He also shows that the qfh pretopology is $\stackrel{fs}{\to}\stackrel{et}{\to}$ refinable where fs is the class of finite surjective morphisms. These facts are the basis for the comparison results in \cite[Section 10]{SV96}, and explicitly recognising this makes the proofs clearer.
 \item The eh pretopology of \cite{Gei06} is $\stackrel{et}{\to}\stackrel{cdp}{\to}$ refinable (it is the same proof as for the cdh pretopology mentioned above).
\end{enumerate}
\end{exam}

\begin{prop} \label{prop:sigmaTauAcyclic}
Suppose that $\tau$ is a Grothendieck pretopology, that $\rho, \sigma$ are Grothendieck subpretopologies of $\tau$, and that $\tau$ is $\stackrel{\rho}{\to}\stackrel{\sigma}{\to}$ refinable. Let $F$ be a presheaf on $\C$ such that
\[ \check{H}^n(U/X, F) = \left \{ \begin{array}{ll} 0 & n > 0 \\ F(X) & n = 0 \end{array} \right . \]
for every $\rho$ cover $U \to X$. Then if $F$ is $\sigma$ acyclic, it is also $\tau$ acyclic.
\end{prop}

\begin{proof}
Because vanishing of \v{C}ech cohomology in degrees $n > 0$ is equivalent to a presheaf being acyclic (Lemma~\ref{lemm:acyclicIsCechAcyclic}), it is sufficient to show that $\check{H}^n_\sigma(X, F) = \check{H}^n_\tau(X, F)$. To calculate the $\tau$ \v{C}ech cohomology we can restrict to covers of the form $V \stackrel{g}{\to} U \stackrel{f}{\to} X$ with $f$ a $\sigma$ cover and $g$ a $\rho$ cover. By our hypothesis, $\check{H}^{q}(V^{\times_X p} / U^{\times_X p}, F) = 0$ for $q > 0$ and $\check{H}^{0}(V^{\times_X p} / U^{\times_X p}, F) = F(U^{\times_X p})$. Hence the spectral sequence of Lemma~\ref{lemm:SS} collapses to give the isomorphism $\check{H}^n(U / X, F) = \check{H}^n(V / X, F)$. Passing to the limit over covers of the form $V \to U \to X$ gives the result.
\end{proof}

The final proposition of this subsection is independent of the rest of this subsection. It collects some elementary properties of refinable pretopologies that we will need later.

\begin{prop} \label{prop:compatSheaf}
Suppose that $\tau$ is a Grothendieck pretopology, that $\rho, \sigma$ are Grothendieck subpretopologies of $\tau$, and that $\tau$ is $\stackrel{\rho}{\to}\stackrel{\sigma}{\to}$ refinable.
\begin{enumerate}
 \item If $F$ is $\rho$ separated then $F_\sigma$ is $\rho$ separated (and hence $\tau$ separated).
 \item If $F$ is a $\rho$ sheaf then $\check{H}^0_\tau(-, F) = \check{H}^0_\sigma(-, F)$.
 \item If $F$ is a $\rho$ sheaf that is $\sigma$ separated then $F_\sigma = F_\tau$. In particular, if $F$ is a $\rho$ sheaf and a $\sigma$ sheaf, then it is also a $\tau$ sheaf.
\end{enumerate}
\end{prop}

\begin{proof}
\begin{enumerate}
 \item We have to show that $F_\sigma \to (F_\sigma)_\rho$ is a monomorphism. That is, for every section $s \in F_\sigma(X)$ sent to zero in $(F_\sigma)_\rho(X)$, we want to show that $s$ is zero. For every section $s \in F_\sigma(X)$ there is a $\sigma$ cover $U \to X$ such that $s|_U$ is in the image of $F \to F_\sigma$, and so it is enough to consider elements in the image of $F \to F_\sigma$. It is clear enough that an element $s \in F(X)$ is sent to zero in $(F_\sigma)_\rho$ if and only if there exists a $\rho$ cover $U \to X$ and a $\sigma$ cover $V \to U$ such that $s|_V = 0$. But $F$ is $\rho$ separated, and so $s|_U = 0$. Since $U \to X$ is a $\sigma$ cover, this implies that $s$ is zero in $F_\sigma(X)$.

 \item By our hypothesis, the $\tau$ \v{C}ech cohomology can be calculated using covers of the form $\{ V_i \stackrel{g_i}{\to} U_i \stackrel{f_i}{\to} X \}$ such that $f_i \in \sigma$ and $g_i \in \rho$. For simplicity we assume that each family has a single element. We have the following morphism of exact sequences
\[ \xymatrix{
0 \ar[r] & \check{H}^0(V / X, F) \ar[r] \ar[d]^\alpha & F(V) \ar[r] \ar[d] & F(V \times_X V) \ar[d] \\
0 \ar[r] & F(U) \ar[r] & F(V) \ar[r] & F(V \times_U V)
} \]
and the morphism $\alpha$ can be inserted into the following morphism of short exact sequences.
\[ \xymatrix{
0 \ar[r] & \check{H}^0(U / X, F) \ar[r] \ar[d]_\beta & F(U) \ar[r] \ar[d] & F(U \times_X U) \ar[d] \\
0 \ar[r] & \check{H}^0(V / X, F) \ar[r] \ar[ur]^\alpha & F(V) \ar[r] & F(V \times_X V)
} \]
Since $F$ is $\rho$ separated, all vertical morphisms are monomorphims, and consequently, the diagram is commutative, and $\alpha$ lifts to give an inverse to $\beta$. Taking the limit over covers of this form gives the result.

 \item This follows immediately from the previous part since for separated presheaves, the zeroth \v{C}ech cohomology calculates the sheafification.
\end{enumerate}
\end{proof}

\begin{rema}
In the third part, we suspect that the assumption that $F$ is $\sigma$ is separated is necessary if we want the result in this level of generality. If this necessity were not the case, this chapter would be considerably shorter.
\end{rema}



\subsection{Comparison of cdh and \texorpdfstring{$\ldh$}{ldh} cohomology}

In this subsection as in Chapter~\ref{chap:cycles} we work with $Sch(S)$ the category of separated schemes of finite type over a separated noetherian base scheme $S$. Lemma~\ref{lemm:fpslAcyclic} is true in any category of schemes in which the $U \times_X \dots \times_X U$ exist.

\begin{lemm} \label{lemm:fpslAcyclic}
Suppose that $F$ is a presheaf of $\zll$ modules with traces. Then for any $\fpsl$ morphism $U \to X$ of degree $d$ the sequence
\[ 0 \to F(X) \to F(U) \to F(U \times_X U) \to \dots \]
is exact. In particular, $F$ is an $\fpsl$ sheaf, $F$ is $\fpsl$ acyclic, and if $F$ happened to be a cdh sheaf of $\zll$ modules with traces, then it is a sheaf for the $\ldh$ pretopology.
\end{lemm}

\begin{rema}
The analogous result is true for any $\fpsl$ hypercover using (almost) the same proof. This is due to the fact that in any $\tau$ hypercover $\mathcal{U}_\bullet \to X$ the face morphisms $d_i : \mathcal{U}_n \to \mathcal{U}_{n - 1}$ are $\tau$-covers as well as the canonical morphisms $\mathcal{U}_{n + 1} \to \mathcal{U}_n \times_{\mathcal{U}_{n - 1}} \mathcal{U}_n$ induced by the identity $d_jd_i = d_{i - 1}d_j$ for $0 \leq i < j \leq n + 1$.

In particular, this implies the $\zll$-linear $\fpsl$ version of Theorem~\ref{theo:CD}(1). The $\zll$-linear $\fpsl$ version of Theorem~\ref{theo:CD}(2) is trivial because every presheaf with transfers is a presheaf with traces, and therefore if we are working $\zll$-linearly, an $\fpsl$ sheaf. 

We also point out that if we are only interested in the $0 \to F(X) \to F(U) \to F(U \times_X U)$ part of this sequence, then the proof takes one line (cf. the last equation of the proof).
\end{rema}


\begin{proof}
We will show that the sequence is exact. It then follows that $F$ is $\fpsl$ acyclic (Lemma~\ref{lemm:acyclicIsCechAcyclic}), and that if $F$ is also a cdh sheaf of $\zll$ modules, then it is an $\ldh$ sheaf (Proposition \ref{prop:compatSheaf}(3)).

If $d_k: U^{\times_X n} \to U^{\times_X (n - 1)}$ is the projection that loses the $k$th coordinate, then the squares
\[ \xymatrix{
U^{\times_X (n + 1)} \ar[r]^{d_i} \ar[d]_{d_j} & U^{\times_X n} \ar[d]^{d_j} \\
U^{\times_X n} \ar[r]_{d_{i - 1}} & U^{\times_X (n - 1)}
} \]
are cartesian for all $j < i$. It follows from (CdB) that we have $F(d_{i - 1})\Tr_{d_{j}} = \Tr_{d_{j}}F(d_i)$ for $j < i$. We claim that $\Tr_{d_{0}}$ is a chain homotopy between zero and $d$ times the identity (in degree zero we take $\Tr_p$). We have
\[ \Tr_{d_{0}} \sum_{i = 0}^{n + 1}(-1)^i F(d_i) = \Tr_{d_{0}} F(d_0) + \sum_{i = 1}^{n + 1}(-1)^i F(d_{i - 1}) \Tr_{d_{0}} = d \cdot \id - \sum_{i = 0}^n (-1)^i F(d_i) \Tr_{d_{0}} \]
in degrees $n = 1, 2, 3, \dots$. In degree zero we have
\[ \Tr_{d_{0}} (F(d_0) - F(d_1)) = d \cdot \id - \Tr_{d_{0}} F(d_1) = d \cdot \id - F(f) \Tr_{f}. \]
Since $d \cdot \id$ is an isomorphism, the complex is acyclic.
\end{proof}

\begin{coro} \label{coro:equivalencecdhltransfers}
The canonical functor $Shv_{\ldh}(Cor(S), \zll) \to Shv_{cdh}(Cor(S), \zll)$ is an equivalence.
\end{coro}

\begin{proof}
Follows immediately from Lemma~\ref{lemm:transfersImpliesTraces} and Lemma~\ref{lemm:fpslAcyclic}.
\end{proof}

We recall the following theorem from \cite{CD}.

\begin{theo}[{\cite{CD}}]\ \label{theo:CD}
\begin{enumerate}
 \item \cite[Definition 9.3.2]{CD}, \cite[Corollary 9.3.16]{CD}, \cite[Proposition 9.4.8]{CD}. \label{theo:CD:enum:1} For any cdh hypercover $\mathcal{U}_\bullet \to X$ the associated sequence of cdh sheaves with transfers is exact
\[ \cdots \to L_{cdh}(\mathcal{U}_2) \to L_{cdh}(\mathcal{U}_1) \to L_{cdh}(\mathcal{U}_0) \to L_{cdh}(X) \to 0 \]
 \item \cite[Lemma 9.3.7]{CD} The inclusion $Shv_{cdh}(Cor(S)) \to PShv(Cor(S))$ admits a left adjoint $a: PShv(Cor(S)) \to Shv_{cdh}(Cor(S))$ such that the diagram
\[ \xymatrix{
PShv(Cor(S)) \ar[r]^a \ar[d]_{Oub} & Shv_{cdh}(Cor(S)) \ar[d]^{Oub} \\
PShv(Sch(S)) \ar[r]_{(-)_{cdh}} & Shv_{cdh}(Sch(S))
} \]
commutes where $Oub$ are the forgetful functors, i.e., precomposition with the graph functor $Sch(S) \to Cor(S)$.
\end{enumerate}
\end{theo}


\begin{lemm}
The category $Shv_{cdh}(Cor(S))$ is a Grothendieck abelian category and hence has enough injectives.
\end{lemm}

\begin{proof}
Every category of presheaves on an essentially small category is a Grothendieck abelian category. Moreover, if $\mathscr{A}$ is a Grothendieck abelian category and $R: \mathscr{B} \to \mathscr{A}$ is a fully faithful functor with a left adjoint then $\mathscr{B}$ is Grothendieck abelian. Grothendieck abelian categories have enough injectives.
\end{proof}

\begin{prop}[{\cite[3.1.8]{Voev00}}] \label{prop:318}
Let $X \in Sch(S)$ and $F$ a cdh sheaf with transfers. Then for any $n$ there is a canonical isomorphism
\[ Ext^n_{Shv_{cdh}(Cor(S))}(L_{cdh}(X), F) \cong H_{cdh}^n(X, F). \]
\end{prop}

\begin{proof}
We follow the proof of \cite[3.1.8]{Voev00}. Let $F \to I^\bullet$ be an injective resolution of $F$ in $Shv_{cdh}(Cor(S))$. The $I^n$ are not necessarily injective in $Shv_{cdh}(Sch(S))$ but if they are acyclic then we can use them to calculate the cohomology groups on the right. It is enough to show that their \v{C}ech cohomology vanishes in positive degrees (Lemma~\ref{lemm:acyclicIsCechAcyclic}). This follows immediately from the adjunction and the exact sequence in Theorem~\ref{theo:CD} and Yoneda.
\end{proof}

\begin{prop} \label{prop:318l}
Let $X \in Sch(S)$ and $F$ a sheaf of $\zll$ modules for the $\ldh$ pretopology equipped with a structure of transfers. Then for any $n$ there is a canonical isomorphism
\[ Ext^n_{Shv_{\ldh}(Cor(S), \zll)}(L_{\ldh}(X), F) \cong H_{\ldh}^n(X, F). \]
\end{prop}

\begin{proof}
As in the Proposition~\ref{prop:318} we need to show that if $I$ is an injective object of $Shv_{\ldh}(Cor(S))$ then it is $\ldh$ acyclic. After the equivalence (Corollary~\ref{coro:equivalencecdhltransfers}) and Proposition~\ref{prop:318} it is cdh acyclic. It has transfers so the higher \v{C}ech cohomology of every $\fpsl$ cover is zero (Lemma~\ref{lemm:fpslAcyclic}). Hence (Proposition~\ref{prop:sigmaTauAcyclic}) it is $\ldh$ acyclic.
\end{proof}

\begin{theo} \label{theo:lcdhcohomologyAgree}
Let $F$ be a cdh sheaf of $\zll$ modules with transfers. Then the canonical morphism $H_{cdh}^n(-, F) \to H_{\ldh}^n(-, F)$ is an isomorphism. Moreover, these functors have a canonical structure of presheaves with transfers.
\end{theo}

\begin{proof}
The first statement is a direct consequence of Proposition~\ref{prop:318}, Proposition~\ref{prop:318l}, and Corollary~\ref{coro:equivalencecdhltransfers}. For the second, recall the isomorphism
\[ Ext^n_{Shv_{cdh}(Cor(S))}(L_{cdh}(X), F) \cong H_{cdh}^n(X, F) \]
 and note that $Ext^n_{Shv_{cdh}(Cor(S))}(L_{cdh}(-), F)$ is functorial with respect to transfers.
\end{proof}



\section{The completely decomposed discrete topology} \label{sec:cd}

The goal of this section is to prove a criterion for a structure of traces on $F$ to pass to a structure of traces on the associated cdh separated presheaf $\im(F \to F_{cdh})$ (see Proposition~\ref{prop:tauSepHasTraces}). To obtain this criterion, we introduce the completely decomposed discrete topology (Definition~\ref{defi:cd}), which we will abbreviate to cdd topology. We show first that a structure of traces on $F$ passes to the sheafification $F_{\dis}$ for this cdd topology  (Theorem~\ref{theo:tracesForcd}). In a further subsection we also show that for a presheaf that is separated for the cdd topology, the properties (Tri1) and (Tri2) are implied by (Tri1)$_{\leq 0}$ (Proposition~\ref{prop:sepImpliesTrin}).

\begin{defi} \label{defi:cd}
The \emph{completely decomposed discrete pretopology} or \emph{cdd pretopology} has as covers families of morphisms (not necessarily of finite type) $\{ U_i \to X\}$ such that for each point $x \in X$ there exists an $i$ and a point $u \in U_i$ over $x$ such that $[k(u):k(x)] = 1$.
\end{defi}

\begin{rema}
The name is motivated by two ideas. Firstly that in the ``discrete'' pretopology every jointly surjective family of morphisms should be a cover. Secondly, adding the adjective ``completely decomposed'' to a pretopology should add the requirement that every cover of the spectrum of a field should have a section.
\end{rema}

We collect here some easy properties of the cdd pretopology. We assume that we are using a category of schemes such that for every scheme $X$ and every point $x \in X$ the morphism $x \to X$ is also in our category.

\begin{lemm} \label{lemm:cdProp}
\begin{enumerate}
 \item For every scheme $X$, every cover for the cdd pretopology admits a refinement by the cover $\{ x \to X \}_{x \in X}$.
 \item If $F$ is a presheaf and $F_{\dis}$ the associated cdd sheaf, there is a canonical isomorphism of presheaves $F_{\dis}(X) \cong \prod_{x \in X} F(x)$.
 \item A presheaf is separated for the cdd topology if and only if for every scheme $X$ the morphism $F(X) \to \prod_{x \in X} F(x)$ is injective.
 \item A presheaf is a cdd sheaf if and only if for every scheme $X$ the morphism $F(X) \to \prod_{x \in X} F(x)$ is an isomorphism.
 \item If $F$ is $\fpsl$ separated on schemes of dimension zero, then $F_{\dis}$ is $\ldh$ separated.
\end{enumerate}
\end{lemm}

\begin{rema} \label{rema:whyCdd}
The cdd pretopology arises quite naturally for us in the following way. The two main classes of presheaves with traces that we are interested in studying - homotopy invariant Nisnevich sheaves with transfers, and the Nisnevich sheafification of algebraic $K$-theory - both satisfy a Gersten exact sequence for regular schemes (in the appropriate categories of schemes). Notably, for any such connected regular scheme $X$ with generic point $\eta$ we have $F(X) \subseteq F(\eta)$. We do not hope to have this property for non-regular schemes, but if $F$ is separated for some topology $\tau$, and a non-regular scheme $X$ admits a regular $\tau$-cover $X' \to X$, then we will have $F(X) \subseteq F(X') \subseteq \prod F(\eta_i)$ where the $\eta_i$ are the generic points of $X'$. In case we assume resolution of singularities and use $\tau = cdh$, this line of reasoning leads to $F(X) \subseteq \prod_{x \in X} F(x)$. That is, $F(X) \subseteq F_{\dis}(X)$. A similar phenomena occurs if we have traces, are $\zll$-linear, and use the $\ldh$-pretopology.
\end{rema}


\subsection{Traces on \texorpdfstring{$F_{\dis}$}{Fcdd}}

In this subsection we show that a structure of traces on a presheaf $F$ passes to a canonical structure of traces on the associated cdd sheaf $F_{\dis}$ (Theorem~\ref{theo:tracesForcd}).

In this subsection we work with a category of schemes that is closed under fibre products, and such that for every scheme $X$ in the category, and every point $x$ of $X$, the morphism $x \to X$ is also in the category.


\begin{theo} \label{theo:tracesForcd}
Suppose that $F$ is a presheaf with traces that satisfies (Tri1)$_{\leq 0}$. Then there is a unique structure of traces on $F_{\dis}$ such that $F \to F_{\dis}$ is a morphism of presheaves with traces. This structure also satisfies (Tri1)$_{\leq 0}$.
\end{theo}

\begin{proof}
We will use the canonical isomorphisms $F_{\dis}(X) \cong \prod_{x \in X} F(x)$ to take $\prod_{x \in X} F(x)$ as the definition of $F_{\dis}(X)$. Let $f: Y \to X$ be a finite flat surjective morphism. We define a morphism
\[ \Tr^{\dis}_f: \prod_{y \in Y} F(y) \to \prod_{x \in X} F(x), \qquad \qquad (s_y) \mapsto (t_x) \]
where
\[ t_x = \sum_{y \in x \times_X Y}  \length \OO_{x \times_X Y, y} \Tr_{f|_{y/x}}(s_y) \]
and we have used $f|_{y/x}$ for the induced morphisms $y \to x$. We claim that these morphisms satisfy (Deg), (CdB), are functorial, and are compatible with the morphism $F \to F_{\dis}$.

\emph{The degree formula.} The axiom (Deg) is straightforward and needs only to be checked in the case where $X$ has a unique reduced point $x$. In this case $X = Spec(k)$ and $Y = Spec(A)$ where $k$ is a field and $A$ is a finite $k$-algebra we have $A = \oplus_i \OO_{Y, y_i}$ where the sum is over the points $y_i$ of $Y$. We can express $\deg f$ as
\[ \deg f = \sum_{y_i \in Y} \dim_k \OO_{Y, y_i} \stackrel{(\ref{lemm:degLength})}{=} \sum_{y_i \in Y} \length \OO_{Y, y_i} \dim_k k(y_i) = \sum_{y_i \in Y} \length \OO_{Y, y_i} \deg f|_{y/x} \]
where $f|_{y/x}$ is still the induced morphism $y \to x$. We then have 
\[ \begin{split}
\Tr_{f}^{\dis} F_{\dis}(f) &= \sum_{y_i \in Y} \length \OO_{Y, y_i} \Tr_{(f|_{y/x})} F(f|_{y/x}) \\
&= \sum_{y_i \in Y} \length \OO_{Y, y_i} \deg (f|_{y/x}) \cdot \id_x \\
&= \deg f \cdot \id_x. 
\end{split} \]

\emph{The change of base formula.} Now consider a cartesian square (\ref{equa:cartSquare}). If $W$ is a point of $X$, say $W = x \in X$, then we have (CdB) by our definition of $\Tr^{\dis}_f$ and $\Tr^{\dis}_g$. To check that two sections in $F_{\dis}(W)$ agree it is sufficient to check them on each point $w \in W$ and so to prove (CdB) in general, it suffices to consider the case when $W$ is an integral scheme of dimension zero. In this case, $W \to X$ factors through the inclusion of a reduced point of $X$, and so we reduce to the case where $W$ and $X$ are both integral dimension zero schemes. Write $W = w$ and $X = x$. By additivity we can assume that $Y$ is connected.

Suppose for the moment that $Y$ is integral and write $Y = y$. Let $z_i$ be the points of $y \times_x w$. To have (CdB) with these assumptions we must show that
\[ F(p) Tr_f = \sum_{z_i \in y \times_x w} \length \OO_{y \times_x w, z_i} Tr_{g|_{z_i/w}} F(q|_{z_i/y}). \]
This follows from (CdB), (Fon), and applying (Tri1)$_{\leq 0}$ to the triangles obtained from $z_i \to y \times_x w \to w$.

Now we return to the case where $Y$ is not necessarily reduced but has a unique point $y$. We can use (Tri1)$_{\leq 0}$ on the triangle $y \to Y \to x$, and (CdB) on the cartesian square having lower row $y \to x$, and so it suffices to show that $\Tr_g^{\dis} = \length \OO_{Y, y} \Tr_h^{\dis} F_{\dis}(\iota)$ where $\iota$ and $h$ are the morphisms $h: y \times_x w \stackrel{\iota}{\to} Y \times_x w \stackrel{g}{\to} w$. Consulting our definition of the $\Tr^{\dis}$, we see that we must show that for each point $z \in Y \times_x w$ we have
\[ \length \OO_{Y \times_x w, z} = \length \OO_{Y, y} \length \OO_{y \times_x w, z}. \]
This follows from Lemma~\ref{lemm:CDBalgebra}.

\emph{Functoriality.} We need to show that if $W \stackrel{g}{\to} Y \stackrel{f}{\to} x$ are finite flat surjective morphisms, and $x$ is a integral dimension zero scheme, then 
\[ \sum_{w \in W} \length \OO_{W, w} = \sum_{y \in Y} \length \OO_{Y, y} \sum_{w \in y \times_Y W} \length \OO_{y \times_Y W, w}. \]
Clearly it suffices to consider the case when $Y$ and $W$ are connected. The result follows now from Lemma~\ref{lemm:fonalgebra}.

\emph{Compatibility with $F \to F_{\dis}$.} By (CdB) it suffices to consider morphisms $Y \to x$ where $x$ is a integral scheme of dimension zero. We can also assume that $Y$ is connected by additivity. Clearly $F \to F_{\dis}$ is compatible with traces when $Y$ is also reduced. We have assumed that $F$ satisfies (Tri1) on dimension zero schemes and we have already noticed that the trace morphisms we have defined on $F_{\dis}$ satisfy (Tri1), so we are done.
\end{proof}

\begin{prop} \label{prop:tauSepHasTraces}
Suppose that $F$ is a presheaf with traces that satisfies (Tri1)$_{\leq 0}$ and suppose that $F_{cdh} \to F_{\dis}$ is a monomorphism. Then $\im(F \to F_{cdh})$ has a unique structure of traces such that $F \to \im(F \to F_{cdh})$ is a morphism of presheaves with traces.
\end{prop}

\begin{proof}
The kernel $K$ of the epimorphism of presheaves $F(X) \to \im(F \to F_{cdh})(X)$ is $K(X) = \{ s \in F(X) $ such that $s|_{X'} = 0$ for some cdh cover $X' \to X \}$. It is enough to show that the trace morphisms of $F$ preserve $K$. After Theorem~\ref{theo:tracesForcd}, for any finite flat surjective morphism $f: Y \to X$ we have a commutative diagram
\[ \xymatrix{
K(Y) \ar[r] & F(Y) \ar[rrr] \ar[d] & & & F_{\dis}(Y) \ar[d] \\
& F(X) \ar[r] & \im(F \to F_{cdh})(X) \ar[r] & F_{cdh}(X) \ar[r] & F_{\dis}(X)
} \]
and so the result follows from the injectivity of $F_{cdh}(X) \to F_{\dis}(X)$.
\end{proof}


\subsection{For cdd separated presheaves (Tri)\texorpdfstring{$_{\leq 0}$}{leq0} implies (Tri)\texorpdfstring{$_{\leq n}$}{leqn}}

In this section we show the following proposition. As in the proof of Theorem~\ref{theo:tracesForcd} we end up chasing multiplicities around and this is done in Lemma~\ref{lemm:triPullbackPoint}.

In this subsection we continue to work with a category of schemes that is closed under fibre products, and such that for every scheme $X$ in the category, and every point $x$ of $X$, the morphism $x \to X$ is also in the category.

\begin{prop} \label{prop:sepImpliesTrin}
Suppose that $F$ is a presheaf with traces such that $F \to F_{\dis}$ is a monomorphism of presheaves. If $F$ satisfies (Tri1)$_{\leq 0}$ then it also satisfies (Tri1)$_{\leq n}$ for all $n$. Moreover, for such a presheaf (Tri2)$_{\leq n}$ is always satisfied for all $n$.
\end{prop}

\begin{proof}
Recall that (Tri2)$_{\leq 0}$ is always satisfied (Lemma~\ref{lemm:Tr2dimZero}). We will show that (Tri1)$_{\leq 0}$ implies (Tri1)$_{\leq n}$ (resp. (Tri2)$_{\leq 0}$ implies (Tri2)$_{\leq n}$) under the assumption that the morphism $F \to F_{\dis}$ is a monomorphism. Since $F \to F_{\dis}$ is a monomorphism, it is sufficient to show that for every (not necessarily closed) point $\iota: x \to X$ of $X$ and every triangle (\ref{equa:tri}) satisfying the hypotheses of (Tri1) (resp. (Tri2)) we have $F(\iota)F(f) = F(\iota) \sum m_k \Tr_{g_{k}}F(h_k)$ (resp. $\tfrac{\deg g}{\deg f} F(\iota)\Tr_f = F(\iota)\Tr_g F(h)$) . By (CdB) it is enough to show that we have $\Tr_{(f \times_X x)} = \sum m_k \Tr_{(g_k \times_X x)} F(h_k \times_X x)$ (resp. $\tfrac{\deg g}{\deg f} \Tr_{(f \times_X x)} = \Tr_{(g \times_X x)} F(h \times_X x)$) for every point $x \in X$. Furthermore, since everything is of dimension zero now, by additivity it suffices to consider the restrictions of these morphisms to each point $y \in Y$ over $x$. Let $y_{k\ell}' \in Y'$ be the points of the $k$th connected component of $Y'$ that lie over $y$, let $W$ be the connected component of $Y \times_X x$ containing $y$ and $W_{k\ell}$ the connected component of $Y' \times_X x$ containing $y_{k\ell}$ so we have the following commutative diagrams.
\[ \xymatrix{
& y_i \ar[rr]^{\eta'} \ar[d]_{\iota'} && y \ar[d]^\iota \\
 \ar[drr]_{g \times_X x} Y' \times_X x \supseteq & W_{k\ell} \ar[rr]^{\eta} \ar[dr]^\gamma && W \ar[dl]_\phi & \subseteq Y \times_X x \ar[dll]^{f \times_X x} \\
&& x
} \]
The calculation for (Tri1) is
\[ \begin{split}
\Tr_{\phi} &\stackrel{(Tri1)_{\leq 0}}{=} \length \OO_{Y \times_X x, y} \Tr_{\phi \iota} F(\iota) \\
&\stackrel{(*)}{=} \length \OO_{Y \times_X x, y} \Tr_{\phi \iota \eta'} F(\iota \eta') \\
&= \length \OO_{Y \times_X x, y} \Tr_{\gamma \iota'} F(\iota') F(\eta) \\
&\stackrel{(**)}{=} \sum_k m_k \sum_\ell \length \OO_{Y' \times_X x, y_{k\ell}} \Tr_{\gamma \iota'} F(\iota') F(\eta) \\
&\stackrel{(Tri1)_{\leq 0}}{=} \sum_k m_k \sum_\ell \Tr_{\gamma} F(\eta)
\end{split} \]
where in the step (*) we have used the hypothesis of (Tri1) that $Y'$ is the disjoint union of the integral components of $Y$ (so $\eta'$ is an isomorphism), and step (**) is the Lemma~\ref{lemm:triPullbackPoint}. Notice that the hypotheses of Lemma~\ref{lemm:triPullbackPoint} include the two cases of (Tri1) and (Tri2) (see Remark~\ref{rema:triPullbackPoint}). The calculation for (Tri2) is similar (there is no $k$ because $Y'$ is connected in the hypotheses of (Tri2)).
\[ \begin{split}
d \cdot \Tr_{\phi} &\stackrel{(Tri1)_{\leq 0}}{=} d \length \OO_{Y \times_X x, y} \Tr_{\phi \iota} F(\iota) \\
&\stackrel{(**)}{=} \sum_\ell \length \OO_{Y' \times_X x, y_{\ell}} [k(y'_\ell):k(y)] \Tr_{\phi \iota} F(\iota) \\
&\stackrel{(Tri2)_{\leq 0}}{=} \sum_\ell \length \OO_{Y' \times_X x, y_{\ell}} \Tr_{\phi \iota \eta'} F(\eta' \iota) \\
&= \sum_\ell \length \OO_{Y' \times_X x, y_{\ell}} \Tr_{\gamma \iota'} F(\iota' \eta) \\
&\stackrel{(Tri1)_{\leq 0}}{=} \sum_\ell \Tr_{\gamma} F(\eta)
\end{split} \]
Again, in the step (**) we have used Lemma~\ref{lemm:triPullbackPoint}.
\end{proof}


\section{Gersten presheaves} \label{sec:gersten}

We have two goals in this section. The first is to find a condition on a presheaf of $\zll$ modules with traces which will imply that $F_{\ldh} \to F_{\dis}$ is injective. We are interested in this because if $F_{\ldh} \to F_{\dis}$ is injective for such a presheaf $F$ then $F_{cdh} \to F_{\ldh}$ is an isomorphism. Our second goal is to promote the structure of traces of $F$ to a structure of traces on $F_{cdh}$.

\subsection{Comparison of the cdh and \texorpdfstring{$\ldh$}{ldh} sheafifications}

In this subsection we continue to work with a category of schemes that is closed under fibre products, and such that for every scheme $X$ in the category, and every point $x$ of $X$, the morphism $x \to X$ is also in the category. We add the hypothesis that the cdh and $\ldh$ pretopologies are defined on our category, and that every $\ldh$ cover is refinable by a regular $\ldh$ cover.

\begin{lemm} \label{lemm:ellcdmono}
Suppose that $F$ is a presheaf of $\zll$ modules with traces that satisfies (Tri1)$_{\leq 0}$ and for every regular scheme $X$ the morphism $F(X) \to \prod_{x \in X^{(0)}} F(x)$ is injective. Then the canonical morphism $F_{\ldh} \to F_{\dis}$ is a monomorphism.
\end{lemm}

\begin{rema} \label{rema:ellcdmono}
The canonical morphism $F_{\ldh} \to F_{\dis}$ is the one obtained from the observation that $F_{\dis}$ is an $\ldh$ sheaf. This is so because the cdd topology is finer than the cdh topology, and $F_{\dis}$ has a structure of traces (Theorem~\ref{theo:tracesForcd}, Lemma~\ref{lemm:fpslAcyclic}).
\end{rema}

\begin{proof}
The injectivity is straightforward: for any scheme $X$ and section $s \in ker(F_{\ldh}(X) \to F_{\dis}(X))$, there exists an $\ldh$ cover $X' \to X$ such that $s|_{X'}$ is in the image of $F(X') \to F_{\ldh}(X')$. By hypothesis on our category of schemes (in practice this we be true via Corollary~\ref{coro:regularlCover}) we can assume that $X'$ is regular. In this case, by hypothesis, the morphism $F(X') \to \prod_{x \in X^{(0)}} F(x)$ is injective and so $s|_{X'} = 0$, hence $s = 0$.
\end{proof}

\begin{coro} \label{coro:cdhliso}
Suppose that $F$ is a presheaf of $\zll$ modules with traces and for every regular scheme $X$ the morphism $F(X) \to \prod_{x \in X^{(0)}} F(x)$ is injective. Then $F_{cdh} \to F_{\ldh}$ is an isomorphism.
\end{coro}

\begin{proof}
The presheaf $F$ is $\fpsl$ separated due to the structure of traces (Lemma~\ref{lemm:fpslAcyclic}) and so after Proposition~\ref{prop:compatSheaf}(1) $F_{cdh}$ is $\fpsl$ separated as well. That is, $F_{cdh} \to F_{\ldh}$ is a monomorphism of presheaves. We have just seen (Lemma~\ref{lemm:ellcdmono}) that $F_{\ldh} \to F_{\dis}$ is injective, and hence $F_{cdh} \to F_{\dis}$ is injective. It now follows from Proposition~\ref{prop:tauSepHasTraces} that its associated cdh separated presheaf $\im(F \to F_{cdh})$ has a structure of traces compatible with that of $F$. Hence, $\im(F \to F_{cdh})$ is an $\fpsl$ sheaf (Lemma~\ref{lemm:fpslAcyclic}) and therefore $F_{cdh}$ is as well (Proposition~\ref{prop:compatSheaf}(3)). That is, the canonical morphism $F_{cdh} \to F_{\ldh}$ is an isomorphism.
\[ \xymatrix{
\ar@<1ex>@{}[rr]^-{\textrm{(i) has traces}} & F \ar[r] \ar@<1ex>@{}[rr]^-{\textrm{(v) has traces}} & \im(F \to F_{cdh}) \ar[r] &  F_{cdh} \ar@<-1ex>@/_12pt/[rr]_{\textrm{(iv) monic}} \ar[r]_{\textrm{(ii) monic}}^{\textrm{(vi) iso}} & F_{\ldh} \ar[r]^{\textrm{(iii) monic}} & F_{\dis}
} \]
The above diagram gives a summary of the argument: (i) implies (ii); then (ii) + (iii) implies (iv) which implies (v) which implies (vi).
\end{proof}

\subsection{Gersten presheaves}

In this subsection we introduce the notion of a Gersten presheaf (Definition~\ref{defi:gersten}). This is a property satisfied by homotopy invariant Nisnevich sheaves with transfers $F$ on the category $Sm(k)$ of separated smooth schemes of finite type with $k$ a perfect field \cite[4.37]{Vctpt}. It conjecturally satisfied by the Zariski sheafification of algebraic $K$-theory for all regular schemes \cite[5.10]{Qui73}.

In the previous section exactness of $0 \to F(X) \to \oplus_{x \in X^{(0)}} F(x)$ allowed us to prove that $F_{\ldh} \to F_{cdd}$ is a monomorphism. Having exactness of $F(X) \to \oplus_{x \in X^{(0)}} F(x) \to \oplus_{x \in X^{(1)}} F_{-1}(x)$ will allow us to recognise the image of $F_{\ldh} \to F_{\dis}$ and enable us to pass a structure of traces on $F$ to a structure of traces on $F_{cdh}$ (Proposition~\ref{prop:tracesForcdh}).


The following definition is inspired by Gersten's conjecture in $K$-theory \cite[5.10]{Qui73}. Recall that for a scheme $X$ we denote the set of points of codimension $n$ by $X^{(n)}$.

\begin{defi} \label{defi:gersten}
Let $F$ be a presheaf on a category of schemes such that for every scheme $X$, and every point $x \in X$ of codimension $\leq 1$ the morphism $x \to X$ is also in the category. We will call $F$ a \emph{Gersten presheaf} if it is equipped with
\begin{enumerate}
 \item an abelian group $F_{-1}(x)$ for every scheme of dimension zero,
 \item a morphism $\partial_{(x_0, x_1)}: F(x_0) \to F_{-1}(x_1)$ for every pair $(x_0, x_1) \in X^{(0)} \times X^{(1)}$ with $x_1 \in \overline{\{x_0\}}$,
\end{enumerate}
such that for each \emph{regular} scheme $X$ the following sequence is exact
\[ 0 \to F(X) \to \prod_{x_0 \in X^{(0)}} F(x_0) \stackrel{\partial_{(x_0, x_1)}}{\to} \prod_{x_1 \in X^{(1)}} F_{-1}(x_1). \]

If $\tau$ is a Grothendieck topology then a \emph{$\tau$ Gersten sheaf} is just a Gersten presheaf that is also a $\tau$ sheaf.
\end{defi}

\begin{rema}
The notation $F_{-1}(x)$ is very suggestive but at the moment we haven't asked for anything more than these be a class of groups. We don't ask that they are functorial, or that they are related to $F$ in any way other than via the $\partial_{(x_0, x_1)}$.

\end{rema}

\begin{exam}
\begin{enumerate}
 \item Homotopy invariant Nisnevich sheaves with transfers are Gersten presheaves with traces on the category of separated schemes essentially of finite type over a perfect field \cite[4.37]{Vctpt}.
 \item Gersten's conjecture (\cite[5.10]{Qui73} or \cite{Ger73}) implies that the Zariski sheafification of algebraic $K$-theory is a Gersten presheaf for all regular schemes. This is known to be true in certain cases, including the case of equicharacteristic schemes \cite{Pan03}.
\end{enumerate}
\end{exam}

\begin{lemm} \label{lemm:gerTri}
Consider a triangle (\ref{equa:tri}) such that $X$ is regular and let $F$ be a Gersten presheaf. Then (Tri2) is satisfied. If moreover (Tri1)$_{\leq 0}$ is satisfied then (Tri1)$_{\leq n}$ is satisfied too, for all $n$.
\end{lemm}

\begin{proof}
The pullback along the generic point of $X$ is injective, and so we reduce to the dimension zero case, which is Lemma~\ref{lemm:Tr2dimZero}.
\end{proof}


\subsection{Nisnevich Gersten sheaves with traces on regular curves}

In this subsection we show that for Nisnevich Gersten sheaves of $\zll$ modules with traces, regular schemes of dimension $\leq 1$ behave like points for the $\ldh$ topology, in the sense that $F \to F_{\ldh}$ is an isomorphism on such schemes (Proposition~\ref{prop:isoDimensionOne}). We use this in the proof of Proposition~\ref{prop:tracesForcdh} to recognise the image of $F_{\ldh} \to F_{\dis}$.

In this subsection we denote by $Sch(S)$ the category of separated schemes essentially of finite type over a noetherian quasi-excellent base $S$.

We begin with three lemmata.

\begin{lemm} \label{lemm:fpslRefine}
If $X$ is an integral noetherian scheme then every $\fpsl$ cover has a refinement of the form $\stackrel{g}{\to} \stackrel{f}{\to}$ where $f$ is a blowup with nowhere dense centre and $g$ is $\fpsl$ with an integral source.
\end{lemm}

\begin{proof}
Let $f: U \to X$ be an $\fpsl$ cover with $X$ integral. If $X$ is zero dimensional then it is $Spec(k)$ with $k$ a field and $U = Spec(A)$ with $A$ a finite $k$ algebra. In particular, if $\m_i$ are the primes of $A$ then $\deg f = \sum [A / \m_i : k] \length A_{\m_i}$. Since $\ell$ doesn't divide $\deg f$, there is some $i$ for which it doesn't divide $[A / \m_i : k] \length A_{\m_i}$ and hence doesn't divide $[A / \m_i : k]$. We take $Spec(A / \m_i)$ as our refinement.

If $X$ is of dimension greater than zero, then by the platification theorem (Theorem~\ref{theo:platification}) there is a blowup with nowhere dense centre $\tilde{X} \to X$ such that the integral components of the proper transform (which is the pull-back in this case due to $f$ being flat) are flat over $X$. By the dimension zero case, there is one of them for which the generic point is $\fpsl$ over the generic point of $X$, and so this integral component gives us the desired refinement.
\end{proof}

\begin{lemm} \label{lemm:cdpPointfpslRefine}
Suppose $X$ is a regular noetherian quasi-excellent scheme of dimension one and $U \to X$ is a morphism which is a composition of Nisnevich and $\fpsl$ covers. Then there exists a refinement of the form $V_1 \to V_0 \to X$ such that $V_0 \to X$ is Nisnevich, $V_1 \to V_0$ is $\fpsl$, the schemes $V_1$ and $V_0$ are regular, and each integral component of $V_0$ has a unique integral component of $V_1$ over it.
\end{lemm}

\begin{proof}
It suffices to consider separately the cases where $U \to X$ is either $\fpsl$ or Nisnevich. In the Nisnevich case $U$ is already regular (\cite[I.9.2]{SGA1}) so only the $\fpsl$ case remains.

By Lemma~\ref{lemm:fpslRefine} we can assume that $U$ is integral.\footnote{Since $X$ is regular of dimension one, every local ring is either a field or a discrete valuation ring i.e., a principal ideal domain, and hence, every blowup is trivial.} Since $X$ is quasi-excellent, the normalisation $\tilde{U} \to U$ is a finite morphism \cite[Theorem 78]{Mat}. Since $X$ is regular of dimension one, flatness is equivalent to every generic point being sent to a generic point and so $\tilde{U} \to X$ is finite, flat, surjective, and of degree prime to $\ell$ (the latter because it is true generically, and the morphism is finite and flat).
\end{proof}

\begin{lemm} \label{lemm:dimOneSequence}
Let $F$ be a Gersten presheaf of $\zll$ modules with traces on $Sch(S)$ that satisfies (Tri1)$_{\leq 0}$, let $Y \to X$ be a finite flat surjective morphism of constant degree prime to $\ell$ between regular integral schemes of dimension one. Let $\widetilde{Y \times_X Y}$ be the normalisation of $Y \times_X Y$. Then the sequence
\[ 0 \to F(X) \to F(Y) \to F(\widetilde{Y \times_X Y}) \]
is exact, where the last morphism is the difference of the morphisms induced by the two projections.
\end{lemm}

\begin{proof}
We already know that $F(X) \to F(Y)$ is injective, so it remains to show exactness at $Y$. Let $Y_i$ be the integral components of $Y \times_X Y$ and $\tilde{Y}_i$ their normalisations. Let $p_1, p_2: Y \times_X Y \to Y$ be the projections, let $n: \widetilde{Y \times_X Y} \to Y \times_X Y$ be the canonical morphism, $n_{i}: \tilde{Y}_i \to \widetilde{Y \times_X Y}$ and $p_{i1}, p_{i2}: \tilde{Y}_{i} \to Y$ be the induced morphisms, and let $m_{i}$ be the multiplicities of the generic point of $Y_{i}$ in $Y \times_X Y$. By Lemma~\ref{lemm:gerTri} we have (Tri1) and (Tri2) for triangles with base $Y$. Moreover, since $Y$ is regular of dimension one, and $Y_i$ are integral, the morphisms $Y_i \to Y$ are flat (hence, finite flat surjective). 
\[ \xymatrix{
\widetilde{Y \times_X Y} = \amalg Y_{i} \ar[rdd]_{\sum p_{i1}} \ar[rrd]^{\sum p_{i2}} \ar[dr]^n \\
&  Y \times_X Y \ar[r]^{p_2} \ar[d]_{p_1} & Y \ar[d]^f \\
& Y \ar[r]_f & X
} \]
Then if $s \in ker(F(Y) \to F(\widetilde{Y \times_X Y}))$ we have 
\[ \begin{split}
 \tfrac{1}{d}F(f)\Tr_f(s) &\stackrel{(CdB)}{=} \tfrac{1}{d}\Tr_{p_{1}} F(p_2)s \\
&\stackrel{(Tri1)}{=} \sum \tfrac{1}{d} (m_{i} \Tr_{p_{i1}} F(n_{i})) F(p_2)s \\
&= \sum \tfrac{1}{d} m_{i} \Tr_{p_{i1}} F(p_{i2}) s \\
&\stackrel{s \in ker}{=} \sum \tfrac{1}{d} m_{i} \Tr_{p_{i1}} F(p_{i1}) s \\
&\stackrel{}{=} \sum \tfrac{1}{d} m_{i} \Tr_{p_{i1}} F(n_{i}) F(p_1) s \\
&\stackrel{(Tri1)}{=} \tfrac{1}{d} \Tr_{p_{1}} F(p_1) s \\
&\stackrel{(Deg)}{=} s
\end{split} \]
So $s$ is in the image of $F(X) \to F(Y)$.
\end{proof}

\begin{prop} \label{prop:isoDimensionOne}
Suppose $F$ is a Nisnevich Gersten sheaf of $\zll$ modules with traces on $Sch(S)$ that satisfies (Tri1)$_{\leq 0}$, and $X$ is a regular noetherian scheme of dimension $\leq 1$. Then $F(X) \to F_{\ldh}(X)$ is an isomorphism.
\end{prop}

\begin{proof}
We claim that for $X$ regular of dimension $\leq 1$ every proper cdh cover $X' \to X$ splits. Indeed, choose a lifting $\eta' \in X'$ of the generic point $\eta$ of $X$, consider its closure $\overline{\eta}'$ and normalise this. The resulting refinement $(\overline{\eta}')^\sim \to X$ is a birational proper morphism between regular schemes noetherian of dimension $\leq 1$. Consequently, it is an isomorphism.

Hence, every $\ldh$ cover is refinable by a cover of the form $V \stackrel{g}{\to} U \stackrel{f}{\to} X$ where $f$ is a Nisnevich cover and $g$ is an $\fpsl$ cover (see Example~\ref{exam:refinable}(3) for the cdh part). Since $F$ is separated with respect to these classes of covers, the morphism $F(X) \to F_{\ldh}(X)$ is injective.

For each $s \in F_{\ldh}(X)$ there exists a cover for which the restriction of $s$ is in the image of $F \to F_{\ldh}$ and we can assume that it has the form mentioned above. We can even assume that $V$ and $U$ are regular schemes of dimension one, and that each integral component of $U$ has a unique integral component of $V$ over it (Lemma~\ref{lemm:cdpPointfpslRefine}). Suppose that $t \in F(V)$ is a lifting of $s|_V$. The section $s|_V$ is in the kernel of $F_{\ldh}(V) \to F_{\ldh}(\widetilde{V \times_U V})$, but we have just seen that $F(\widetilde{V \times_U V}) \to F_{\ldh}(\widetilde{V \times_U V})$ is injective, and so $t$ lifts to a section $t' \in F(U)$ (Lemma~\ref{lemm:dimOneSequence}), which clearly, is a lift of $s|_U \in F_{\ldh}(U)$. The same argument lifts $t'$ to a section of $F(X)$: the section $s|_U$ is in the kernel of $F_{\ldh}(U) \to F_{\ldh}(U \times_X U)$ and the scheme $U \times_X U$ are regular of dimension one, and so since $F \to F_{\ldh}$ is injective on such schemes, $t'$ is in the kernel of $F(U) \to F(U \times_X U)$; since $F$ is a Nisnevich sheaf, we find a section $t'' \in F(X)$ sent to $s$.
\end{proof}

\subsection{Traces on \texorpdfstring{$F_{cdh}$}{Fcdh}}

In this subsection we continue to denote by $Sch(S)$ the category of separated schemes essentially of finite type over a noetherian quasi-excellent base $S$.

\begin{prop} \label{prop:tracesForcdh}
Suppose that $F$ is a Nisnevich Gersten sheaf of $\zll$ modules with traces on $Sch(S)$ such that (Tri1)$_{\leq 0}$ is satisfied. Then there is a unique structure of traces on $F_{cdh}$ such that $F \to F_{cdh}$ is a morphism of presheaves with traces. This structure satisfies (Tri1) and (Tri2).

Moreover, if $X$ is regular, then the canonical morphism $F(X) \to F_{cdh}(X)$ is an isomorphism.
\end{prop}

\begin{proof}
Recall that with these hypotheses the canonical morphism $F_{cdh} \to F_{\ldh}$ is an isomorphism (Corollary~\ref{coro:cdhliso}), and the canonical morphism $F_{\ldh} \to F_{\dis}$ is a monomorphism (Lemma~\ref{lemm:ellcdmono}). The plan is to find a criterion for a section in $F_{\dis}(X)$ to be in the image of $F_{\ldh} \to F_{\dis}$, and show that the trace morphisms of $F_{\dis}$ (Theorem~\ref{theo:tracesForcd}) preserve this criterion. Our criterion also shows that $F(X) \cong F_{cdh}(X)$ for $X$ regular. This proves the proposition.

For each scheme $X$ let $\rs F(X)$ be the set of sections $s \in F_{\dis}(X)$ such that
\begin{enumerate}
 \item[] for every morphism $i: T \to X$ from a semi-local regular scheme of dimension one there exists a section $s_T \in F(T)$ such that $s$ agrees with $s_T$ in $F_{\dis}(T)$.
\end{enumerate}
\[ \xymatrix{
& F(T) \ar[d] \\
F_{\dis}(X) \ar[r] & F_{\dis}(T)
} \]

We make and prove the following claims.

\emph{For every scheme $X$, the image of $F_{\ldh}(X)$ is contained in $\rs F(X)$ and the groups $\rs F(X)$ form a subpresheaf of $F_{\dis}$.} The first statement follows directly from the square on the left.
\[ \xymatrix{
F_{\lhd}(X) \ar[d] \ar[r] & F_{\ldh}(T) \stackrel{\ref{prop:isoDimensionOne}}{\cong} F(T) \ar[d] 
&& F_{\dis}(X) \ar[d] \ar[dr] &  F(T) \ar[d] \\
F_{\dis}(X) \ar[r] & F_{\dis}(T)
&& F_{\dis}(Y) \ar[r] & F_{\dis}(T)
} \]
For the second, suppose that $Y \to X$ is a morphism and $T$ is a regular scheme of dimension one and $T \to Y$ a morphism. Then the commutativity of the diagram on the right implies the second statement.

\emph{For $X$ regular, $\rs F(X)$ is precisely the image of $F(X)$ (and hence $F_{\ldh}(X)$ as well) in $F_{\dis}(X))$.}
Consider now a section $(s_x) \in \rs F(X)$ (where $s_x \in F(x)$) that we want to lift. For every point of codimension one $x \in X$, the localisation $\OO_{X, x}$ is a discrete valuation ring. Let $s_{\OO_{X, x}}$ be the section of $F(Spec(\OO_{X, x}))$ obtained via the criterion of $\rs F(X)$ and let $\eta$ be the generic point of $\OO_{X, x}$. By the exact sequence
\[ 0 \to F(Spec(\OO_{X, x})) \to F(\eta) \stackrel{\partial_{(\eta, x)}}{\to} F_{-1}(x) \]
since $s_\eta$ lifts we have $\partial_{(\eta, x)}s = 0$. This is true for every pair $(\eta, x) \in X^{(0)} \times X^{(1)}$ and so by the exact sequence
\[ 0 \to F(X) \to \prod_{x \in X^{(0)}} F(x) \to \prod_{x \in X^{(1)}} F_{-1}(x) \]
the section $(s_\eta)_{\eta \in X^{(0)}}$ lifts to a section $s \in F(X)$ such that $s|_\eta = s_\eta$ for each generic point. We claim that $s|_x = s_x$ for all points of $X$ and we prove it by induction on the codimension.

Suppose that it is true for points of codimension less than $n$ and let $x$ be a point of codimension $n$. Then as a result of the regularity of the local ring $\OO_{X, x}$ there exists a discrete valuation ring $R$ and a morphism $Spec(R) \to X$ such that the image of the closed point is $x$ and the image of the open point is a point $y$ of codimension $n - 1$ (in fact due to the existence of a regular sequence in $\OO_{X, x}$ we can choose $R$ such that the morphism induces an isomorphism on residue fields). By the criterion of $\rs F$ there is a section $s_R \in F(Spec(R))$ whose restrictions to $y$ and $x$ agree with $s_y$ and $s_x$. Hence, the restriction of $s|_{Spec(R)}$ to $y$ agrees with $s_R|_y$. But $F(Spec(R)) \to F(y)$ is injective by the Gersten sequence, and so $s|_{Spec(R)} = s_R$. But this implies that $s|_x =s_R|_x$ and by construction this is $s_x$.

At this point we have shown $F(X) \cong F_{cdh}(X)$ for regular $X$, since $F(X) \to \rs F(X)$ is injective by the Gersten sequence.

\emph{For $X$ any scheme, $\rs F(X)$ is precisely the image of $F_{\ldh}(X)$ in $F_{\dis}(X))$.}
Let $X$ be a scheme and $X' \to X$ an $\ldh$ cover with $X'$ regular, and $X'' \to X' \times_X X'$ an $\ldh$ cover with $X''$ regular (Corollary~\ref{coro:regularlCover}). We have the following diagram.
\[ \xymatrix{
F_{\ldh}(X) \ar[r] \ar[d]_\cong & \rs F(X) \ar[d] \\
ker( F_{\ldh}(X') \to F_{\ldh}(X'')) \ar[r] & ker( \rs F(X') \to \rs F(X''))
} \]
We have seen that $F_{\ldh} \to \rs F$ is an isomorphism on regular schemes and so the lower horizontal morphism in the square is an isomorphism. So for any section $s$ in $\rs F(X)$ there exists a section $t$ in $F_{\ldh}(X)$ which agrees with $s$ in $\rs F(X')$. Now $\rs F$ is a subpresheaf of $F_{\dis}$ and $F_{\dis}$ is $\ldh$ separated (Remark~\ref{rema:ellcdmono}) and hence $\rs F$ is $\ldh$ separated so $t$ agrees with $s$ in $\rs F(X)$.

\emph{The trace morphisms on $F_{\dis}$ preserve the subgroups $\rs F(X)$.}
It suffices to show that for $f: Y \to X$ a finite flat surjective morphism of schemes and $s \in \rs F(Y) \subseteq F_{\dis}(Y)$ the image $\Tr_f(s) \in F_{\dis}(X)$ is in the subgroup $\rs F(X) \subseteq F_{\dis}(X)$. Let $T$ be a regular integral semi-local scheme of dimension one and $T \to X$ a morphism. We must find a section in $F(T)$ that agrees with $\Tr_f(s)$ in $F_{\dis}(T)$.

Let $\widetilde{T \times_X Y}$ be the normalisation of $T \times_X Y$, let $T_i \to T \times_X Y$ be the inclusions of the integral components of $T \times_X Y$. Since $T$ is regular and integral of dimension one, the induced morphisms $T_i \to T$ are finite flat surjective. Let $\widetilde{T_i}$ be their normalisations and $n_i: \widetilde{T}_i \to T \times_X Y$ and $f_i: \widetilde{T}_i \to T$ the canonical morphisms. Since $F_{\dis}$ satisfies (Tri1) and (Tri2) (Proposition~\ref{prop:sepImpliesTrin}) we have
\[ \Tr_{T \times_X f} = \sum m_i \Tr_{f_i} F_{\dis}(n_i). \]

Now consider the following diagram.
\[ \xymatrix@!=12pt{
F(\widetilde{T \times_X Y}) \ar[drr] \ar@/_6pt/[ddrr] && &&  \\
&& F_{\dis}(\widetilde{T \times_X Y}) \ar@/_6pt/[ddrr] && F_{\dis}(T \times_X Y) \ar[dd] \ar[ll] && F_{\dis}(Y) \ar[dd] \ar[ll] \\
&& F(T) \ar[drr]&&&& \\
&&&& F_{\dis}(T) && \ar[ll] F_{\dis}(X)
} \]
Since $s \in \rs F(Y)$, we can find a section $(t_i) \in \oplus F(\widetilde{T}_i) = F(\widetilde{T \times_X Y})$ that agrees with $s$ in $F_{\dis}(\widetilde{T \times_X Y})$. Write $(s_i) \in \oplus F_{\dis}(\widetilde{T}_i) = F_{\dis}(\widetilde{T \times_X Y})$ for the image of $s$ in this group. Due to (Tri1) and (Tri2) for $F_{\dis}$ the image of $s$ in $F_{\dis}(T)$ is equal to the image of $(m_i s_i)$ in $F_{\dis}(T)$. Hence, $\sum m_i \Tr_{f_i}(t_i) \in F(T)$ is a section which agrees with the image of $s$ in $F_{\dis}(T)$.
\end{proof}


\section{From traces to transfers} \label{sec:tracesToTransfers}

\subsection{Statement and strategy}

In this section $Sch(S)$ denotes the category of separated schemes of finite type over a separated noetherian base $S$. This section is independent from the rest of the chapter; we only use properties of the category $Cor(S)$ (and Theorem~\ref{theo:CD} which just cites a result from \cite{CD}).

In this section we prove the following theorem.

\begin{theo} \label{theo:cdhTriImpliesTransfers}
Consider the functor $Shv_{cdh}(Cor(S)) \to PreShvTra(Sch(S))$ that takes a cdh sheaf with transfers to its corresponding presheaf with traces (Lemma~\ref{lemm:transfersImpliesTraces}). This functor is fully faithful, and its essential image is the full subcategory of cdh sheaves with traces which satisfy (Tri1) and (Tri2).
\end{theo}

The idea is that up to cdh refinement, each correspondence\footnote{A correspondence is by definition a morphism in the category $Cor(S)$. That is, an element of $c_{equi}(X \times_S Y / X, 0)$ for some schemes $X, Y \in Sch(S)$ (cf. Definition~\ref{defi:categoryOfCorrespondances}).} $\alpha \in \hom_{Cor(S)}([X], [Y])$ is of the following form.

\begin{defi} \label{defi:fn}
We say that a correspondence $\alpha \in \hom_{Cor(S)}([X], [Y])$ is of the form (FN) if
\begin{enumerate}
 \item[(FN)] there exists integers $n_i$, and closed integral subschemes $Z_i$ of $X \times Y$, such that the morphisms $g_i: Z_i \to X$ induced by the first projection are flat and finite, and 
\[ \alpha = \sum n_i [f_i] \circ {^t[g_i]} \]
where $f_i: Z_i \to Y$ are the morphisms induced by the second projection.
\end{enumerate}
\end{defi}

In this section the brackets $[-]$ and the composition sign $\circ$ will quickly become tiresome and so we omit them. So for example, if $p: X \to Y$ is a morphism of schemes and $\alpha \in \hom_{Cor(S)}([Y], [W])$ a correspondence, instead of $\alpha \circ [p] \in \hom_{Cor(S)}([X], [W])$ we will write $\alpha p \in \hom_{Cor(S)}(X, W)$. We will also use the notation $\alpha: X \corr Y$ to indicate that $\alpha \in \hom_{Cor(S)}(X, Y)$.

The strategy is the following.

\begin{enumerate}
 \item The definition:
 \begin{enumerate}
 \item (Definition~\ref{defi:corrfn}) If $\alpha: X \corr Y$ is of the form (FN) then we define $F(\alpha): F(Y) \to F(X)$ as $\sum n_i \Tr_{g_{i}} F(f_i)$.
 \item (Lemma~\ref{lemm:defiUnique}) In general, for a correspondence $\alpha: X \corr Y$ we define $F(\alpha): F(Y) \to F(X)$ as the unique morphism such that: for every cdh cover $p: X' \to X$ such that $\alpha p$ is of the form (FN) we have $F(p)F(\alpha) = F(\alpha p)$.
 \end{enumerate}
 \item We then need to show (Proposition~\ref{prop:mainProp}): If $X \stackrel{\alpha}{\corr} Y \stackrel{\beta}{\corr} Z$ is a pair of composable correspondences then $F(\alpha)F(\beta) = F(\beta\alpha)$.
 \item To do this, by the definition we need to put $\beta$, $\alpha$ and $\beta\alpha$ in the form (FN). Once we have the appropriate commutative diagram in the category of correspondences (Diagram~\ref{equa:commDiagramCorr}), we show $F(\alpha)F(\beta) = F(\beta\alpha)$ using the properties:
\begin{enumerate}
 \item (Lemma~\ref{lemm:postCompWithAMorphism}) For $\alpha: X \corr Y$ a correspondence and $f: X' \to X$ any morphism of schemes, we have $F(f)F(\alpha) = F(\alpha f)$.
 \item (Lemma~\ref{lemm:preCompWithAMorphism}) For $\alpha: X \corr Y$ a correspondence of the form (FN) and $g: Y \to Y'$ a morphism of schemes, we have $F(\alpha)F(g) = F(g \alpha)$.
 \item (Lemma~\ref{lemm:compFN}) If $X \stackrel{\alpha}{\corr} Y \stackrel{\beta}{\corr} Z$ is a pair of composable correspondences such that $\alpha, \beta, \beta\alpha,$ and $\Gamma_{\beta} \alpha$ (see Definition~\ref{defi:graph}) are of the form (FN) then $F(\alpha)F(\beta) = F(\beta\alpha)$.
\end{enumerate}
 \item (Lemma~\ref{lemm:fullyFaithfulCorrTraces}) Showing fully faithfulness is a straightforward reduction to correspondences of the form (FN).
\end{enumerate}

\subsection{Proof}

\begin{defi} \label{defi:corrfn}
Suppose that $F$ is a presheaf with traces and $\alpha: X \corr Y$ is a correspondence of the form (FN). We define $F(\alpha)$ as $\sum n_i \Tr_{g_{i}} F(f_i)$.
\end{defi}

\begin{lemm} \label{lemm:platificationCorr}
For every correspondence $\alpha: X \corr Y$ there exists a cdh covering $p: X' \to X$ such that $\alpha p$ is of the form (FN). Moreover, if $\alpha_i: X \corr Y_i$ is a finite family of correspondences, we can find a cdh cover $p: X' \to X$ such that each $\alpha_i p$ is of the form (FN).
\end{lemm}

\begin{proof}
Suppose that $\alpha_i = \sum n_{ij} z_{ij}$. By the platification theorem (Theorem~\ref{theo:platification}) there exists a blowup with nowhere dense centre $X_0 \to X_{red}$ such that the proper transform of $\amalg \overline{\{z_{ij}\}} \to X_{red}$ is flat over $X'$ where $\overline{\{z_{ij}\}}$ is the closure of $z_{ij}$ in $X \times_S Y_i$. Let $p$ be the composition $p: X_0 \to X_{red} \to X$. Hence (Lemma~\ref{lemm:agreesWithOldPullbacks}) each $\alpha_i p$ is of the form (FN). We let $i_0: W_0 \to X$ be a closed subscheme such that $W_0 \amalg X_0 \to X$ is a cdh cover and then repeat with $\alpha_i i_0$. Eventually we end up with a reduced $W_n$ of dimension zero and every correspondence with source a reduced scheme of dimension zero is of the form (FN). So by induction on the dimension we are done.
\end{proof}

\begin{lemm} \label{lemm:alphaptwoFN}
Let $F$ be a presheaf with traces that satisfies (Tri1). Let $\alpha: X \corr Y$ be a correspondence of the form (FN) and $p: X' \to X$ a morphism such that $\alpha p$ is also of the form (FN). Then $F(p)F(\alpha) = F(\alpha p)$.
\end{lemm}

\begin{proof}
Let $\alpha = \sum n_i z_i$. Since $\alpha$ is of the form (FN) we have $\alpha p = \sum n_i m_{ij} w_{ij}$ where $w_{ij}$ are the generic points of $X' \times_X \overline{\{z_i\}}$ and $m_{ij}$ the lengths of their local rings. We have diagrams such as
\[ \xymatrix{
\overline{\{w_{ij}\}} \ar[r]^{p_{ij}''} \ar[dr]_{g_{ij}''} & X' \times_X \overline{\{z_i\}} \ar[r]^{p'} \ar[d]^{g_i'} & \overline{\{z_i\}} \ar[r]^{f_i} \ar[d]^{g_i} & Y \\
& X' \ar[r]_p & X
} \]
where $g_i, g_i', g_{ij}''$ are flat. We then have
\[ \begin{split}
F(p)F(\alpha) &= F(p)\sum n_i \Tr_{g_{i}}F(f_i) \\
&= \sum n_i F(p)\Tr_{g_{i}}F(f_i) \\
&\stackrel{(CdB)}{=} \sum n_i \Tr_{g_{i}}' F(p')F(f_i) \\
&\stackrel{(Tri1)}{=} \sum n_i \biggl ( \sum m_{ij} \Tr_{g_{ij}''}F(p_{ij}'') \biggr ) F(p')F(f_i) \\
&= \sum n_i m_{ij} \Tr_{g_{ij}''} F(f_ip'p_{ij}'') \\
&= F(\alpha p) 
\end{split}  \]
\end{proof}

\begin{lemm} \label{lemm:defiUnique}
Suppose that $F$ is a cdh sheaf with traces that satisfies (Tri1) and $\alpha: X \corr Y$ a correspondence. There exists a unique morphism $F(\alpha): F(Y) \to F(X)$ such that: for every cdh cover $f: X' \to X$ such that $\alpha f$ is of the form (FN) we have $F(f)F(\alpha) = F(\alpha f)$.
\end{lemm}

\begin{proof}
There always exists such an $f$ (Lemma~\ref{lemm:platificationCorr}). Chose one. Let $p, q: X'' = X' \times_X X' \to X'$ be the two projections. We have $\alpha f p = \alpha f q$ and we chose another cdh cover $g: W \to X''$ such that $\alpha f p g$ (and hence $\alpha f q g$ as well) is of the form (FN). Lemma~\ref{lemm:alphaptwoFN} tells us then that $F(pg)F(\alpha f) = F(\alpha f p g) = F(\alpha f p g) = F(qg)F(\alpha f)$ and hence $F(g)F(p) F(\alpha f) = F(g)F(q) F(\alpha f)$. Since $F$ is a cdh sheaf and $g$ is a cdh cover, it follows that $F(p)F(\alpha f) = F(q)F(\alpha f)$. Again, $F$ is a cdh sheaf and so this implies that the morphism $F(\alpha f): F(Y) \to F(X')$ factors uniquely as $F(Y) \to F(X) \stackrel{F(f)}{\to} F(X')$. So we have found our $F(\alpha): F(Y) \to F(X)$ and it remains to show that it is independent of the choice of $f$.

Every two covers such as $f$ that put $\alpha$ in the form (FN) are dominated by a third one (Lemma~\ref{lemm:platificationCorr}) and so to show the independence it suffices to consider the case of two covers $f_0, f_1$ with a factorisation $f_1: X_1' \stackrel{h}{\to} X_0' \stackrel{f_0}{\to} X$. We will write $F(\alpha)_{f_0}$ and $F(\alpha)_{f_1}$ momentarily for the two morphisms $F(Y) \to F(X)$ induced by $f_0$ and $f_1$ respectively. We have
\[ F(f_1) F(\alpha)_{f_0} = F(h) F(f_0) F(\alpha)_{f_0} = F(h) F(\alpha f_0) \stackrel{(\ref{lemm:alphaptwoFN})}{=} F(\alpha f_0 h) = F(\alpha f_1) \]
and hence uniqueness of the factorisation in the definition of $F(\alpha)_{f_1}$ implies that $F(\alpha)_{f_1} = F(\alpha)_{f_0}$.
\end{proof}

\begin{lemm} \label{lemm:postCompWithAMorphism}
The morphism defined in Lemma~\ref{lemm:defiUnique} associated to a correspondence $\alpha: X \corr Y$ satisfies: for any morphism of schemes $f: X' \to X$ we have $F(\alpha f) = F(f)F(\alpha)$.
\end{lemm}

\begin{proof}
Choose a cdh cover $p: U \to X$ such that $\alpha p$ is of the form (FN). We consider the cartesian square
\[ \xymatrix{
U' \ar[r]^{p'} \ar[d]_{f'} & X' \ar[d]^f \\
U \ar[r]_p & X
} \]
and find a cdh cover $q: U'' \to U'$ such that $\alpha f p' q$ is of the form (FN). Note that the cdh cover $p'q$ puts $\alpha f$ in the form (FN). We have the following commutative diagram
\[ \xymatrix{
F(Y) \ar[d]_{F(\alpha)} \ar[dr]^{F(\alpha p)} \\
F(X) \ar[r]_{F(p)} \ar[d]_{F(f)} & F(U) \ar[d]^{F(f')} \ar[dr]^{F(f'q)} \\
F(X') \ar[r]_{F(p')} & F(U') \ar[r]_{F(q)} & F(U'')
} \]
It follows that
\[ \begin{split}
F(p'q)F(f)F(\alpha) &= F(q)F(p')F(f)F(\alpha) = F(q)F(f')F(p)F(\alpha) = F(q)F(f')F(\alpha p) \\
&= F(f'q)F(\alpha p) \stackrel{\ref{lemm:alphaptwoFN}}{=} F(\alpha p f'q) = F(\alpha f p'q)  
\end{split} \]
Since by definition $F(\alpha f)$ is the unique morphism that satisfies $F(p'q)F(\alpha f) = F(\alpha f p' q)$ it follows that $F(f)F(\alpha) = F(\alpha f)$.
\end{proof}

\begin{lemm} \label{lemm:preCompWithAMorphism}
If $F$ satisfies (Tri2) as well, the morphisms $F(\alpha)$ from Lemma~\ref{lemm:defiUnique} satisfy: if $\alpha: X \corr Y$ is of the form (FN) and $f: Y \to W$ is a morphism of schemes we have $F(\alpha)F(f) = F(f \alpha)$.
\end{lemm}

\begin{proof}
By the definition of $F(f \alpha)$, we must show that $F(f\alpha p) = F(p)F(\alpha)F(f)$ where $p: X' \to X$ is a covering such that $f\alpha p$ is of the form (FN).

Let $\alpha = \sum n_i z_i$. Since the $\overline{\{z_i\}} \to X$ are flat, $\alpha p$ is the correspondence $\sum n_i m_{ij} w_{ij}$ where the $w_{ij}$ are the generic points of $\overline{\{z_i\}} \times_X X'$ and $m_{ij}$ the lengths of their local rings. Let $f': X' \times Y \to X \times W$ be the morphism induced by $f$ and $d_{ij} = [k(w_{ij}) : k(f'w_{ij})]$ so the correspondence $f \alpha p$ is $\sum n_i m_{ij} d_{ij} f'w_{ij}$. We have diagrams
\[ \xymatrix{ 
\overline{\{f'w_{ij}\}} \ar[r]^{h_{ij}'} \ar@/_12ex/[ddr]_{s_{ij}} & X' \times W \ar[rr]^{pr} & & W \\
\overline{\{w_{ij}\}} \ar[r]^{h_{ij}} \ar[dr]_{r_{ij}} \ar[u]_{f_{ij}'} & \overline{\{z_i\}} \times_X X' \ar[r]^{p_i} \ar[d]^{q_i'} \ar[u]^{f'} & \overline{\{z_i\}} \ar[r]^{g_i} \ar[d]^{q_i} & Y \ar[u]_f \\
& X' \ar[r]_p & X
} \]
and
\[ \begin{split}
F(p)F(\alpha)F(f) &= \sum n_i F(p)\Tr_{q_{i}}F(g_i)F(f) \\
&\stackrel{(CdB)}{=} \sum n_i \Tr_{q_{i}'}F(p_i)F(g_i)F(f) \\
&\stackrel{(Tri1)}{=} \sum n_i m_{ij}\Tr_{r_{ij}}F(h_{ij})F(p_i)F(g_i)F(f) \\
&= \sum n_i m_{ij}\Tr_{r_{ij}}F(f_{ij}')F(h_{ij}')F(pr) \\
&\stackrel{(Tri2)}{=} \sum n_i m_{ij} d_{ij} \Tr_{s_{ij}}F(h_{ij}') F(pr) \\
&= F(f \alpha p).
\end{split} \]
\end{proof}

We will use the following definition to apply Lemma~\ref{lemm:platificationCorr} in the proof of Proposition~\ref{prop:mainProp}.

\begin{defi} \label{defi:graph}
Suppose that $\alpha: X \corr Y$ is a correspondence of the form (FN) and $\alpha = \sum n_i z_i$. The canonical morphism $\delta: X \times Y \to X \times X \times Y$ is a closed immersion and since the $\overline{\{z_i\}} \to X$ are flat, $\sum n_i \delta(z_i)$ defines a correspondence $X \corr X \times Y$. We will denote this correspondence by
\[ \Gamma_\alpha = \sum n_i \delta(z_i) : X \corr X \times Y. \]
\end{defi}

\begin{rema}
The reader familiar with Voevodsky and Suslin's theory of relative cycles will recognise $\Gamma_\alpha$ as the external product of $[\id]$ and $\alpha$. We actually don't need the condition that $\alpha$ is of the form (FN) but in our application of this definition we will have it.
\end{rema}

\begin{lemm} \label{lemm:compFN}
Let $F$ be a cdh sheaf with traces that satisfies (Tri1) and (Tri2). Suppose that $X \stackrel{\alpha}{\corr} Y \stackrel{\beta}{\corr} Z$ are a pair of composable correspondences and suppose further that $\alpha, \beta, \beta\alpha$, and $\Gamma_\beta\alpha$ (see Definition~\ref{defi:graph}) are all of the form (FN).

Then $F(\alpha)F(\beta) = F(\beta\alpha)$.
\end{lemm}

\begin{proof}
It is enough to consider the case where $\beta = w$ and $\alpha = z$ are formal sums consisting of a single point with multiplicity one. Our diagram is
\[ \xymatrix{
\overline{\{v_{i}\}} \ar[r] \ar[d] \ar@/_12pt/[dd] & \overline{\{z\}} \times_Y \overline{\{w\}} \ar[r] \ar[d] & \overline{\{w\}} \ar[r]_\iota \ar[d] & Z \\
\overline{\{v_i'\}} \ar[rrru] \ar[d] & \overline{\{z\}} \ar[r] \ar[dl] & Y \\
X
} \]
where $v_i$ are the generic points of $\overline{\{z\}} \times_Y \overline{\{w\}}$ and the $v_i'$ their images in $X \times Z$. Lets say that $\ell_i = length(\OO_{\overline{\{z\}} \times_Y \times \overline{\{w\}}, v_{i}})$ and $d_{i} = [k(v_{i}) : k(v_{i}')]$. Since $\alpha$ and $\beta$ are of the form (FN) we have $\beta\alpha = \sum \ell_i d_i v_i'$ and $\Gamma_\beta\alpha = \sum \ell_i v_i$ so by hypothesis the schemes $\overline{\{v_i\}}$ and $\overline{\{v'_i\}}$ are flat over $X$. It now suffices to apply (CdB), functoriality, and (Tri1) and (Tri2) to see that $F(\alpha)F(\beta) = F(\beta\alpha)$.
\end{proof}

The following proposition is the cdh analogue of \cite[3.1.5]{Voev00} (which incidentally follows from \cite[3.1.3]{Voev00} via the same argument we use here).

\begin{prop} \label{prop:cdhpullback}
Let $\alpha: X \corr Y$ be a correspondence and $p: Y' \to Y$ a cdh cover. Then there exists a correspondence $\alpha' : X' \corr Y'$ and a cdh cover $p': X' \to X$ such that the square
\[ \xymatrix{
X' \ar[r]^{\alpha'} \ar[d]_{p'} & Y' \ar[d]^p \\
X \ar[r]_\alpha & Y
} \]
commutes in $Cor(S)$.
\end{prop}

\begin{proof}
This follows from Theorem~\ref{theo:CD}(\ref{theo:CD:enum:1}).

More explicitely, we have the following three elementary facts which hold for any Grothendieck pretopology $\tau$.
\begin{enumerate}
 \item A morphism of $\tau$ sheaves $G \to F$ is surjective as a morphism of $\tau$ sheaves if and only if: for every object $X$ and every section $s \in F(X)$ there exists a cover $X' \to X$ and a section $t \in G(X')$ such that the images of $t$ and $s$ agree in $F(X')$ via the morphisms in the obvious commutative square.
 \item For any presheaf $F$, object $X$, and section $s \in F_\tau(X)$ there exists a $\tau$ cover $U \to X$ such that $s|_U$ is in the image of $F(U) \to F_\tau(U)$.
 \item For any presheaf $F$ and object $X$, two sections $s, t \in F(X)$ are sent to the same section of $F_\tau(X)$ if and only if there exists a $\tau$ cover $U \to X$ such that the two restrictions $s|_U, t|_U$ are equal in $F(U)$.
\end{enumerate}

The result Theorem~\ref{theo:CD}(\ref{theo:CD:enum:1}) implies that $L_{cdh}(Y') \to L_{cdh}(Y)$ is a surjective morphism of cdh sheaves. The correspondence $\alpha: X \corr Y$ gives a section in $L(Y)(X)$ and hence in $L_{cdh}(Y)(X)$. By the first fact mentioned above, we find a cdh cover $W \to X$ and an element $t \in L_{cdh}(Y')(W)$ such that the images of $t$ and $\alpha$ agree in $L_{cdh}(Y)(W)$. By the second fact mentioned above, we find a cdh cover $V \to W$ and a section $u \in L(Y')(V)$ whose image in $L_{cdh}(Y')(V)$ agrees with the restriction of $t$. Finally, by the third fact, there is a cdh cover $X' \to V$ such that the restriction $\alpha|_{X'}$ of $\alpha$ in $L(Y)(X')$ agrees with the restriction $L(p)u|_{X'}$ of the image $L(p)u$ of $u$. That is, we have found a section $\alpha' = u|_{X'} \in L(Y')(X')$ whose image $L(p)(\alpha')$ in $L(Y)(X')$ agrees with the restriction $\alpha|_{X'} \in L(Y)(X')$ where $X' \to X$ is a cdh cover. This is equivalent to the desired commutative square.
\end{proof}

\begin{prop} \label{prop:mainProp}
Let $F$ be a cdh sheaf with traces that satisfies (Tri1) and (Tri2). Then there exists a unique structure of presheaf with transfers on $F$ such that for the correspondences $\alpha$ of the form (FN) the morphism $F(\alpha)$ is that described in Definition~\ref{defi:corrfn}.
\end{prop}

\begin{proof}
Let $X \stackrel{\alpha}{\corr} Y \stackrel{\beta}{\corr} Z$ be a pair of composable correspondences. To prove that $F(\alpha)F(\beta) = F(\beta\alpha)$ we need to show that if $q: U \to X$ is a cdh cover such that $\beta \alpha q$ is of the form (FN), then $F(q)F(\alpha)F(\beta) = F(\beta\alpha q)$.

Suppose that $p: Y' \to Y$ is a cdh cover such that $\beta p$ is of the form (FN), and suppose that
\[ \xymatrix{
X' \ar[r]^{\alpha'} \ar[d]_{p'} & Y' \ar[d]^p \\
X \ar[r]_\alpha & Y
} \]
is a commutative square as in Proposition~\ref{prop:cdhpullback}. By composing with a further cdh cover $X'' \to X'$ we can assume (Lemma~\ref{lemm:platificationCorr}) that $\alpha', \beta p \alpha'$ and $\Gamma_{\beta p} \alpha'$ are of the form (FN). The commutative diagram is the following.
\begin{equation} \label{equa:commDiagramCorr}
\xymatrix{
X' \ar[r]^{\alpha'} \ar[d]_{p'} & Y' \ar[d]^p \ar[dr]^{\beta p}\\
X \ar[r]_\alpha & Y \ar[r]_\beta & Z
} \end{equation}
We now have
\[ \begin{split}
F(p')F(\alpha)F(\beta) &\stackrel{\footnotesize\ref{lemm:postCompWithAMorphism}}{=} F(\alpha p') F(\beta)\\
&= F(p\alpha')F(\beta) \\
&\stackrel{\footnotesize\ref{lemm:preCompWithAMorphism}}{=} F(\alpha')F(p)F(\beta) \\
&\stackrel{\footnotesize\ref{lemm:postCompWithAMorphism}}{=} F(\alpha')F(\beta p) \\
&\stackrel{\footnotesize\ref{lemm:compFN}}{=} F(\beta p\alpha') \\
&= F(\beta \alpha p').
\end{split} \]
\end{proof}

\begin{lemm} \label{lemm:fullyFaithfulCorrTraces}
The functor that associates a presheaf with traces to a cdh sheaf with transfers $Shv_{cdh}(Cor(S)) \to PreShvTra(Sch(S))$ is fully faithful.
\end{lemm}

\begin{proof}
Suppose that $F$ and $G$ are two cdh sheaves with transfers and $\phi: F \to G$ is a morphism of presheaves with traces. For a correspondence $\alpha: X \corr Y$ of the form (FN) it is clear that we have $\phi_X F(\alpha) = G(\alpha) \phi_Y$. If $\alpha$ is not of the form (FN) then there exists a cdh cover $p: X' \to X$ such that $\alpha p$ is of the form (FN) (Lemma~\ref{lemm:platificationCorr}). Now the commutativity of the outside rectangle, the rightmost square, and the injectivity of $G(X) \to G(X')$ implies the commutativity of the leftmost square in the following diagram
\[ \xymatrix{
F(Y) \ar[r]^{F(\alpha)} \ar[d]_{\phi_Y} & F(X) \ar[d]^{\phi_X} \ar[r]^{F(p)} & F(X') \ar[d]^{\phi_{X'}} \\
G(Y) \ar[r]_{G(\alpha)} & G(X) \ar[r]_{G(p)} & G(X')
} \]
\end{proof}

\section{Summary}

In this last section we collect the main results of this chapter. Depending on the context we will want them in various different forms. In this last section
\begin{enumerate}
 \item[] $S$ is a quasi-excellent separated noetherian scheme,
 \item[] $Sch(S)$ is the category of separated $S$-schemes of finite type,
 \item[] $Sm(S)$ is the full subcategory of smooth schemes in $Sch(S)$,
 \item[] $Reg(S)$ is the full subcategory of regular schemes in $Sch(S)$, and
 \item[] $EssSch(S)$ is the category of schemes which are inverse limits of left filtering systems in $Sch(S)$ for which each of the transition morphisms are affine open immersions.
\end{enumerate}
We remind the reader that if $j^{ess}: PreShv(Sch(S)) \to PreShv(EssSch(S))$ is the left Kan extension along $j: Sch(S)^{op} \to EssSch(S)^{op}$ then for any presheaf with traces $F$, the arguments in \cite[Section 8]{EGAIV3} give a canonical structure of presheaf with traces on $j^{ess}F$ such that $F \to (j^{ess}F) \circ j$ is a morphism of presheaves with traces.

The following theorem is the most general collection of the results in this chapter.

\begin{theo} \label{theo:bigComparisonTheorem}
Let $\ell$ be a prime invertible on $S$. Suppose that $F$ is a Nisnevich Gersten sheaf of $\zll$ modules with traces on $EssSch(S)$ that satisfies (Tri1)$_{\leq 0}$. Then
\begin{enumerate}
 \item the canonical morphism $F_{cdh} \to F_{\ldh}$ is an isomorphism and on regular schemes $X$ in $EssSch(S)$ we have $F(X) = F_{cdh}(X) = F_{\ldh}(X)$,
 \item for every $n \in \ZZ_{\geq 0}$ and $X \in Sch(S)$ the canonical morphism $H_{cdh}^n(X, F_{cdh}) \to H_{\ldh}^n(X, F_{\ldh})$ is an isomorphism, and 
 \item each of the presheaves $F|_{Reg(S)}$, $F_{\ldh}|_{Sch(S)}$, and $H_{\ldh}^n(-, F_{\ldh})|_{Sch(S)}$ has a canonical structure of presheaf with transfers.
\end{enumerate}
\end{theo}

\begin{proof}
The first statement is just Corollary~\ref{coro:cdhliso} and the second part of Proposition~\ref{prop:tracesForcdh}. Now the first part of Proposition~\ref{prop:tracesForcdh} says that $F_{cdh}$ has a structure of traces satisfying (Tri1) and (Tri2) and so we can apply Theorem~\ref{theo:cdhTriImpliesTransfers} to get a structure of transfers. Now that we have a structure of transfers, Theorem~\ref{theo:lcdhcohomologyAgree} says that the cohomologies agree. Finally, the structure of transfers on $F_{cdh} = F_{\ldh}$ has already been mentioned. The structure of transfers on $F|_{Reg(S)}$ comes from the isomorphism $F|_{Reg(S)} = F_{\ldh}|_{Reg(S)}$ and the structure of transfers on the cohomology is part of Theorem~\ref{theo:lcdhcohomologyAgree}.
\end{proof}

The following theorem is designed to be applied to the homotopy presheaves of an oriented object in the Morel-Voevodsky stable homotopy category $\SH(k)$.

\begin{theo} \label{theo:cdhlprimeComparison}
Let $k$ be a perfect field and $\ell$ a prime that is invertible in $k$. Let $F$ be a presheaf of $\zll$-modules with traces on $Sch(k)$, such that
\begin{enumerate}
 \item $F(X) \to F(X_{red})$ is an isomorphism for every $X \in Sch(k)$,
 \item $F(X) \to F(\AA^1_X)$ is an isomorphism for every $X \in Sm(k)$, and
 \item $F|_{Sm(k)}$ has a structure of transfers,
\end{enumerate}
then for every $n \in \ZZ_{\geq 0}$ and every $X \in Sch(S)$, the canonical morphism
\[ H^n_{cdh}(X, F_{cdh}) \to H^n_{\ldh}(X, F_{\ldh}) \]
is an isomorphism.
\end{theo}

\begin{proof}
After \cite[Theorem 3.1.12]{Voev00}, Proposition~\ref{prop:tracesForNis}, and \cite[Theorem 18.1.2]{EGAIV4} we can assume that $F$ is a Nisnevich sheaf. This implies in particular, that it is a Gersten sheaf after we extend it to $EssSch(k)$ (\cite[4.37]{Vctpt}). So $F_{cdh} \cong F_{\ldh}$ (Corollary~\ref{coro:cdhliso}). Then, as before, $F_{cdh}$ has a structure of traces satisfying (Tri1) and (Tri2) (Proposition~\ref{prop:tracesForcdh}) and so we can apply Theorem~\ref{theo:cdhTriImpliesTransfers} to get a structure of transfers. Theorem~\ref{theo:lcdhcohomologyAgree} then tells us that the cohomologies agree.
\end{proof}

\begin{rema}
We can avoid Theorem~\ref{theo:cdhTriImpliesTransfers} in the proof of Theorem~\ref{theo:cdhlprimeComparison} by considering directly the left Kan extension of $F|_{Sm(k)}$ along $SmCor(k)^{op} \to Cor(k)^{op}$. We will see this in the proof of Proposition~\ref{prop:nisellcomparison}. In fact, if we use this technique, then we only need the $0 \to F(X) \to \prod_{x \in X^{(0)}} F(x)$ part of the Gersten sequence. So in fact, we could replace the assumption that $F$ is homotopy invariant, with the assumption that $F_{Nis}(X) \to \prod_{x \in X^{(0)}} F(x)$ is injective. However, we can't avoid the assumption that $F$ has traces on $Sch(k)$ because we need this to get $F_{cdh} \cong F_{\ldh}$.
\end{rema}

Finally we have the following proposition which is useful for working with Voevodsky motives.

\begin{prop} \label{prop:nisellcomparison}
Suppose that $k$ is a perfect field and $\ell$ a prime invertible in $k$. Let $F$ be a homotopy invariant presheaf with transfers of $\zll$ modules on $Sm(k)$. Then the canonical morphism $F_{Nis} \to F_{\ldh}$ is an isomorphism.

In particular, every homotopy invariant Nisnevich sheaf of $\zll$ modules with transfers is an $\ldh$ sheaf.
\end{prop}

\begin{rema}
Note that Voevodsky's work tells us that $F_{Nis}$ is a homotopy invariant presheaf with transfers \cite[Theorem 3.1.12]{Voev00}.
\end{rema}

\begin{proof}
We can assume that $F$ is a Nisnevich sheaf after \cite[Theorem 3.1.12]{Voev00}. This implies in particular that $F$ is a Gersten presheaf (\cite[4.37]{Vctpt}). Let $i: SmCor(k) \to Cor(k)$ be the canonical morphism and consider $i^*F$ the left Kan extension of $F$ along $i$. Then since $F = (i^*F)\circ i$, the presheaf $i^*F$ is still a Gersten presheaf. Moreover, it is by definition a presheaf with transfers and therefore a presheaf with traces that satisfies (Tri1)$_{\leq 0}$ (Lemma~\ref{lemm:transfersImpliesTraces}). We then apply Corollary~\ref{coro:cdhliso} and Proposition~\ref{prop:tracesForcdh}.
\end{proof}

\chapter{Traces and the slice filtration} \label{chap:slices}

\section{Introduction}

\lettrine[lines=2, findent=3pt, nindent=0pt]{R}{ecall} Theorem~\ref{theo:cdhlprimeComparison} from the last chapter that we want to apply to homotopy presheaves of oriented objects in the Morel-Voevodsky stable homotopy category $\SH(k)$. The piece that we don't have is a structure of traces on these presheaves for non-smooth schemes. The goal of this chapter is to address this. In particular, we want to have a structure of non-smooth traces on the homotopy presheaves of $\HZl$ where $\HZl$ is the object that represents motivic cohomology with $\zll$-coefficients and $\ell$ is invertible in the perfect base field $k$. We do this via the slice filtration.

The goal of Section~\ref{sec:slice} is to develop the necessary material to have the isomorphisms $s_0f_*f^*E \cong f_*f^*s_0E$ that we will use to transfer a structure of traces on $\KH$ -- the object that represents algebraic $K$-theory -- to a structure of traces on $\HZ$. We begin by recalling the definition of the slice filtration on the Morel-Voevodsky stable homotopy category. We do this in the context of Ayoub's stable homotopy 2-functors as this language makes it much easier to discuss how the slice filtration interacts with the functors $f_\#, f^*,$ and  $f_*$. We state and prove a theorem of Pelaez using Ayoub's language which gives a criterion for a functor to behave well with respect to the slice filtration (Theorem~\ref{theo:criteria}). This is no more than a translation; the theorem and proof belong to Pelaez. We then apply this theorem to various functors in Lemma~\ref{coro:smoothPullback}. Some of these functors do not appear \cite{Pel12}, notably the functors associated to closed immersions.

We then introduce the notion of an object with a weak structure of smooth traces (Definition~\ref{defi:wsst}). This is essentially a structure of traces on smooth schemes, where the only axiom we ask for is (Deg). This is to apply the resolution of singularities argument in \cite{Pel12} using the theorem of Gabber on alterations, which we do in Theorem~\ref{theo:pullbackAll}. Finally, we want to apply the theorem of Pelaez (Theorem~\ref{theo:criteria}) to the functors $f_*$ and $f^*$ where $f$ is a finite flat surjective morphism between non-smooth schemes. The case when $f$ is {\'e}tale is already taken care of but the radicial case is trickier. It turns out to be easier if we consider instead the composition $f_*f^*$. We attain this in Proposition~\ref{prop:ffslices} after some lemmas.

In Section~\ref{sec:traces} we define what it means for an object $E \in \SH(S)$ to have a structure of traces. We show that a structure of traces on an object $E$ induces a structure of traces in the sense of Definition~\ref{defi:traces} on each of its homotopy presheaves (Lemma~\ref{lemm:tracesImpliesTraces}). We show that if we work $\zpi$-linearly, a structure of traces on an object induces a canonical structure of traces on the connective covers $f_qE$ and the slices $s_qE$ (Proposition~\ref{prop:tracesOnSlices}).

\section{On the functoriality of the slice filtration} \label{sec:slice}

\subsection{Preliminaries} \label{sec:prelim}

The material developed in the \cite{Pel12} holds in greater generality than it is presented in that article. We develop it in the setting of Ayoub's stable homotopy 2-functors \cite{Ayo07}. This makes the proofs involving the functoriality cleaner. We don't recall all the axioms of a stable homotopy 2-functor but we will recall the properties that we need, as we need them.

As such we will consider 2-functors $\H^*: Sch(S) \to TriCat$ to the 2-category of triangulated categories. For each scheme $X \in Sch(S)$ we set $\H(X) = \H^*(X)$ and for a morphism $f$ we set $f^* = H^*(f)$. We recall briefly that included in the definition of a 2-functor are 2-isomorphisms $c^*(f, g): (gf)^* \stackrel{\sim}{\to} f^*g^*$ for every two composable morphisms $f, g$ and these satisfy an appropriate coherency condition.

\begin{defi} \label{defi:sliceable2Functor}
We define a \emph{sliceable 2-functor} to be a quadruple $(S, \H^*, \Sigma, \mathcal{G})$ such that $S$ is a separated noetherian scheme,
\[ \H^*: Sch(S) \to TriCat \]
is a contravariant 2-functor to the 2-category of triangulated categories, $\Sigma: \H^* \to \H^*$ is an autoequivalence, and $\mathcal{G}$ is a set of compact objects in $\H^*(S)$. We require that:
\begin{enumerate}
 \item \label{item:smoothLeftAdjoint} If $f: Y \to X \in Sch(S)$ happens to be smooth then $f^*$ has a left adjoint 
\[ f_\#: \H(Y) \to \H(X). \]
 \item Each of the triangulated categories $\H(X)$ admits all small sums, and each of the functors $f^*$ preserves compact objects and small sums.
 \item For every $a: X \to S$ in $Sch(S)$ the category $\H(X)$ is generated by the set of compact objects $\{ \Sigma^n f_\#f^* a^*A : n \in \ZZ, Y \stackrel{f}{\to} X $ smooth, $A \in \mathcal{G} \}$ .
\end{enumerate}
\end{defi}

\begin{defi} \label{defi:sliceFiltration}
Let $(S, \H^*, \Sigma, \mathcal{G})$ be a sliceable 2-functor. For any scheme $X \stackrel{a}{\to} S \in Sch(S)$ we define
\[ \H(X)\eff \]
as the smallest full triangulated subcategory of $\H(X)$ containing all small sums and the set of objects 
\[ \{\Sigma^n f_\#f^* a^*A : n \in \ZZ_{\geq 0},  Y \stackrel{f}{\to} X \textrm{ smooth}, A \in \mathcal{G} \}. \]
More generally, we define $\Sigma^q\H(X)\eff$ as the smallest full triangulated subcategory of $\H(X)$ containing all small sums and the set of objects 
\[ \{\Sigma^n f_\#f^* a^*A : n \geq q,  Y \stackrel{f}{\to} X \textrm{ smooth}, A \in \mathcal{G} \}. \]
If $X$ is clear from the context we will write $\H\eff$ and $\Sigma^q\H\eff$. We obtain in this way a sequence of compactly generated triangulated subcategories, each containing small sums:
\[ \dots \hookleftarrow \Sigma^{q - 1}\H(X)\eff \hookleftarrow \Sigma^{q}\H(X)\eff \hookleftarrow \Sigma^{q + 1}\H(X)\eff \hookleftarrow \dots \]
\end{defi}

\begin{rema}
The question of how the subcategories $\Sigma^{q}\H(X)\eff$ behave with respect to the morphisms $f_\#$ and $f^*$, as well as various other functors is the main topic of this chapter.

We don't need the hypothesis (3) for any of the definitions but we will use it often and so we include it for expositional reasons.
\end{rema}

\begin{rema} \label{rema:exampleSH}
In practice, $\H$ will be the Morel-Voevodsky stable homotopy category $\SH$, the base scheme $S$ will be noetherian and $\S$ will be the category of schemes of finite type over $S$. The autoequivalence will be $(\PP^1, \infty) \wedge -$ smash with the projective line pointed at infinity and the set $\mathcal{G}$ will consist of a single object, $\un_S$ the unit in $\H(S)$ for the smash product.

If $\H = \SH$ is the Morel-Voevodsky stable homotopy category and $\un_S \in \SH(S)$ is the unit for the smash product then for any smooth $S$-scheme $f: X \to S$ there is a canonical isomorphism $f_\#f^*\un_S \cong \Sigma^\infty(X_+)$ functorial in $X$. In this case Definition~\ref{defi:sliceFiltration} is Voevodsky's original definition of the slice filtration.
\end{rema}

It is straightforward from properties of adjunctions that $\Sigma$ preserves compact objects and if $f: Y \to X$ is a smooth morphism then $f_\#$ preserves compact objects. In particular, if $A \in \H(S)$ is a compact object, $Y \stackrel{f}{\to} X$ a smooth morphism in $S$, and $X \stackrel{a}{\to} S$ the structural morphism of $X$, then $\Sigma^n f_\# f^* a^*A$ is also compact for all $n \in \ZZ$.


We have the following theorem of Neeman.

\begin{theo}[{\cite[Theorem 4.1]{Nee96}}] \label{theo:neeman}
Suppose that $\T$ is a compactly generated triangulated category, $\T'$ any other triangulated category, and $\T \to \T'$ a triangulated functor that preserves coproducts. Then $\T \to \T'$ has a right adjoint.
\end{theo}

Since each $\Sigma^{q}\H(X)\eff$ contains small sums and is generated by compact objects, by Theorem~\ref{theo:neeman} the inclusion $i_q: \Sigma^{q}\H(X)\eff \hookrightarrow \H(X)$ admits a right adjoint $r_q: \H(X) \to \Sigma^{q}\H(X)\eff$.

\begin{defi}
Let $(S, \H^*, \Sigma, \mathcal{G})$ be a sliceable 2-functor. For any scheme $X \in Sch(S)$ we define the endofunctor 
\[ f_q: \H(X) \to \H(X) \]
as the composition $f_q = i_qr_q$ where $r_q$ is the right adjoint to the inclusion $i_q: \Sigma^{q}\H(X)\eff \hookrightarrow \H(X)$. The counit of the adjunction $(i_q, r_q)$ gives us a natural transformation $f_q \to \id$.
\end{defi}

Evidently, for any object $E \in \H(X)$, the object $f_qE$ is in the subcategory $\Sigma^{q}\H(X)\eff$ and the morphism $f_qE \to E$ is determined up to unique isomorphism by the property that it induces an isomorphism $\hom(F, f_qE) \to \hom(F, E)$ for any object $F \in \Sigma^{q}\H(X)\eff$.

Due to our definitions, for any $q' < q$ the adjunction $(i_q, r_q)$ factors through the adjunction $(i_{q'}, r_{q'})$
\[ \Sigma^{q}\H(X)\eff \rightleftarrows \Sigma^{q'}\H(X)\eff \rightleftarrows \H(X) \]
and so since $r_{q'}i_{q'}$ is the identity (again due to our definitions) the morphism $f_qf_{q'} \to f_q$ is invertible. So we obtain a canonical morphism $f_q \to f_{q'}$. In fact, for each $q$, these form a functor $\ZZ_{\geq 0} \to End(\H(X))$ from the category associated to the totally ordered set $\ZZ_{\geq 0}$ to the category of endomorphisms of $\H(X)$ which sends $n \in \ZZ_{\geq 0}$ to $f_{q - n}$. This functor is equipped with a morphism of diagrams towards the constant diagram with value the identity endofunctor. The following lemmata are straightforward.

\begin{lemm}
Let $(S, \H^*, \Sigma, \mathcal{G})$ be a \emph{sliceable 2-functor}. Let $X \in Sch(S)$. For any $q \in \ZZ$, and any $E \in \H(X)$ the canonical morphism
\[ \hocolim_{{q'} \leq q} f_{q'}E  \to E \]
is an isomorphism.
\end{lemm}

\begin{proof}
To show that this morphism is an isomorphism it suffices to show it is an isomorphism after evaluating on $\hom(G, -)$ for each $G$ in
$\{ \Sigma^n f_\#f^* a^*A : n \in \ZZ, Y \stackrel{f}{\to} X $ smooth, $A \in \mathcal{G} \}$
By assumption, each $G$ is contained in some $\Sigma^p\H(X)\eff$  and so by the universal property of the morphism $f_pE' \to E'$ mentioned above, it suffices to show that this is an isomorphism after applying $\hom(G, f_p-)$ for a suitable $p$. The functor $f_p$ preserves small sums\footnote{The functor $i_p$ preserves small sums because it is a left adjoint, and $r_p$ preserves small sums because its left adjoint sends each object in a set of compact generating objects to a compact object.} and therefore it preserves homotopy colimits. So we have reduced to showing that
\[ \hom(G, \hocolim_{{q'} \leq q} f_pf_{q'}E) \to \hom(G, f_pE) \]
is an isomorphism which is clear since $f_pf_{q'} = f_p$ for all $q'$ sufficiently small.
\end{proof}

\begin{lemm} \label{lemm:defsq}
Let $(S, \H^*, \Sigma, \mathcal{G})$ be a \emph{sliceable 2-functor}. Let $X \in Sch(S)$. There exists for each $q$ a unique endofunctor $s_q$ together with morphisms $f_q \to s_q \to f_{q + 1}[1]$ such that for any object $E \in \H(X)$ the triangle
\[ f_{q + 1}E \to f_qE \to s_qE \to f_{q + 1}E[1] \]
is distinguished. For any object $F \in \Sigma^r\H(X)\eff$ with $r > q$ the group $\hom(F, s_{q}E)$ is zero.
\end{lemm}

\begin{proof}
Using the usual triangulated category techniques,\footnote{i.e., for every object $E$ we chose a cone of $f_{q + 1}E \to f_qE$ and for every morphism we chose a morphism between the cones and then to show that this defines a functor we show that the morphisms we have chosen are unique.} it suffices to show that for any pair of objects $E, F$ there are no non-zero morphisms from $f_{q + 1}E[1]$ to any cone of $f_{q + 1}F \to f_qF$. Now $f_{q + 1}E[1] \in \Sigma^{q + 1}\H(X)\eff$ (and also $\Sigma^{q}\H(X)\eff$) and so by the universal property of $f_q \to \id$ and $f_{q + 1} \to \id$ the two vertical morphisms in the triangle
\[ \xymatrix{
\hom(f_{q + 1}E[1], f_{q + 1}F) \ar[rr] \ar[dr] && \hom(f_{q + 1}E[1], f_{q}F) \ar[dl] \\
& \hom(f_{q + 1}E[1], F)
} \]
are isomorphisms. Hence, the third one is an isomorphism as well which implies that there are no non-zero morphisms from $f_{q + 1}E[1]$ to any cone of $f_{q + 1}E \to f_qE$.
\end{proof}

\begin{lemm}
Let $(S, \H^*, \Sigma, \mathcal{G})$ be a \emph{sliceable 2-functor}. Let $X \in Sch(S)$.  There exists for each $q$ a unique endofunctor $s_{< q}$ together with morphisms $\id \to s_{<q} \to f_q[1]$ such that for any object $E \in \H(X)$ the triangle
\[ f_q E \to E \to s_{<q}E \to f_qE[1] \]
is distinguished. For any object $F \in \Sigma^r\H(X)\eff$ with $r \geq q$ the group $\hom(F, s_{<q}E)$ is zero.
\end{lemm}

\begin{proof}
The same proof as for Lemma~\ref{lemm:defsq} works.
\end{proof}

\begin{defi}
Let $(S, \H^*, \Sigma, \mathcal{G})$ be a \emph{sliceable 2-functor}. We define
\[ \H(X)^\perp(q) \]
as the full triangulated subcategory of $\H(X)$ whose objects $E$ satisfy $\hom(F, E) = 0$ for all $F \in \Sigma^q\H(X)\eff$. If $X$ is clear from the context we will write $\H^\perp(q)$.
\end{defi}

We finish this preliminary subsection with some properties of stable homotopy 2-functors that we will use in developing and applying Theorem~\ref{theo:criteria} of Pelaez. The following theorem should really just be a reference to \cite[Chapter 1]{Ayo07} which contains all the properties we need. However, the material there is for a 2-functor on the category of quasi-projective schemes. As such, we give some indication of the easy generalisation of the properties we want to the category of all schemes of finite type.

\begin{theo} \label{theo:ayo}
Let $(S, \H^*, \Sigma, G)$ be a sliceable 2-functor that satisfies:
\begin{enumerate}[(I)]
 \item \label{item:smoothBaseChange}(Smooth base change) For every cartesian square
\begin{equation} \label{equa:cartesianSquare}
\xymatrix{
Y \times_X W \ar[r]^-g \ar[d]_q & W \ar[d]^p \\
Y \ar[r]_f & X
} \end{equation}
in $Sch(S)$ with $f$ smooth the comparison exchange 2-morphism $f^*p_* \stackrel{\sim}{\to} q_*g^*$ is invertible.\footnote{The right adjoints $p_*, q_*$ to the functors $p^*, q^*$ exist for any sliceable 2-functor. See the proof below.}
 \item \label{item:zariskiSeparated}(Zariski separated) For every Zariski cover $\{j_i: U_i \to X\}$ in $Sch(S)$ the family of functors $\{j_i^*\}$ is conservative.
 \item \label{item:stability}(Stability) We have $\Sigma = p_\#s_*$ where $s$ is the zero section of the canonical projection $p: \AA^1_X \to X$ for each $X \in Sch(S)$
\end{enumerate}

If the restriction of $\H^*$ to the category $QProj(S)$ of quasi-projective $S$-schemes is a stable homotopy 2-functor (\cite[Definition 1.4.1]{Ayo07}). Then $\H^*$ has the following properties.
\begin{enumerate}
 \item Adjoints.
  \begin{enumerate}
 \item \label{item:rightAdjoint}(Right adjoint) For every morphism $f: Y \to X$ in $Sch(S)$ the 1-functor $f^*$ has a right adjoint 
\[ f_*: \H(Y) \to \H(X). \]
  \item \label{item:projRightAdjoint}For every projective morphism $f: Y \to X$ in $Sch(S)$ the 1-functor $f_*$ has a right adjoint
\[ f^!: \H(X) \to \H(Y). \]
  \item \label{item:finiteEtale}If $f$ is a finite {\'e}tale morphism in $Sch(S)$ then $f_\#$ is canonically isomorphic to $f_*$.
  \item \label{item:nilpotent}If $i: Z \to X$ is a nilpotent immersion in $Sch(S)$ then $i^*$, and hence $i_*$, is an equivalence of categories.
 \end{enumerate}

 \item \label{item:otherStability} Tate twists. The auto-equivalences $\Sigma$ form an auto-equivalence of $\H^*$. That is, for any morphism $f: Y \to X$ in $Sch(S)$ we have 2-isomorphisms $\phi: \Sigma f^* \stackrel{\sim}{\to} f^* \Sigma$, and so by adjunction a 2-isomorphism $\psi: f_* \Sigma \stackrel{\sim}{\to} \Sigma f_*$, and if $f$ is smooth a 2-isomorphism $\chi: f_\# \Sigma \stackrel{\sim}{\to} \Sigma f_\#$.

 \item \label{item:localisation}Localisation. Suppose that $j: U \to X$ is an open immersion in $Sch(S)$ and $i: Z \to X$ a complementary closed immersion. There exist unique 2-morphisms $\phi, \psi$ such that
\[ j_\#j^* \to \id_{\H(X)} \to i_*i^* \stackrel{\phi}{\to} j_\#j^*[1] \]
and 
\[ i_*i^! \to \id_{\H(X)} \to j_*j^* \stackrel{\psi}{\to} i_*i^![1] \]
are distinguished triangles where the other morphisms are the units and counits of the adjunctions.

 \item Base change.
 \begin{enumerate}
 \item  \label{item:properBaseChange}For every cartesian square (\ref{equa:cartesianSquare}) in $Sch(S)$ the exchange 2-morphism $f^*p_* \stackrel{\sim}{\to} q_*g^*$ is invertible if $f$ is proper.
 \item \label{item:otherSmoothBaseChange}If $f$ is smooth then the exchange 2-morphism $g_\#q^* \stackrel{\sim}{\to} p^*f_\#$ is invertible.
 \item \label{item:smoothClosedBaseChange}If $f$ is smooth and $p$ a closed immersion then the exchange 2-morphism $f_\#q_* \stackrel{\sim}{\to} p_*g_\#$ is invertible.
 \end{enumerate}

 \item \label{item:mayerVietoris}Mayer-Vietoris. Suppose that $\{ j_U: U \to X, j_V: V \to X \}$ is an open Zariski cover and $j_{U \cap V}: U \cap V \to X$ the intersection. Then there exists a distinguished triangle
\[ {j_{U \cap V}}_\#j_{U \cap V}^* \to {j_{U}}_\#j_{U}^* \oplus {j_{V}}_\#j_{V}^* \to id \to {j_{U \cap V}}_\#j_{U \cap V}^*[1]. \]

 \item \label{item:homotopyInvariance}Homotopy invariance. For any scheme $X \in Sch(S)$ if $p: \AA^1_X \to X$ is the canonical projection then the unit of the adjunction $\id \to p_*p^*$ is an isomorphism. Equivalently, $p_\#p^* \to \id$ is an isomorphism.
\end{enumerate}
\end{theo}

\begin{rema}
In the case $\H^* = \SH$, if we denote by $T \in \SH(S)$ the Tate object (i.e., either the projective line pointed at infinity or the homotopy cokernel of $\Sigma^\infty((\AA^1_S - 0)_+) \to \Sigma^\infty({\AA^1_S}_+)$) then there is a canonical isomorphism between the endomorphism $\Sigma$ defined in this theorem and the endomorphism $T \wedge -$.
\end{rema}

\begin{rema} \label{rema:stableHomotopy2functor}
To aid the reader who is familiar with the theory but unable to recall the precise definition of a stable homotopy 2-functor, we recall that \cite[Definition 1.4.1]{Ayo07} asks that on the category of quasi-projective $S$-schemes we have: $\H(\varnothing) = 0$, smooth left adjoints (Definition~\ref{defi:sliceable2Functor}(\ref{item:smoothLeftAdjoint})), right adjoints (\ref{item:rightAdjoint}), smooth base change (\ref{item:smoothBaseChange}), stability (\ref{item:stability}) in the form ``each $p_\#s_*$ is an equivalence'', homotopy invariance (\ref{item:homotopyInvariance}), and finally, for a closed immersion $i$ with open compliment $j$ the pair $(i^*, j^*)$ is conservative, and $i^*i_* = \id$.
\end{rema}

\begin{proof}
\begin{enumerate}
 \item
  \begin{enumerate}
 \item Since each $\H(X)$ is compactly generated by a theorem of Neeman \cite[Theorem 4.1]{Nee96} it suffices to show that $f^*$ preserves small sums. This is one of our assumptions.
 \item If the morphism $f$ is a projective morphism in $QProj(S)$ then this is \cite[Proposition 1.6.46, Theorem 1.7.17]{Ayo07}. Suppose that $f$ is projective but not in $QProj(S)$. By what we have just mentioned it suffices to show that $f_*$ preserves small sums.  That is, the canonical morphism $\sum_i f_*E_i \to f_*\sum_i E_i$ is an isomorphism. There exists a Zariski cover $\{j_i : U_i \to X \}$ of $X$ such that each $U_i$ is in $QProj(S)$.  Then smooth base change and the quasi-projective case gives the result. 
 \item Notice that under our hypotheses, the restriction of $\H^*$ to $QProj(X)$ is a stable homotopy 2-functor for any $X \in Sch(S)$. Hence, by replacing $S$ with $X$ we can assume that $f$ is in $QProj(S)$. This case is \cite[Section 1.5.3, Theorem 1.7.17]{Ayo07}.
 \item Again, we can assume that $X = S$. In this case it follows from localisation and the identity $i^*i_* = \id$ (see Remark~\ref{rema:stableHomotopy2functor}).
\end{enumerate}

 \item This is assumed in the definition of a sliceable 2-functor.

 \item Localisation in the quasi-projective case is \cite[Lemma 1.4.6]{Ayo07}. Assuming $X = S$ puts us in the quasi-projective case.

 \item
 \begin{enumerate}
 \item For the case when $f$ is projective and the square is in $QProj(S)$ this is \cite[Corollary 1.7.18]{Ayo07}. The generalisation follows from Chow's Lemma and is detailed in \cite[Proposition 2.3.11]{CD}.
 \item This follows by adjunction directly from the smooth base change we have assumed.
 \item The morphism is defined in the usual way using smooth base change $f_\#q_* \to p_*p^* f_\# q^* \cong p_*g_\# q^*q^* \to p_*g_\#$. If our square is in $QProj(S)$ then this is \cite[Corollary 1.4.18]{Ayo07}. Replacing $S$ with $X$ it is also true for any square for which $f$ is a quasi-projective morphism. In the general case, let $\{U_i \to X\}_{i = 1, \dots, n}$ be a finite Zariski cover of $Y$ such that each $U_i \to X$ is quasi-projective. Notice that this implies that for each non-empty subset $I \subset \{ 1, \dots, n\}$ the scheme $U_I = \cap_{i \in I} U_i$ is also quasi-projective over $X$. For each such $I$ the cartesian square
\[ \xymatrix{
U_I \times_X W \ar[r] \ar[d] & Y \times_X W  \ar[d] \\
U_I \ar[r] & Y
} \]
satisfies the property we want, and so the natural transformation $f_\#q_* \to p_*g_\#$ evaluated on any object $E$ in the image of $(U_I \times_X W \to Y \times_X W)_\#$ is an isomorphism. We will show by induction on the size of a subset $J \subseteq \{ 1, \dots, n\}$ the natural transformation $f_\#q_* \to p_*g_\#$ is an isomorphism when evaluated on any object $E$ in the image of $(\cup_{i \in J} U_i \times_X W \to Y \times_X W)_\#$. Let $\phi_J: \cup_{i \in J} U_i \times_X W \to Y \times_X W$ denote the morphism. If $J$ is empty, the object is necessarily zero, and so we clearly get an isomorphism. If not, then there exists two subsets of smaller size $J'$ and $J''$ such that $J' \cup J'' = J$ and by the Mayer-Vietoris triangle
\[ {\phi_{J' \cap J''}}_\#\phi_{J' \cap J''}^* \to {\phi_{J'}}_\#\phi_{J'}^* \oplus {\phi_{J''}}_\#\phi_{J''}^* \to {\phi_{J}}_\#\phi_{J}^* \to {\phi_{J' \cap J''}}_\#\phi_{J' \cap J''}^*[1]. \]
and the inductive assumption, we are done. This proves that the natural transformation  $f_\#q_* \to p_*g_\#$ is an isomorphism as $\phi_{\{1, \dots, n\}} = \id_{Y \times_X W}$.
 \end{enumerate}

 \item Mayer-Vietoris is an immediate consequence of Zariski separatedness.
 \item Replacing $S$ with $X$ we can assume that $X$ is quasi-projective (over $X$). In this case it is one of the axioms of a stable homotopy 2-functor.
\end{enumerate}
\end{proof}

\begin{defi}
If $S$, $\H^*$, $\Sigma$ and $\mathcal{G}$ are as in Theorem~\ref{theo:ayo} we will refer to $\H^*$ as a \emph{stable homotopy 2-functor}. We take $S$ and $\mathcal{G}$ to be implicit in the definition of such a stable homotopy 2-functor. Of course $\Sigma$ is defined in Theorem~\ref{theo:ayo}(\ref{item:stability}).
\end{defi}

\begin{rema}
This is a mild abuse of the terminology as Ayoub's stable homotopy 2-functors are defined on the category of quasi-projective schemes over $S$ and we have further asked for $\H^*$ to be what he calls ``compactly generated by the base'' \cite[Definition 2.1.155]{Ayo07}.
\end{rema}

Suppose that $Y, X$ are two schemes and $\Phi: \H(X) \to \H(Y)$ is a functor. We will say that \emph{$\Phi$ preserves $\Sigma^q\H\eff$ (resp. $\H^\perp(q)$)} if for every object $E \in \Sigma^q\H(X)\eff$ (resp. $E \in \H(X)^\perp(q)$) the object $\Phi E$ is in $\Sigma^q\H(Y)\eff$ (resp. $\H(Y)^\perp(q)$).

\begin{lemm} \label{lemm:generate}
Let $\H^*$ be a stable homotopy 2-functor. Let $a: X \to S$ be a scheme in $Sch(S)$ and $\{U_i \to X\}_{i = 1, \dots, N}$ a Zariski cover. Then for every $q$ the category $\Sigma^q\H(X)\eff$ is the smallest full triangulated subcategory of $\H(X)$ containing all sums and the objects $\Sigma^n f_\#f^* a^*A$ where $n \geq q$, $A \in \mathcal{G}$, and $f:W \to X$ is a smooth morphism from an affine scheme $W$ whose image is contained in some $U_i$.
\end{lemm}

\begin{proof}
Let $T$ be the smallest full triangulated subcategory of $\H(X)$ containing the objects of the form described in the statement. Suppose that $f:W \to X$ is a smooth morphism with source an affine scheme whose image is contained in some $U_i$. We claim that for every open subscheme $j: W' \to W$ of $W$ and $n \geq q$, the object $\Sigma^n (fj)_\#(fj)^* a^*\un_S$ is in $T$. Indeed, this is obviously true if $W'$ is also affine. Now every open subscheme of $W$ can be covered by finitely many basic open affine subschemes (i.e., affine subschemes of the form $Spec(A_f)$ where $W = Spec(A))$. We work by induction on the smallest number $r$ of such subschemes it takes to cover $W'$. If  $r = 1$ there is nothing to show since in this case $W'$ is affine. Suppose it is true for $i < r$ and $W'$ can be covered by $r$ basic affine open subschemes of $W$. Then in particular, there is a cover of the form $\{ W'' \subset W', W''' \subset W'\}$ where $W''$ is a basic affine open and both $W'''$ and $W'' \cap W'''$ can be covered by $r - 1$ basic affine opens (since the intersection of two basic affine opens is abasic affine open). It then follows from Mayer-Vietoris (Theorem~\ref{theo:ayo}(\ref{item:mayerVietoris})) that the object $\Sigma^n (fj)_\#(fj)^* a^*\un_S$ corresponding to $W'$ is in $T$.

We use the same argument twice more.

Let $f: W \to X$ be a smooth morphism whose image is contained in some $U_i$. We claim that for every open subscheme $j: W' \to W$ of $W$ and $n \geq q$, the object $\Sigma^n (fj)_\#(fj)^* a^*\un_S$ is in $T$. We have just seen that this is true if $W'$ is contained in an affine open subscheme of $W$. We use the same argument as above with ``basic open affine'' replaced by ``an open subscheme that is contained in an open affine of $W$''.

Let $f: W \to X$ be any smooth morphism. We claim that for every open subscheme $j: W' \to W$ of $W$ and $n \geq q$, the object $\Sigma^n (fj)_\#(fj)^* a^*\un_S$ is in $T$. Indeed, every such $W'$  can be covered by a finite number of open subschemes whose images are contained in some $U_i$. We use the same induction argument again.

We have shown that $T$ contains the generators for $\Sigma^q\H(X)\eff$. Since it is a triangulated subcategory with small sums, this is enough to conclude.
\end{proof}

\begin{lemm} \label{lemm:preserve}
Let $\H^*$ be a stable homotopy 2-functor. The following functors preserve the following categories.
\begin{enumerate}
  \item For any morphism $f: Y \to X$ in $Sch(S)$
 \begin{enumerate}
  \item $f^*$ preserves $\Sigma^q\H\eff$, and
  \item $f_*$ preserves $\H^\perp(q)$.
 \end{enumerate}
  \item For a smooth morphism $f: Y \to X$ in $Sch(S)$
 \begin{enumerate}
  \item $f_\#$ preserves $\Sigma^q\H\eff$, and
  \item $f^*$ preserves $\H^\perp(q)$.
 \end{enumerate}
 \item For $i: Z \to X$ a closed immersion between quasi-projective $S$ schemes, $i_*i^*$ and $i_*$ both preserve $\Sigma^q\H\eff$.
 \item For $i: Z \to X$ a nilpotent immersion between quasi-projective $S$ schemes, both $i^*$ and $i_*$ preserve both $\Sigma^q\H\eff$ and $\H^\perp(q)$.
 \item For $f: Y \to X$ a finite {\'e}tale morphism between quasi-projective $S$ schemes, each of $f_\#, f^*, f_*$ preserve both $\Sigma^q\H\eff$ and $\H^\perp(q)$.
\end{enumerate}
\end{lemm}

\begin{proof}
In the first two (b) follows from (a) by adjunction. Suppose that $a: X \to S$ is the structural morphism of $X$ and $p: W \to X$ a smooth morphism. For every $r$ there are canonical isomorphisms $f^*\Sigma^rp_\#p^*a^*\un_S \stackrel{\ref{theo:ayo}}{\cong} \Sigma^rf^*p_\#p^*a^*\un_S \stackrel{\ref{theo:ayo}(\ref{item:smoothBaseChange})}{\cong} \Sigma^rp'_\#p'^*f^*a^*\un_S \stackrel{}{\cong} \Sigma^rp'_\#p'^*(af)^*\un_S$ where $p' = Y \times_X p$. Hence, $f^*$ sends generators of $\Sigma^q\H(X)\eff$ to $\Sigma^q\H(Y)\eff$. Since $f^*$ is triangulated and preserves homotopy colimits (because it is a left adjoint), this is enough to conclude that $f^*$ preserves $\Sigma^q\H\eff$. The same argument works for $f_\#$: Suppose that $p: W \to Y$ is a smooth morphism. We have isomorphisms $f_\#\Sigma^rp_\#p^*(af)^*\un_S \stackrel{}{\cong} f_\#\Sigma^rp_\#p^*f^*a^*\un_S \stackrel{\ref{theo:ayo}(\ref{item:otherStability})}{\cong} \Sigma^r f_\#p_\#p^*f^*a^*\un_S \cong \Sigma^r (fp)_\#(fp)^*a^*\un_S$.

Consider the case of a closed immersion. The functor $i_*i^*$ is straight-forward. Let $j: X - Z \to X$ be the complementary open immersion. We have a localisation distinguished triangle $j_\#j^* \to \id \to i_*i^* \to j_\#j^*[1]$ (Theorem~\ref{theo:ayo}(\ref{item:localisation})) and since $j$ is smooth, since $j_\#$ and $j^*$ both preserve $\Sigma^q\H\eff$ it follows that $i_*i^*$ also preserves $\Sigma^q\H\eff$.

Now we consider the functor $i_*$. Let $T$ be the full subcategory of $\Sigma^q\H(Z)\eff$ consisting of objects $E$ such that $i_*E \in \Sigma^q\H(X)\eff$. The triangulated category $\Sigma^q\H(X)\eff$ has small sums and $i_*$ commutes with small sums (as it is a left adjoint (Theorem~\ref{theo:ayo}(\ref{item:projRightAdjoint}))) and so $T$ is a triangulated category with small sums. We will show that $T$ contains a generating family for $\Sigma^q\H(Z)\eff$, which then implies it contains all of $\Sigma^q\H(Z)\eff$. Suppose $\{ U_i \to X\}_{i = 1, \dots, N}$ is a Zariski cover of $X$ by affine schemes. We consider the generating family of $\Sigma^q\H(Z)\eff$ described in Lemma~\ref{lemm:generate} associated to the cover $\{ U_i \cap Z \to Z \}$. Suppose $\Sigma^n f_\#f^* a^*\un_S$ is a member of this generating family where $a: Z \to S$ is the structural morphism, $n \geq q$ and $f:W \to Z$ is a smooth morphism from an affine scheme $W$ whose image is contained in some $U_i \cap Z$. Recall a theorem of Arabia \cite{Ara01} that says that in general if $Z' \to X'$ is a closed immersion of affine schemes and $W' \to Z'$ is a smooth morphism then there exists a smooth morphism $V' \to X'$ and a $Z'$-isomorphism $Z' \times_{X'} V' \cong W'$. In our case, this gives us a smooth morphism $g: V \to X$ and a $Z$-isomorphism $W \cong Z \times_X V$. Let $a': X \to S$ be the structural morphism of $X$ so that $a = a'i$. By the appropriate parts of Theorem~\ref{theo:ayo} we find an isomorphism $\Sigma^n f_\#f^* a^*\un_S \cong \Sigma^n f_\#f^* i^*a'^*\un_S \cong i^*\Sigma^n g_\#g^* a'^*\un_S$. That is, our object is in the image of $i^*: \Sigma^q\H(X)\eff \to \Sigma^q\H(Z)\eff$. Now we have $i_*(\Sigma^n f_\#f^* a^*\un_S) \cong i_*(i^*\Sigma^n g_\#g^* a'^*\un_S)$. So $\Sigma^n f_\#f^* a^*\un_S$ is in $T$ because $i_*i^*$ preserves $\Sigma^q\H\eff$. So $i_*$ preserves $\Sigma^q\H\eff$.

Suppose $i: Z \to X$ is a nilpotent immersion. It remains only to see that $i^*$ preserves $\H^\perp(q)$. Since $i$ is a nilpotent immersion, the functors $i^*$ and $i_*$ are equivalences of categories, each inverse to the other. In particular, $i^*$ is now a right adjoint of $i_*$. We have seen $i_*$ preserves $\Sigma^q\H\eff$ and so it follows by adjunction that $i^*$ preserves $\H^\perp(q)$.

Lastly, in the finite {\'e}tale case, we have already seen above that $f_\#$ and $f^*$ preserve $\Sigma^q\H\eff$, and $f^*$ and $f_*$ preserve $\H^\perp(q)$. But $f_\#$ is isomorphic to $f_*$ (Theorem~\ref{theo:ayo}(\ref{item:finiteEtale})) and so $f_\#$ also preserves $\H^\perp(q)$ and $f_*$ also preserves $\Sigma^q\H\eff$.
\end{proof}


\subsection{After Pelaez}

In this subsection we continue with $S, \H^*$ and $\mathcal{G}$ as in Theorem~\ref{theo:ayo}. Recall that for any object $E \in \H(Y)$ we have $f_qE \in \Sigma^q\H(Y)\eff$ (by definition) and $s_qE \in \H(Y)^\perp(q + 1)$.

\begin{defi}
For the rest of this section, we will have $\Phi: \H(Y) \to \H(X)$ a triangulated functor, $E \in \H(Y)$ an object, and $q \in \ZZ$ an integer. We will be considering whether  the following conditions hold.
\begin{enumerate}
 \item[{(Pel0)$_q$}] $\Phi \hocolim_{p \leq q} f_pE = \hocolim_{p \leq q} \Phi f_pE$.
 \item[{(Pel1)$_q$}] $\Phi f_qE \in \Sigma^q\H(X)\eff$.
 \item[{(Pel2)$_q$}] $\Phi s_qE \in \H(X)^\perp(q + 1)$.
\end{enumerate}
\end{defi}

\begin{rema}
If the object $E$ is not clear from the context we will write (Pel$i$)$_q$($E$) with $i = 0, 1,$ or $2$. In particular note that (Pel1)$_q$($E$) implies (Pel1)$_q$($f_rE$) for all $r \leq q$. We will also use the notation (Pel$i$)$_I$ for $I \subseteq \ZZ$ to indicate that (Pel$i$)$_q$ is true for all $q \in I$ and (Pel$i$) for (Pel$i$)$_\ZZ$.
\end{rema}

We collect here some functors that are known to satisfy some of these conditions.

\begin{lemm} \label{lemm:table}
Let $\H^*$ be a stable homotopy 2-functor. Then the following functors satisfy the following conditions for all objects.
\begin{center}
\begin{tabular}{cccc}
Functor & (Pel0) & (Pel1) & (Pel2)  \\  \hline
$f^*$ &  Yes &  Yes & $-$ \\
$f_*$ & Yes & $\times$ &  Yes \\
$f_\#$ with $f$ smooth &  Yes & Yes & $\times$ \\
$f^*$ with $f$ smooth &  Yes &  Yes & Yes \\
$i^*$ with $i$ a nilpotent immersion &  Yes & Yes & Yes \\
$i_*$ with $i$ a closed immersion &  Yes & Yes & Yes \\
$i_*i^*$ with $i$ a closed immersion &  Yes & Yes & $-$ \\
$f_\#, f^*, f_*$ with $f$ finite and {\'e}tale &  Yes &  Yes &  Yes
\end{tabular}
\end{center}
``$\times$'' indicates ``not in general'' (there are certainly examples where the property is satisfied, for example $f = \id$), and ``$-$'' indicates ``unknown''.
\end{lemm}

\begin{rema}
Theorem~\ref{theo:pullbackAll} gives conditions under which $f^*$ preserves (Pel2). This is a version of a theorem of Pelaez. His theorem has fewer restrictions but assumes resolution of singularities. Given what we know about $i_*$ this applies then to $i_*i^*$ as well.

While our counter-examples for $f_*$ and $f_\#$ show that they don't preserve the slice filtration, they suggest that they ``shift'' it in a suitable sense, at least in certain cases (cf. \cite[Theorem 4.4]{Pel11}).
\end{rema}

\begin{proof}
The columns (Pel1) and (Pel2) follow directly from Lemma~\ref{lemm:preserve}. The column (Pel0), apart from $f_*$, follows from the functors in question being left adjoints (cf. Theorem~\ref{theo:ayo}). For (Pel0) for $f_*$, we note that $f^*$ preserves a set of compact generators (due to them being compatible with $\Sigma$ and smooth base change) and therefore its right adjoint preserves small sums.

For a counter example to $f_\#$ satisfying (Pel2) suppose that $s_qE \neq 0$ and consider the canonical projection of the affine line $p: \AA^1_S \to S$. Let $s$ be the zero section. If $p_\#$ satisfies (Pel2), then $p_\#s_* = \Sigma$ would satisfy (Pel2) as well. Now $s_qE' \in \Sigma^q\H\eff \cap \H^\perp(q + 1)$ for every object $E'$ and so $\Sigma s_q E \in \Sigma^{q + 1}\H\eff$. But if (Pel2) is satisfied then we also have $\Sigma s_q E \in \H^\perp(q + 1)$. Hence, the identity morphism of $\Sigma s_q E$ is zero, and therefore $s_qE$ is zero. So in this case, $p_\#$ does not satisfy (Pel2).

A similar phenomena gives a counter example to $f_*$ satisfying (Pel1). Suppose that $s_qE \neq 0$. Let $p: \PP^1_S \to S$ be the projection, $s$ the section at infinity, and $j: \AA^1 \to \PP^1$ the complimentary affine line.  If $j: \AA^1_S \to \PP^1_S$ is the open immersion then we have the localisation distinguished triangle
\[ s_*s^!  \to \id \to j_*j^* \to s_*s^![1]. \]
Evaluating this triangle on $p^*$ and applying $p_*$ gives
\[ s^!p^*  \to p_*p^* \to \id \to s^!p^*[1] \]
where we have used homotopy invariance to obtain the $\id$. So if $p_*$ satisfies (Pel1) then so does $p_*p^*$ and $s^!p^*$. This latter is the right adjoint to $p_\#s_*$ which is isomorphic to $a_\#s_*$ where $a: \AA^1_S \to S$ is the projection.\footnote{To see this use base change Theorem~\ref{theo:ayo}(\ref{item:smoothClosedBaseChange}) on $\AA^1_S \to \PP^1_S \leftarrow S$ with the latter the embedding at zero.} Hence, $p_\#s_* = \Sigma$ and $s^!p^* = \Sigma^{-1}$. So $\Sigma^{-1}$ would satisfy (Pel1) in this case. But then we would have  $\Sigma^{-1} s_q E \in \Sigma^q\H\eff$. However, $s_q E \in \H^\perp(q + 1)$ and so $\Sigma^{-1} s_q E \in \H^\perp(q)$ leading to $s_q E = 0$ as before. Hence, $p_*$ does not satisfy (Pel1).
\end{proof}

\begin{defi} \label{defi:comparisonMorphismssf}
Let $(S, \H^*, \Sigma, \mathcal{G})$ be a sliceable 2-functor. Let $\Phi: \H(Y) \to \H(X)$ be a functor, $E \in \H(Y)$ an object, and $q \in \ZZ$ an integer. We consider the canonical morphisms
\[ \Phi f_q \leftarrow f_q \Phi f_q \to f_q\Phi. \]
If $\Phi f_q E \leftarrow f_q \Phi f_q E$ is an isomorphism (for example if (Pel1)$_q(E)$ is satisfied) we denote the resulting canonical morphism by
\[ \alpha_q(E): \Phi f_q E \to f_q \Phi E \]
or $\alpha_q$ or $\alpha$ if $E$ and $q$ are clear from the context.
\end{defi}

\begin{lemm} \label{lemm:comparisonMorphismssf}
Let $(S, \H^*, \Sigma, \mathcal{G})$ be a sliceable 2-functor. Let $\Phi: \H(Y) \to \H(X)$ be a functor, $E \in \H(Y)$ an object, and $q \in \ZZ$ an integer. If the two morphisms
\[  \Phi f_q E \leftarrow f_q \Phi f_q E \qquad \textrm{ and } \qquad \Phi f_{q + 1} E \leftarrow f_{q + 1} \Phi f_{q + 1} E \]
are isomorphisms (for example if (Pel1)$_q(E)$ and (Pel1)$_{q + 1}(E)$ are satisfied) then there is a unique morphism
\[ \beta_q(E): \Phi s_q E \to s_q \Phi E \]
such that the following diagram is commutative.
\[ \xymatrix{
\Phi(f_{q + 1}E) \ar[r] \ar[d]^{\alpha_{q + 1}} & \Phi(f_qE) \ar[r] \ar[d]^{\alpha_{q}} & \Phi(s_qE) \ar[r] \ar[d]^{\beta_{q}} & \Phi(f_{q + 1}E) \ar[d]^{\alpha_{q + 1}[1]} \\
f_{q + 1}\Phi(E) \ar[r] & f_q\Phi(E) \ar[r] & s_q\Phi(E) \ar[r] & f_{q + 1}\Phi(E)
} \]
The morphisms $\alpha_p(E)$ and $\beta_p(E)$ are functorial in $\Phi$ in two senses:
\begin{enumerate}
 \item If $\eta: \Phi \to \Psi$ is a natural transformation between functors the appropriate $\alpha$'s are defined then the square
\[ \xymatrix{
\Phi(f_qE) \ar[r]^\alpha \ar[d]_\eta & f_q \Phi(E) \ar[d]^\eta \\
\Psi(f_qE) \ar[r]_\alpha & f_q \Psi(E) 
} \]
commutes (and similarly for $s_q$ if the $\beta$'s are defined).
 \item If $\Phi: \H(Y) \to \H(X)$ and $\Psi: \H(W) \to \H(Y)$ are triangulated functors such that the appropriate $\alpha$'s are defined then the triangle
\[ \xymatrix{
\Psi \Phi (f_q E) \ar[r]_\alpha \ar@/^12pt/[rr]^\alpha & \Psi f_q \Phi (E) \ar[r]_\alpha & f_q \Psi \Phi(E) 
} \]
commutes (as well as the analogous statement for $s_q$).
\end{enumerate}
\end{lemm}

\begin{proof}
There certainly exists such a morphism $\Phi s_q E \to  s_q\Phi E$ since the two triangles $\Phi f_{q + 1}E \to \Phi f_qE \to \Phi s_qE \to \Phi f_{q + 1}E[1]$ and $f_{q + 1}\Phi E \to f_q\Phi E \to s_q\Phi E \to f_{q + 1} \Phi E[1]$ are distinguished. Uniqueness comes from the fact that $\hom(\Phi f_{q + 1}E[1], s_q\Phi E) = 0$. This latter is a consequence of the facts that $\Phi f_{q + 1}E[1] \cong f_{q + 1}\Phi f_{q + 1}E[1]$ is in $\Sigma^{q + 1}\H(Y)\eff$ and $s_q\Phi E$ is in $\H(Y)^\perp(q + 1)$.

The functoriality for the $\alpha$'s is clear from the appropriate functoriality of the $f_q$'s. The functoriality for the $\beta$'s is again a consequence of the fact that there are no non-zero morphisms from $\Sigma^{q + 1}\H\eff$ to $\H^\perp(q + 1)$.
\end{proof}

\begin{lemm}[Pelaez] \label{lemm:pel1}
Let $(S, \H^*, \Sigma, \mathcal{G})$ be a sliceable 2-functor. Let $\Phi: \H(Y) \to \H(X)$ a functor, $E \in \H(Y)$ an object, and $q \in \ZZ$ an integer. Suppose that (Pel1)$_{q + 1}$, (Pel1)$_q$ and (Pel2)$_q$ are satisfied. Then 
the two morphisms
\[ \alpha_{q + 1}(f_qE) : \Phi(f_{q + 1}f_qE) \to f_{q + 1}(\Phi(f_qE)) \]
\[ \beta_{q}(f_qE) : \Phi(s_qf_qE) \to s_q(\Phi(f_qE)) \] 
are isomorphisms in $\H(X)$.
\end{lemm}

\begin{proof}[Proof (Pelaez).]
We have the commutative diagram associated to $f_q$
\[ \xymatrix{
\Phi(f_{q + 1}f_qE) \ar[r] \ar[d]^{\alpha_{q + 1}} & \Phi(f_qf_qE) \ar[r] \ar[d]^{\alpha_{q}} & \Phi(s_qf_qE) \ar[r] \ar[d]^{\beta_{q}} & \Phi(f_{q + 1}f_qE) \ar[d]^{\alpha_{q + 1}[1]} \\
f_{q + 1}\Phi(f_qE) \ar[r] & f_q\Phi(f_qE) \ar[r] & s_q\Phi(f_qE) \ar[r] & f_{q + 1}\Phi(f_qE)
} \]
The property (Pel1) implies that $\alpha_q(f_qE)$ is an isomorphism. Using the octahedral axiom we have a commutative diagram where all the rows and columns are distinguished triangles
\[ \xymatrix{
\Phi(f_{q + 1}f_qE) \ar[r] \ar[d]^{\alpha_{q + 1}} & \Phi(f_{q}f_qE) \ar[r] \ar[d]^{\alpha_q} & \Phi(s_qf_qE) \ar[r] \ar[d]^{\beta_q} & \Phi(f_{q + 1}f_qE)[1] \ar[d] \\
f_{q + 1}\Phi(f_qE) \ar[r] \ar[d] & f_q\Phi(f_qE) \ar[r] \ar[d] & s_q\Phi(f_qE) \ar[r] \ar[d] & f_{q + 1}\Phi(f_qE)[1] \ar[d] \\
A \ar[r] & 0 \ar[r] & A[1] \ar@{=}[r] & A[1]
} \]
and it now suffices to show that $A = 0$. We note that $A$ is in $\Sigma^{q + 1}\H(X)\eff$ since both $f_{q + 1}\Phi(f_qE)$ and $\Phi(f_{q + 1}f_qE) = \Phi(f_{q + 1}E)$ are.

On the other hand, $\Phi(s_qE) \cong \Phi(s_qf_qE)$ is in $\H(X)^\perp(q + 1)$ by hypothesis and $s_q\Phi(f_qE)$ is in $\H(X)^\perp(q + 1)$ (as $s_q$ always is) so $A[1]$ is also in $\H(X)^\perp(q + 1)$. Since $\H(X)^\perp(q + 1)$ is a triangulated subcategory, $A$ is also in $\H(X)^\perp(q + 1)$.

Now there are no nonzero morphisms from $\Sigma^{q + 1}\H(X)\eff$ to $\H(X)^\perp(q + 1)$ and so the identity of $A$ is zero, hence $A$ is isomorphic to zero.
\end{proof}

\begin{theo}[Pelaez] \label{theo:criteria}
Let $(S, \H^*, \Sigma, \mathcal{G})$ be a sliceable 2-functor. Let $\Phi: \H(Y) \to \H(X)$ be a functor, $E \in \H(Y)$ an object, and $q \in \ZZ$ an integer. Suppose that (Pel0)$_q$, (Pel1)$_{\leq q + 1}$, and (Pel2)$_{\leq q}$ are satisfied. Then the morphisms
\[ \alpha_r(E): \Phi(f_r E) \to  f_r \Phi(E) \]
\[ \beta_r(E): \Phi(s_r E) \to  s_r \Phi(E) \]
are isomorphisms for all $r \leq q$.
\end{theo}

\begin{proof}[Proof (Pelaez).]
The hypotheses are stable by lowering $q$ and so it suffices to prove that $\alpha_q(E)$ and $\beta_q(E)$ are isomorphisms. The same proof works for both, and we will give the proof for $\beta$ but the reader can check that the proof remains valid with $\beta$ replaced with $\alpha$ everywhere (and $s_r$ replaced with $f_r$ where appropriate).

For any fixed integer $N$ we have $E \cong hocolim_{p \leq N} f_p E$ and so since $\Phi$ and $s_q$ commute with homotopy colimits $\beta_q(E) = hocolim_{p \leq N} \beta_q(f_p E)$ hence it suffices to show that each $\beta_q(f_p E)$ is an isomorphism for all $p \leq N$ for some $N$. We chose $N = q$. This way, Lemma~\ref{lemm:pel1} implies that $\beta_q(f_q E)$ is an isomorphism. We now proceed by induction.

Suppose that $\beta_q(f_rE)$ is an isomorphism for some $r \leq q$. We must show that $\beta_q(f_{r - 1}E)$ is an isomorphism. We have a commutative diagram
\[ \xymatrix{
\Phi(s_qf_rE) \ar[r]^{\beta_q(f_rE)} \ar[d]_{\Phi s_q\rho} & s_q\Phi f_rE \ar[d]^{s_q\Phi \rho} \\
\Phi(s_qf_{r - 1}E) \ar[r]_-{\beta_q(f_{r - 1}E)} & s_q\Phi f_{r - 1}E
} \]
where $\rho: f_r \to f_{r - 1}$ is the canonical natural transformation. The inductive hypothesis says that the upper morphism is an isomorphism, and we have that $s_qf_r = s_q$ and $s_qf_{r - 1} = s_q$, hence $s_q\rho$ is an isomorphism by construction of the slice filtration. Hence, it suffices to show that the morphism on the right is an isomorphism.

We have another commutative square
\[ \xymatrix{
s_q\Phi f_rE \ar[r] \ar[d]_{s_q\Phi\rho} & s_q\Phi f_rf_{r - 1}E \ar[d]^{\alpha_r} \\
s_q\Phi f_{r - 1}E \ar[r] & s_qf_r\Phi f_{r - 1}E
} \]
with the horizontal morphisms isomorphisms. The right morphism is an isomorphism by Lemma~\ref{lemm:pel1} above and the third hypothesis. Hence, the morphism on the left is an isomorphism as desired.
\end{proof}

\subsection{Applications of Pelaez's Theorem}

\begin{coro} \label{coro:smoothPullback}
Let $\H^*$ be a stable homotopy 2-functor. Suppose $f:Y \to X$ is a morphism in $Sch(S)$, and $E \in \H(X)$ and $F \in \H(Y)$ are any objects. The canonical morphisms
\[ \alpha_r(E): f^*(f_r E) \to  f_r f^* (E) \]
\[ \beta_r(E): f^*(s_r E) \to  s_r f^*(E) \]
are isomorphisms for all $r$ if $f$ is smooth, or if $f$ is a nilpotent immersion. Similarly, the canonical morphisms
\[ \alpha_r(E): f_*(f_r F) \to  f_r f_* (F) \]
\[ \beta_r(E): f_*(s_r F) \to  s_r f_*(F) \]
are isomorphisms for all $r$ if $f$ is a closed immersion, or if $f$ is a finite {\'e}tale morphism.
\end{coro}

\begin{proof}
This follows from Theorem~\ref{theo:criteria} and Lemma~\ref{lemm:table}.
\end{proof}

We now discuss some consequences of Gabber's theorem (Theorem~\ref{theo:gabberGlobal}).

\begin{defi} \label{defi:wsst}
We will say that an object $E \in \H(S)$ has a \emph{weak structure of smooth traces} if for every $Y \stackrel{f}{\to} X \stackrel{a}{\to} S$ in $Sch(S)$ with $f$ a finite flat surjective morphism between \textbf{smooth} $S$-schemes $X, Y$ such that $f_*\OO_Y$ is a globally free $\OO_X$-module, we are given a morphism $\Tr_f: f_*f^*a^*E \to a^*E$ in $\H(X)$ such that the composition with $a^*E \to f_*f^*a^*E$ is $\deg f \cdot id_{a^*E}$.
\end{defi}

\begin{defi} \label{defi:zlllocal}
Suppose that $\Lambda \subseteq \QQ$ is a subring of the rational numbers. We will say that an object $E$ in an additive category is \emph{$\Lambda$-local} if $\hom(E, E)$ is a $\Lambda$-module. It is equivalent to ask that for every integer $n$ that is invertible in $\Lambda$ the endomorphism $n \cdot \id_E$ is an isomorphism.
\end{defi}

There is some material on $\Lambda$-local objects in Section~\ref{sec:invertTriangle}.

\begin{theo} \label{theo:pullbackAll}
Let $\H^*$ be a stable homotopy 2-functor, suppose $S$ is the spectrum of a perfect field $k$ of exponential characteristic $p$, let $E \in \H(k)$ be a $\zpi$-local object (Definition~\ref{defi:zlllocal}), and $q \in \ZZ$ an integer. If $s_rE$ has a weak structure of smooth traces (Definition~\ref{defi:wsst}) for every $r \leq q$ then for any separated $k$-scheme of finite type $a: X \to Spec\ k$ the morphisms
\[ \beta_q: a^*s_qE \to s_qa^*E \]
\[ \alpha_q: a^*f_qE \to f_qa^*E \]
are isomorphisms in $\H(X)$.
\end{theo}

\begin{proof}

As discussed in Section~\ref{sec:invertTriangle}, it suffices that the statement is true when $E$ is a $\zll$ local object, for every prime $\ell$ different from $p$.

Note that $a^*$ satisfies (Pel0) and (Pel1) for all objects (Lemma~\ref{lemm:table}) so by Theorem~\ref{theo:criteria} it suffices to verify that $a^*$ satisfies (Pel2)$_{\leq q}$. That is, we wish to see that $a^*s_rE \in \H(X)^\perp(r + 1)$ for all $r \leq q$. The hypotheses of the theorem are stable under lowering $q$ and so it suffices to consider the case $r = q$. The proof is by Noetherian induction.

For the morphisms $i: X_{red} \to X$ we have seen that $i^*$ and $i_*$ are inverse equivalences of categories that both preserve $\H^\perp(q + 1)$ and so we can assume that $X$ is reduced. Let $p: X' \to X$ be a proper morphism furnished by Gabber's Theorem (\ref{theo:gabberGlobal}) with $X'$ connected, quasi-projective and smooth, and $j: U \to X$ an non-empty open subset such that $X' \times_X U \to U$ is a finite flat surjective morphism of constant degree prime to $\ell$. Let $Z$ be a closed compliment to $U$. Our diagram is the following:
\[ \xymatrix{
Z \times_X X' \ar[d]^{\tilde{p}} \ar[r]^{\tilde{i}} & X' \ar[d]^p & \ar[l]_{\tilde{j}} X' \times_X U \ar[d]^h \\ 
Z \ar[r]_i & X & \ar[l]^j U
} \]
By the inductive hypothesis and the localisation distinguished triangle $j_\#j^* \to id \to i_*i^* \to j_\#j^*[1]$ it suffices to show that $j_\#j^*a^*s_qE \in \H(X)^\perp(q + 1)$ (Lemma~\ref{lemm:table}). The weak structure of smooth traces on $s_qE$ and the fact that we are working $\zll$-locally, implies that $j_\#(j^*a^*s_qE) \to j_\#(h_*h^*)(j^*a^*s_qE)$ is a monomorphism. Since $\H(X)^\perp(q + 1)$ is idempotent complete, it now suffices to show that $j_\#h_*h^*j^*a^*s_qE$ is in $\H(X)^\perp(q + 1)$.

The base change properties in Theorem~\ref{theo:ayo} give isomorphisms $h_*h^*j^* \cong h_*\tilde{j}^*p^* \cong j^*p_*p^*$ and so it now suffices to show that $(j_\#j^*)(p_*p^*a^*s_qE) \in \H(X)^\perp(q + 1)$. Since $(ap): X' \to Spec(k)$ is smooth $p^*a^*s_qE \in \H(X)^\perp(q + 1)$ (Corollary~\ref{coro:smoothPullback}) and we have seen that $p_*$ preserves $\H^\perp(q + 1)$ (Lemma~\ref{lemm:table}) so $p_*p^*a^*s_qE \in \H(X)^\perp(q + 1)$. Using again the localisation distinguished triangle $j_\#j^* \to id \to i_*i^* \to j_\#j^*[1]$, it suffices to show that $(i_*i^*)(p_*p^*a^*s_qE) \in \H(X)^\perp(q + 1)$. But now by base change (Theorem~\ref{theo:ayo}) we have an isomorphism $i_*(i^*p_*)p^*a^*s_qE \cong i_*(\tilde{p}_*\tilde{i}^*)p^*a^*s_qE$ and by the inductive hypothesis $\tilde{i}^*p^*a^*s_qE \in  \H(X)^\perp(q + 1)$ and so since $i_*\tilde{p}_*$ preserves $ \H^\perp(q + 1)$ (Lemma~\ref{lemm:table}) the proof is complete.
\end{proof}

The remainder of this section is devoted to the proof of Proposition~\ref{prop:ffslices}, and it is much more enjoyable if read in reverse. That is, in the order Proposition~\ref{prop:ffslices}, Proposition~\ref{prop:ffPreserve}, Lemma~\ref{lemm:sqffEffective}, Lemma~\ref{lemm:phiphi}, and then Lemma~\ref{lemm:radicialIso}.

\begin{defi}
If $\H^*$ is a 2-functor on $Sch(S)$ and $E \in \H(S)$ an object, for each scheme $a: X \to S$ in $Sch(S)$ we denote by $E_X$ the object $a^*E$. Note that for any morphism $f: Y \to X$ there is a canonical isomorphism $f^*E_X = E_Y$.
\end{defi}

\begin{lemm} \label{lemm:radicialIso}
Let $\H^*$ be a stable homotopy 2-functor, $f: Y \to X$ a radicial finite flat surjective morphism of degree $d$ between smooth $S$-schemes, and $E \in \H(S)$ a $\ZZ[\tfrac{1}{d}]$-local object with a weak structure of smooth traces. Then
\[ E_X \to f_*f^*E_X \]
is an isomorphism in $\H(X)$.
\end{lemm}

\begin{proof}
First we make a general observation. Suppose $\mathcal{A}$ is an additive category, $\Phi$ an additive endomorphism, $\eta: id \to \Phi$ a natural transformation of additive endofunctors (i.e., $\eta_{(A \oplus B)} = \eta_A \oplus \eta_B$), $A$ an object of $\mathcal{A}$, and suppose that $A$ is a direct summand of $\Phi A$ via the morphism $A \to \Phi A$. In this situation, if $\Phi A \to \Phi \Phi A$ is an isomorphism, then $A \to \Phi A $ is an isomorphism. In effect, writing $\nu: A \oplus B \stackrel{\sim}{\to} \Phi A$ we have a commutative square
\[ \xymatrix{
A \oplus B \ar[d]_{\eta_A \oplus \eta_B} \ar[r]_\cong^\nu & \Phi A \ar[d]^{\eta_{(\Phi A)}}_\cong \\
\Phi A \oplus \Phi B \ar[r]_{\Phi \nu}^\cong  & \Phi \Phi A
} \]
We will apply this to our situation with $\Phi = f_*f^*$ and $A = E_X$. Due to the invertibility of $d$ and the trace morphism, the morphism $E_X \to f_*f^*E_X$ is a monomorphism, and every monomorphism in a triangulated category splits. So $E_X$ is a direct summand of $f_*f^*E_X$. To prove the lemma, it suffices then to see that $f_*f^*E_X \to f_*f^*f_*f^*E_X$ is an isomorphism. We make the cartesian square
\[ \xymatrix{
Y' \ar[r]^q \ar[d]_p & Y \ar[d]^f \\
Y \ar[r]_f & X
} \]
By projective base change (Theorem~\ref{theo:ayo}(\ref{item:properBaseChange})) we have an isomorphism $f_*(f^*f_*)f^* \cong f_*p_*q^*f^*$ and since the square is commutative an isomorphism $f_*p_*q^*f^* \cong f_*p_*p^*f^*$. Now $p: Y' \to Y$ admits a section which is a closed immersion (since all our schemes are separated) and surjective (since $p$ is radicial). Consequently, $p^*$ is an equivalence of categories (Theorem~\ref{theo:ayo}(\ref{item:nilpotent})), and it follows that $id \to p_*p^*$ is an isomorphism. So we have reduced to showing the commutativity of the following square
\[ \xymatrix{
f_*f^* \ar[r]_\cong \ar[d]_{\eta_{(f_*f^*)}} & f_*p_*p^*f^* \ar[d]^{\alpha}_\cong \\
f_*f^*f_*f^* \ar[r]^\cong & f_*p_*q^*f^*
} \]
where $\alpha$ is the comparison $p^*f^* \cong q^*f^*$. The commutativity of this square follows from the commutativity of the following diagram since the lower row is precisely the morphism which projective base change (Theorem~\ref{theo:ayo}(\ref{item:properBaseChange})) states is an isomorphism
\[ \xymatrix{
f_*f^* \ar[r] \ar[d]_{\eta_{(f_*f^*)}} & f_*p_*p^*f^* \ar[d]_{\eta_{(f_*p_*p^*f^*)}} \ar@{=}[r]^\alpha & f_*p_*q^*f^* \ar@{=}[dr]^{id} \ar[d]_{\eta_{(f_*p_*q^*f^*)}} \\
f_*f^*(f_*f^*) \ar[r] & f_*f^*(f_*p_*p^*f^*) \ar@{=}[r]_\alpha & f_*f^*(f_*p_*q^*f^*) \ar[r]_{f_* \epsilon_{p_*q^*f^*}} & f_*p_*q^*f^*
} \]
We have used $\eta$ for the units of adjunction and $\epsilon$ for the counit. The commutativity of the squares is just the naturality of the transformations $\eta$, and the commutativity of the triangle is from the definition of adjunction: $(f_* \epsilon_A) \circ \eta_{f_*A} = id_{f_*A}$.
\end{proof}

\begin{lemm} \label{lemm:phiphi}
Let $\H^*$ be a sliceable 2-functor, $X \in Sch(S)$, $E \in \H(X)$, and $q \in \ZZ$. Suppose we have an endofunctor $\Phi: \H(X) \to \H(X)$ that preserves colimits and is equipped with a natural transformation $id \to \Phi$ such that the morphisms
\[ E \to \Phi E,  \qquad \textrm{ and } \qquad s_r E \to \Phi s_rE \]
are isomorphisms for all $r < q$. Then the morphism
\[ f_qE \to \Phi f_qE \] 
is an isomorphism as well.
\end{lemm}

\begin{proof}
We have a morphism of distinguished triangles
\[ \xymatrix{
f_qE \ar[r] \ar[d] & E \ar[r] \ar[d] & s_{< q}E \ar[r] \ar[d] & f_qE[1] \ar[d] \\
\Phi f_qE \ar[r]  & \Phi E \ar[r] & \Phi s_{< q}E \ar[r] & \Phi f_qE[1]
} \]
from which we see that $f_q E \to \Phi f_q E$ is an isomorphism if and only if $s_{< q} E \to \Phi s_{< q} E$ is an isomorphism. We will prove the latter. Recall that there is a canonical isomorphism $E \cong hocolim_{r < q} f_r E$. Since all the functors in question commute with colimits, it suffices to prove that
\[  \xymatrix{ s_{< q} f_rE \ar[r]^{\eta_r} & \Phi s_{< q} f_rE } \]
is an isomorphism for all $r < q$. We do this by induction.

In the case $r = q - 1$ we find the following commutative square
\[ \xymatrix{
s_{q - 1}E \ar[r] \ar[d] & \Phi s_{q - 1}E \ar[d] \\
s_{< q} f_{q - 1}E \ar[r]_{\eta_{q - 1}} &  \Phi s_{< q} f_{q - 1}E
} \]
where the vertical morphisms are isomorphisms\footnote{One can see this by considering the distinguished triangle $s_{q - 1}f_{q - 1} \to s_{< q} f_{q - 1} \to s_{< q - 1} f_{q - 1} \to s_{q - 1}f_{q - 1}[1]$.} and the upper morphism is an isomorphism by assumption. So assume that our inductive hypothesis is true for $r + 1$. We have the following morphism of distinguished triangles
\[ \xymatrix{
s_{< q} f_{r + 1}E \ar[r] \ar[d]_{\eta_{r + 1}} & s_{< q} f_rE \ar[r] \ar[d]_{\eta_{r}} & s_{< q}s_r E \ar[r] \ar[d] & s_{< q} f_{r + 1}E[1] \ar[d]^{\eta_{r + 1}[1]} \\
\Phi s_{< q} f_{r + 1}E \ar[r]  & \Phi s_{< q} f_rE \ar[r] & \Phi s_{< q} s_rE \ar[r] & \Phi s_{< q} f_{r + 1}E[1]
} \]
and due to the inductive hypothesis and the fact that the natural transformation $s_{< q} s_r \cong s_r$ is an isomorphism, the result is proven.
\end{proof}

\begin{rema} \label{rema:genericDecomposition}
If $Y \to X$ is the morphism given in Theorem~\ref{theo:gabberGlobal}, notice that there exists a non-empty open subscheme $U \subset X$ such that the induced morphism $f: Y \times_X U \to U$ satisfies:
\begin{enumerate}
 \item[(*)] the morphism $f_{red}$ is a composition $\stackrel{r}{\to} \stackrel{e}{\to}$ where $r$ is a radicial finite flat surjective morphism and $e$ is an {\'e}tale finite surjective morphism, and both are morphisms between smooth $k$-schemes.
\end{enumerate}
\end{rema}

For the following results we use the following hypotheses:

\begin{enumerate}
 \item[(**)] Let $\H^*$ be a stable homotopy 2-functor, $S$ the spectrum of a $k$ a perfect field of exponential characteristic $p$, and $q \in \ZZ$. Suppose that $E \in \H(k)$ is a $\zpi$-local object such that $s_rE$ and $E$ have a weak structure of smooth traces for every $r \leq q + 1$. Let $f: Y \to X$ be a finite flat surjective morphism in $Sch(k)$ and $a: X \to k$ the structural morphism.
\end{enumerate}

\begin{lemm} \label{lemm:sqffEffective}
Assume the hypotheses (**). If $f$ satisfies the condition (*) of Remark~\ref{rema:genericDecomposition}. Then
\[ f_*f^*(f_qE_X) \in \Sigma^q\H(X)^{eff}. \]
\end{lemm}

\begin{proof}
Let $\tilde{i}: Y_{red} \to Y$ and $i: X_{red} \to X$ be the canonical closed immersions and $Y_{red} \stackrel{r}{\to} W \stackrel{e}{\to} X_{red}$ the factorisation.
\[ \xymatrix{
Y_{red} \ar[r]^r \ar[d]_{\tilde{i}} & W \ar[r]^e & X_{red} \ar[d]^i \\
Y \ar[rr]_f && X
} \]
The canonical natural transformation $id \to \tilde{i}_*\tilde{i}^*$ is a natural isomorphism (Theorem~\ref{theo:ayo}(\ref{item:nilpotent})) and so the canonical morphism $f_*f^* (f_q E_X) \to f_*\tilde{i}_*\tilde{i}^*f^* (f_qE_X)$ is an isomorphism. So it suffices to show that $f_*\tilde{i}_*\tilde{i}^*f^* (f_qE_X) \in \Sigma^q\H(X)^{eff}$.The morphism $f \tilde{i} $ also factors as $i e r$. Now both $i^*$ and $e^*$ commute with $f_q$ (Corollary~\ref{coro:smoothPullback}) and $i_*$ and $e_*$ preserve $\Sigma^q\H^{eff}$ (Lemma~\ref{lemm:preserve}), so it suffices to show that $r_*r^* (f_qE_W) \in \Sigma^q\H(W)^{eff}$. By additivity, it suffices to consider the case when $W$ and $Y_{red}$ are connected. That is, we assume that $r: Y_{red} \to W$ is a radicial finite flat surjective morphism between \textbf{connected} smooth $k$-schemes in $Sch(k)$.

We claim that $f_qE_W \to r_*r^*f_qE_W$ is actually an isomorphism. By Lemma~\ref{lemm:phiphi} to prove this claim it suffices to show that $E_W \to r_*r^*E_W$ and $s_rE_W \to r_*r^*s_rE_W$ are isomorphisms for all $r < q$ (to see that $r_*r^*$ preserves colimits, notice that it has a right adjoint $r_*r^!$ by Theorem~\ref{theo:ayo}(\ref{item:projRightAdjoint})). To show that these are isomorphisms, by Lemma~\ref{lemm:radicialIso} it suffices to show that $E \in \H(k)$ is $\ZZ[\tfrac{1}{d}]$-local where $d = \deg(Y_{red} \to W)$. By assumption $E$ is $\zpi$-local and so it suffices to show that $d$ is a power of $p$. Using the assumption that $Y_{red}$ and $W$ are connected, we have $d = [k(Y_{red}):k(W)]$ and since this is radicial, its degree must be a power of $p$, and we are done.
\end{proof}

\begin{prop} \label{prop:ffPreserve}
Assume the hypotheses (**). For all $r \leq q$
\[ f_*f^*(f_rE_X) \in \Sigma^r\H(X)^{eff}. \]
\end{prop}

\begin{proof}
It suffices to consider the case $q = r$ because the hypotheses are stable under lowering $q$. We use induction on the dimension of $X$. Suppose that $f: Y \to X$ is a finite flat surjective morphism. Since $f$ satisfies the property (*) of Remark~\ref{rema:genericDecomposition} generically, there exists a dense open $U$ of $X$ such that $U \times_X f$ satisfies the property (*). We form the following cartesian squares
\[ \xymatrix{
Z' \ar[d]_{\tilde{f}} \ar[r] & Y \ar[d]^f & U' \ar[d]^{g} \ar[l] \\
Z \ar[r]_i & X & \ar[l]^j U 
} \]
Consider the exact triangle $j_!j^! \to id \to i_*i^* \to j_!j^![1]$ evaluated on the object $f_*f^*(f_qE_X)$:
\[ j_!j^!(f_*f^*(f_qE_X)) \to f_*f^*(f_qE_X) \to i_*i^*(f_*f^*(f_qE_X)) \to j_!j^!(f_*f^*(f_qE_X))[1]. \]
By projective base change (Theorem~\ref{theo:ayo}(\ref{item:properBaseChange})) this triangle is isomorphic to the triangle
\[ (j_!g_*g^*j^*)(f_qE_X) \to f_*f^*(f_qE_X) \to (i_*\tilde{f}_*\tilde{f}^*i^*)(f_qE_X) \to (j_!g_*g^*j^*)(f_qE_X)[1]. \]
We will show that $(j_!g_*g^*j^*)(f_qE_X)$ and $(i_*\tilde{f}_*\tilde{f}^*i^*)(f_qE_X)$ are in $\Sigma^q\H(X)^{eff}$ and the result will follow since $\Sigma^q\H(X)\eff$ is triangulated.

By Theorem~\ref{theo:pullbackAll} we have an isomorphism $i^*f_qE_X \cong f_qi^*E_X$, by definition $f_qi^*E_X = f_qE_Z$, and so by induction $\tilde{f}_*\tilde{f}^*(f_qE_Z) \cong {f}_*\tilde{f}^*i^*(f_qE_X) \in \Sigma^q\H(Z)\eff$. Lastly, by Lemma~\ref{lemm:preserve} $i_*$ preserves $\Sigma^q\H\eff$ and so $(i_*\tilde{f}_*\tilde{f}^*i^*)(f_qE_X) \in \Sigma^q\H(X)\eff$.

For the other corner of the triangle, by Theorem~\ref{theo:pullbackAll} we have an isomorphism $j^*f_qE_X \cong f_qj^*E_X$, by definition $f_qj^*E_X = f_qE_U$, and so by Lemma~\ref{lemm:sqffEffective} $g_*g^*f_qE_U \cong g_*g^*j^*(f_qE_X) \in \Sigma^q\H(U)\eff$. Lastly, by Lemma~\ref{lemm:preserve} $j_!$ preserves $\Sigma^q\H$ and so $j_!g_*g^*j^*(f_qE_X) \in \Sigma^q\H(X)\eff$.
\end{proof}

\begin{prop} \label{prop:ffslices}
Assume the hypotheses (**). The functor $f_*f^*$ on $\H(X)$ satisfies (Pel0)$_q$, (Pel1)$_{\leq q + 1}$, and (Pel2)$_{\leq q}$ for $E_X$. Consequently, the morphisms
\[ \alpha_r(E): f_*f^*(f_r E_X) \to  f_r f_*f^*E_X \]
\[ \beta_r(E): f_*f^*(s_r E_X) \to  s_r f_*f^*E_X \]
are isomorphisms for all $r \leq q$.
\end{prop}

\begin{proof}
Proposition~\ref{prop:ffPreserve} says precisely that (Pel1)$_{\leq q + 1}$ is satisfied. Both $f_*$ and $f^*$ are right adjoints (Theorem~\ref{theo:ayo}) and so (Pel0) is satisfied. Consider $s_rE_X$ for some $r \leq q$. By Theorem~\ref{theo:pullbackAll} there is a canonical isomorphism $f^*s_rE_X \cong s_rf^*E_X$ and so since $f_*$ preserves $SH^\perp$ it follows that (Pel2)$_{\leq q}$ is satisfied. For the stated isomorphisms we need only to recall Theorem~\ref{theo:criteria}.
\end{proof}

\begin{rema} \label{rema:functorialitybetaflowerstar}
We note some consequences of Proposition~\ref{prop:ffslices}. We keep the notation used in the statement. Combining this proposition with Theorem~\ref{theo:pullbackAll} we have canonical induced isomorphisms $\beta: f_*s_rf^*E_X \stackrel{\sim}{\to} s_rf_*f^*E_X$ that fit into diagrams 
\[ \xymatrix{ s_qf_*f^*E_X \ar[r]_\beta \ar@/^12pt/[rr] & f_*s_qf^*E_X \ar[r] & f_*f^*s_qE_X  } \]
Consequently, these $\beta$ satisfy the same functoriality as those mentioned in Lemma~\ref{lemm:comparisonMorphismssf}. The same applies to isomorphisms $\beta: f_*s_r(af)^*E \stackrel{\sim}{\to} s_rf_*(af)^*E$ and the analogous $\alpha$ with $f_r$.
\end{rema}

\section{Traces in the context of a stable homotopy 2-functor} \label{sec:traces}

In this section we develop a notion of an object $E$ of $\H(S)$ having a structure of traces. We show that this induces a structure of traces on the slices $s_qE$ (and the same proof shows that there is an induced structure of traces on the connective covers $f_qE$. 

\subsection{Definition}

We make the following definition.

\begin{defi} \label{defi:sectionWithTraces}
Let $\H_*$ be a covariant 2-functor assigning to every object $X \in Sch(S)$ an additive category $\H(X)$, and each morphism $f: Y \to X$ in $Sch(S)$ an additive functor $f_*: \H(Y) \to \H(X)$. Let $E_-$ be a section of $\H_*$. That is, for each scheme $X$ we are given an object $E_X \in \H(X)$, for each morphism $f: Y \to X$ of schemes we have a morphism $c_f: E_X \to f_*E_Y$ and these morphisms satisfy a suitable coherency condition.

A \emph{structure of traces} on the section $E_-$ is the data of a morphism $Tr_f: f_*E_Y \to E_X$ in $\H(X)$ for each finite flat surjective morphism $f: Y \to X$ in $Sch(S)$ and these morphisms are required to satisfy the following axioms.

\begin{enumerate}
 \item[(Fon)] If we have $W \stackrel{g}{\to} Y \stackrel{f}{\to} X$ in $Sch(S)$ with $f$ and $g$ finite flat surjective then $\Tr_{fg} = \Tr_f \circ f_*\Tr_g$. That is, the following diagram commutes.
\begin{equation} \label{equa:fonDef}
\xymatrix{
f_*g_*E_W \ar[d]_{f_*\Tr_g} \ar[r]^\cong & (fg)_*E_W \ar[d]^{\Tr_{fg}} \\
f_*E_Y \ar[r]_{\Tr_f} & E_X
} \end{equation}
where the isomorphism is the connection isomorphism $f_*g_* \stackrel{\cong}{\to} (fg)_*$.

 \item[(CdB)] Suppose that (\ref{equa:cartSquareDefiTraces}) is a cartesian square in $Sch(S)$ with $f$ finite flat surjective. Then $c_p \circ \Tr_f = p_*\Tr_g  \circ f_*c_q$. That is, the following diagram commutes
\begin{equation} \label{equa:cdbDef}
\xymatrix{
p_*g_*E_{Y \times_X W} \ar[r]^{p_*\Tr_g} & p_*E_W \\
f_*q_*E_{Y \times_X W} \ar[u]^\cong \\
f_*E_Y \ar[r]_{\Tr_f} \ar[u]^{f_*c_q} & E_X \ar[uu]_{c_p}
} \end{equation}
where the isomorphism is built out of the connection morphisms of the 2-functor $H_*$.

 \item[(Deg)] If we have $Y \stackrel{f}{\to} X$ in $Sch(S)$ with $f$ a finite flat surjective morphism of constant degree $d$ then the composition of $\Tr_f: f_*E_Y \to E_X$ with the connection morphism $c_f: E_X \to f_*E_Y$ is $d$ times the identity. That is, we have
\[ \Tr_f c_f = d \cdot id_{E_X}. \]
\end{enumerate}
\end{defi}

\begin{lemm} \label{lemm:cdbdash}
Continuing with the assumptions and notation of of Definition~\ref{defi:sectionWithTraces} suppose that for every morphism $p: W \to X$ in $Sch(S)$ the functor $p_*$ has a left adjoint $p^*: \H(X) \to \H(W)$. Then (CdB) is equivalent to:
 \item[(CdB$'$)] The following diagram commutes
\begin{equation} \label{equa:cdbDef}
\xymatrix{
g_*E_{Y \times_X W} \ar[r]^{\Tr_g} & E_W \\
g_*q^*E_{Y} \ar[u]^{g_*c_q'} \\
p^*f_*E_Y \ar[r]_{p^*\Tr_f} \ar[u] & p^*E_X \ar[uu]_{c_p'}
} \end{equation}
where the $c'$ are the adjoints to the $c$ and the unlabelled morphism is the canonical comparison morphism built from adjunctions $p^*f_* \to p^*f_*q_*q^* = p^*p_*g_*q^* \to g_*q^*$.
\end{lemm}

\begin{proof}
This is an exercise in adjunctions that is left to the reader.
\end{proof}

\begin{defi}
In the notation and assumptions of Lemma~\ref{lemm:cdbdash}, suppose we are given an object $E \in \H(S)$ over the base scheme. A \emph{structure of traces on $E$} is a structure of traces on the canonical section that associates to $a: X \to S$ the object $a^*E$, and to a morphism $f: Y \to X$ the unit of the adjunction $c_f: a^*E \to f_*f^*a^*E = f_*(af)^*E$. That is, for every $Y \stackrel{f}{\to} X \stackrel{a}{\to} S$ with $f$ finite flat surjective, we have a morphism
\[ \Tr_f: f_*(af)^*E \to a^*E \]
and these morphisms satisfy the appropriate axioms.
\end{defi}

The following two lemmata are clear from the definitions.

\begin{lemm} \label{lemm:tracesImpliesTraces}
In the notation of of Definition~\ref{defi:sectionWithTraces}, let $E_-$ be a section of $H_*$. The presheaf $F: Sch(S) \to \H(S)$ that sends an $S$-scheme $a: X \to S$ to the object $a_*E_X$ and a morphism $f: Y \to X$ to the morphism $a_*E_X \to a_*f_*E_Y \cong (af)_*E_Y$ has a canonical structure of presheaf with traces in the sense of Definition~\ref{defi:traces}.
\end{lemm}

\begin{rema} \label{rema:tracesImpliesTraces}
An immediate consequence of this lemma is that for every object $E' \in \H(S)$, the presheaf of abelian groups sending an $S$-scheme $a: X \to S$ to the abelian group $\hom_{\H(S)}(E', a_*E_X)$ also have a structure of presheaf with traces. In particular, if $\H$ is the Morel-Voevodsky stable homotopy category and $E \in \SH(S)$ is an object with traces, then for each $p, q \in \ZZ$ the presheaf on $Sch(S)$ that takes $a: X \to S$ to $\hom_{\SH(S)}(\Sigma^{-q}\un_S[2q - p], a_*a^*E)$ has a canonical structure of traces. Due to the adjunction $(a_\#a^*, a_*a^*)$ when $a$ is smooth, the restriction of this presheaf to $Sm(S)$ agrees with the cohomology sheaf $E^{p,q}(-)$ of $E$ defined in \cite[Section 6]{Voev98}.
\end{rema}

\begin{lemm} \label{lemm:localisationOfTraces}
Suppose that $\H^1_*$ and $\H^2_*$ are two 2-functors as in Definition~\ref{defi:sectionWithTraces} and $\phi: H^1_* \to H^2_*$ is a morphism between them. Let $E_-$ be a section of $H^1_*$. If $E_-$ has a structure of traces, then there is a canonical induced structure of traces on the canonical section $\phi E_-$ of $H^2_*$.
\end{lemm}

\subsection{Traces on slices}

\begin{prop} \label{prop:tracesOnSlices}
Let $\H^*$ be a stable homotopy 2-functor, suppose $S$ is the spectrum of a perfect field $k$, and $p$ its exponential characteristic. Suppose that $E \in \H(k)$ is a $\zpi$-local object with a structure of traces, and such that $s_rE$ has a weak structure of smooth traces for all $r \leq q + 1$. Then $f_qE$ and $s_qE$ both have canonical structures of traces.
\end{prop}

\begin{proof}
The proof for $f_qE$ and $s_qE$ is the same. We give the proof for $s_qE$.

Consider morphisms of schemes $Y \stackrel{f}{\to} X \stackrel{a}{\to} k$ with $f$ finite flat surjective. After Proposition~\ref{prop:ffslices} and Theorem~\ref{theo:pullbackAll} we have canonical isomorphisms $f_*(af)^*s_qE \cong s_qf_*(af)^*E$ and $a^*s_qE \cong s_qa^*E$ which are functorial in an appropriate way (see Lemma~\ref{lemm:comparisonMorphismssf} for the details). These give rise to candidate trace morphisms
\[ f_*(af)^*s_qE \cong s_qf_*f^*a^*E \stackrel{s_q\Tr_f}{\to} s_qa^*E \cong a^*s_qE \]
induced by the trace morphisms $\Tr_f$ of $E$. We will label these new morphisms $\Tr^{s_q}_f$. The diagrams that we wish to prove commute are the following.

Functoriality:
\[ \xymatrix{
f_*g_*(afg)^*s_qE \ar[r]^\cong \ar[d]_{f_*\Tr^{s_q}_g} & (fg)_*(afg)^*s_qE \ar[d]^{\Tr^{s_q}_{fg}} \\
f_*(af)^*s_qE \ar[r]_{\Tr^s_f} & a^*s_qE
} \]

Base-change (Lemma~\ref{lemm:cdbdash}):
\[ \xymatrix{
g_*(apg)^*s_qE \ar[r]^{\Tr_g^{s_q}} & (ap)^*s_qE \\
g_*q^*(af)^*s_qE \ar[u] \\
p^*f_*(af)^*s_qE \ar[r]_{p^*\Tr_f^{s_q}} \ar[u] & p^*a^*s_qE \ar[uu]
} \]

Degree:
\[ \xymatrix{
a^*s_qE \ar@/_18pt/[rr]_{d \cdot id_{a^*s_qE}} \ar[r] & f_*(af)^*s_qE \ar[r]^-{\Tr^{s_q}_f} & a^*s_qE
} \]

Each of these diagrams arises in the following way. We begin with a 2-category $I$ of a special form: there exists some positive integer $n$ and a 2-functor $I \to \{0, \dots, n\}$ sending each object of $I$ to a unique object of the totally ordered set $\{0, \dots, n\}$ considered as a 2-category with no non-identity 2-morphisms. We identify the objects of $I$ with the objects of $\{0, \dots, n\}$. Then we have a 2-diagram $F: I \to Cat$ such that there exists a (not necessarily unique) scheme $X_i$ for each object $i$ of $I$ such that $F(i) = \H(X_i)$. We have an object $E' \in \H(X_0)$ and consequently, an induced diagram $F_{E'}$ in $\H(X_n)$ indexed by the 1-category $\hom_I(0, n)$.

For example, in the case of (Deg), $E' = s_q E$, and we could take $S, X, Y, X$ to be the sequence $X_0, \dots, X_n$ of schemes. The 1-functors involved are the $a^*, f^*, f_*, (af)^*, d \cdot id_{\H(X)}$, and their various compositions such as $f^*a^*, f_*f^*a^*, f_*(af)^*, $etc. The two functors are made from the various connection isomorphisms such as $(af)^* \cong f^*a^*$ and their horizontal and vertical compositions. What we would like is that the $\beta$ of Lemma~\ref{lemm:comparisonMorphismssf} induce an isomorphism of diagrams between the diagram $F_{s_qE}$ just described, and the diagram $s_qF_E$ obtained in the same way, but starting with $E$ and and applying $s_q$ at the end. The functoriality described in Lemma~\ref{lemm:comparisonMorphismssf} says precisely that this is true.
\end{proof}

\subsection{Traces on products}

Now we continue with a covariant 2-functor $\H_*$ as Definition~\ref{defi:sectionWithTraces} but we further assume that it factors through the 2-category of additive monoidal categories with lax functors. That is, each of the categories $\H(X)$ is equipped with a product $\otimes$ and for every morphism $f: Y \to X$ of $S$-schemes we have a binatural transformation $f_*(-) \otimes f_*(-) \to f_*(- \otimes -)$ which are not required to be isomorphisms (and in practice they won't be). These binatural transformations are required to be compatible with the isomorphisms $(gf)_* = f_*g_*$ in the obvious way.

Given such a structure, the category of sections of $\H_*$ has an obvious product structure where the product of two sections $E_-, F_-$ associates to a scheme $X$ the object $E_X \otimes F_X$ and to a morphism $f: Y \to X$ the composition $E_X \otimes F_X \to f_*E_Y \otimes f_*F_Y \to f_*(E_Y \otimes F_Y)$.

\begin{exam}
If $\H$ a unitary monoidal stable homotopy 2-functor in the sense of \cite[Definition 2.3.1]{Ayo07} then all the above assumptions are satisfied. The functors $f_*$ are lax monoidal due to the $f^*$ being strong monoidal (i.e., the $f^*(- \otimes -) \to f^*- \otimes f^*-$ are isomorphisms). In particular, this applies to the Morel-Voevodsky stable homotopy category $\SH$ as well as the stable homotopy 2-functor obtained from a ring spectrum in $\SH(S)$ by taking the homotopy category of its category of modules.
\end{exam}

\begin{defi} \label{defi:cartesianSection}
With the assumptions and notation just established, we will say that a section $E_-$ is \emph{cartesian} if for every section $F_-$ and every projective morphism $f: Y \to X$ the morphism $E_X \otimes f_*F_Y \to f_*(E_Y \otimes F_Y)$ is an isomorphism.
\end{defi}

\begin{exam}
In the example of a unitary monoidal stable homotopy 2-functor mentioned above, if we have a section $E_-$ which is cartesian in the sense that the connection morphisms $f^*E_X \to E_Y$ are isomorphisms, then $E_-$ is cartesian in the sense of Definition~\ref{defi:cartesianSection} (see \cite[Theorem 2.3.40, Theorem 1.7.17]{Ayo07}).
\end{exam}

\begin{prop} \label{prop:tracesOnProducts}
Let $\H_*$ be a covariant 2-functor of additive monoidal categories with lax functors as described above. Suppose that $E_-$ and $F_-$ are two sections. If $E_-$ is cartesian (Definition~\ref{defi:cartesianSection}) and $F_-$ has a structure of traces, then there is a canonical structure of traces on the product $(E \otimes F)_-$.
\end{prop}

\begin{proof}
Suppose $f: Y \to X$ is a finite flat surjective $S$-morphism. To define trace morphisms $f_*(E_Y \otimes F_Y) \to E_X \otimes F_X$ we use the isomorphism $f_*(E_Y \otimes F_Y) \stackrel{\sim}{\leftarrow} E_X \otimes f_*F_Y$ coming from the assumption that $E_-$ is cartesian, composed with the traces on $F_-$. We will denote these morphisms by $\Tr_f^\otimes$.

Each of the axioms are satisfied as a result of the functoriality and compatibility conditions that we have asked for. Here are the diagrams. 

Functoriality:
\[ \xymatrix{
f_*g_*(E_W \otimes F_W) \ar@{=}[rr] && (fg)_*(E_W \otimes F_W) \\
f_*(E_Y \otimes g_* F_W) \ar[d] \ar[u]^\cong & \ar[l]_\cong E_X \otimes f_*g_*F_W \ar@{=}[r] \ar[d] & E_X \otimes (fg)_*F_W \ar[d] \ar[u]_\cong  \\
f_*(E_Y \otimes F_Y) & E_X \otimes f_*F_Y \ar[r] \ar[l]_\cong & E_X \otimes F_X
} \]

Base-change:
\[ \xymatrix{
p_*g_*(E_{Y \times_X W} \otimes F_{Y \times_X W}) & \ar[l]_-\cong p_*(E_Y \otimes g_*F_Y) \ar[r] & p_*(E_Y \otimes F_Y) \\
& E_X \otimes p_*g_*F_{Y \times_X W} \ar[r] \ar[u]  & E_X \otimes p_*F_{Y \times_X W} \ar[u]  \\
f_*q_*(E_{Y \times_X W} \otimes F_{Y \times_X W}) \ar@{=}[uu] & \ar[l] E_X \otimes f_*q_*F_{Y \times_X W} \ar@{=}[u] & \\
f_*(E_Y \otimes F_Y) \ar[u] & \ar[l]_\cong E_X \otimes f_*F_Y \ar[r] \ar[u] & E_X \otimes F_X \ar[uu]
} \]

Degree:
\[ \xymatrix{
& E_X \otimes f_*F_Y \ar[d] \ar[dr]_\cong \ar[r]^{E_X \otimes \Tr_f} & E_X \otimes F_X \\
E_X \otimes F_X \ar[r] \ar[ur]^{E_X \otimes c_f} & f_*E_Y \otimes f_*F_Y \ar[r] & f_*(E_Y \otimes F_Y) \ar[u]_{\Tr^\otimes_f}
} \]
\end{proof}

\chapter{Motivic applications} \label{chap:applications}

\section{Introduction}

\lettrine[lines=2, findent=3pt, nindent=0pt]{I}{n} this chapter we use the previous material to give a proof of Theorem~\ref{theo:MainResultRos} which is our main technical result. We then demonstrate how this theorem may be applied to obtain $\zpi$ linear versions of results that previously assumed resolution of singularties.

In Section~\ref{sec:objectsWithTrace} we show that the object $\HZ$ representing motivic cohomology in the Morel-Voevodsky stable homotopy category has a weak structure of smooth trace (Definition~\ref{defi:wsst}, Proposition~\ref{prop:weakTracesOnHZ}) and that the object representing algebraic $K$-theory has a structure of traces (Definition~\ref{defi:sectionWithTraces}, Proposition~\ref{prop:ktheoryTraces}). Applying the material of the previous chapter and a theorem of Levine (\cite[Theorems 6.4.2 and 9.0.3]{Lev08}), this implies that $\HZpi$ has a structure of traces (Corollary~\ref{coro:tracesOnHZmod}).

In Section~\ref{sec:RoS} we prove Theorem~\ref{theo:MainResultRos}. The main technical results that we use are Corollary~\ref{coro:tracesOnHZmod}, Theorem~\ref{theo:cdhlprimeComparison}, and a result of Cisinski applying a theorem of Ayoub that says that every object in the Morel-Voevodsky stable homotopy category satisfies cdh descent.

In Section~\ref{sec:bivCycCoh} and Section~\ref{sec:traCatMotFie} we show how Theorem~\ref{theo:MainResultRos} implies $\zpi$-linear versions of all the results in \cite{FV} and \cite{Voev00} without having to use resolution of singularities. We show in Section~\ref{sec:higChoGro} how this works for \cite{Sus00}.

In Section~\ref{sec:kTheory} we use the $\ldh$ topology and the theorem of Gabber to give a partial answer to a conjecture of Weibel about vanishing of algebraic $K$-theory (Theorem~\ref{theo:KtheoryVanishing}).

\section{Some objects of $\SH(k)$ with traces} \label{sec:objectsWithTrace}

In this section we show that the object representing motivic cohomology in $\SH(S)$ has a weak structure of smooth traces, and the object representing algebraic $K$-theory has a structure of traces.

\begin{prop} \label{prop:weakTracesOnHZ}
Suppose $S$ is a noetherian scheme. The object $\HZ \in \mathsf{SH}(S)$ that represents motivic cohomology (\cite[Section 6.1]{Voev98}) has a weak structure of smooth traces (Definition~\ref{defi:wsst}).
\end{prop}

\begin{rema} \label{rema:twoHZs}
We can construct by hand a structure of traces on the section $\HZ_-$ which assigns to each scheme $X \in Sch(S)$ the object $\HZ_X$ representing motivic cohomology defined by Voevodsky. This is a consequence of the component terms of each spectrum $\HZ_X$ being presheaves with transfers on $Sch(S)$. However, for our purposes we need a structure of traces on the section $(-)^*\HZ_k$ determined by the object $\HZ_k \in \SH(k)$, and for non-smooth schemes $a: X \to S$ it is an open conjecture (\cite[Conjecture 17]{Voe02}) whether $a^*\HZ_S$ is isomorphic to $\HZ_X$.
\end{rema}

\begin{proof}
Let $Comp(Shv_{Nis}(SmCor(S)))$ denote the category of unbounded chain complexes in the abelian category $Shv_{Nis}(SmCor(S))$, and denote by $D(Comp(Shv_{Nis}(SmCor(S))))$ its associated derived category (obtained by localising at quasi-isomorphisms). The category $DM\eff(S)$ is by definition the localisation of $D(Comp(Shv_{Nis}(SmCor(S))))$ at the class of morphisms $L_{Nis}(\AA^1_X) \to L_{Nis}(X)$ for all $X \in Sm(S)$ (recall the notation from Definition~\ref{defi:representableSheafWithTransfers}).

First we claim that every object of $Comp(Shv_{Nis}(SmCor(S)))$ has a weak structure of smooth traces as a section of $Comp(Shv_{Nis}(SmCor(-)))$. Consider $Shv_{Nis}(SmCor(-))$ as a 2-functor  on $Sch(S)$. For any smooth scheme $a: X \to S$ the functor $a^*$ is just restriction $(-)|_{SmCor(X)}: Shv_{Nis}(SmCor(S)) \to Shv_{Nis}(SmCor(X))$ and for any morphism $f: X \to S$ the functor $f_*$ is composition with $X \times_S -: SmCor(S) \to SmCor(X)$. Let $a: X \to S$ be a smooth morphism and $f: Y \to X$ a finite flat surjective morphism with $af$ smooth as well. Note that since $f$ is finite, $f_*: PreShv(Y) \to PreShv(X)$ is exact and preserves Nisnevich sheaves. We can explicitely describe the functors $f_* (af)^*$ and $a^*$ on $Comp(Shv_{Nis}(SmCor(S)))$ by evaluating them on a sheaf $F \in Shv_{Nis}(SmCor(S))$, and describing the two resulting sheaves $f_* (af)^*F$ and $a^*F$ in $Shv_{Nis}(SmCor(X))$ by evaluating them on an object $U \in SmCor(X)$. We have $(f_*(af)^*F)(U) = F(Y \times_X U)$ and $((af)^*F)(U) = F(U)$, and the correspondence $[^tf]: [X] \to [Y]$ in $SmCor(S)$ (Definition~\ref{defi:graphAndTranspose}) gives us a morphism between these two groups. Since (CdB) is satisfied in $SmCor(S)$ (Proposition~\ref{prop:corHasTraces}) these morphisms are functorial in the appropriate way and we obtain a canonical natural transformation $f_*(af)^* \to a^*$. Moreover, since (Deg) is satisfied in $SmCor(S)$ (Proposition~\ref{prop:corHasTraces}), the composition $a^* \to f_*(af)^* \to a^*$ is $d$ times the identity when $f$ is of constant degree $d$. Hence, the claim.

Let $\mathbb{L}$ denote the cokernel of the morphism $L(s): L_{Nis}(S) \to L_{Nis}(\PP^1_S)$ given by the section $s: S \to \PP^1_S$ at infinity. To obtain the category $DM(S)$ we formally adjoint a tensor inverse to $\mathbb{L}$. That is, we consider the category $Sp_{\mathbb{L}}(Comp(Shv_{Nis}(SmCor(S))))$ of $\mathbb{L}$-spectra in $Comp(Shv_{Nis}(SmCor(S)))$. Such a spectrum is a sequence $(K_0, K_1, \dots )$ of objects of $Comp(Shv_{Nis}(SmCor(S)))$ together with connection morphisms $K_n \to \ihom(\mathbb{L}, K_{n + 1})$. Let $p: \PP^1_S \to S$ be the canonical projection. We will use the same notation $p$ for bases other than $S$ as well. Let $\Omega_S = ker(p_*p^*(-) \to (-))$ where the morphism is induced by the unit of the adjunction $id \to s_*s^*$ and the identity $ps = id$. There is a canonical isomorphism $\Omega_S \cong \ihom(\mathbb{L}, -)$. 

To show that the trace morphisms we defined above pass to $\mathbb{L}$-spectra, we must show that the following square is commutative.
\[ \xymatrix{
f_*(af)^*K_n \ar[r] \ar[d] &  \Omega_X f_*(af)^*K_{n + 1} \ar[d] \\
a^*K_n \ar[r] & \Omega_X a^*K_{n + 1} 
} \]
We can see this immediately by evaluating on an object $U \in Sm(X)$ as we obtain the following square.
\[ \xymatrix{
K_n(Y \times_X U) \ar[r] \ar[d] & \ker \biggl (K_{n +1}(Y \times_X U \times_X \PP^1_X) \to K_{n +1}(Y \times_X U) \biggr ) \ar[d] \\
K_n(U) \ar[r] & \ker \biggl (K_{n +1}(U \times_X \PP^1_X) \to K_{n +1}(U) \biggr )
} \]
For a morphism $f$ and a smooth morphism $a$ the functors $a^*$ and $f_*$ preserve $\AA^1$-local objects as they have left adjoints which preserve representables. When $f$ is finite they are also both exact. Consequently, we have shown that every object of $DM(S)$ has a weak structure of smooth traces.

Finally, we observe that $\HZ_S \in \SH(S)$ is by definition the image of the object represented by $S$ in $DM(S)$ and that for smooth morphisms $a: X \to S$ we have $a^*\HZ_S \cong \HZ_X$. It now follows from Lemma~\ref{lemm:localisationOfTraces} that $\HZ_S$ (and indeed, any object in the image of $DM(S) \to \SH(S)$) has a weak structure of smooth traces.
\end{proof}




We now turn our attention to algebraic $K$-theory.  See \cite{Wei89} for homotopy invariant algebraic $K$-theory and \cite{Cis13} for its representability in the Morel-Voevodsky stable homotopy category. We recall one construction of the object $\KH$ in $\SH(S)$ that represents homotopy invariant algebraic $K$-theory. For a category $\C$ we denote by $Sp_{S^1}(\C)$ the category of presheaves of $S^1$-spectra on $\C$. When $\C = Sm(X)$ for some scheme $X \in Sch(S)$ we denote by $Sp_{\PP^1}(Sp_{S^1}(Sm(X)))$ the category of $\PP^1$-spectra in $Sp_{S^1}(Sm(X))$ where $\PP^1$ is pointed at infinity. By definition, $\SH(S)$ is the homotopy category of $Sp_{\PP^1}Sp_{S^1}(Sm(X))$, where $Sp_{S^1}(Sm(X))$ is given the model category structure that is the Bousfield localisation with respect to $\AA^1$ invariance and Nisnevich descent. The notation $f^*, f_*$ will be overused, sometimes referring to inverse image and direct image of $\OO_X$-modules, and sometimes referring to inverse image and direct image of presheaves of $S^1$-spectra, or $\PP^1$-spectra. It should be clear from the context which is intended.

We will end up discussing four different incarnations of $K$-theory: a presheaf of $S^1$-spectra on $Sch(S)$, a section of the 2-functor $Sp_{S^1}(Sm(-))$, a section of the 2-functor $Sp_{\PP^1}Sp_{S^1}(Sm(-))$, and a section of the 2-functor $SH(-)$.

\begin{enumerate}
 \item \emph{$K$, a presheaf of $S^1$-spectra on $Sch(S)$.} Following \cite[3.1]{TT90} (cf. \cite[Definition 1.5.3]{TT90} and \cite[Lemma 3.5]{TT90} as well) we denote by $K(X)$ the $S^1$-spectra associated to the biWaldhausen category of perfect complexes on the scheme $X \in Sch(S)$. In order to end up with an actual presheaf of $S^1$-spectra (instead of just a a lax functor), when we say perfect complex, we mean a presheaf on $Sch(X)$ (as opposed to the small Zariski site of $X$) with the appropriate structure and properties (see \cite[Section C4]{FS02}).

 \item \emph{$K_-$, a section of the 2-functor $Sp_{S^1}(Sm(-))$.} For a scheme $X \in Sch(S)$ we define $K_X = K|_{Sm(X)}$ as the restriction of $K$ to smooth schemes over $X$ for $X \in Sch(S)$. For a morphism of $S$-schemes $f: Y \to X$, we have a corresponding adjunction of presheaves of $S^1$-spectra
\[ f^*: Sp_{S^1}(Sm(X)) \rightleftarrows Sp_{S^1}(Sm(Y)): f_* \]
with the right adjoint given by $f_*E(-) = E(Y \times_X -)$. Hence, there is a canonical morphism $K_X \to f_*K_Y$. The $K_X$ together with these canonical morphisms give a section $K_-$ of the 2-functor $Sp_{S^1}(Sm(X))$.

 \item \emph{$\underline{K}_X$, a section of the 2-functor $Sp_{\PP^1}Sp_{S^1}(Sm(-))$.} For each scheme $X \in Sch(S)$ define $F(X) = \mathrm{hofib}(K(\PP^1_X) \stackrel{K(\infty)}{\to} K(X))$ where $\infty: X \to \PP^1_X$ is the closed embedding at infinity. These $F(X)$ form a presheaf of $S^1$-spectra. As with $K$, define $F_X = F|_{Sm(X)}$ as the restriction of $F$ to smooth schemes over $X$ for $X \in Sch(S)$.

On $\PP^1_\ZZ$ choose a global section of $\OO(1)$ whose fibre at infinity is invertible. There is a corresponding morphism $\OO \to \OO(1)$ which can be regarded as a perfect complex concentrated in (cohomological) degrees $0$ and $1$. Its pullback to $Spec(\ZZ)$ along $\infty$ is acyclic. We will denote this complex by $u$. Inverse image gives us a corresponding complex on $\PP^1_X$ for every scheme $X$ which we will denote by $u_X$. Let $p: \PP^1_X \to X$ be the canonical projection. We consider the map $u_X \otimes p^*-: Perf(X) \to Perf(\PP^1_X)$. Notice that as we are using big vector bundles (\cite[Section C4]{FS02}) this is natural in $X$. Notice also that this is exact as $p$ is flat and $u_X$ is a complex of vector bundles. Denote the corresponding map of $K$-theory spectra by $b: K(X) \to K(\PP^1_X)$, also natural in $X$. The composition $\infty^* (u_X \otimes p^*)$ is $(\infty^*u_X) \otimes -$, tensor with an acyclic complex of vector bundles. Hence, $b$ gives rise to a map $\beta: K(X) \to F(X) = \mathrm{hofib}(K(\PP^1_X) \stackrel{K(\infty)}{\to} K(X))$, natural in $X$. That is, we have a map of presheaves of $S^1$-spectra
\[ \beta: K \to F. \]

It follows from our definitions and the fact that $p: \PP^1_X \to X$ is smooth that there is a canonical isomorphism $F_X \cong \ihom(\PP^1_X, K_X)$ in $Ho(Sp_{S^1}(Sm(X)))$ where $\PP^1_X$ is pointed at infinity. Via this canonical morphism, the morphisms $\beta$ give rise to a $\PP^1$-spectrum $(K_X, K_X, K_X, \dots)$ in $Sp_{S^1}(Sm(X))$ which we call $\underline{K}_X$.

 \item \emph{$\KH$, the object representing homotopy invariant $K$-theory in $\SH(X)$.} Finally, the localisation $\KH_X$ of each $\underline{K}_X$ in $Sp_{\PP^1}Sp_{S^1}(Sm(X))$ with respect to Nisnevich descent and $\AA^1$-homotopy (that is, a fibrant replacement for the localised model category structure) gives the object in $\SH(X)$ representing homotopy algebraic $K$-theory (see \cite{Cis13}).
\end{enumerate}

\begin{prop} \label{prop:ktheoryTraces}
The object $\KH \in \mathsf{SH}(S)$ that represents homotopy invariant algebraic $K$-theory has a structure of traces (Definition~\ref{defi:sectionWithTraces}).
\end{prop}


\begin{proof}
Each of the four incarnations of algebraic $K$-theory mentioned above will have traces in their own sense, and each one induces the traces on the next. For the ``trace'' morphisms that we will associate with $K$ (resp. $K_-$, $\underline{K}_-, \KH$) we will use $\Tr_f^K$ (resp. $\Tr^{S^1}_f, \Tr^{\PP^1}_f, \Tr_f$).

We begin with traces on $K$ and the properties we need. The construction of $K$ is functorial in complicial biWaldhausen categories. Notably, for each finite flat surjective morphism $f: Y \to X$ we obtain a corresponding exact functor $f_*: Perf(Y) \to Perf(X)$ between the corresponding biWaldhausen categories of perfect complexes. Hence, there are  morphisms $\Tr_f^K: K(Y) \to K(X)$. Due to the functoriality and the standard properties of $\OO_X$-modules we have the following properties. For a morphism $f: Y \to X$, we denote by $K(f): K(X) \to K(Y)$ the morphism of spectra induced by inverse image $f^*: Perf(X) \to Perf(Y)$ (discussed in \cite[3.14]{TT90}).

\emph{Functoriality.} (cf. \cite[1.5.4]{TT90}) If $W \stackrel{g}{\to} Y \stackrel{f}{\to} X$ are finite flat surjective then we have a homotopy $\Tr_f^K\Tr_g^K \cong \Tr_{fg}^K$.

\emph{Base-change.} (cf. \cite[3.18]{TT90}) If we have a cartesian square (\ref{equa:cartSquare}), then there is a homotopy $K(p)\Tr_f^K \cong \Tr_g^KK(q)$.

\emph{Degree.} (cf. \cite[1.7.3.2]{TT90}) If $f: Y \to X$ is finite flat surjective and there is an isomorphism $f_*\OO_Y \cong \OO_X^d$ then there is a homotopy of maps of $S^1$-spectra $\Tr_f^K K(f) \cong d \cdot K(id_X)$.

Suppose that $f: Y \to X$ is a finite flat surjective morphism. To ease the notation we use $f$ for the induced morphism $\PP^1_Y \to \PP^1_X$ as well, and $p$ for both projections $\PP^1_X \to X$ and $\PP^1_Y \to Y$. Recall that above we have defined $F(X) = \mathrm{hofib}(K(\PP^1_X) \stackrel{K(\infty)}{\to} K(X))$.

After the base-change property, the morphisms $\Tr^K$ induce a morphism $\Tr^F_f: F(Y) \to F(X)$. We claim that the following square on the left is commutative up to homotopy. This follows from the commutativity of the square in the middle up to homotopy, which is a consequence of the commutativity of the square on the right up to natural isomorphism. This latter commutativity is a consequence of the projection formula $f_*(u_Y \otimes p^*-) = f_*(f^*u_X \otimes p^*-) \cong u_X \otimes f_*p^*(-)$, and base change $p^*f_* \cong f_*p^*$.
\[ \xymatrix@C=15pt{
K(Y) \ar[d]_{\Tr^K_f} \ar[r]^{\beta} & F(Y) \ar[d]^{\Tr^F_f} && K(Y) \ar[d]_{\Tr^K_f} \ar[r]^b & K(\PP^1_Y) \ar[d]^{\Tr^K_f}  && Perf(Y) \ar[r]^{u \otimes p^*-} \ar[d]_{f_*} & Perf(\PP^1_Y) \ar[d]^{f_*} \\
K(X) \ar[r]_{\beta} & F(X) && K(X) \ar[r]_b & K(\PP^1_X) && Perf(X) \ar[r]_{u \otimes p^*-} & Perf(\PP^1_X)
} \]
That is, $\beta: K \to F$ (defined above) is a ``morphism of presheaves of $S^1$-spectra with traces'' (although we haven't formally defined what that means).

Now we pass to the sections $K_-$ of the 2-functor $Ho(Sp_{S^1}(Sm(-)))$. We will use the above properties to show that $K_-$ has a structure of traces as a section. Recall that for $f: Y \to X$ a morphism of $S$-schemes, we have a corresponding adjunction of presheaves of $S^1$-spectra
\[ f^*: Sp_{S^1}(Sm(X)) \rightleftarrows Sp_{S^1}(Sm(Y)): f_*. \]
Due to the base-change mentioned above, if $f$ is finite flat surjective we have a morphism of presheaves of $S^1$-spectra $\Tr_f^{S^1}: f_*K_Y \to K_X$ induced by the morphisms $\Tr^K$.

\emph{Functoriality, Base-change, and Degree.} These follow immediately from the corresponding properties of the $\Tr^K$ and the description of $f_*$ as $f_*E(-) = E(Y \times_X -)$.

\emph{Periodicity.} Via the canonical isomorphism $F_X \cong \ihom(\PP^1_X, K_X)$ in $Ho(Sp_{S^1}(Sm(X)))$, the traces $\Tr^F_f$ that we have defined in $F_-$ correspond to the morphisms $\ihom(\PP^1_X, f_*K_Y) \to \ihom(\PP^1_X, K_X)$ induced by the traces $\Tr^K_f$ of $K_-$. Hence commutative diagrams
\[ \xymatrix{
f_*K_Y \ar[r]^-{f_*\beta'} \ar[d]_{\Tr^K_f} & f_*\ihom(\PP^1_Y, K_Y) \ar@{=}[r] & \ihom(\PP^1_X, f_*K_Y) \ar[d]^{\ihom(\PP^1_X, \Tr^K_f)} \\
K_X \ar[rr]_-{\beta'} && \ihom(\PP^1_X, K_X)
} \]
where the $\beta': K_X \to \ihom(\PP^1_X, K_X)$ are the morphisms corresponding to the $\beta: K_X \to F_X$.

Recall that for any morphism $f: Y \to X$ of schemes there is an adjunction
\[ f^*: Sp_{\PP^1}Sp_{S^1}(Sm(X)) \rightleftarrows Sp_{\PP^1}Sp_{S^1}(Sm(Y)): f_* \]
with the right adjoint given by $f_*(E_0, E_1, \dots) = (f_*E_0, f_*E_1, \dots)$ and the new structural morphisms are the compositions $f_*E_n \to f_*\ihom(\PP^1_Y, E_{n + 1}) \cong \ihom(\PP^1_X, f_*E_{n + 1})$. It follows from our remarks on periodicity that when $f$ is finite flat surjective we have induced trace morphisms $\Tr^{\PP^1}_f: f_*\underline{K}_Y \to \underline{K}_X$.

\emph{Functoriality, Base-change, Degree.} These follow immediately from the corresponding properties of the $\Tr^{S^1}$ and the description of $f_*$ that we have given. Hence, the section $\underline{K}_-$ has a structure of traces.

Now we consider $\KH_X$. The category $\SH(X)$ can be presented as the localisation of the homotopy category $Ho(Sp_{\PP^1}(Sm(X)))$ with respect to $\AA^1$-localisation and Nisnevich descent. Inverse image preserves Nisnevich hypercovers, and the projections $\AA^1_U \to U$ so the class of morphisms that we are localising with respect to is preserved. Consequently, direct image preserves local objects. That is, the localisation functors $Ho(Sp_{\PP^1}(Sm(X))) \to SH(X)$ satisfy the properties required to apply Lemma~\ref{lemm:localisationOfTraces} to the section $\underline{K}_-$ of the 2-functor $Ho(Sp_{\PP^1}(Sm(-)))$.
\end{proof}

\begin{coro} \label{coro:tracesOnHZmod}
Suppose $k$ is a perfect field of exponential characteristic $p$. Then for any object $M$ of $\SH(k)$, the object $\HZpi_k\wedge M$ has a structure of traces.
\end{coro}


\begin{proof}
We have seen that $\KH$ has a structure of traces (Proposition~\ref{prop:ktheoryTraces}) and after work of Levine we know that the zero slice of $\KH$ is $\HZ$ (\cite[Theorems 6.4.2 and 9.0.3]{Lev08}). Hence, $\HZ$ has a structure of traces (Proposition~\ref{prop:tracesOnSlices}). So applying Proposition~\ref{prop:tracesOnProducts} shows that $\HZpi_k\wedge M$ has a structure of traces.
\end{proof}






\section{Resolution of singularities for relative cycles} \label{sec:RoS}

In this section our goal is to prove the following theorem.

\begin{theo} \label{theo:MainResultRos}
Let $k$ be a perfect field of exponential characteristic $p$. Suppose that $F$ is a presheaf with transfers on $Sch(k)$ such that $F_{cdh} \otimes \zpi = 0$. Then $\underline{C}_*(F|_{Sm(k)})_{Nis} \otimes \zpi$ is quasi-isomorphic to zero.
\end{theo}

Recall that the 2-functor $X \mapsto \SH(X)$ factors through a 2-functor $X \mapsto \M(X)$ where
\begin{enumerate}
 \item for each $X \in Sch(S)$ the category $\M(X)$ is a stable model category (hence enriched in symmetric $S^1$-spectra \cite{Dug06}) which is combinatorial\footnote{Combinatorial categories were introduced by Jeff Smith. The definition can be found in \cite[Section 2]{Dug01}.} and cellular,
 \item for each $f: Y \to X$ in $Sch(S)$ the functor $f^*: \M(X) \to \M(Y)$ is a left Quillen functor,
 \item for each smooth $f: Y \to X$ in $Sch(S)$ the functor $f^*$ is has a left adjoint $f_\#$ which is a left Quillen functor.
 \item for each cartesian square (\ref{equa:cartSquare}) with $f$ smooth, the corresponding natural transformations $g_\#q^*\to p^*f_\#$ are isomorphisms. 
\end{enumerate}
In other words, $\M$ is a stable $Sm$-fibred combinatorial model category (\cite[Definitions 1.1.2, 1.1.9, 1.3.2, 1.3.20]{CD}). Moreover, $\SH$ is obtained by passing to the homotopy categories of $\M$. That is, $\SH$ is the homotopy $Sm$-fibred category associated with $\M$ (\cite[1.3.23]{CD}). These statements follow directly from the construction of $\SH$ given in \cite{Ayo07}.

As each $\M(S)$ is enriched in symmetric $S^1$-spectra \cite{Dug06}, for any pair of objects $\E, \F \in \M(S)$ we can associate a presheaf of $S^1$-spectra that sends a scheme $a: X \to S$ to the $S^1$-spectrum
\[ \underline{(\F, \E)}(X) \stackrel{def}{=} \ihom(\F, a_*a^*\E). \]

In the following theorem, ``descent'' is in the sense of \cite[Definition 3.2.5]{CD}. In the case where $\M$ is the stable $Sm$-fibred model category that associates to each scheme $X \in Sch(S)$ the corresponding category of presheaves of $S^1$-spectra $Sp_{S^1}(Sm(-))$, this definition of agrees with that of Jardine-Thomason (see \cite{Mit97} for a civilised discussion of this notion of descent).

\begin{theo}[{\cite[Corollary 3.2.18]{CD}}] \label{theo:descent}
Suppose that $\M$ is a stable $Sm$-fibred combinatorial model category over $Sch(S)$ and $\E \in \M(S)$. Let $\tau$ be a Grothendieck topology and $\mathcal{G}$ a set of generators for $Ho(\M(S))$. Then $\E$ satisfies $\tau$-descent if and only if for every $\F \in \mathcal{G}$ the presheaf of $S^1$-spectra $\underline{(\F, \E)}$ satisfies $\tau$-descent.
\end{theo}

\begin{rema}
The statement in \cite{CD} is for all $\F$, not just a set of generators, but a glance at the proof of \cite[Corollary 3.2.17]{CD} shows that it suffices to consider generators.
\end{rema}

Now for any $S$ the triangulated category $\SH(S)$ is compactly generated by objects of the form $\Sigma^{-q} f_\#f^*\un_S$ for $f: Y \to S$ a smooth morphism and $q > 0$. If $\F$ in Theorem~\ref{theo:descent} is of this form and $a: X \to S$ is also smooth, then we have the following canonical isomorphisms
\[ \begin{split}
\pi_n\underline{(\F, \E)}(X) &\cong \hom(\un_S[n], \quad \ihom(f_\#f^*(\Sigma^{-q} \un_S), a_*a^*\E)) \\
&\cong \hom(f_\#f^*(\Sigma^{-q} \un_S[n]), \quad  a_*a^*\E) \\
&\cong \hom(a_\#a^*f_\#f^*(\Sigma^{-q} \un_S[n]), \quad  \E) \\
&\cong \hom(\Sigma^{-q} \Sigma^\infty X \times_S Y_+ [n]), \quad  \E)
\end{split} \]
This group is denoted by $\E^{2q-n,q}(X \times_S Y)$ in \cite[Section 6]{Voev98} and $\pi_{n - q}(\E)(X \times_S Y)_{q}$ in \cite{Mor04}.

\begin{defi}
We introduce the notation
\[ \uE^{q, Y}(X) \stackrel{def}{=} \ihom(\Sigma^{-q} f_\#f^*\un_S, a_*a^*\E). \]
\end{defi}

The following corollary is a summary of what we have just discussed.

\begin{coro} \label{coro:descentSummary}
Let $S$ be a noetherian scheme and suppose that $\tau$ is a Grothendieck topology in $Sch(S)$. Then an object $\E \in \SH(S)$ satisfies $\tau$-descent if and only if for every $q > 0$ and every smooth $S$-scheme $Y \to S$ the presheaf of $S^1$-spectra $\uE^{q, Y}$ satisfies $\tau$-descent.
\end{coro}

Due to the isomorphisms mentioned above, after work of D{\'e}glise, if $\E$ is oriented then the Nisnevich sheaf associated to the presheaf $\pi_n\uE^{q, Y}$ on $Sm(S)$ has a structure of Nisnevich sheaf with transfers (\cite{Deg11}).

\begin{prop}[{\cite{Deg11}}] \label{prop:orientHomMod}
Let $k$ be a perfect field and $\E \in \SH(k)$ an oriented object (in the sense of Morel \cite[Definition 2.1]{Vez01}). Then for any $n, q \in \ZZ$ and smooth $Y \to k$, the Nisnevich sheaf $(\pi_n\uE^{q, Y})_{Nis}$ associated to the presheaf of homotopy groups $\pi_n\uE^{q, Y}$ has a structure of transfers on $Sm(k)$.
\end{prop}


\begin{proof}
The presheaves $\pi_n \uE^{q, Y}(-)$ and $\pi_n \uE^{q, S}(Y \times_S -)$ are canonically isomorphic. The functor $Y \times_k - : Sm(k) \to Sm(k)$ lifts to a functor $Y \times_k -: SmCor(k) \to SmCor(k)$ compatible with the inclusion $Sm(k) \to SmCor(k)$ and so it suffices to show that the Nisnevich sheaf associated to the presheaf $\pi_n \uE^{q, Spec(k)}( -)$ has transfers on $Sm(k)$.

First we claim that if $\E$ is orientable, then it is weakly orientable (\cite[Definition 4.2.3]{Deg11}). Recall that the Hopf map $\eta: \Sigma^\infty (\GG_m, 1) \to \un$ is defined as the map in $\SH(k)$ induced by $\AA^2 - \{0\} \to \PP^1$ after applying $\Sigma^{-1}$ \cite[1.2.6]{Deg11}. As $\E$ is oriented, the projective bundle formula holds (\cite[Proposition 2.4(ii)]{Vez01} attributes this to Morel) and so for any smooth scheme $W$ and any linear embedding $\PP^1 \to \PP^2$ the induced morphism 
\[ \hom_{\SH(k)}(\Sigma^i\Sigma^\infty (\PP^2 \times W)_+, \E) \to \hom_{\SH(k)}(\Sigma^i\Sigma^\infty (\PP^1 \times W)_+, \E) \]
is split surjective for all $i \in \ZZ$. After the homotopy exact sequence \cite[6.2.1]{Mor04}\[ \Sigma^\infty (\AA^1 - \{0\})_+ \stackrel{\Sigma \eta}{\to} \Sigma^\infty(\PP^1)_+ \to \Sigma^\infty(\PP^2)_+ \to \Sigma^\infty(\AA^1 - \{0\})_+[1] \]
this implies that the morphism $\hom_{\SH(k)}((\Sigma \eta) \wedge \Sigma^i\Sigma^\infty W_+, \E)$ is zero for all $i \in \ZZ$ and smooth $W$. Equivalently, the morphism
\[ \hom_{\SH(k)}(\Sigma^{i + 1} \Sigma^\infty W_+, \ihom(\eta, \E)) \] is zero for all $i \in \ZZ$ and smooth $W$ where $\ihom$ is the internal hom in the monoidal category $\SH(k)$. As the $\Sigma^{i + 1} \Sigma^\infty W_+$ form a compact generating family for the triangulated category $\SH(k)$ this implies that the morphism $\ihom(\eta, \E)$ is zero. That is, $\E$ is weakly orientable \cite[Definition 4.2.3, Lemma 4.2.2(ii')]{Deg11}. 

An equivalent condition for $\E$ to be weakly orientable is that the associated homotopy modules $\underline{\pi}_m(\E)_*$ (\cite[1.1.2, Definition 1.2.2, 1.2.3]{Deg11}) are orientable (\cite[Definition 1.2.7]{Deg11}) for each $m$ (\cite[Lemma 4.2.2(i)]{Deg11}). One of the main results of \cite{Deg11} is that orientable homotopy modules are precisely those homotopy modules which admit a structure of transfers on $Sm(k)$ (\cite[Corollary 4.1.5(2)(i), Corollary 4.1.5(ii)]{Deg11}). By definition, the Nisnevich sheaves $\underline{\pi}_m(\E)_i$ on $Sm(k)$ are the sheaves associated to the presheaves $\hom_{\SH(k)}(\Sigma^{-i}(-)_+[i + m], \E)$ on $Sm(k)$. That is, the presheaves $\pi_{i + m}\uE^{i, Spec(k)}$. Hence, the Nisnevich sheaf associated to the presheaf $\pi_n \uE^{q, Spec(k)}( -)$ has transfers on $Sm(k)$ for each $n, q \in \ZZ$.
\end{proof}

We can deduce now the following theorem.


\begin{theo} \label{theo:ldescentInSH}
Suppose $k$ is a perfect field and $\ell$ a prime that is invertible in $k$. Let $\E$ be an oriented $\zll$-local object (Definition~\ref{defi:zlllocal}) of $\SH(k)$ with a structure of traces (Definition~\ref{defi:sectionWithTraces}). Then any direct factor of $\E$ satisfies $\ldh$-descent. In particular, any module in $\SH(k)$ over the ring spectrum $\HZl$ satisfies $\ldh$-descent.
\end{theo}

\begin{proof}
It is immediate from the definition (\cite[Definition 3.2.5]{CD}) that any direct factor of an object satisfying descent also satisfies descent. Let $\E \in \SH(k)$ be a $\zll$-local object with a structure of traces. After Corollary~\ref{coro:descentSummary} it suffices to show that $\uE^{q, Y}$ satisfies $\ldh$-descent for every $q > 0$ and every smooth $Y \to k$. Let $\uE^{q, Y} \to \uE' = \HH_{\ldh}(-, \uE^{q, Y})$ be the Godement-Thomason construction \cite[1.33]{Tho85}. The morphism of associated $\ldh$ sheaves $(\pi_n\uE^{q, Y})_{\ldh} \to (\pi_n\uE')_{\ldh}$ is an isomorphism for all $n$. If we can show that $\pi_n\uE^{q, Y} \to \pi_n\uE'$ is an isomorphism of presheaves for every $n$ then $\uE^{q, Y} \to \uE'$ is a weak equivalence of presheaves of $S^1$-spectra, so $\uE^{q, Y}$ satisfies $\ldh$ descent, and we are done.

We know that $\E$ and hence $\uE^{q, Y}$ (Theorem~\ref{theo:descent}) satisfies cdh descent, since every object of $\SH(k)$ satisfies cdh descent (\cite[3.7]{Cis13}). This gives us a cdh descent spectral sequence for $\uE^{q, Y}$
\[ E_2^{s,t} = H_{cdh}^s(X, (\pi_{-t}\uE^{q, Y})_{cdh}) \implies \pi_{-s-t}\uE^{q, Y}(X) \]
together with a morphism towards the $\ldh$ descent spectral sequence for $\uE'$
\[ E_2^{s, t} = H_{\ldh}^s(X, (\pi_{-t}\uE')_{\ldh}) \implies \pi_{-s-t}\uE'(X). \]
The first spectral sequence converges since the cdh topology has finite cohomological dimension \cite[Theorem 12.5]{SV00}. As we know that $(\pi_n\uE)_{\ldh} \to (\pi_n\uE')_{\ldh}$ is an isomorphism for all $n$, it suffices to show that
\begin{equation} \label{equa:uEcdhlprime}
H_{cdh}^s(X, (\pi_{-t}\uE^{q, Y})_{cdh}) \to H_{\ldh}^s(X, (\pi_{-t}\uE^{q, Y})_{\ldh})
\end{equation}
is an isomorphism for all $s, t$. This will imply that the morphism of spectral sequences is an isomorphism, and therefore give the convergence of the second spectral sequence, and an isomorphism $\pi_{n}\uE^{q, Y}(X) \stackrel{\sim}{\to} \pi_{n}\uE'(X)$ for all $n$.

That the morphism (\ref{equa:uEcdhlprime}) is an isomorphism will follow from Theorem~\ref{theo:cdhlprimeComparison}. To apply this theorem, we must show that $\pi_{n}\underline{\E}^{q, Y}$ is a homotopy invariant presheaf of $\zll$-modules with traces such that $\pi_{n}\underline{\E}^{q, Y}|_{Sm(k)}$ has a structure of presheaf with transfers, and $\pi_n\uE^{q, Y}(U) \to \pi_{n}\underline{\E}^{q, Y}(U_{red})$ is an isomorphism for all $U \in Sch(k)$. Recall that for $a: U \to k$ in $Sch(k)$ we have a canonical isomorphism
\[ \pi_{n}\underline{\E}^{q, Y}(U) \cong \hom_{\SH(k)}(\Sigma^{-q} f_\#f^*\un_S[n], a_*a^*\E). \]
This presheaf is homotopy invariant and doesn't see nilpotents because the same is true of the functor $Sch(k) \to End(\SH(k))$ defined by $(a:U \to k) \mapsto a_*a^*$. It is a presheaf of $\zll$-modules by our hypothesis that $\E$ is $\zll$-local. It has traces as a result of $\E$ having a structure of traces (Lemma~\ref{lemm:tracesImpliesTraces}). Finally, the hypothesis that $\E$ is oriented implies that $\pi_n\uE^{q, Y}|_{Sm(k)}$ is a presheaf with transfers after a theorem of D{\'e}glise (Proposition~\ref{prop:orientHomMod}).

For the last assertion in the statement, it suffices to notice that any $\HZl$ module $M$ is canonically a direct factor of $\HZl \otimes M$, and the latter has a structure of traces after Corollary~\ref{coro:tracesOnHZmod}.
\end{proof}

\begin{coro} \label{coro:lHypercoversInDM}
Let $k$ be a perfect field and $\ell$ a prime that is invertible in $k$. Suppose that $\X \to X$ is a smooth $\ldh$-hypercover in $Sm(k)$ of a smooth scheme $a: X \to k$ in $Sm(k)$. Then the corresponding morphism $M(\X) \to M(X)$ is an isomorphism in $DM\eff(k, \zll)$.
\end{coro}

\begin{proof}
By Voevodsky's Cancellation Theorem \cite{Voe10} it suffices to show that this morphism is an isomorphism in $DM(k, \zll)$. There is a canonical adjunction
\[ \ZZ : \SH(k) \rightleftarrows DM(k, \zll) : U \]
with $\ZZ$ symmetric monoidal such that for any simplicial smooth scheme ${q: \Y \to k}$ we have $\ZZ(q_\#q^*\un) = M(\Y)$. Moreover, for any object $E \in DM(k, \zll)$, the spectrum $U(E)$ has a structure of $\HZl$-module. Let $p: \X \to k$ and $a : X \to k$ be the structural morphisms. Due to the last assertion of Theorem~\ref{theo:ldescentInSH}, it is sufficient to prove that, for any object $\E$ of $\SH(k)$ which satisfies $\ldh$ descent, the map
\[ \hom(a_\# a^* \un, \E) \to \hom(p_\# p^* \un, \E) \]
is bijective. By adjunction it is equivalent to show that the map 
\[ \hom(\un, a_*a^*\E) \to \hom(\un, p_*p^*\E) \]
is a bijection, and this holds since $\E$ satisfies $\ldh$ descent.
\end{proof}

We are now in a position to prove Theorem~\ref{theo:MainResultRos}.

\begin{proof}[Proof of Theorem~\ref{theo:MainResultRos}.]
It is enough to show that $\underline{C}_*(F)_{Nis} \otimes \zll$ quasi-isomorphic to zero for each prime $\ell \neq p$. Note that our assumptions imply that ${F_{\ldh} \otimes \zll = 0}$ for each $\ell \neq p$. 

Corollary~\ref{coro:lHypercoversInDM} is precisely the condition \cite[Proposition 5.2.10]{CD}(i), and \cite[Proposition 5.2.10]{CD}(ii$'$) applied to $F$ is the condition that ${\underline{C}_*(F)_{Nis} \otimes \zll = 0}$ since this is the image of $F$ in $DM\eff(k, \zll)$ under the canonical morphism
\[ D(PreShv(SmCor(k), \zll)) \to DM\eff(k, \zll). \]
Hence, after \cite[Proposition 5.2.10]{CD}, the former implies the latter.

We can be a bit more verbose. For $\tau = Nis, \ldh$ we have canonical equivalences
\[ D(PreShv(SmCor(k), \zll)) / \mathscr{L}_\tau \cong D(Shv_\tau(SmCor(k), \zll)) \]
where $\mathscr{L}_\tau$ is the class of cones of morphisms of the form $L(\X) \to L(X)$ with $\X \to X$ a $\tau$-hypercovering. In the light of these equivalences, Corollary~\ref{coro:lHypercoversInDM} implies that we have a commutative triangle
\begin{equation} \label{equa:dmNislTri}
\xymatrix{
D(Shv_{Nis}(SmCor(k), \zll)) \ar[rr] \ar[dr] && D(Shv_{\ldh}(SmCor(k), \zll)) \ar[dl] \\
& DM\eff(k, \zll)
}
\end{equation}
As before, $\underline{C}_*(F)_{Nis} \otimes \zll$ is the image of $F \otimes \zll$ in the lower category, and $F_{\ldh} \otimes \zll$ is its image in the upper right category. It follows that if $F_{\ldh} \otimes \zll$ is zero, then $\underline{C}_*(F)_{Nis} \otimes \zll$ is zero.
\end{proof}

We have the following easy consequence of our main theorem.

\begin{coro} \label{coro:a1localisations}
Let $k$ be a perfect field of exponential characteristic $p$ and $\ell$ a prime difference from $p$. Then there are canonical functors
\[ D(Shv_{\ldh}(SmCor(k), \zll)) \to DM\eff(k, \zll) \]
\[ D(Shv_{cdh}(Cor(k), \zpi)) \to DM\eff(k, \zpi) \]
which identify the targets as the localisations of the sources with respect to morphisms of the form $L_\tau(\AA^1_X) \to L_\tau(X)$ (where $\tau = cdh$ or $\ldh$ as applicable) and $X \in Sm(k)$.
\end{coro}

\begin{proof}
The existance of the first functor has already been seen in the proof of  Theorem~\ref{theo:MainResultRos} and considering the categories in question as localisations of the category $D(PreShv(SmCor(k), \zll))$ as discussed in that proof leads to the universal property which identifies $DM\eff(k, \zll)$ as the appropriate localisation. For the second functor consider the following commutative diagram of functors:
\[ \xymatrix{
D(Shv_{Nis}(Cor(k), \zll)) \ar[r] \ar[d] & D(Shv_{Nis}(SmCor(k), \zll)) \ar[dd] \\
D(Shv_{cdh}(Cor(k), \zll)) \ar[d]_{(1)} & \\
D(Shv_{\ldh}(Cor(k), \zll)) \ar[r]_{(2)} & D(Shv_{\ldh}(SmCor(k), \zll))
} \]
As a consequence of the theorem of Gabber giving smooth $\ldh$ covers (Corollary~\ref{coro:regularlCover}) the functor (2) is an equivalence. Since cdh sheaves with transfers are already $\ldh$ sheaves (Corollary~\ref{coro:equivalencecdhltransfers}) the functor (1) is an equality. Hence the desired functor, at least with $\zll$ coefficients exists by the existence of the first functor in the statement, and moreover, it is identified with the localisation with respect to morphisms of the form $L_{cdh}(\AA^1_X) \to L_{cdh}(X)$. Now we have commutative squares
\[ \xymatrix{
D(Shv_{cdh}(Cor(k), \zpi)) \ar[r] \ar[d] & D(Shv_{cdh}(Cor(k), \zll)) \ar[d] \\
D(Shv_{cdh}(Cor(k), \zpi)) / \langle L_{cdh}(\AA^1_X) \to L_{cdh}(X) \rangle \ar[r] & DM\eff(k, \zpi)
} \]
and the result follows from Section~\ref{sec:invertTriangle}.
\end{proof}

\begin{defi}
For $X \in Sch(k)$ we will denote by $\mp{X}$ (resp. $\mcp{X}$) the image of $c_{equi}(X /k, 0)$ (resp. $z_{equi}(X / k, 0)$) in $DM\eff(k, \zpi)$ under the functor
\[ D(Shv_{cdh}(Cor(k), \zpi)) \to DM\eff(k, \zpi). \]
\end{defi}

\begin{prop}[{cf. \cite[Theorem 4.1.10]{Voev00}}] \label{prop:V4110}
Let $k$ be a perfect field of exponential characteristic $p$. Suppose that $X$ is a scheme of finite type and $F$ a presheaf with transfers on $Sch(k)$. Then there is a canonical isomorphism
\[ \hom_{DM\eff(k, \zpi)}(\mp{X}, \underline{C}_*(F|_{Sm(k)})\ozpi) \cong \HH^i_{cdh}(X, \underline{C}_*(F)_{cdh})\ozpi. \]
\end{prop}

\begin{proof}
Use the description of $DM\eff(k, \zpi)$ as a localisation of $D(Shv_{cdh}(Cor(k), \zpi))$ together with the analogue of \cite[Proposition 3.1.9]{Voev00}.
\end{proof}

\section{Bivariant cycle cohomology --- After Friedlander, Voevodsky} \label{sec:bivCycCoh}

In this section we collect some results of \cite{FV} for which Theorem~\ref{theo:MainResultRos} allows us to remove the resolution of singularities assumption.

\subsection{Bivariant cycle cohomology}

Recall the following definition.

\begin{defi}[{\cite[Definition 4.3]{FV}}]
Let $X, Y$ be schemes of finite type over a field $k$ and $r \geq 0$ be an integer. The bivariant cycle cohomology groups of $Y$ with coefficients in cycles on $X$ are the groups
\[ A_{r, i}(Y, X) = \HH^{-i}_{cdh}(Y, (\underline{C}_*(z_{equi}(X, r))_{cdh}) \]
The notation $A_{r,i}(X)$ is also used for the groups $A_{r, i}(Spec(k), X)$.
\end{defi}

\begin{theo}[{cf. \cite[Theorem 5.5]{FV}}] \label{theo:FV55}
Let $k$ be a perfect field of exponential characteristic $p$, let $\ell$ be a prime different from $p$, and suppose $F$ is a presheaf with transfers on $Sch(k)$. 
\begin{enumerate}
 \item For any smooth scheme $U$ and $n \geq 0$ there are canonical isomorphisms
\[ \HH^n_{cdh}(U, \underline{C}_*(F)_{cdh})\ozpi \cong \HH^n_{Zar}(U, \underline{C}_*(F|_{Sm(k)})_{Zar})\ozpi, \]
\[ \HH^n_{cdh}(U, \underline{C}_*(F)_{cdh})\otimes \zll \cong \HH^n_{\ldh}(U, \underline{C}_*(F)_{\ldh})\otimes \zll. \]
 \item For any separated scheme of finite type $X$ over $k$, and any $n \geq 0$ the projection $X \times \AA^1 \to X$ induces isomorphisms
\[ \HH_{cdh}^n(X, \underline{C}_*(F)_{cdh})\ozpi \cong \HH_{cdh}^n(X \times \AA^1, \underline{C}_*(F)_{cdh})\ozpi. \]
\end{enumerate}
\end{theo}

\begin{proof}
Due to the fact that we can calculate hypercohomology as hom groups in derived categories of sheaves with transfers, Corollary~\ref{coro:a1localisations} gives us 
\[ \HH^n_{cdh}(U, \underline{C}_*(F)_{cdh})\ozpi = \HH^n_{Nis}(U, \underline{C}_*(F|_{Sm(k)})_{Nis})\ozpi. \]
That the Nisnevich and Zariski hypercohomology are the same follow from the hypercohomology spectral sequence and \cite[Theorem 5.1(2)]{FV}. The second equality also follows from the hypercohomology spectral sequence and Theorem~\ref{theo:lcdhcohomologyAgree}. The third equality also follows from Corollary~\ref{coro:a1localisations} and calculating hypercohomology using $\hom$'s in the derived categories of sheaves with transfers.
\end{proof}

\begin{prop}[{cf. \cite[Proposition 5.9]{FV}}] \label{prop:FV59}
Let $k$ be a perfect field of exponential characteristic $p$ and let $X, Y \in Sch(k)$. Then for all $r, i$ the homomorphisms
\[ A_{r,i}(Y, X)\ozpi \to A_{r,i}(Y \times \AA^1, X)\ozpi \]
induced by the projection are isomorphisms.
\end{prop}

\begin{proof}
This is a special case of Theorem~\ref{theo:FV55} with $F = z_{equi}(X, r)$.
\end{proof}

\begin{theo}[{cf. \cite[Theorem 5.11]{FV}}] \label{theo:FV511}
Let $k$ be a perfect field of exponential characteristic $p$ and let $X \in Sch(k)$. Let $Y \subset X$ be a closed subscheme of $X$, and let $U_1, U_2$ be Zariski open subsets with $X = U_1 \cup U_2$. Then there are canonical exact triangles (in the derived category of complexes of sheaves on $Sm(k)_{Zar}$) of the form
\[ \begin{split}
\underline{C}_*(z_{equi}(Y, r))_{Zar}\ozpi &\to \underline{C}_*(z_{equi}(X, r))_{Zar}\ozpi  \\
&\to \underline{C}_*(z_{equi}(X - Y, r))_{Zar}\ozpi \to \underline{C}_*(z_{equi}(Y, r))_{Zar}\ozpi[1] 
\end{split} \]
and
\[ \begin{split}
\underline{C}_*(z_{equi}(X, r))_{Zar}\ozpi  &\to \underline{C}_*(z_{equi}(U_1, r))_{Zar}\ozpi \oplus \underline{C}_*(z_{equi}(U_2, r))_{Zar}\ozpi \\
&\to \underline{C}_*(z_{equi}(U_1 \cap U_2, r))_{Zar}\ozpi \to \underline{C}_*(z_{equi}(X, r))_{Zar}\ozpi[1]. \end{split} \]
\end{theo}

\begin{proof}
We have a the sequence
\[ 0 \to z_{equi}(Y, r) \to z_{equi}(X, r) \to z_{equi}(X - Y, r) \]
where the right-most morphism becomes surjective after taking the associated cdh sheaves (\cite[Theorem 4.2.9]{SV}, \cite[Theorem 4.3.1]{SV}). Hence, by Theorem~\ref{theo:MainResultRos} after applying $\underline{C}_*(-)\ozpi$ we get a short exact sequence of complexes of Nisnevich sheaves on $Sm(S)$. That this is also a short exact sequence of complexes of Zariski sheaves is \cite[Lemma 4.1]{FV} and \cite[Theorem 5.1]{FV}.

The proof for the second sequence is the same using \cite[Corollary 4.3.2]{SV} instead of \cite[Theorem 4.3.1]{SV}.
\end{proof}

\begin{coro}[{cf. \cite[Corollary 5.12]{FV}}] \label{coro:512}
With the notation and assumptions of Theorem~\ref{theo:FV511} for any scheme $U \in Sch(k)$ there are long exact sequences
\[ \dots A_{r,i}(U, Y)\ozpi \to A_{r,i}(U, X)\ozpi \to A_{r,i}(U, X - Y)\ozpi \to A_{r, i - 1}(U, Y)\ozpi \to \dots \]
and 
\[ \begin{split}
\dots A_{r,i}(U, X)\ozpi &\to A_{r,i}(U, U_1)\ozpi \oplus A_{r,i}(U, U_2)\ozpi \to A_{r, i}(U, U_1 \cap U_2)\ozpi \\
&\to A_{r,i - 1}(U, X)\ozpi \to \dots 
\end{split} \]
\end{coro}

\begin{proof}
Use Theorem~\ref{theo:FV511} Corollary~\ref{coro:a1localisations}, and the fact that we can calculate hypercohomology in the derived category of sheaves with transfers.
\end{proof}

\begin{theo}[{cf. \cite[Theorem 5.13]{FV}}] \label{theo:FV513}
Let $k$ be a perfect field of exponential characteristic $p$ and let $X \in Sch(k)$. Let $Z \subset X$ be a closed immersion and $X' \to X$ a proper morphism in $Sch(k)$ such that $X' \to X$ is an isomorphism outside of $Z$. Let $Z' = Z \times_X X'$. Then there is a canonical exact triangle (in the derived category of complexes of sheaves on $Sm(k)_{Zar}$)
\[ \begin{split}
\underline{C}_*(z_{equi}(Z', r))_{Zar}\ozpi &\to \underline{C}_*(z_{equi}(Z, r))_{Zar}\ozpi \oplus \underline{C}_*(z_{equi}(X', r))_{Zar}\ozpi \\
&\to \underline{C}_*(z_{equi}(X, r))_{Zar}\ozpi \to \underline{C}_*(z_{equi}(Z', r))_{Zar}\ozpi [1].
\end{split} \]
\end{theo}

\begin{proof}
Exactly the same as for Theorem~\ref{theo:FV511} using \cite[Proposition 4.3.3]{SV} instead of \cite[Theorem 4.3.1]{SV}.
\end{proof}

\begin{coro}[{cf. \cite[Corollary 5.14]{FV}}]
With the notation and assumptions of Theorem~\ref{theo:FV513}, for any scheme $U \in Sch(k)$ there is a canonical long exact sequence of the form
\[ \begin{split}
\dots A_{r,i}(U, Z')\ozpi &\to A_{r,i}(U, Z)\ozpi \oplus A_{r,i}(U, X')\ozpi \to A_{r,i}(U, X)\ozpi \\
&\to A_{r,i - 1}(U, Z')\ozpi \to \dots
\end{split} \]
\end{coro}

\begin{proof}
As for Corollary~\ref{coro:512}.
\end{proof}

\subsection{Duality}

We now turn to the section on Duality.

\begin{rema}
We recall that all the material in the subsection ``The moving lemma'' \cite[Section 6]{FV} apply to varieties over an arbitrary field $k$. This is pointed out in the first paragraph of that section. This is also true of \cite[Theorem 7.1]{FV}. The assumption that the base field admits resolution of singularities is said to resume between \cite[Theorem 7.1]{FV} and \cite[Lemma 7.2]{FV}, but the latter doesn't use it (if we take the smoothness of $\overline {U}$ as an assumption). It is needed for \cite[Proposition 7.3]{FV} and the material which follows it.
\end{rema}

\begin{defi}[{cf. \cite[after Proposition 2.1]{FV}}]
For $X, U \in Sch(S)$ and $r \geq 0$ the presheaf $z_{equi}(U, X, r)$ on $Sch(S)$ is defined as
\[ z_{equi}(U, X, r)(-) = z_{equi}(X / S, r)(- \times_S U). \]
That is, it is the composition of $z_{equi}(X, r)$ with the endomorphism $- \times_S U$ of the category $Cor(S)$.

Recall that the correspondence homomorphisms \cite[Section 3.7]{SV} induce a morphism of presheaves \cite[Corollary 3.7.5]{SV}
\[ Cor(-, -): z_{equi}(U, X, r) \otimes z_{equi}(U / S, n) \to z_{equi}(X \times_S U / S, n). \]
If $U \in Sch(k)$ is flat and equidimensional over $S$ of dimension $n$, then $U$ determines an element $cycl_{U/S}(U)$ in $z_{equi}(U/S, n)$. That is, a global section of the presheaf $z_{equi}(U / S, n)$. Evaluating $Cor(-, -)$ on this section defines a morphism of presheaves
\[\dua: z_{equi}(U, X, r) \to z_{equi}(X \times_S U, r + n). \]
\end{defi}

\begin{lemm} \label{lemm:functorialDua}
The morphism $\dua$ is always injective. Furthermore, it is covariantly functorial in $X$ for proper morphisms via the proper push-forward, contravariantly functorial in $X$ for flat equidimensional morphisms ($r$ obviously increases by the relative dimension of the morphism), and contravariantly functorial in $U$ with respect to flat equidimensional morphisms (with the appropriate change in $n$)
\end{lemm}

\begin{proof}
For the injectivity we recall the definition of $Cor(-, -)$. Given a cycle $\beta = \sum n_i z_i \in z_{equi}(U/S, n)(S)$ with $\iota_i: \overline{z}_i \to U$ the canonical closed immersions, and a cycle $\alpha \in z_{equi}(X/S, r)(U)$ we obtain cycles $\iota_i^\circledast \alpha \in z_{equi}(X/S, r)(\overline{z}_i)$ for each $i$. These are formal sums of points of $\overline{z}_i \times_S X$, which we can also consider as formal sums of points of $U \times_S X$. The definition of $Cor(\alpha, \beta)$ is $Cor(\alpha, \beta) = \sum n_i \iota_i^\circledast \alpha $ considered as a formal sum of points in $U \times_S X$. Now if $\beta$ is of the form $cycl_{U / S}(U)$ then the morphism $\amalg \overline{z}_i \to U$ is birational and so $\oplus \iota_i^\circledast$ is injective. Since we are dealing with free abelian groups, $\oplus n_i \iota_i^\circledast$ is also injective, and finally, for each $i$, the points in the formal sum $\iota_i^\circledast \alpha$ considered as points in $U \times_S X$ lie over the generic point $z_i$ of $U$. Hence, each of the formal sums $\iota_i^\circledast \alpha$ contains distinct points. So $\dua$ is injective when evaluated on $S$. To see that it is injective on every scheme $f: V \to S$ in $Sch(S)$ we just replace $S$ with $V$, $U$ with $U \times_S V$ and $X$ with $X \times_S V$. Since $U$ is flat over $S$ we have $f^\circledast cycl_{U / S}(U) = cycl_{V \times_S U / V}(V \times_S U)$.

The functoriality in $X$ is an immediate consequence of \cite[Proposition 3.6.2]{SV} and \cite[Lemma 3.6.4]{SV}. For the contravariance in $U$ suppose that $p: U' \to U$ is a flat equidimensional $S$-morphism of relative dimension $m$. We have an induced morphism or presheaves
\[ z_{equi}(U, X, r) \to z_{equi}(U', X, r) \]
given on $V$ by the appropriate $(U' \times_S V  \to U \times_S V)^\circledast$ and morphism of presheaves
\[ z_{equi}(U / S, r) \to z_{equi}(U'/ S, r + m) \]
\[ z_{equi}(X \times_S U / S, r) \to z_{equi}(X \times_S U'/ S, r + m) \]
given by the flat pullbacks \cite[Lemma 3.6.4]{SV}. By \cite[Lemma 3.7.2]{SV} these fit into a commutative square 
\[ \xymatrix{
z_{equi}(U, X, r) \otimes z_{equi}(U / S, n) \ar[r]^-{Cor(-, -)} \ar[d]_{p^\circledast \otimes p^*} &  z_{equi}(X \times_S U / S, n) \ar[d]^{p^*} \\
z_{equi}(U', X, r) \otimes z_{equi}(U' / S, n + m) \ar[r]_-{Cor(-, -)} &  z_{equi}(X \times_S U' / S, n + m)
} \]
Now since $p^* cycl_{U / S}(U) = cycl_{U' / S}(U')$ we are done.
\end{proof}

For the convenience of the reader we reproduce the following theorem.

\begin{theo}[{\cite[Theorem 7.1]{FV}}]
Let $X, Y$ be smooth projective equidimensional schemes over a field $k$. Then the embedding of presheaves
\[ \dua: z_{equi}(Y, X, r) \to z_{equi}(X \times_k Y, r + n) \]
induces a quasi-isomorphism of presheaves on $Sm(k)$ after applying $\underline{C}_*(-)$.
\end{theo}

We wish to extend this theorem to non-smooth non-projective quasi-projective schemes. To do this we use the presheaves $z\eff(X, n)$ and $z\eff(U, X, n)$ (with $U, X \in Sch(k), n \geq 0$) which are the subpresheaves of $z(X, n)$ and $z(U, X, n)$ consisting of those cycles of the form $\sum n_i z_i$ with all $n_i \geq 0$.


\begin{defi} \label{defi:phiDuality}
Suppose that
\begin{enumerate}
 \item $k$ is a perfect field of exponential characteristic $p$,
 \item $\overline{U}, \overline{X}, \overline{Y}$ are proper schemes in $Sch(k)$, with $\overline{U}$ equidimensional of dimension $n$,
 \item $U \to \overline{U}$, $X \to \overline{X}$ are open immersions,
 \item $\overline{Y} \to \overline{X}$ is a proper morphism,
\end{enumerate}
Given $U, X$ we can always find a suitable $\overline{U}, \overline{X}$ \cite{Nag62}. We define $\alpha\eff_{\overline{Y}}$ as the morphism of presheaves of abelian monoids
\[ \alpha\eff_{\overline{Y}}: z\eff(\overline{U} \times_k \overline{Y}, r + n) \to z\eff(U \times_k X, r + n) \]
which is the composition of the proper push-forward
\[ z\eff(\overline{U} \times_k \overline{Y}, r + n) \to z\eff(\overline{U} \times_k \overline{X}, r + n) \]
and the flat pullback
\[ z\eff(\overline{U} \times_k \overline{X}, r + n) \to z\eff(U \times_k X, r + n). \]
We also define the corresponding morphism of presheaves of abelian groups
\[ \alpha_{\overline{Y}}: z_{equi}(\overline{U} \times_k \overline{Y}, r + n) \to z_{equi}(U \times_k X, r + n). \]
The presheaf of abelian monoids $\Phi\eff_{\overline{Y}}$ is the presheaf which fits into the following cartesian diagram
\[ \xymatrix{
 \Phi\eff_{\overline{Y}} \ar[r] \ar[d] & z\eff_{equi}(U, X, r) \ar[d] \\
z_{equi}(\overline{U} \times \overline{Y}, r + n) \ar[r]_{\alpha} & z_{equi}(X \times_k U, r + n) 
} \]
and the subpresheaf of abelian groups
\[ \delta_{\overline{Y}}: \Phi_{\overline{Y}} \to z_{equi}(\overline{U} \times_k \overline{Y}, r + n) \]
is defined to be the subpresheaf of abelian groups generated by the subpresheaf of abelian monoids $ \Phi\eff_{\overline{Y}}$.
\end{defi}

Hence, we have a corresponding commutative square
\begin{equation} \label{equa:dualityLemma}
\xymatrix{
 \Phi_{\overline{Y}} \ar[r] \ar[d]_{\delta_{\overline{Y}}} & z_{equi}(U, X, r) \ar[d]^{\dua} \\
z_{equi}(\overline{U} \times \overline{Y}, r + n) \ar[r]_\alpha & z_{equi}(X \times_k U, r + n) 
} \end{equation}
of presheaves of abelian groups. Voevodsky-Friedlander warn us that this is not in general cartesian \cite[before Lemma 7.2]{FV} but that it is in the case that $\overline{Y} = \overline{X}$.

\begin{exam}
Consider the case $U = \overline{U}, X = \overline{X}$ and $\overline{Y} = \overline{X} \amalg \overline{X}$. Let $\alpha$ be any cycle in $ z\eff(X \times_k U, r + n) $ that is not in $z\eff(U, X, r)$.  Then $(\alpha, -\alpha)$ is in the pullback of the square  (\ref{equa:dualityLemma}) but not in $\Phi_{\overline{Y}}$.
\end{exam}


\begin{lemm}[{cf. \cite[proof of  Theorem 7.4]{FV}}] \label{lemm:iffDuality}
With the notation and assumptions of Definition~\ref{defi:phiDuality} the morphism $\delta_{\overline{X}}$ induces a quasi-isomorphism of complexes of abelian groups after applying $\underline{C}_*(-)(k)\ozpi$ if and only if $\dua$ does.
\end{lemm}

\begin{proof}
We use the following diagram of morphisms of presheaves 
\[ \xymatrix{
0 \ar[r] & ker' \ar[r] \ar[d]^a & \Phi_{\overline{X}} \ar[r] \ar[d]^{\delta_{\overline{X}}} & z_{equi}(U, X) \ar[r] \ar[d]^{\dua} & coker_1 \ar[r] \ar[d]^c & 0 \\
0 \ar[r] & ker(\alpha_{\overline{X}}) \ar[r] & z_{equi}(\overline{U} \times_k \overline{X}) \ar[r]_{\alpha_{\overline{X}}} & z_{equi}(U \times_k X) \ar[r] & coker_2 \ar[r] & 0
} \]
where we have used the abbreviations
\[ z_{equi}(U, X) = z_{equi}(U, X, r) \]
\[ z_{equi}(\overline{U} \times_k \overline{X}) = z_{equi}(\overline{U} \times_k \overline{X}, r + n) \]
\[ z_{equi}(U \times_k X) = z_{equi}(U \times_k X, r + n) \]
Since $\dua$ is a monomorphism and the square involving $\delta_{\overline{X}}$ and $\dua$ is cartesian, the morphism $a$ is an isomorphism and the other vertical morphisms are all monomorphisms. Hence, it suffices to show that $\underline{C}_*(coker_i)(k)\ozpi$ is acyclic for $i = 1, 2$. By Theorem~\ref{theo:MainResultRos} it suffices to show that $(coker_i)_{cdh} = 0$ for $i = 1, 2$ and since $c$ is a monomorphism we can restrict our attention to $i = 2$. This is a standard application of the platification theorem (Theorem~\ref{theo:platification}) as described in \cite[Theorem 4.3.1]{SV} and \cite[Theorem 4.2.9]{SV}.
\end{proof}

\begin{prop}[{cf. \cite[Proposition 7.3]{FV}}] \label{prop:smoothDuality}
With the notation and assumptions of Definition~\ref{defi:phiDuality} suppose further that $\overline{Y}$ and $\overline{U}$ are smooth. Then the morphism $\delta_{\overline{Y}}$ induces a quasi-isomorphism of complexes of abelian groups after applying $\underline{C}_*(-)(k)$.
\end{prop}

\begin{proof}
This is the first case treated in the proof of \cite[Proposition 7.3]{FV}.
\end{proof}

\begin{theo}[{cf. \cite[Theorem 7.4]{FV}}] \label{theo:firstDuality}
Suppose that $k$ is a perfect field, $p$ its exponential characteristic, $U$ a reduced quasi-projective equidimensional scheme of dimension $n$ over $k$, and $X$ a scheme of finite type over $k$. Then for any $r \geq 0$ the embedding
\[ \dua: z_{equi}(U, X, r) \to z_{equi}(X \times_k U, r + n) \]
induces a quasi-isomorphism of complexes of abelian groups after applying $\underline{C}_*(-)(k)\ozpi$.
\end{theo}

\begin{proof}
We can assume that $X$ is reduced as the canonical morphism $X_{red} \to X$ induce isomorphisms of all the presheaves involved.

Choose embeddings of $U$ and $X$ as open subschemes of proper $k$-schemes $U \to \overline{U}$, $X \to \overline{X}$ \cite{Nag62} so that we are in the situation of Definition~\ref{defi:phiDuality}. By Lemma~\ref{lemm:iffDuality} it suffices to show that $\delta_{\overline{X}}$ induces a quasi-isomorphism after applying $\underline{C}_*(-)(k)\ozpi$. We will show that we obtain a quasi-isomorphism after applying $\underline{C}_*(-)(k) \otimes \zll$ for each $\ell \neq p$.

Let $\overline{U}' \to \overline{U}$ and $\overline{X}' \to \overline{X}$ be morphisms given by Theorem~\ref{theo:gabberGlobal} and let $V \to \overline{X}$ (resp. $W\to \overline{U}$) be an open immersion such that the induced morphism $V \times_{\overline{X}} \overline{X}' \to V$ (resp. $W \times_{\overline{U}} \overline{U}' \to W$) is finite flat surjective locally (on the target) of degree prime to $\ell$. Define $V' = V \times_{\overline{X}} \overline{X}'$ and $W' = W \times_{\overline{U}} \overline{U}'$.

Replacing $U$ and $X$ with $W$ and $V$ and using Lemma~\ref{lemm:iffDuality} again if suffices to show that
\[ \dua \otimes \zll: z_{equi}(W, V, r) \otimes \zll \to z_{equi}(V \times_k W, r + n) \otimes \zll \]
induces a quasi-isomorphism after applying $\underline{C}_*(-)(k)$. Now $\dua$ is functorial with respect to flat pullback and proper push-forward (Lemma~\ref{lemm:functorialDua}), and so since the degrees of our flat finite surjective morphisms are invertible in $\zll$, this $\dua$ just mentioned is a retraction of 
\[ \dua \otimes \zll : z_{equi}(W', V', r) \otimes \zll \to z_{equi}(V' \times_k W', r + n) \otimes \zll \]
(cf. \cite[Lemma 2.3.5]{SV}). So it suffices to show that this latter induces a quasi-isomorphism after applying $\underline{C}_*(-)(k)$. Now $W'$ and $V'$ are open subschemes of $\overline{U}'$ and $\overline{X}'$ respectively which are both smooth and proper over $k$. Using Lemma~\ref{lemm:iffDuality} a final time we reduce to showing that 
\[ \delta_{\overline{X}'} \otimes \zll : \Phi_{\overline{X}'} \otimes \zll \to z_{equi}(\overline{X}' \times_k \overline{U}', r + n) \otimes \zll \]
induces a quasi-isomorphism after applying $\underline{C}_*(-)(k)$. This is given by Proposition~\ref{prop:smoothDuality}.
\end{proof}

\begin{coro} \label{coro:firstDualityPresheaf}
Suppose that $k$ is a perfect field, $p$ its exponential characteristic, $U$ a reduced quasi-projective equidimensional scheme of dimension $n$ over $k$, and $X$ a scheme of finite type over $k$. Then for any $r \geq 0$ the embedding
\[ \dua: z_{equi}(U, X, r) \to z_{equi}(X \times_k U, r + n) \]
induces a quasi-isomorphism after applying $\underline{C}_*(-)\ozpi$ of complexes of presheaves on the category of quasi-projective smooth $k$-schemes.
\end{coro}

\begin{proof}
For any smooth scheme $V$ of dimension $m$ we have
\[ \xymatrix{
\underline{h}_i(z_{equi}(U, X, r))(V)\ozpi \ar@{=}[d]  \ar[r] & \underline{h}_i(z_{equi}(X \times_k U, r + n))(V)\ozpi \ar@{=}[d] \\
\underline{h}_i(z_{equi}(X, r))(V \times_k U)\ozpi \ar@{=}[d] & \underline{h}_i(z_{equi}(V , X \times_k U, r + n))(k)\ozpi \ar[d]^\cong\\
\underline{h}_i(z_{equi}(U \times_k V, X, r))(k)\ozpi \ar[r]_-\cong & \underline{h}_i(z_{equi}(X \times_k U \times_k V, r + n + m))(k)\ozpi
} \]
where the isomorphisms are given by Theorem~\ref{theo:firstDuality}.
\end{proof}

\subsection{Properties}

\begin{defi}[{cf. \cite[Beginning of Section 4]{FV}}]
If $F$ is a presheaf on $Sch(k)$ or $Sm(k)$ recall that \cite{FV} denote by $\underline{h}_i(F)$ the homology presheaves of the complex of presheaves $\underline{C}_*(F)$.
\end{defi}

\begin{theo}[{cf. \cite[Theorem 8.1]{FV}}] \label{theo:FV81}
Let $k$ be a perfect field of exponential characteristic $p$, let $U$ be a smooth quasi-projective scheme over $k$ and $X$ a separated scheme of finite type over $k$. Then the natural homomorphisms of abelian groups
\[ \underline{h}_i(z_{equi}(X, r))(U)\ozpi \to A_{r,i}(U, X)\ozpi \]
are isomorphisms for all $i \in \ZZ$.
\end{theo}

\begin{proof}
The proof of \cite[Theorem 8.1]{FV} works fine after applying $(-)\ozpi$ to everything.
\end{proof}

\begin{theo}[{cf. \cite[Theorem 8.2]{FV}}] \label{theo:FV82}
Let $k$ be a perfect field of exponential characteristic $p$, let $U$ be a smooth scheme of pure dimension $n$ over $k$, and let $X, Y$ be separated schemes of finite type over $k$. Then there are canonical isomorphisms
\[ A_{r,i}(Y \times U, X) \ozpi \stackrel{\sim}{\to} A_{r + n, i}(Y, X \times U)\ozpi. \]
\end{theo}

\begin{proof}
We begin with canonical morphisms
\[ \HH^{-i}_{cdh}(Y, \underline{C}_*(z_{equi}(U, X, r))_{cdh}) \to \HH^{-i}_{cdh}(Y \times_k U, \underline{C}_*(z_{equi}(X, r))_{cdh}). \]
Let $p: U \to k$ denote the projection. We have two canonical left exact functors
\[ Shv_{cdh}(Sch(U)) \stackrel{p_*}{\to} Shv_{cdh}(Sch(k)) \stackrel{\Gamma(Y, -)}{\to} Ab \]
where the first is composition with the functor $- \times_k U: Sch(k) \to Sch(U)$ and the second is evaluation at $Y$. Let $G = p_*, F = \Gamma(Y, -)$, and $K = \underline{C}_*(z_{equi}(X, r))_{cdh}$. Then the hypercohomology groups on the left are
\[ H^{i-1}RF(G (K)) \]
and the hypercohomology groups on the right are 
\[ H^{i-1}R(FG)(K) = H^{i-1}RF(RG (K)). \]
Our morphism is induced by the canonical morphism $G(K) \to RG (K)$. Notice that these are natural in $Y$.

To show that these canonical morphisms are isomorphisms (after $(-)\ozpi$), it suffices to do so after $-\otimes \zll$ for each prime $\ell$ different from $p$. As we are dealing with presheaves with transfers, we can replace the cdh topology with the $\ldh$ topology, and so we are now trying to show that the canonical morphisms 
\[ \HH^{-i}_{\ldh}(Y, \underline{C}_*(z_{equi}(U, X, r))_{\ldh}) \otimes \zll \to \HH^{-i}_{\ldh}(Y \times_k U, \underline{C}_*(z_{equi}(X, r))_{\ldh}) \otimes \zll \]
are isomorphisms. Due to the theorem of Gabber (Corollary~\ref{coro:regularlCover}) it suffices to consider the case when $Y$ is smooth and quasi-projective, and indeed we can also replace $Sch(k)$ by $Sm(k)$. In this case, we claim that, in the notation used above, $G(K) \to RG (K)$ is an isomorphism in the derived category of complexes of $\ldh$ sheaves. For this morphism to be an isomorphism it suffices that it induces isomorphisms after applying $\HH^{-i}(V, -)$ for each smooth quasi-projective $V$. That is, we have returned to the following version of our initial morphism
\[ \HH^{-i}_{\ldh}(Y, \underline{C}_*(z_{equi}(U, X, r))_{\ldh}) \otimes \zll \to \HH^{-i}_{\ldh}(Y \times_k U, \underline{C}_*(z_{equi}(X, r))_{\ldh}) \otimes \zll \]
but now assuming that $Y$ is smooth and quasi-projective and that we are on the $\ldh$ site $Sm(k)$. We have the morphism
\[ \dua: z_{equi}(U, X, r) \to z_{equi}(X \times_k U, r + n) \]
and Corollary~\ref{coro:firstDualityPresheaf} tells us that this induces a quasi-isomorphism of presheaves on $Sm(k)$ after applying $\underline{C}_*(-) \otimes \zll$ and hence it suffices to show that the induced morphism 
\[ \HH^{-i}_{\ldh}(Y, \underline{C}_*(z_{equi}(X \times_k U, r + n))_{\ldh}) \otimes \zll \to \HH^{-i}_{\ldh}(Y \times_k U, \underline{C}_*(z_{equi}(X, r))_{\ldh}) \otimes \zll \]
is an isomorphism. By definition, this morphism is the morphism
\[ A_{r + n, i}(Y, X \times_k U) \otimes \zll \to A_{r, i}(Y \times_k U, X) \otimes \zll. \]
Our result now follows from Theorem~\ref{theo:FV81} using the following diagram.
\[ \xymatrix{
A_{r + n, i}(Y, X \times_k U) \otimes \zll \ar@{=}[dd]_{\ref{theo:FV81}}  \ar[r] & A_{r, i}(Y \times_k U, X) \otimes \zll \ar@{=}[d]^{\ref{theo:FV81}} \\
 & \underline{h}_i(z_{equi}(X, r))(Y \times_k U) \ar@{=}[d] \\
\underline{h}_i(z_{equi}(X \times_k U, r + n))(Y) & \ar@{=}[l]^-{\ref{coro:firstDualityPresheaf}} \underline{h}_i(z_{equi}(U, X, r))(Y)
} \]
\end{proof}

\begin{theo}[{cf. \cite[Theorem 8.3]{FV}}] \label{theo:FV83}
Let $k$ be a perfect field of exponential characteristic $p$ and let $X, Y$ be separated schemes of finite type over $k$.
\begin{enumerate}
 \item (Homotopy invariance) The pull-back homomorphism $z_{equi}(X, r) \to z_{equi}(X \times \AA^1, r + 1)$ induces for any $i \in \ZZ$ an isomorphism
\[ A_{r,i}(Y, X)\ozpi \to A_{r + 1, i}(Y, X \times \AA^1)\ozpi. \]
 \item (Suspension) Let
\[ p: X \times \PP^1 \to X \]
\[ i: X \to X \times \PP^1 \]
be the natural projection and closed embedding. Then the morphism
\[ i_* \oplus p^*: z_{equi}(X, r + 1) \oplus z_{equi}(X, r) \to z_{equi}(X \times \PP^1, r + 1) \]
induces an isomorphism
\[ A_{r + 1, i}(Y, X)\ozpi \oplus A_{r, i}(Y, X)\ozpi \to A_{r + 1, i}(Y, X \times \PP^1)\ozpi. \]
 \item (Cosuspension) There are canonical isomorphisms:
\[ A_{r,i}(Y \times \PP^1, X)\ozpi \stackrel{\sim}{\to} A_{r + 1, i}(Y, X)\ozpi \oplus A_{r,i}(Y, X)\ozpi. \]
 \item (Gysin) Let $Z \subset U$ be a closed immersion of smooth schemes everywhere of codimension $c$ in $U$. Then there is a canonical long exact sequence of abelian groups of the form
\[ \dots A_{r + c, i}(Z, X)\ozpi \to A_{r,i}(U, X)\ozpi \to A_{r, i}(U - Z, X)\ozpi \]
\[ \to A_{r + c, i - 1}(Z, X)\ozpi \to \dots \]
\end{enumerate}
\end{theo}

\begin{proof}
\begin{enumerate}
 \item Follows from Theorem~\ref{theo:FV81} and homotopy invariance (Proposition~\ref{prop:FV59}).

 \item Follows from the localisation sequence (Corollary~\ref{coro:512}) and the first part.
 \item Follows from Theorem~\ref{theo:FV81} and the second part.
 \item Follows from Theorem~\ref{theo:FV81} and (Corollary~\ref{coro:512}).
\end{enumerate}
\end{proof}








\section{Triangulated categories of motives over a field --- After Voevodsky} \label{sec:traCatMotFie}

In this section we show how Theorem~\ref{theo:MainResultRos} can be used to lift the assumption of resolution of singularities on all of the results in \cite{Voev00}, if we are willing to work $\zpi$-linearly. The principle is that every time Voevodsky assumes the existence of a smooth cdh cover we can use the existence of a smooth $\ldh$-cover, and the only other way he uses resolution of singularities is via \cite[Theorem 4.1.2]{Voev00} which we replace with Theorem~\ref{theo:MainResultRos}.

\begin{defi}
We define $DM\eff_{gm}(k, \zpi)$ as the full triangulated subcategory of compact objects in $DM\eff(k, \zpi)$. The category $DM_{gm}(k, \zpi)$ is obtained by formally adjoining a tensor inverse to $\zpi(1)$ as is done for Chow motives.
\end{defi}

\begin{lemm}[{cf. \cite[Corollary 4.1.4]{Voev00}}] \label{lemm:geometricMotive}
Let $k$ be a perfect field of exponential characteristic $k$. Then $DM\eff\gm(k, \zpi)$ contains $\mp{X}$ for any scheme $X$ of finite type over $k$.
\end{lemm}

\begin{proof}
Follows immediately from Corollary~\ref{coro:a1localisations}.
\end{proof}

\begin{prop}[{cf. \cite[Corollary 3.5.5]{Voev00}}] \label{prop:generatedBysmproj}
Let $k$ be a perfect field of exponential characteristic $p$. Then $DM\eff\gm(k, \zpi)$ is generated as a pseudo-abelian triangulated category by objects of the form $\mp{X}$ for smooth projective varieties $X$ over $k$.
\end{prop}

\begin{proof}
We will show that the image of the family $SP = \{ \mp{X} : X $ is smooth and projective $\}$ in $DM\eff(k, \zll)$ is a compact generating family for every prime $\ell$ different from $p$. It then follows from Lemma~\ref{lemm:generatingLocalTri} that the smallest triangulated category of $DM\eff(k, \zpi)$ containing $SP$ is in fact the full subcategory of compact objects.

Let $\mathscr{P}$ denote the smallest pseudo-abelian triangulated subcategory of $DM\eff(k, \zll)$ containing objects of the form $\ml{X}$ for smooth projective varieties $X$ over $k$. As $\mathscr{L} = \{ \ml{X} : X$ is smooth $\}$ is a compact generating family for $DM\eff(k, \zll)$, the smallest triangulated category containing $\mathscr{L}$ is the full subcategory of compact objects of $DM\eff(k, \zll)$. So it suffices to show that every $\ml{X}$ with $X$ smooth (but not necessarily projective) is contained in $\mathscr{P}$. We will do so by induction on the dimension of $X$. Suppose it is true for all smooth schemes of dimension strictly less than $d$ and let $X$ be a smooth scheme of dimension $d$. Due to the Mayer-Vietoris distinguished triangles \cite[Lemma 2.1.2]{Voev00} it suffices to consider $X$ quasi-projective. Let $X \to \overline{X}$ be a compactification of $X$ and $\overline{Y} \to \overline{X}$ a morphism given by Theorem~\ref{theo:gabberGlobal}. So $\overline{Y}$ is a smooth projective variety and there exists a dense open subscheme $U \subset X$ such that $U \times_{\overline{X}} \overline{Y} \to U$ is finite flat surjective of degree prime to $\ell$.

We claim that for any dense open embedding $V \to V'$ of smooth schemes, $M\gm(V)$ is in $\mathscr{P}$ if and only if $M\gm(V')$ is in $\mathscr{P}$. Assuming this claim we proceed as follows. By definition of $\mathscr{P}$ it contains $\overline{Y}$. Due to our claim, $M\gm(U \times_{\overline{X}} \overline{Y})$ is in $\mathscr{P}$. By the degree formula for correspondences  and the fact that we are working $\zll$-linearly, $M\gm(U)$ is a retract of $M\gm(U \times_{\overline{X}} \overline{Y})$. Since $\mathscr{P}$ is pseudo-abelian, this implies that $M\gm(U)$ is in $\mathscr{P}$. Finally, by using our claim again, this implies that $M\gm(X)$ is in $\mathscr{P}$.

It remains to prove our claim. Let $V \to V'$ be a dense open embedding of smooth schemes. Since the base field is perfect, every reduced scheme contains an open dense smooth scheme. Consequently, there exists a sequence $V = V_0 \subset V_1 \subset \dots \subset V_n = V'$ of dense open immersions such that each $V_i - V_{i - 1}$ is smooth and everywhere of codimension $c_i$ for some $c_i$. Then our claim follows from the inductive hypothesis and the triangles
\[ M\gm(V_{i - 1}) \to M\gm(V_i) \to M\gm(V_i - V_{i - 1})(c_i)[2c_i] \to M\gm(V_{i - 1})[1] \]
given by \cite[Proposition 3.5.4]{Voev00}.
\end{proof}

\begin{prop}[{cf. \cite[Proposition 4.1.3]{Voev00}}] \label{prop:cdhTriangleDM}
Consider a cartesian square of morphisms of schemes of finite type over $k$ of the form
\begin{equation} \label{equi:blowupSquare}
\xymatrix{
Z' \ar[r] \ar[d] & X' \ar[d]^p \\
Z \ar[r]_i & X
} \end{equation}
such that $p$ is proper, $i$ is a closed immersion, and $p$ is an isomorphism over $X - Z$. Then there is a canonical distinguished triangle in $DM\eff(k, \zpi)$ of the form
\[ \mp{Z'} \to \mp{Z} \oplus \mp{X'} \to \mp{X} \to \mp{Z'}[1]. \]
\end{prop}

\begin{proof}
Follows from the following short exact sequence of cdh sheaves (\cite[Theorem 4.2.9]{SV}, \cite[Proposition 4.3.3]{SV}) and our definitions.
\[ 0 \to c_{equi}(Z'/k, 0)_{cdh} \to c_{equi}(Z/k, 0)_{cdh} \oplus c_{equi}(X'/k, 0)_{cdh} \to c_{equi}(X/k, 0)_{cdh} \to 0. \]
\end{proof}

\begin{prop}[{cf. \cite[Proposition 4.1.5]{Voev00}}] \label{prop:compactLocalisationTriangle}
Let $k$ be a perfect field of exponential characteristic $p$ and $X$ a scheme of finite type over $k$. Let $Z$ be a closed subscheme of $X$. Then there is a canonical distinguished triangle in $DM\eff(k, \zpi)$ of the form
\[ \mcp{Z} \to \mcp{X} \to \mcp{X - Z} \to \mcp{Z}[1]. \]
If $X$ is proper then there is a canonical isomorphism $\mcp{X} \cong \mp{X}$.
\end{prop}

\begin{proof}
The second statement follows from the equality $c_{equi}(X/k, 0) = z_{equi}(X/k, 0)$ when $X$ is proper. The proof of the first is the same as that of Proposition~\ref{prop:cdhTriangleDM} with $z_{equi}, z$ replacing $c_{equi}, c$ and the short exact sequence \cite[Theorem 4.3.1]{SV} replacing \cite[Theorem 4.3.3]{SV}.
\end{proof}

\begin{lemm}[{cf. \cite[Corollary 4.1.6]{Voev00}}] \label{lemm:geometricCompactMotive}
Let $k$ be a perfect field of exponential characteristic $k$. Then $DM\eff\gm(k, \zpi)$ contains $\mcp{X}$ for any scheme $X$ of finite type over $k$.
\end{lemm}

\begin{proof}
Follows immediately from Lemma~\ref{lemm:geometricMotive} and Proposition~\ref{prop:compactLocalisationTriangle} using a compactification \cite{Nag62}.
\end{proof}



\begin{prop} \label{prop:reductionGabber}
Let $k$ be a perfect field of exponential characteristic $p$.
Suppose that $M_1, M_2: Cor(k)^{op} \to \T$ are functors to a $\zpi$-linear triangulated category. Let $\eta: M_1 \to M_2$ be a natural transformation between these functors. Suppose further, that for every cartesian square of the form (\ref{equi:blowupSquare}) such that $p$ is proper, $i$ is a closed immersion, and $p$ is an isomorphism over $X - Z$, there exist morphisms $\mu_i: M_i(X) \to M_i(Z')[1]$ for $i = 1, 2$ such that the triangles
\[ M_i(Z') \to M_i(Z) \oplus M_i(X') \to M_i(X) \stackrel{\mu_i}{\to} M_i(Z')[1]. \]
are distinguished for $i = 1, 2$ and the squares
\[ \xymatrix{
M_1(X) \ar[r] \ar[d] & M_1(Z')[1] \ar[d] \\
M_2(X) \ar[r] & M_2(Z')[1]
} \]
are commutative. Then if $\eta$ is an isomorphism for every smooth scheme, it is an isomorphism for every scheme.

\end{prop}

\begin{proof}
After Lemma~\ref{lemm:conservativeLocal} it suffices to show that the statement is true for every $\zll$-linear triangulated category for each prime $\ell$ different to $p$. We will work by induction on the dimension of $X$. The natural transformation $\eta$ can be seen to be an isomorphism in dimension zero by considering the distinguished triangle with $Z = X' = X_{red}$. So suppose that $\eta$ is an isomorphism for all schemes of dimension less than $n$, and let $X$ be a scheme of dimension $n$. Let $Y \to X$ be a morphism given by Theorem~\ref{theo:gabberGlobal} and $\tilde{X} \to X$ a blow-up such that the proper transform $\tilde{Y} \to \tilde{X}$ is flat (Theorem~\ref{theo:platification}). Let $Z, W$ be closed subschemes of $X, Y$ such that the dimensions of $Z, W, Z \times_X \tilde{X}$ and $W \times_Y \tilde{Y}$ are less than $n$, and $\tilde{X} \to X$ (resp. $\tilde{Y} \to Y$) is an isomorphism over $X - Z$ (resp. $Y - W$). Then considering the associated distinguished triangles, to show that $M_1(X) \to M_2(X)$ (resp. $M_1(\tilde{Y}) \to M_2(\tilde{Y})$) is an isomorphism, due to the induction hypothesis, it suffices to show that $M_1(\tilde{X}) \to M_2(\tilde{X})$ (resp. $M_1(Y) \to M_2(Y)$) is an isomorphism. Since $Y$ is smooth, it follows that $M_1(\tilde{Y}) \to M_2(\tilde{Y})$ is an isomorphism. Now since $Y \to X$ is generically finite surjective of degree prime to $\ell$, the flat morphism $\tilde{Y} \to \tilde{X}$ is globally finite surjective and locally of degree prime to $\ell$. In particular, as $M_1, M_2$ are natural in $Cor(k)$ and we are working $\zll$ linearly, the morphism $M_1(\tilde{X}) \to M_2(\tilde{X})$ is a retract of $M_1(\tilde{Y}) \to M_2(\tilde{Y})$ (Proposition~\ref{prop:corHasTraces} - traces). Hence, $M_1(\tilde{X}) \to M_2(\tilde{X})$  is an isomorphism, and therefore $M_1(X) \to M_2(X)$  is an isomorphism.

\end{proof}

\begin{prop}[{cf. \cite[Proposition 4.1.7]{Voev00}}] \label{prop:productIsomorphisms}
Let $k$ be a perfect field of exponential characteristic $p$. Let $X, Y$ be schemes of finite type over $k$. In $DM\eff(k, \zpi)$ there are canonical isomorphisms:
\[ \mp{X} \otimes \mp{Y} \cong \mp{X \times_k Y} \textrm{ and } \]
\[ \mcp{X} \otimes \mcp{Y} \cong \mcp{X \times_k Y}. \]
\end{prop}

\begin{proof}
Due to Proposition~\ref{prop:reductionGabber}, for the first isomorphism it is sufficient to give an isomorphism
\[ c_{equi}(X / k, 0)_{cdh} \stackrel{tr}{\otimes} c_{equi}(Y / k, 0)_{cdh} \to c_{equi}(X \times_k Y / k, 0)_{cdh} \]
in $Shv_{cdh}(Cor(k))$. This follows immediately from the definition of $\stackrel{tr}{\otimes}$ as the functor induced by the left Kan extension along $Cor(k) \times Cor(k) \to Cor(k)$.

For the second isomorphism, use compactifications \cite{Nag62} and Proposition~\ref{prop:compactLocalisationTriangle}.
\end{proof}

\begin{coro}[{cf. \cite[Corollary 4.1.8]{Voev00}}]
Let $k$ be a perfect field of exponential characteristic $p$. For any scheme $X$ of finite type over $k$ one has canonical isomorphisms
\[ \mp{X \times \AA^1} \cong \mp{X} \textrm{ and } \]
\[ \mcp{X \times \AA^1} \cong \mcp{X}(1)[2]. \]
In particular, we have
\[ \mcp{\AA^n} \cong \zpi(1)[2]. \]
\end{coro}

\begin{proof}
After Proposition~\ref{prop:productIsomorphisms} it is sufficient to show that
\[ \mp{\AA^1} \cong \zpi, \textrm{ and } \]
\[ \mcp{\AA^1} \cong \zpi(1)[2]. \]
The first follows from the definition of $DM\eff(k, \zpi)$ as we have inverted $\AA^1_k \to k$. The second follows from the definition of the Tate object and Proposition~\ref{prop:compactLocalisationTriangle}.
\end{proof}

\begin{coro}[{cf. \cite[Corollary 4.1.11]{Voev00}}]
Let $k$ be a perfect field of exponential characteristic $p$ and $X$ a scheme of finite type over $k$. Let $\E$ be a vector bundle on $X$. Denote by $p: \PP(\E) \to X$ the projective bundle over $X$ associated with $\E$. Then one has a canonical isomorphism in $DM\eff(k, \zpi)$ of the form
\[ \mp{\PP(\E)} \cong \oplus_{n = 0}^{\dim \E - 1} \mp{X}(n)[2n]. \]
\end{coro}

\begin{proof}
The proof of \cite[Corollary 4.1.11]{Voev00} works fine with the usual adjustments to use the theorem of Gabber as done in Proposition~\ref{prop:productIsomorphisms}.
\end{proof}

\begin{prop}[{cf. \cite[Proposition 4.2.3, Corollaries 4.2.4, 4.2.5, 4.2.6, 4.2.7, Theorem 4.3.2]{Voev00}}] \label{prop:V423}
Let $k$ be a perfect field of exponential characteristic $p$ and $X, Y$ schemes of finite type over $k$.
\begin{enumerate}
 \item \label{item:V423}For any $r \geq 0$ there are canonical isomorphisms
\[ \hom_{DM\eff(k, \zpi)}(\mp{Y}(r)[2r + i], \mcp{X}) \cong A_{r, i}(Y, X)\ozpi. \]

 \item If $f: X \to Y$ is a flat equidimensional morphism of relative dimension $n$ then there is a canonical morphism in $DM\eff(k, \zpi)$ of the form
\[ f^*: \mcp{Y}(n)[2n] \to \mcp{X} \]
and these morphisms satisfy the standard properties of the contravariant functoriality of algebraic cycles.

 \item \label{item:V425}If $X$ happens to be smooth, and we denote by $A^i(X)$ the group of cycles of codimension $i$ on $X$ modulo rational equivalence, then there is a canonical isomorphism
\[ A^i(X)\ozpi \cong \hom_{DM\eff(k, \zpi)}(\mp{X}, \zpi(i)[2i]). \]

 \item \label{prop:V423:comparisonChowGroups} If $X, Y$ are smooth and proper then one has
\[ \hom_{DM\eff(k, \zpi)}(\mp{X}, \mp{Y}) \cong A_{\dim(X)}(X \times_k Y)\ozpi \]
\[ \hom_{DM\eff(k, \zpi)}(\mp{X}, \mp{Y}[i]) = 0 \textrm{ for } i > 0. \]

 \item If $X$ is smooth then there is a canonical isomorphism 
\[ \ihom(\mp{X}, \mcp{Y}) \cong \underline{C}_*(z_{equi}(X, Y, 0)). \]

 \item Suppose that $k$ is a perfect field of exponential characteristic $p$. Let $X$ be a smooth proper scheme of dimension $n$ over $k$. The morphism
\[ \mp{X} \to \ihom(\mp{X}, \zpi(n)[2n]) \]
induced by the diagonal $X \to X \times X$ (\ref{prop:V423}(\ref{item:V425})) is an isomorphism.
\end{enumerate}
\end{prop}

\begin{proof}
The proofs in \cite{Voev00} go through without changes.
\end{proof}

We recall Voevodsky's Cancellation Theorem.

\begin{theo}[{\cite[Corollary 4.10]{Voe10}}] \label{theo:cancellationTheorem}
Suppose that $k$ is a perfect field. Then for any $K, L \in DM\eff(k)$ the map
\[ \hom(K, L) \to \hom(K(1), L(1)) \]
is a bijection.
\end{theo}

\begin{prop}[{cf. \cite[Proposition 4.3.3]{Voev00}}]
Let $k$ be a perfect field of exponential characteristic $p$. Let $X$ be a scheme of finite type of dimension $n$ over $k$. Then for any $n, r \geq 0$ the morphism
\[ \ihom\biggl(\mp{X}, \zpi(n)\biggr)(r) \to \ihom\biggl(\mp{X}, \zpi(n + r) \biggr ) \]
is an isomorphism.
\end{prop}

\begin{proof}
For any $A$ we have the following isomorphisms given by Theorem~\ref{theo:cancellationTheorem}:
\[ \xymatrix{
\hom(A(r),\ \  \ihom \biggl(\mp{X},\zpi(n)\biggr)(r)) \ar@{=}[dd] &\hom(A(r), \ \ihom \biggl (\mp{X}, \zpi(n + r) \biggr )) \ar@{=}[d] \\
& \hom(A(r) \otimes \mp{X}, \ \  \zpi(n + r)) \ar@{=}[d] \\
\hom(A, \ \ \ihom \biggl (\mp{X}, \zpi(n) \biggr )) \ar@{=}[r] & \hom(A \otimes \mp{X},\   \zpi(n))
} \]
Taking $A = \ihom(\mp{X}, \zpi(n))$, the identity $\id_{A(r)}$ induces the desired morphism, and moreover, this morphism is natural in $X$. We then use Proposition~\ref{prop:reductionGabber} (together with Proposition~\ref{prop:cdhTriangleDM}) to reduce to the case when $X$ is smooth, in which case it follows from Proposition~\ref{prop:V423}(\ref{item:V423}) and Theorem~\ref{theo:FV82}.
\end{proof}


\begin{theo}[{cf. \cite[Theorem 4.3.7]{Voev00}}] \label{theo:dualityDMp}
Let $k$ be a perfect field of exponential characteristic $p$. Then $DM\gm(k, \zpi)$ has an internal hom. Setting $A^* = \ihom_{DM\gm(k, \zpi)}(A, \zpi)$ one has:
\begin{enumerate}
 \item For any object $A$ in $DM\gm(k, \zpi)$ the canonical morphism $A \to (A^*)^*$ is an isomorphism.
 \item For any pair of objects $A, B$ of $DM\gm(k, \zpi)$ there are canonical morphisms
\[ (A \otimes B)^* = A^* \otimes B^* \]
\[ \ihom(A, B) = A^* \otimes B. \]
 \item For a smooth scheme $X$ of pure dimension $n$ over $k$ one has canonical isomorphisms
\[ \mp{X}^* \cong \mcp{X}(-n)[-2n] \]
\[ \mcp{X}^* \cong \mp{X}(-n)[-2n]. \]
\end{enumerate}
\end{theo}

As in \cite[before Proposition 2.1.4]{Voev00}, let $Chow\eff(k)$ be the pseudo-abelian envelope of the category $\mathcal{C}_0$, whose objects are smooth projective $k$ schemes, and morphisms are $\hom_{\mathcal{C}_0}(X, Y) = \oplus_{X_i} A_{\dim(X_i)}(X_i \times_k Y)$ where the $X_i$ are the connected components of $X$ and $A_d(-)$ is the group of cycles of dimension $d$ modulo rational equivalence. We also write
\[ Chow : SmProj(k) \to Chow^{eff}(k) \]
for the corresponding functor.

\begin{prop}[{cf. \cite[Theorem 2.1.4]{Voev00}, \cite[Remark after Corollary 2.15]{Voev00}}] \label{theo:embeddingOfChowMotives}
Let $k$ be a perfect field of exponential characteristic $p$. There exists a functor $Chow\eff(k) \to DM\gm\eff(k)$ such that the following diagram commutes
\[ \xymatrix{
SmProj(k) \ar[r] \ar[d]_{Chow} & Sm(k) \ar[d]^{M\gm} \\
Chow\eff(k) \ar[r] & DM\gm\eff(k)
} \]
Furthermore, the induced functors
\begin{equation} \label{equa:effChowEmbedding}
Chow\eff(k)\ozpi \to DM\gm\eff(k)\ozpi
\end{equation}
and 
\begin{equation} \label{equa:chowEmbedding}
Chow(k)\ozpi \to DM\gm(k)\ozpi
\end{equation}
are fully faithful.
\end{prop}

\begin{proof}
The existance of the functor is already given in \cite[Theorem 2.1.4]{Voev00}. To show that (\ref{equa:effChowEmbedding}) is fully faithful, it suffices to show that for each pair $X, Y$ if smooth projective varieties, the morphism
\[ \hom_{Chow\eff(k)\ozpi}(X, Y) \to \hom_{DM\eff(k, \zpi)}(X, Y) \]
is an isomorphism. This is precisely what Proposition~\ref{prop:V423}(\ref{prop:V423:comparisonChowGroups}) says. To show that (\ref{equa:chowEmbedding}) is an isomorphism, it now suffices to recall that
\[ Chow\eff(k) \to Chow(k) \]
is fully faithful, and so is
\[ DM\eff(k) \to DM(k) \]
(Theorem~\ref{theo:cancellationTheorem}).
\end{proof}


\section{Higher Chow Groups and {\'e}‰tale Cohomology --- After Suslin} \label{sec:higChoGro}

Our goal in this section is to relate Bloch's higher Chow groups of varieties over an perfect field to {\'e}tale cohomology. We follow Suslin's article \cite{Sus00} very closely, but we replace the theorem of Voevodsky \cite[Theorem 3.1]{Sus00} he cites with a $\zpi$-version that uses Gabber's theorem on alterations instead of resolution of singularities.

For the rest we follow his strategy to the letter. In Section 1 and  Section 2 of \cite{Sus00} Suslin shows that the higher Chow groups of an affine equidimensional separated scheme of finite type over a field can be calculated using equidimensional cycles. This is valid with integral coefficients and no restrictions on the base field. In Section 3 he generalises this, showing that the higher Chow groups of any quasi-projective scheme $X$ of characteristic zero can be calculated as the Suslin homology of the presheaves $z_{equi}(X/Spec(k), -)$. This is proven using induction on the dimension, a localisation long exact sequence, the result for affine varieties, and the theorem \cite[Theorem 3.1]{Sus00} of Voevodsky that we will replace.

Voevodsky's theorem assumes resolution of singularities, and this is the only place Suslin's proof assumes the base field is of characteristic zero. Replacing this with our theorem that uses Gabber's theorem on alterations, permits us to have this result in characteristic $p$ if we use $\zpi$-coefficients. That is, the higher Chow groups with $\zpi$-coefficients of any quasi-projective scheme $X$ of characteristic $p$ can be calculated as the Suslin homology of the presheaves $z_{equi}(X/Spec(k), -)\ozpi$.

In Section 4 of \cite{Sus00} Suslin goes on to use the main result of \cite{SV96} to show that the higher Chow groups of codimension $d = \dim X$ are dual to $\Ext^*_{qfh}(z_{equi}(X/Spec(k), 0), \ZZ/m)$ if $X$ is an equidimensional quasi-projective variety over an algebraically closed field $k$ of characteristic zero. Having removed the reliance on resolution of singularities, we now have this result over algebraically closed fields of positive characteristic when $m$ is prime to the characteristic of the field. This latter implies the following theorem.

\begin{theo} \label{theo:suslin}
Let $X$ be an equidimensional quasi-projective scheme over an algebraically closed field $k$. Let $i \geq d = \dim X$ and suppose that $m$ is prime to the characteristic of $k$. Then
\[ CH^i(X, n; \ZZ/m) \cong H_c^{2(d - i) + n}(X, \ZZ/m(d - i))^\# \]
where $H_c$ is {\'e}tale cohomology with compact supports. If the scheme $X$ is smooth then this formula simplifies to $CH^i(X, n; \ZZ/m) \cong H_{\acute{e}t}^{2i - n}(X, \ZZ/m(i))$.
\end{theo}

In what follows we reproduce the argument used in Section 3 of \cite{Sus00}, with the appropriate adjustments, fixing some small mistakes along the way.

We make a final remark. We have claimed that Voevodsky's theorem \cite[Theorem 3.1]{Sus00} is the only place that Suslin assumes resolution of singularities. This is not strictly true, as \cite{SV96}, published in 1996, assumes resolution of singularities. However, de Jong's theorem on alterations \cite{deJ96}, published that same year, is sufficient for the purposes of \cite{SV96}. See \cite{Gei00} for a discussion of this fact. 

Denote by $\Delta^n$ the linear subvarieties of $\AA^{n + 1}$ given by the equation $t_0 + \dots + t_n = 1$. Any non-decreasing morphism $\phi: \{0, \dots, n\} \to \{0, \dots, m\}$ induces a canonical morphism $\Delta^n \to \Delta^m$ and these morphisms give $\Delta^\bullet$ the structure of a cosimplicial scheme. If $\phi$ is injective, the image of the corresponding morphism $\Delta^n \to \Delta^m$ is called a face. 

Suppose that $X \in Sch(k)$ is equidimensional. Let $z^i(X, n)$ denote the free abelian group generated by codimension $i$ subvarieties $V \subset X \times \Delta^n$ which intersect $X \times \Delta^m$ properly for every face $\Delta^m \to \Delta^n$. Using a suitable definition of intersection, as outlined at the beginning of Section 2 of \cite{Sus00}, we obtain the structure of a simplicial abelian group on $z^i(X, -)$ for each $i$. The $n$th homotopy group of this simplicial abelian group is denoted $CH^i(X, n)$. These groups were introduced in \cite{Blo86}.

Now suppose that $i \leq d = \dim X$. Denote by $z^i_{equi}(X, n)$ the free abelian group generated by the closed subvarieties $V$ in $X \times \Delta^n$ such that the projection $V \to \Delta^i$ is an equidimensional morphism of relative dimension $d - i$. It can be shown that $z^i_{equi}(X, n)$ is a subgroup of $z^i(X, n)$ -- see the discussion in \cite{Sus00} before \cite[Theorem 2.1]{Sus00}. We have the following theorem.

\begin{theo}[{\cite[Theorem 2.1]{Sus00}}] \label{theo:sus21}
Assume that $X$ is an affine equidimensional scheme and $i \leq d = \dim X$. Then the embedding of complexes $z^i_{equi}(X, n) \hookrightarrow z^i(X, n)$ is a quasi-isomorphism.
\end{theo}

Recall that to a presheaf $F$ on $Sch(k)$ we associate the presheaves $\underline{C}_n(F)(-) = F(\Delta^n \times -)$. The face maps give the $\underline{C}_*(F)(-)$ the structure of a simplicial abelian presheaf, and taking the alternating sums of the face morphisms, the presheaves $\underline{C}_*(F)$ gain the structure of a complex of presheaves. It is immediate that if $F$ is a Nisnevich sheaf, then so are the $\underline{C}_n(F)$. The (co)homology groups of the complexes $\underline{C}_*(F), \underline{C}_*(F) \otimes^L \ZZ / m$, and $RHom(\underline{C}_*(F), \ZZ / m)$ (i.e., the (co)homology sheaves evaluated on the base field) are written as $H_*^\sing(F), H_*^\sing(F, \ZZ/m)$, and $H^*_\sing(F, \ZZ / m)$ respectively.

We replace \cite[Theorem 3.1]{Sus00} with the following theorem, which is an immediate consequence of Theorem~\ref{theo:MainResultRos}.

\begin{theo} \label{theo:ros}
Let $F$ be a presheaf with transfers on $Sch(k)$. If $F_{cdh} \otimes \zpi = 0$ then $H_*^\sing(F) \otimes \zpi = 0$.
\end{theo}

Recall that by the definition of the presheaves $z_{equi}(-/Spec(k), n)$ and the functor $\underline{C}_*(-)$, each presheaf $\underline{C}_n(z_{equi}(X/Spec(k), d - i))$ is a subgroup of $z^i_{equi}(X, n)$.

\begin{theo}[{cf. \cite[Theorem 3.2]{Sus00}}]
Let $k$ be a perfect field of characteristic $p$, let $X \in Sch(k)$ be an equidimensional quasi-projective scheme, and let $i \leq d = \dim X$. Then the composition 
\[ \underline{C}_*(z_{equi}(X/Spec(k), d - i)) \hookrightarrow z^i_{equi}(X, -) \hookrightarrow z^i(X, -) \]
is a quasi-isomorphism.
\end{theo}

\begin{proof}
The proof is by noetherian induction. Clearly none of the presheaves in question see nilpotents and so we can assume that $X$ is reduced. We can find an open affine subscheme $U \to X$ which is regular. Let $Y$ be a closed complement. The sequence of presheaves
\[ 0 \to z_{equi}(Y/Spec(k), d - i) \to z_{equi}(X/Spec(k), d - i) \to z_{equi}(U/Spec(k), d - i) \]
is exact and the cdh sheaf associated to the quotient
\[ z_{equi}(U/Spec(k), d - i) / z_{equi}(X/Spec(k), d - i) \]
is trivial. This is because the cdh sheafifications of the $z_{equi}(-/Spec(k), d - i)$ are isomorphic to $z(-/Spec(k), d - i)$ (\cite[Theorem 4.2.9]{SV00}) and the corresponding sequence for the $z(-/Spec(k), d - i)$ is exact \cite[Theorem 4.3.1]{SV00}. Thus, applying $\underline{C}_*(-)$ and using Theorem~\ref{theo:ros} we get a long exact sequence associated to the homology groups of the $\underline{C}_*(z_{equi}(-/Spec(k), d - i))$. The inclusion $\underline{C}_*(z_{equi}(X/Spec(k), d - i)) \hookrightarrow z^i_{equi}(X, -) \hookrightarrow z^i(X, -)$ gives a morphism between this long exact sequence and the localisation sequence for higher Chow groups \cite[Theorem 3.1]{Blo86}. By the inductive hypothesis, this is an isomorphism on the terms containing $Y$. On the terms containing $U$, it is an isomorphism because $\underline{C}_*(z_{equi}(U/Spec(k), d - i)) \hookrightarrow z^i_{equi}(U, -)$ is an equality for regular $U$ (\cite[Corollary 3.4.5]{SV00}) and Theorem~\ref{theo:sus21} says that $z^i_{equi}(X, -) \hookrightarrow z^i(X, -)$ is a quasi-isomorphism for affine equidimensional $U$. Hence, by the five lemma, the morphism of long exact sequences is an isomorphism.
\end{proof}

\section{Vanishing of negative $K$-theory} \label{sec:kTheory}

In \cite[2.9]{Wei80} Weibel asks if $K_n(X) = 0$ for $n < - \dim X$ for every noetherian scheme $X$ where $K_n$ is the $K$-theory of Bass-Thomason-Trobaugh. This question was answered in the affirmative in \cite{CHSW} for schemes essentially of finite type over a field of characteristic zero. Assuming strong resolution of singularities, it is also answered in the affirmative in \cite{GH10} for schemes essentially of finite type over a field of positive characteristic. Both of these proofs compare $K$-theory with cyclic homology, and then use a cdh descent argument.

In this section we will give a partial answer to Weibel's conjecture. The proof is actually very short, and uses almost none of the machinery we have developed. Its key is a theorem of Cisinski which says that the homotopy invariant $K$-theory presheaf of $S^1$-spectra $KH$ satisfies cdh descent (\cite[3.7]{Cis13}). To prove this, he proves that $KH$ is representable in the Morel-Voevodsky stable homotopy category, and then applies Ayoub's projective base change result \cite{Ayo07} that we reproduced as Theorem~\ref{theo:ayo}(4a). It is perhaps needless to say that Voevodsky's reduction of cdh descent to the statement that appropriate squares are homotopy cartesian is used as well.

It is possible to give a self-contained account of the proof (that is, without referring to the previous sections) in a few pages. The references we make to earlier results are Proposition~\ref{prop:ktheoryTraces}, Lemma~\ref{lemm:tracesImpliesTraces}, Proposition~\ref{prop:compatSheaf}, Example~\ref{exam:refinable}, and Corollary~\ref{coro:regularlCover}).

Proposition~\ref{prop:ktheoryTraces} says that $K$-theory has a structure of object with traces but we actually only use (Deg), and this is straight-forward. Remark~\ref{rema:tracesImpliesTraces} says that traces on an object in $\SH$ imply traces on its homotopy presheaves but we don't need this if we show directly that the groups $KH_n$ have trace morphisms which satisfy (Deg). Proposition~\ref{prop:compatSheaf} is an elementary property of refinable Grothendieck pretopologies, Example~\ref{exam:refinable} notices that the $\ldh$-pretopology is $\stackrel{\fpsl}{\to} \stackrel{cdh}{\to}$ refinable, and Corollary~\ref{coro:regularlCover} is a statement of Gabber's theorem on alterations.

\begin{theo} \label{theo:KtheoryVanishing}
Let $X$ be a quasi-separated quasi-excellent noetherian scheme and $p$ a prime that is nilpotent on $X$. Then $K^B_n(X) \otimes \zpi = 0$ for $n < - \dim X$.
\end{theo}

\begin{proof}
Since $p$ is nilpotent on $X$ the canonical morphism $K^B_n \otimes \zpi \to KH_n \otimes \zpi$ is an isomorphism \cite[9.6]{TT90}. Hence it suffices to prove that \mbox{$KH_n(X) \otimes \zll$} vanishes for every prime $\ell \neq p$ and $n < - \dim X$. Since $KH$ satisfies cdh descent (\cite[3.7]{Cis13}) we have a spectral sequence (cf. \cite[1.36]{Tho85}, \cite[Corollary 5.2]{Wei89})
\[ E_2^{p, q} = H^{p}_{cdh}(X, KH_{-q}(-)_{cdh}) \implies KH_{- p - q}(X) \]
which converges due to the cdh cohomological dimension being bounded by $\dim X$ \cite[12.5]{SV00}. Furthermore, we see that the $E_2$ sheet is zero outside of $0 \leq p \leq \dim X$. We tensor this spectral spectral sequence with $\zll$ to obtain a second spectral sequence 
\[ E_2^{p, q} = H^{p}_{cdh}(X, KH_{-q}(-)_{cdh}) \otimes \zll \implies KH_{- p - q}(X) \otimes \zll \]
Due to the vanishing of $E_2$ terms already mentioned, we have reduced to showing that $KH_{-q}(-)_{cdh} \otimes \zll = 0$ for $q > 0$. The presheaves $KH_{-q}(-)$ have a structure of traces (Proposition~\ref{prop:ktheoryTraces}, Lemma~\ref{lemm:tracesImpliesTraces}, Remark~\ref{rema:tracesImpliesTraces}) and so their cdh associated sheaves are $\ldh$ separated (Lemma~\ref{lemm:fpslAcyclic}, Proposition~\ref{prop:compatSheaf}, Example~\ref{exam:refinable}). Now every scheme admits an $\ldh$ cover $\{U_i \to X \}$ with $U_i$ regular (Corollary~\ref{coro:regularlCover}) and $KH_{-q}(U) \otimes \zpi$ vanishes for $U$ regular and $q > 0$ (\cite[Proposition 6.8]{TT90}). It follows that $KH_{-q}(-)_{cdh} \otimes \zll = 0$ for $q > 0$. That is, the $E_2$ terms of the spectral sequence vanish unless $q \leq 0$ and $0 \leq p \leq \dim X$. This implies that $KH_n(X) \otimes \zll = 0$ for $n < - \dim X$.
\end{proof}

\appendix
\chapter{Appendix}
\section{Some local algebra}

\begin{lemm} \label{lemm:degLength} \label{lemm:appendix1}
Suppose that $\phi: A \to B$ is a finite flat $A$ algebra with $A$ and $B$ Artinian local rings. Then
\[ \length A \cdot \deg \phi = \length B \cdot [B / \m_B : A / \m_A]. \]
\end{lemm}

\begin{proof}
Both sides are equal to $\length_A B$. Alternatively (and more explicitely), define $L = A / \m_A$. Every indecomposable $A$-module is isomorphic to $L$. So since $A \to B$ is flat, applying $- \otimes_A B$ to a composition series for the $A$-module $A$ shows that
\[ \length_B B = \length_A A \length_B (L \otimes_A B), \]
since $\length_B$ is additive. Now $L \otimes_A B$ is a finite local $L$-algebra with residue field $B / \m_B$, and so we have
\[ \dim_{L}(L \otimes_A B) = \length_{L \otimes_A B} (L \otimes_A B) \cdot [B / \m_B : L]. \]
Finally, recognising that $\length_B (L \otimes_A B) = \length_{L \otimes_A B} (L \otimes_A B)$ and putting everything together, we find the desired equality.
\end{proof}

\begin{lemm} \label{lemm:fonalgebra}
Suppose that $k \stackrel{\phi}{\to} A \stackrel{\psi}{\to} B$ are finite flat morphisms of local rings with $k$ a field. Then
\[ \length (B / \m_A B) \cdot \length A = \length B. \]
\end{lemm}

\begin{proof}
We apply Lemma~\ref{lemm:degLength} to the three morphisms $\phi$, $\psi\phi$, and $A / \m_A \to B / \m_A B$ to obtain the following equalities:
\[ \deg \phi = \length A \cdot [A / \m_A : k] \]
\[ \deg \phi \deg \psi = \length B \cdot [B / \m_B : k] \]
\[  \deg \psi = \length (B / \m_A B) \cdot [B / \m_B : A / \m_A] \]
Multiplying the first and last together and comparing with the second gives the desired result.
\end{proof}

\begin{lemm} \label{lemm:CDBalgebra} \label{lemm:appendix3}
Let $K / k$ be a field extension, $A$ a finite local $k$ algebra. Let $Spec(B_i)$ be a connected component of $Spec(A \otimes_k K)$ and $Spec(C_i)$ the corresponding connected component of $K \otimes_k (A / \m_A)$. Then we have
\[ \length B_i = \length C_i \cdot \length A. \]
\end{lemm}

\begin{proof}
We have the following cartesian square
\[ \xymatrix{
Spec(C_i) \ar[r] \ar[d] & Spec(A / \m_A) \ar[d] \\
Spec(B_i) \ar[r] & Spec(A)
} \]
with the lower, and hence upper, horizontal morphism being flat. Every indecomposable $A$-module is isomorphic to $A / \m_A$ and so applying $B_i \otimes_A -$ to a composition series for $A$ we find that $\length_{B_i} B_i = \length_A A \cdot \length_{B_i} (B_i \otimes_A (A / \m_A))$. But $\length_{B_i} (B_i \otimes_A (A / \m_A)) = \length_{C_i} C_i$, hence the result.
\end{proof}

\begin{lemm} \label{lemm:triPullbackPoint}
Suppose that we have a finite set of commutative triangles
\[ \xymatrix{
Y'_k \ar[rr]^{h_k}  \ar[dr]_{g_k} && Y \ar[dl]^f \\
& X
} \]
($k = 1, \dots, n$ say) such that 
\begin{enumerate}
 \item all of $f_k, g_k$ are finite flat surjective,
 \item each $Y'_k$ is generically reduced,
 \item there exists an integer $d$ such that for each generic point $\eta$ of $Y$ we have $d = \sum \deg (h_k \times_Y \eta)$,
 \item there exist integers $m_k$ such that for each generic point $\eta \in h_k(Y'_k)$ we have $m_k = \length \OO_{Y, \eta}$.
\end{enumerate}
Now suppose that $x$ is a point in $X$, let $y$ be a point of $Y$ over $x$, and let $y_{k\ell}'$ be the points of $Y'_k$ over $y$. Then we have
\[ d \cdot \length \OO_{x \times_X Y, y} = \sum_{k} m_{k} \sum_{\ell} \length \OO_{x \times_X Y', y_{k\ell'}} [k(y_{k\ell}'):k(y)]. \]
\end{lemm}

\begin{rema} \label{rema:triPullbackPoint}
One of the applications of this lemma is to Proposition~\ref{prop:sepImpliesTrin}. In this case, it will be applied to triangles arising from the axioms (Tri1) and (Tri2), both of which are included in the hypotheses of the lemma, and both of which are easier cases. We have combined them because they are both special cases of this more general result. In the case of (Tri1) we have $d = 1$ and the equality becomes
\[ \length \OO_{x \times_X Y, y} = \sum_{k} m_{k} \sum_{\ell} \length \OO_{x \times_X Y', y_{k\ell'}} [k(y_{k\ell}'):k(y)]. \]
In the case of (Tri2) there is a unique $Y'_k$, and $m_k = 1$ so the equality is
\[ d \cdot \length \OO_{x \times_X Y, y} = \sum_{\ell} \length \OO_{x \times_X Y', y_{k\ell'}} [k(y_{k\ell}'):k(y)]. \]
\end{rema}

\begin{proof}
We first consider the case where there is a unique point $y \in Y$ over $x \in X$. We begin with some identities. Choose a generic point $\xi \in X$. We let $\eta_i$ be the generic points of $Y$ that are over $\xi$, and $\eta_{ij}$ the generic points of $\amalg Y'_k$ that are over $\eta_i$ (we don't care which $Y'_k$ they belong to). For each $i$ we let $k_i$ be an index such that $\eta_i \in h_{k_i}(Y_{k_i}')$ so we have $m_{k_i} = \length \OO_{Y, \eta_i}$. We claim that
\begin{equation} \label{equa:triRelationGeneric}
d \cdot \deg f = \sum_k m_k \cdot \deg g_k.
\end{equation}
Notice that since $Y'$ is generically reduced, so is $X$ and consequently, the multiplicity of a generic point of $Y$ in $Y$ is the same as its multiplicity in $Y \times_X \xi$.

Consider the base change by the chosen generic point $\xi \to X$. Due to the hypotheses in the statement, we have
\[ \begin{split}
d \cdot \deg f &= d \cdot \deg (f \times_X \xi) \stackrel{(\ref{lemm:degLength})}{=} d \sum_i [k(\eta_i) : k(\xi) ] m_{k_i} \\
&\stackrel{(\mathrm{Hyp 2}, \mathrm{Hyp 3})}{=} \sum_{ij} [k(\eta_{ij}) : k(\eta_i) ][k(\eta_i) : k(\xi) ] m_{k_i} \\
&= \sum_{ij} [k(\eta_{ij}') : k(\xi) ] m_{k_i} \stackrel{}{=} \sum_k \deg (g_k \times_X \xi) m_k = \sum_k \deg g_k m_k.
\end{split} \]
Hence, our claim (\ref{equa:triRelationGeneric}) is proven.

We now use the same ideas for the point $x$ in $X$. Let $n = \length \OO_{Y \times_X x, y}$, and $r_{k\ell} = \length \OO_{Y' \times_X x, y_{k\ell}'}$.
Using the identity (\ref{equa:triRelationGeneric}) we find that we have
\[ \begin{split}
d \cdot n \cdot [k(y):k(x)] &\stackrel{(\ref{lemm:degLength})}{=} d \cdot \deg\ f
\stackrel{(Equ.({\footnotesize\ref{equa:triRelationGeneric}}))}{=} \sum_k m_k \deg\ g_k \\
&\stackrel{(\ref{lemm:degLength})}{=} \sum_k m_k \sum_{\ell} r_{k\ell} [k(y_{k\ell}') : k(x)] \\
&\stackrel{}{=} \sum_k m_k \sum_{\ell} r_{k\ell} [k(y_{k\ell}') : k(y)][k(y):k(x)]
\end{split} \]
and cancelling $[k(y) : k(x)]$ from either side gives the desired equality.

Now we remove the hypothesis that there is a unique point $y$ over $x$ and consider the general case. The hypotheses on the triangle in the statement are preserved by change of base by {\'e}tale morphisms $X' \to X$ (Lemma~\ref{lemm:baseChangeTrianglePreserved}) and consequently, by base change by the henselisation $\ ^hx$ at the point $x$. After base change we find that in each connected component of $\ ^hx \times_X Y$ there is a unique point of $\ ^hx \times_X Y$ over the closed point $x$ of $\ ^hx$. Since it suffices to verify the formula on each connected component of $Y$, we are done.
\end{proof}

\begin{lemm} \label{lemm:baseChangeTrianglePreserved}
Suppose that (\ref{equa:tri}) is a triangle that satisfies the hypotheses of Lemma~\ref{lemm:triPullbackPoint} and $X' \to X$ is an {\'e}tale morphism. Then
\[ \xymatrix{
Y' \times_X X' \ar[rr] \ar[dr] && Y \times_X X' \ar[dl] \\
& X'
} \]
also satisfies the hypotheses of Lemma~\ref{lemm:triPullbackPoint}.
\end{lemm}

\begin{proof}
Let $f', g', h'$ be the pullbacks of $f, g, $ and $h$ respectively. Clear the first hypothesis is preserved by base change and so we only need concern ourselves with the other three. These all take place within the generic fibers, and so we can assume $X$ (and hence everything else as well) is of dimension zero.

$Y' \times_X X'$ is generically reduced because base change of a field extension by a separable field extension does not introduce nilpotents. The third hypothesis is satisfied because degree is preserved by pullback. Finally the forth hypothesis is satisfied due to Lemma~\ref{lemm:CDBalgebra}: let $K / k$ be $X' \to X$, $A$ the local ring of a generic point of $Y$ and note that since $K / k$ is {\'e}tale, $A_{red} \otimes_k K$ is reduced and so $\length C_i = 1$.
\end{proof}

\section{Inverting integers in triangulated categories} \label{sec:invertTriangle}


In this section we discuss some ways to invert integers in triangulated categories. The goal is to be able to make analogues of familiar statements such as ``an element $a$ of an abelian group $A$ is zero if and only if its image in $A \otimes \zll$ is zero for every prime $\ell$''.

We present two analogues of $-\otimes \zll$ for triangulated categories. The first is to tensor all the hom groups with $\zll$. We show that this category has a canonical structure of a triangulated category (Proposition~\ref{prop:minimalStructureTriangulated}), however the canonical functor $\T \to \T \otimes \zll$ does not behave well with respect to sums and so we don't have the necessary adjunctions. The second is to consider the triangulated subcategory of $\zll$-local objects. For nice triangulated categories this is a localisation (Proposition~\ref{prop:localisationTriLambda}). We compare these two constructions in Corollary~\ref{coro:compareTwoLocalisations}. We then state the ``local-global'' principles that we need (Lemma~\ref{lemm:conservativeLocal}, Lemma~\ref{lemm:generatingLocalTri}). The first says that we can detect isomorphisms via these localisations, and the second says that we can detect compact generating families via these localisations.

This section is a little notation heavy and so we make the following summary.

\begin{enumerate}
 \item $\T \otimes \Lambda$ (Definition~\ref{defi:tensorNaive}) The category obtained by applying $- \otimes \Lambda$ to each hom group.
 \item $\Lambda$-local (Definition~\ref{defi:localObject}).
 \item $\T[S^{-1}]$ (Definition~\ref{defi:localObject}) The subcategory of $\ZZ[S^{-1}]$ local objects.
 \item ${_S\T}$ (Definition~\ref{defi:localObject}) The left orthogonal to $\T[S^{-1}]$.
 \item $(-)[S^{-1}], {_S(-)}$ (Proposition~\ref{prop:localisationTriLambda}) Localisation functors.
 \item $\T_c$ (Definition~\ref{defi:compactObjects}) The full subcategory of compact objects.
 \item $S_{(\ell)} \subset \ZZ$ (Lemma~\ref{lemm:conservativeLocal}) the set of integers coprime to $\ell$.
\end{enumerate}

We also recall for reference the following standard terms from the theory of triangulated categories.

\begin{enumerate}
 \item A \emph{compact} object is an object $F$ such that $\hom(F, \oplus E_\lambda) \cong \oplus \hom(E, F_\lambda)$ for any family of objects whose sum exists.
 \item A \emph{thick} triangulated subcategory is a triangulated subcategory that is closed under direct summands.
 \item A \emph{localising} subcategory of a triangulated category admitting small sums is a triangulated subcategory closed under small sums.
 \item A \emph{compactly generated} triangulated category with small sums is a triangulated category with small sums which itself is the smallest localising subcategory containing its subcategory of compact objects.
\end{enumerate}

\begin{defi} \label{defi:tensorNaive}
Suppose $\Lambda$ is a ring and $\A$ an additive category. Define $\A \otimes \Lambda$ to be the (a postiori additive) category which has the same objects as $\A$ and $\hom_{\A \otimes \Lambda}(A, B) = \hom_{\A}(A, B) \otimes \Lambda$.
\end{defi}

\begin{exam}
Let $\mathcal{A}b$ be the category of abelian groups, and suppose that $\mathcal{A}b \to \mathcal{A}b \otimes \QQ$ preserves sums. Then we would have
\[ \begin{split}
(\prod \ZZ / p) \otimes \QQ &= (\prod \hom_{\mathcal{A}b}(\ZZ / p, \QQ / \ZZ)) \otimes \QQ \\
&= \hom_{\mathcal{A}b}(\oplus \ZZ / p, \QQ / \ZZ) \otimes \QQ \\
&= \hom_{\mathcal{A}b \otimes \QQ}(\oplus \ZZ / p, \QQ / \ZZ) \\
&= \prod \hom_{\mathcal{A}b \otimes \QQ}( \ZZ / p, \QQ / \ZZ) \\
&= \prod (\ZZ / p \otimes \QQ).
\end{split} \]
However, the last groups is zero and the first group is not.

A similar problem is exhibited with the derived category of abelian groups.
\end{exam}

\begin{prop} \label{prop:minimalStructureTriangulated}
Suppose that $\T$ is a triangulated category and $\Lambda \subseteq \QQ$ a subring of the rationals. Then there is a ``minimal'' structure of triangulated category on $\T \otimes \Lambda$ such that the canonical functor $F: \T \to \T \otimes \Lambda$ is exact.
\end{prop}

\begin{rema}
We mean minimal in the sense that given any other structure of triangulated category on $\T \otimes \Lambda$ such that $F: \T \to \T \otimes \Lambda$ is exact, the class of distinguished triangles contains this ``minimal'' class of distinguished triangles.
\end{rema}

\begin{rema}
This is a consequence of every morphism in $\T \otimes \Lambda$ being isomorphic to a morphism in the image of $\T \to \T \otimes \Lambda$ in some kind of nice way. If the category $\T$ happens to be $R$-linear for a ring $R$ then the result holds for any localisation $R[S^{-1}]$ of $R$. The proof does not seem to easily generalise to other kinds of $R$-algebras.
\end{rema}

\begin{proof}
The proof follows via the principle that the diagrams in $\T \otimes \Lambda$ in question are isomorphic to diagrams in the image of $\T$. We will elaborate. Let $S \subset \ZZ$ be the set of integers which are invertible in $\Lambda$, so that $\Lambda = \ZZ[S^{-1}]$.

Clearly, any triangle of $\T \otimes \Lambda$ that is isomorphic to a distinguished triangle in the image of $F$ is necassarily distinguished. We will show that this class of triangles satisfies the axioms for a triangulated category. 

\emph{TR0.} The triangle $X \stackrel{id}{\to} X \to 0 \to X[1]$ is in the image of $F$.

\emph{TR1.} Suppose that $X \stackrel{s^{-1}f}{\to} Y$ is a morphism of $\T \otimes \Lambda$ (where $s \in S$). Completing $f: X \to Y$ to a distinguished triangle in $\T$ we find a commutative diagram
\[ \xymatrix{
X \ar@{=}[d] \ar[r]^{s^{-1}f} & Y \ar[d]^{s \cdot id_Y} \ar@{-->}[r]^{s \cdot g} & Z \ar@{=}[d] \ar@{-->}[r]^h & X[1] \ar@{=}[d] \\
X \ar[r]_f & Y \ar[r]_g & Z \ar[r]_h & X[1]
} \]
where the vertical morphisms are isomorphisms in $\T \otimes \Lambda$ (the inverse to $s \cdot id$ is $s^{-1} \cdot id$. Hence, the upper triangle is a distinguished triangle in $\T \otimes \Lambda$.

\emph{TR2.} It is clear from our definition of distinguished triangles in $\T \otimes \Lambda$ that the class is closed under rotation.

\emph{TR3.} Since every distinguished triangle of $\T \otimes \Lambda$ is isomorphic to a triangle that is the image of a distinguished triangle in $\T$ it suffices to consider the case when the two triangles in question are in the image of $F$ (see the following diagram).
\begin{equation} \label{equ1}
\xymatrix{
\ar[r] \ar[d]^\cong & \ar[r] \ar[d]^\cong & \ar[r] \ar[d]^\cong & \ar[d]^\cong & \not \in Im\ F\\
\ar[r] \ar[d] & \ar[r] \ar[d] & \ar[r] \ar@{-->}[d] & \ar[d] & \in Im\ F \\
\ar[r] \ar[d]^\cong & \ar[r] \ar[d]^\cong & \ar[r] \ar[d]^\cong & \ar[d]^\cong & \in Im\ F \\
\ar[r] & \ar[r] & \ar[r] & & \not \in Im\ F
}
\end{equation}
Suppose we are given such a diagram as follows (where $s, t \in S$).
\[ \xymatrix{
X \ar[r]^f \ar[d]^{s^{-1}a} & X' \ar[r]^{f'} \ar[d]^{t^{-1}b} & X'' \ar[r]^{f''} & X[1] \ar[d]^{s^{-1}a[1]} \\
Y \ar[r]^g  & Y' \ar[r]^{g'}  & Y'' \ar[r]^{g''} & Y[1] \\
} \]
We use (TR3) in $\T$  to find the dashed morphism $c$ in $\T$ making the following diagram commute and then compose with $(st)^{-1} \cdot id$ in $\T \otimes \Lambda$.
\begin{equation} \label{equ2}
\xymatrix{
X \ar[r]^f \ar[d]^{t \cdot a} & X' \ar[d]^{s \cdot b} \ar[r]^{f'} & X'' \ar@{-->}[d]^c \ar[r]^{f''} & X[1] \ar@{=}[d]^{t \cdot a[1]} \\
Y \ar[r]^{g} \ar[d]^{(st)^{-1} \cdot id} & Y' \ar[r]^{g'} \ar[d]^{(st)^{-1} \cdot id} & Y'' \ar[d]^{(st)^{-1} \cdot id} \ar[r]^{g''} & Y[1] \ar[d]^{(st)^{-1} \cdot id[1]} \\
Y \ar[r]^{g} & Y' \ar[r]^{g'}  & Y'' \ar[r]^{g''} & Y[1]  \\
}
\end{equation}

\emph{TR4'.} Instead of proving the octohedral axiom (TR4) we will prove (TR4') (\cite[Definition 1.3.13]{Nee01}). This is equivalent to (TR4) (\cite[Remark 1.4.7]{Nee01}) but is sometimes easier to work with. The statement of (TR4') is as follows: Given a diagram
\[ \xymatrix{
X \ar[r]^u \ar[d]_f & Y \ar[r]^v \ar[d]_g & Z \ar[r]^w & X[1] \ar[d]_{f[1]} \\
X' \ar[r]_{u'} & Y' \ar[r]_{v'} & Z' \ar[r]_{w'} & X'[1]
} \]
where the rows are distinguished triangles, the morphism $h: Z \to Z'$ given by (TR3) may be chosen such that the mapping cone $MapCone(f, g, h)$ is a distinguished triangle. The mapping cone $MapCone(f, g, h)$ is by definition (\cite[Definition 1.3.1]{Nee01}) the diagram
\[ \xymatrix{
Y \oplus X' \ar[r] \ar@<1ex>@{}[r]^{\left ( \begin{array}{cc} -v & 0 \\ g & u' \end{array} \right ) }
 & Z \oplus Y' \ar[r] \ar@<1ex>@{}[r]^{\left ( \begin{array}{cc} -w & 0 \\ h & v' \end{array} \right ) }
 & X[1] \oplus Z' \ar[r] \ar@<1ex>@{}[r]^{\left ( \begin{array}{cc} -u[1] & 0 \\ f[1] & w' \end{array} \right ) }
 & Y[1] \oplus X'[1] 
} \]

By (TR4') in $\T$, the morphism $c$ in the diagram (\ref{equ2}) above may be chosen so that the mapping cone of $(t\cdot a, s \cdot b, c)$ is distinguished in $\T$. By the following lemma, the isomorphism in Diagram~\ref{equ2}, and the remark preceding Diagram~\ref{equ1} this is enough to conclude that we have (TR4') in $\T \otimes \Lambda$. 

\end{proof}

\begin{lemm}
Consider a composable pair of morphisms of triangles
\[ \xymatrix{
X \ar[r]^u \ar[d]^f & Y \ar[r]^v \ar[d]^g & Z \ar[r]^w \ar[d]^h & X[1] \ar[d]^{f[1]} \\
X' \ar[r]^{u'} \ar[d]^{f'} & Y' \ar[r]^{v'} \ar[d]^{g'} & Z' \ar[r]^{w'} \ar[d]^{h'} & X'[1] \ar[d]^{f'[1]} \\
X'' \ar[r]^{u''}   & Y'' \ar[r]^{v''}   & Z'' \ar[r]^{w''}   & X''[1] \\
} \] 
Then the obvious potential morphism $(f', g', h')$ (resp. $(f, g, h)$) between the mapping cones (\cite[Definition 1.3.1]{Nee01}) $MapCone(f, g, h) \to MapCone(f'f, g'g, h'h)$ (resp. $MapCone(f'f, g'g, h'h) \to MapCone(f, g, h)$) is actually a morphism of triangles. That is, all the appropriate squares commute.

In particular, if $(f', g', h')$ (resp. $(f, g, h)$) is an isomorphism of triangles, then there exists an isomorphism $MapCone(f, g, h) \cong MapCone(f'f, g'g, h'h)$ (resp. $MapCone(f'f, g'g, h'h) \cong MapCone(f, g, h)$).
\end{lemm}

\begin{proof}
We write down the first square in each case. The following is the first square for the map $MapCone(f, g, h) \to MapCone(f'f, g'g, h'h)$.
\[ \xymatrix@!=42pt{
Y \oplus X' \ar[r] \ar@<1ex>@{}[r]^{\left ( \begin{array}{cc} -v & 0 \\ g & u' \end{array} \right ) }
\ar[d] \ar@<-1ex>@{}[d]_{\left ( \begin{array}{cc} id_Y & 0 \\ 0 & f' \end{array} \right ) } & 
Z \oplus Y' \ar[d] \ar@<1ex>@{}[d]^{\left ( \begin{array}{cc} id_Z & 0 \\ 0 & g' \end{array} \right ) } \ar[r] \ar[d] & \ar[r] \ar[d] & \ar[d] \\
Y \oplus X'' \ar[r] \ar@<-1ex>@{}[r]_{\left ( \begin{array}{cc} -v & 0 \\ g'g & u'' \end{array} \right ) } &
Z \oplus Y'' \ar[r] & \ar[r] & 
} \]
The two compostions are equal to 
\[ \left ( \begin{array}{cc} -v & 0 \\ g'g & g'u' \end{array} \right )  \textrm{ and } \left ( \begin{array}{cc} -v & 0 \\ g'g & u''f' \end{array} \right )  \]
and these are equal by the commutativity of the appropriate square in the statement $u''f' = g'u'$.

The first square for the map $MapCone(f'f, g'g, h'h) \to MapCone(f, g, h)$ is
\[ \xymatrix@!=42pt{
Y \oplus X'' \ar[r] \ar@<1ex>@{}[r]^{\left ( \begin{array}{cc} -v & 0 \\ g'g & u'' \end{array} \right ) }
\ar[d] \ar@<-1ex>@{}[d]_{\left ( \begin{array}{cc} g & 0 \\ 0 & id_{X''} \end{array} \right ) } & 
Z \oplus Y'' \ar[d] \ar@<1ex>@{}[d]^{\left ( \begin{array}{cc} h & 0 \\ 0 & id_{Y''} \end{array} \right ) } \ar[r] \ar[d] & \ar[r] \ar[d] & \ar[d] \\
Y' \oplus X'' \ar[r] \ar@<-1ex>@{}[r]_{\left ( \begin{array}{cc} -v' & 0 \\ g' & u'' \end{array} \right ) } &
Z' \oplus Y'' \ar[r] & \ar[r] & 
} \]
The two compostions are equal to 
\[ \left ( \begin{array}{cc} -hv & 0 \\ g'g & u'' \end{array} \right )  \textrm{ and } \left ( \begin{array}{cc} -v'g & 0 \\ g'g & u'' \end{array} \right )  \]
and these are equal by the commutativity of the appropriate square in the statement $hv = v'g$.
\end{proof}

\begin{defi} \label{defi:localObject}
Suppose that $\mathcal{A}$ is an additive category, $E \in \mathcal{A}$ an object, and $\Lambda \subseteq \QQ$ a subring of the rationals. We will say $E$ is \emph{$\Lambda$-local} if for every integer $n \in \ZZ$ which is invertible in $\Lambda$, the morphism $n \cdot \id_E$ is an isomorphism. If $S \subseteq \ZZ$ is a multiplicative system\footnote{i.e., $S$ is closed under multiplication.} then we will denote the full subcategory of $\ZZ[S^{-1}]$ local objects by $\mathcal{A}[S^{-1}]$. It will also make life easier to define
\[ {_S\mathcal{A}} = (\mathcal{A}[S^{-1}])^\perp. \]
That is, ${_S\mathcal{A}}$ is the full subcategory of $\mathcal{A}$ of objects $F$ such that $\hom(F, E) = 0$ for every $E \in \mathcal{A}[S^{-1}]$.
\end{defi}

\begin{prop} \label{prop:localisationTriLambda}
Suppose that $\T$ is a compactly generated triangulated category admitting small sums and $S \subset \ZZ$ is a multiplicative system. Then:
\begin{enumerate}
 \item ${_S\T}$ is the smallest localising subcategory of $\T$ containing the objects $Cone(E \stackrel{n \cdot \id_E}{\to} E)$ for $n \in S$ and $E$ compact in $\T$.
 \item The inclusion $i: \T[S^{-1}] \to \T$ has a left adjoint $(-)[S^{-1}]$ and the inclusion $j: {_S\T \to \T}$ has a right adjoint ${_S(-)}$.
\[ \xymatrix{
{_S\T} \ar@<0.5ex>[r]^j & \T \ar@<0.5ex>[r]^{(-)[S^{-1}]} \ar@<0.5ex>[l]^{{_S(-)}} & \ar@<0.5ex>[l]^i \T[S^{-1}] 
} \]
 \item There is a canonical equivalence $\T[S^{-1}] \cong \T / ({_S\T})$.
 \item Both the functor $(-)[S^{-1}]$ and the inclusion $i: \T[S^{-1}] \to \T$ preserve small sums.
 \item The functor $(-)[S^{-1}]$ preserves compact objects.
 \item If $\G$ is a small generating family of compact objects for $\T$ then
\[ \G[S^{-1}] = \{ E[S^{-1}] : E \in \G \} \]
is a generating family of compact objects for $\T[S^{-1}]$.
 \item Both $\T[S^{-1}]$ and ${_S\T}$ are compactly generated triangulated categories admitting small sums.
\end{enumerate}
\end{prop}

\begin{proof}
Let $A$ denote the smallest localising subcategory of $\T$ containing $Cone(E \stackrel{n \cdot \id_E}{\to} E)$ for every compact object $E \in \T$ and every $n \in S$ which is in $S$. The localisation $\T \to \T / A$ always exists \cite[Theorem 2.1.8]{Nee01} and given our assumptions it has a right adjoint if $A$ is compactly generated (\cite[Proposition 8.4.2]{Nee01} and see also the end of \cite[Remark 8.1.7]{Nee01}). But $A$ is compactly generated by definition. The image of the right adjoint $\T / A \to \T$ is canonically identified with $^\perp A$ \cite[Theorem 9.1.16]{Nee01} where $^\perp A$ is the full subcategory of $\T$ whose objects $E$ satisfy $\hom_\T(F, E) = 0$ for all $F \in A$ \cite[Definition 9.1.10]{Nee01}. As $A$ is compactly generated, for $E$ to belong to $^\perp A$ it is sufficient that $\hom_\T(F, E) = 0$ for all $F$ in a compact generating family of $A$. For example, $F$ of the form $Cone(F' \stackrel{n \cdot \id_{F'}}{\to} F')$ for $n \in S$ and $F'$ compact in $\T$. But since $\T$ is compactly generated, for an object $E$ to satisfy this means precisely that $n \cdot \id_E$ is invertible for every $n \in S$. Hence, $^\perp A = \T[S^{-1}]$.
\begin{enumerate}
 \item It now follows from \cite[Corollary 9.1.14]{Nee01} that $A = (^\perp A)^\perp = {_S\T}$.
 \item This is also \cite[Corollary 9.1.14]{Nee01}.
 \item This too is \cite[Corollary 9.1.14]{Nee01}.
 \item $(-)[S^{-1}]$ preserves small sums by virtue of it having a right adjoint. Similarly, for any adjunction $(L, R)$ if $LR = \id$ then $R$ preserves small sums.
 \item In general, a left adjoint preserves compact objects if its right adjoint preserves small sums.
 \item Again, this is true in general whenever $L$ preserves compact objects and $LR = \id$.
 \item $\T[S^{-1}]$ admits small sums via the presentation as $^\perp A = \T[S^{-1}]$ given above. we have just seen that it is compactly generated. We have seen already that $A = {_S\T}$ is compactly generated and localising.
\end{enumerate}
\end{proof}

\begin{defi} \label{defi:compactObjects}
For $\T$ a triangulated category admitting small sums, let $\T_c$ denote the full subcategory of compact objects.
\end{defi}

We will need the following result of Thomason-Neeman.

\begin{lemm}[{\cite[Lemma 4.4.5]{Nee01}}] \label{lemm:NeeCompact}
Suppose that $\T$ is a triangulated category admitting small sums and $\G$ is a small generating family of compact objects. Then the smallest thick subcategory of $\T$ containing $\G$ is the full subcategory of compact objects.
\end{lemm}

\begin{rema}
This is a special case of \cite[Lemma 4.4.5]{Nee01} where $\alpha = \aleph_0 = \beta$. In this case $\T^\alpha$ is the subcategory of compact objects (see the end of \cite[Example 1.10]{Nee01}). It might also be relevant to point out to the conscientious reader unfamiliar with the theory that [TR5] and $\alpha$-localising are defined after \cite[Remark 1.4]{Nee01} and the notation $\T^\alpha, \langle S \rangle, \langle S \rangle^\alpha$ is defined in \cite[Definition 1.12]{Nee01}.
\end{rema}

We now have the following comparison of our two ways of inverting integers.

\begin{coro} \label{coro:compareTwoLocalisations}
With the assumptions and notation of Proposition~\ref{prop:localisationTriLambda} we have a canonical equivalence of categories $\T_c \otimes \ZZ[S^{-1}] \cong \T_c / ({_S\T})_c$. Consequently, the induced functor
\[ \T_c \otimes \ZZ[S^{-1}] \to (\T[S^{-1}])_c \]
is fully faithful.
\end{coro}

\begin{proof}
Consider the unversal property that $\T_c \to \T_c \otimes \ZZ[S^{-1}]$ satisfies. It is the universal exact functor towards $\ZZ[S^{-1}]$ linear triangulated category. It is also the unversal exact functor whose image is a $\ZZ[S^{-1}]$ linear triangulated category. One more way of saying this, is that it is the unversal exact functor sending each $E \stackrel{n \cdot \id_E}{\to} E$ to an isomorphism where $E \in \T_c$ and $n \in S$. That is, if $B$ is the smallest thick triangulated subcategory of $\T_c$ containing each $Cone(E \stackrel{n \cdot \id_E}{\to} E)$ where $E \in \T_c$ and $n \in S$, then we have a canonical equivalence
\[ \T_c \otimes \ZZ[S^{-1}] \cong \T_c / B. \]
Lemma~\ref{lemm:NeeCompact} tells us that in fact we have $B = ({_S\T})_c$.

We are now trying to show that the canonical functor
\[ \T_c / ({_S\T})_c \to (\T / ({_S\T}))_c \]
is fully faithful. Glancing at the statement of \cite[Theorem 2.1.8]{Nee01} it suffices to show that $({_S\T})_c$ consists precisely of the objects of $\T_c$ which are sent to zero under the canonical functor
\[ \T_c \to \T / ({_S\T}). \]
The objects that this functor sends to zero are precisely the objects in the thick subcategory $\T_c \cap ({_S\T})$. Notably, $Cone(E \stackrel{n \cdot \id_E}{\to} E)$ is in this subcategory for every $E \in \T_c$ and $n \in S$. Since these compactly generate ${_S\T}$ it follows that $\T_c \cap ({_S\T}) = ({_S\T})_c$, hence the result by Lemma~\ref{lemm:NeeCompact}.
\end{proof}

\begin{coro} \label{coro:compareHomsLocalTri}
With the notation and assumptions of Proposition~\ref{prop:localisationTriLambda}, for any objects $E, F$ in $\T$ with $F$ compact, the canonical morphism
\[ \hom_\T(F, E) \otimes \ZZ[S^{-1}] \to \hom_{\T[S^{-1}]}(F[S^{-1}], E[S^{-1}]) \]
is an isomorphism.
\end{coro}

\begin{proof}
If $E$ is also compact, then this follows immediately from Corollary~\ref{coro:compareTwoLocalisations}. If not, consider the natural transformation of homological functors
\[ \hom_\T(F, -) \otimes \ZZ[S^{-1}] \to \hom_{\T[S^{-1}]}(F[S^{-1}], -[S^{-1}]). \]
Note also that since $F$ is compact and $(-)[S^{-1}]$ is very nice (Proposition~\ref{prop:localisationTriLambda}) both these functors send small sums to small sums. Hence, the full subcategory of $\T$ on which this natural transformation is an isomorphism is localising, and contains all compact objects. Since $\T$ is compactly generated, this natural transformation is an isomorphism on all of $\T$.
\end{proof}

Now we come to our ``local-global'' principles.

\begin{lemm} \label{lemm:conservativeLocal}
With the assumptions and notation of Proposition~\ref{prop:localisationTriLambda}, a morphism $f: F \to E$ between $\ZZ[S^{-1}]$ local objects is an isomorphism if and only if $f[S_{(\ell)}^{-1}]$ is an isomorphism in $\T[S_{(\ell)}^{-1}]$ for every prime $\ell$ that is not in $S$ where $S_{(\ell)} = (\ZZ - \ell \ZZ)$.
\end{lemm}

\begin{proof}
It is equivalent to prove the following. Let $E$ be a $\ZZ[S^{-1}]$ local object such that $E[S_{(\ell)}^{-1}]$ is zero for every prime $\ell$ that is not in $S$. Then $E$ is zero.

Suppose that $E$ is a $\ZZ[S^{-1}]$ local object such that $E[S_{(\ell)}^{-1}]$ is zero for every prime $\ell$ that is not in $S$. To show that $E$ is zero it suffices to show that $\hom_{\T}(F, E) = 0$ for every compact object $F$ of $\T$ (as $\T$ is compactly generated). Since $E$ is $\ZZ[S^{-1}]$-local, this group is a $\ZZ[S^{-1}]$ module and so it is sufficient to show that $\hom_{\T}(F, E) \otimes \zll = 0$ for every prime $\ell$ that is not in $S$. But now since $F$ is compact we have by Corollary~\ref{coro:compareHomsLocalTri}
\[ \hom_{\T}(F, E) \otimes \zll \cong \hom_{\T[S_{(\ell)}^{-1}]}(F[S_{(\ell)}^{-1}], E[S_{(\ell)}^{-1}]) \]
which vanishes since we are assuming $E[S_{(\ell)}^{-1}] = 0$.
\end{proof}

\begin{lemm} \label{lemm:generatingLocalTri}
With the assumptions and notation of Proposition~\ref{prop:localisationTriLambda} suppose that $\G$ is a small family of compact objects of $\T$. Then the following are equivalent:
\begin{enumerate}
 \item \label{lemm:generatingLocalTri1} The family 
\[ \G[S^{-1}] = \{ E[S^{-1}] : E \in \G\} \]
is a generating family for $\T[S^{-1}]$.
 \item \label{lemm:generatingLocalTri2} The smallest thick triangulated subcategory of $\T[S^{-1}]$ containing $\G[S^{-1}]$ is the full subcategory $\T[S^{-1}]_c$ of compact objects.
 \item \label{lemm:generatingLocalTri3} For each prime $\ell$ which is not in $S$, the family
\[ \G[S_{(\ell)}^{-1}] = \{ E[S_{(\ell)}^{-1}] : E \in \G\} \]
is a generating family of compact objects for $\T[S_{(\ell)}^{-1}]$.
 \item \label{lemm:generatingLocalTri4} The smallest thick triangulated subcategory of $\T[S_{(\ell)}^{-1}]$ containing $\G[S_{(\ell)}^{-1}]$ is the full subcategory $\T[S_{(\ell)}^{-1}]_c$ of compact objects.
\end{enumerate}
\end{lemm}

\begin{proof}
By Proposition~\ref{prop:localisationTriLambda} the categories $\T[S^{-1}]$ also admit small sums and are compactly generated. Moreover, the functor $\T \to \T[S^{-1}]$ preserves compact generating families. So we can apply Lemma~\ref{lemm:NeeCompact} to show that (\ref{lemm:generatingLocalTri1}) $\iff$ (\ref{lemm:generatingLocalTri2}) and (\ref{lemm:generatingLocalTri3}) $\iff$ (\ref{lemm:generatingLocalTri4}). To see that (\ref{lemm:generatingLocalTri1}) implies (\ref{lemm:generatingLocalTri3}) it suffices to note that we can replace $\T$ with $\T[S^{-1}]$ and apply Proposition~\ref{prop:localisationTriLambda} again. Finally, suppose (\ref{lemm:generatingLocalTri3}) is true and consider $E \in \T[S^{-1}]$ such that
\[ \hom_{\T[S^{-1}]}(F[S^{-1}], E) = 0 \]
for every $F \in \G$. We wish to show that $E[S^{-1}] = 0$. As this latter is $\ZZ[S^{-1}]$ local, by Lemma~\ref{lemm:conservativeLocal} it suffices to show that $E[S_{(\ell)}^{-1}] = 0$ for each prime $\ell$ not in $S$. For this, by the assumption that (\ref{lemm:generatingLocalTri3}) is true, it suffices to have 
\[ \hom_{\T[S_{(\ell)}^{-1}]}(F[S_{(\ell)}^{-1}], E[S_{(\ell)}^{-1}]) = 0. \]
For each $F \in \G$. But since $F$ is compact we have
\[ \hom_{\T[S_{(\ell)}^{-1}]}(F[S_{(\ell)}^{-1}], E[S_{(\ell)}^{-1}]) = \hom_{\T}(F, E) \otimes \zll \]
and 
\[ \hom_{\T[S^{-1}]}(F[S^{-1}], E) = \hom_{\T}(F, E) \otimes \ZZ[S^{-1}] \]
and so the required vanishing follows from $\ell$ not being in $S$ together with the assumed vanishing.
\end{proof}


\bibliographystyle{alpha}
\bibliography{bib}

\end{document}